\documentclass[a4paper,11pt]{article}
\usepackage[dvipsnames]{xcolor}
\usepackage{graphicx} 
\usepackage{subcaption}
\usepackage[a4paper,left=2cm,right=2cm,top=3cm,bottom=3cm]{geometry}
\usepackage{amsmath,amssymb,pgfplots}
\usetikzlibrary{arrows.meta}
\pgfplotsset{compat=1.18}
\usepackage{natbib}
\usepackage{amsfonts}
\usepackage{amsthm}
\usepackage{import}
\usepackage{enumerate}
\usepackage{authblk, bm, bbm, soul}
\usepackage{amssymb}
\usepackage{enumitem}
\usepackage{multirow}
\usepackage{capt-of}
\usepackage{booktabs}
\usepackage{hyperref}
\usepackage{cleveref}
\usepackage{natbib}
\hypersetup{colorlinks=true, citecolor=blue, linkcolor = blue, urlcolor = blue}
\usepackage{stackrel}
\usepackage{algorithm}
\usepackage{multirow}
\usepackage{dsfont}
\usepackage{bbm}
\usepackage{algpseudocode}
\usepackage[toc,page]{appendix}

\usepackage{aliascnt}
\newtheorem{theorem}{Theorem}[]
\newtheorem{definition}{Definition}[]

\newtheorem{assumption}{Assumption}
\newtheorem{assumptionprime}{Assumption}

\newtheorem{proposition}[theorem]{Proposition}
\newtheorem{corollary}[theorem]{Corollary}
\newtheorem{lemma}[theorem]{Lemma}
\theoremstyle{remark}
\newtheorem{remark}{Remark}
\newtheorem{example}{Example}[section]
\crefname{example}{example}{examples}
\Crefname{example}{Example}{Examples}
\crefname{assumption}{assumption}{assumptions}
\Crefname{assumption}{Assumption}{Assumptions}
\crefname{assumptionprime}{assumption}{assumptions}
\Crefname{assumptionprime}{Assumption}{Assumptions}

\newcommand{\Reg}{\mathrm{Reg}}

\DeclareMathOperator*{\op}{op}

\DeclareMathOperator*{\argmax}{arg\,max}
\global\long\def\indc{\mathds{1}}
\newcommand{\cX}{\mathcal{X}}

\newcommand{\1}{\mathbf{1}}
\newcommand{\KL}{\mathrm{KL}}

\newcommand{\dist}{\mathrm{dist}}
\newcommand{\supp}{\mathrm{supp}}

\allowdisplaybreaks[4]

\title{Optimal Learning in Fairness-Aware Contextual Bandits}

\author[1]{Gengyu Xue}
\author[2]{Fan Wang}
\author[3]{Feiyu Jiang}
\author[4]{Zifeng Zhao}
\author[1]{Yi Yu}
\affil[1]{Department of Statistics, University of Warwick}
\affil[2]{School of Mathematics and Statistics, University of Melbourne}
\affil[3]{School of Management, Fudan University}
\affil[4]{Mendoza College of Business, University of Notre Dame}
\date{\today}
\begin{document}
\maketitle

\begin{abstract}

Mitigating unfair disparities in sequential decision-making is a growing concern in machine learning.  We study contextual bandits under a global covariate-based Lipschitz fairness constraint, pioneered by  \cite{dwork2012fairness}, requiring individuals with similar contexts to receive similar action distributions.  We show that this constraint induces a randomisation region near the decision boundary and explicitly quantify the population level cost of fairness.  We further establish a canonical representation of the optimal fair policy, which reduces the fair policy design to a lower-dimensional boundary optimisation problem.

We develop a lazy fair greedy algorithm that exploits fairness-induced covariate diversity under continuous contexts, with policy updates computed via finite-dimensional linear programming.  For general context distributions, including discrete and mixed contexts, we propose a lazy
polyhedral upper confidence bound (UCB) algorithm, whose optimistic optimisation reduces to finitely many linear programs over confidence-set vertices.  Both methods satisfy fairness at every round and achieve regret of order $O(\sqrt{T} + T^{(d-2)/(d-1)}{\indc\{d\geq 4\}})$ up to poly-logarithmic factors, where $d$ is the dimension of the contexts.  A matching minimax lower bound for $\ell_1$-based fairness establishes optimality, and reveals an unavoidable nonparametric cost of learning the decision boundary, which dominates the parametric learning cost when  $d \geq 4$. Experiments on synthetic and real data support our theoretical findings and demonstrate the practical effectiveness of our methods.

\end{abstract}

\section{Introduction}

Multi-armed bandits (MABs) provide a fundamental framework for interactive decision-making, in which a policymaker repeatedly selects actions, observes feedback, and adapts its strategy based on past observations. In many modern applications, the policymaker also has access to covariate information that can inform its decisions beyond the feedback collected from previous actions. This naturally motivates the study of contextual MABs, which have become a central framework for adaptive sequential policy designs, with successful applications in healthcare and clinical trials \citep{ameko2020offline}, online recommendation \citep{zhu2023scalable} and finance \citep{cannelli2023hedging}, among others.

In many applications,  designing policies solely to maximise rewards, however, may raise serious ethical or even legal concerns. Empirical evidence from online advertising demonstrates that optimisation-driven delivery systems can produce unequal exposure across demographic groups, even in the absence of explicitly discriminatory targeting \citep{ali2019discrimination}. These concerns are also recognised in broader AI governance frameworks. For instance, the National Institute of Standards and Technology's AI Risk Management Framework of the U.S. emphasises the importance of identifying and mitigating bias and potential harms to individuals, organisations and society through the design and deployment of AI systems \citep{nist2023airmf}. Similarly, China's Provisions on the Administration of Algorithmic Recommendation for Internet Information Services explicitly prohibit the use of algorithmic recommendations to impose unreasonable differential pricing based on consumers' preferences and transaction histories \citep{cac2022algorithm}. In response, algorithmic fairness has emerged as a critical consideration in both research and practice, motivating substantial efforts to mitigate bias in existing machine learning methods and to embed explicit fairness guarantees into sequential decision-making policies.

Fairness in decision-making has been studied from various perspectives, with two prominent approaches focusing on group fairness and individual fairness. Group fairness concerns disparities across groups defined by sensitive attributes, such as race or gender, and typically requires certain statistical measures, such as selection rates or false positive rates, to be comparable across groups \citep[e.g.][]{zeng2026bayes, hou2024finite, fukuchi2023demographic}. Individual fairness instead requires individuals who are similar with respect to the task at hand to receive similar treatment, a property not guaranteed by group-level criteria \citep{dwork2012fairness}.

Individual fairness is particularly relevant to contextual bandit problems, where decisions are made sequentially based on individuals' covariates. In automated lending, for example, an unconstrained  reward-maximising policy may deterministically approve one applicant and reject another with a nearly identical financial profile simply because they lie on opposite sides of a decision boundary. 
This discontinuity underscores   the importance of treating similar individuals similarly, as advocated by the \textit{fairness through awareness} principle of \citet{dwork2012fairness}. 

In this paper, we study contextual bandits under a global covariate-based Lipschitz fairness constraint,
requiring individuals with similar contexts to receive similar action selection probabilities. Imposing this constraint fundamentally changes both the oracle problem structure and the statistical complexity
of bandit learning. Without fairness constraints, the optimal policy in a linear contextual bandit is
determined pointwise by the sign of a linear reward gap between the two arms, yielding a   parametric learning problem
whose minimax regret scales as $\sqrt{T}$ up to logarithmic factors \citep{lattimore2020bandit}. Under Lipschitz fairness, by contrast, a policy cannot switch abruptly across the decision boundary. Instead, it must transition smoothly across a neighbourhood of the boundary, and the optimal shape
of this transition depends on how covariate probability mass is distributed near the boundary. Consequently,
even under a linear reward model, learning an optimal fair policy entails a nonparametric
functional optimisation problem along the $(d-1)$-dimensional decision boundary. We formalise this phenomenon by characterising the geometry and canonical structure of optimal fair policies, quantifying
the population-level cost of fairness, and establishing matching upper and minimax lower regret bounds.

\subsection{Related literature} \label{sec_lit_review}

To highlight the contribution of this paper, we briefly summarise the literature most closely related to fair bandits and interactive decision making in this subsection. 

A prominent strand of the fair bandit literature imposes arm-level allocation constraints, requiring each arm to receive a minimum share of selections. \citet{chen2020fair} study contextual MABs with constraints that guarantee the marginal probability of each arm being pulled is uniformly lower-bounded. \citet{patil2021achieving} consider stochastic Bernoulli bandits with per-arm quota constraints, requiring each arm to be pulled at least a prescribed fraction of rounds. In the combinatorial bandit setting,  where users can choose a combination of arms per round, \citet{li2019combinatorial} enforce asymptotic minimum pull-proportion constraints using virtual queues, while \citet{xu2020combinatorial} impose lower bounds on the empirical selection proportion of each arm throughout the time horizon. Furthermore, \citet{liu2022combinatorial} develop a more general framework for combinatorial bandits with linear constraints, which includes \citet{xu2020combinatorial} as a special case. \citet{wu2024achieving} consider Bernoulli bandits with minimum average reward requirements for each arm, shifting the
requirement from selection frequency to realised reward. Crucially, these formulations enforce fairness
across arms by guaranteeing minimum exposure or reward per action, rather than across individuals
with varying covariates.

Another related line of work  studies individual fairness through meritocratic or similarity-based constraints on the policy. \citet{joseph2016fairness} introduce a meritocratic notion of fairness, requiring that an arm should not be selected with higher probability than another arm unless it has higher expected reward. This meritocratic perspective is further developed by \citet{joseph2018meritocratic}, studying fairness in contextual bandits with infinite convex action sets. In a similar spirit, \citet{wang2021fairness} study fairness of exposure, under which the probability of selecting an arm is required to be proportional to an application-dependent merit function of its expected reward. Moving from merit-based exposure to similarity-based constraints, \citet{liu2017calibrated} study smooth and calibrated fairness in non-contextual MABs, where selection probabilities are constrained by similarities or divergences between reward distributions. \citet{gillen2018online} consider linear contextual bandits with Lipschitz-type individual fairness constraints. Their constraints compare arm-specific contexts within each single round rather than coupling treatment probabilities globally across individuals arriving over the covariate space. 

Other approaches incorporate fairness into the learning objective, rather than imposing it as an explicit constraint. For example, \citet{blanchard2025collaborative} study collaborative bandits using min-max regret as a fairness-aware objective. Similar objective based formulations also include Nash regret in \citet{barman2023fairness} and \citet{sawarni2023nash}, and $\alpha$-fair utility functions in \citet{chaudhary2023alpha}. Group fairness has also been studied in contextual bandits. \citet{grazzi2022group} consider linear contextual bandits with a group meritocratic fairness constraint based on sensitive-group rankings. This notion is closely related to the meritocratic fairness concepts of \citet{joseph2016fairness} and \citet{wang2021fairness}, but applies fairness at the level of sensitive groups rather than individual arms.  These group-level criteria do not generally ensure similar action probabilities for individuals with similar contexts, as required by our individual fairness constraint.

Beyond bandits, individual fairness has also been studied in other interactive decision-making problems. 
\citet{gupta2021individual} introduce a notion named individual fairness in hindsight, comparing current and past decisions through a one-sided similarity constraint, and establish sublinear regret. \citet{salem2026algorithmic} study envy-freeness across groups and time. 
In contextual dynamic pricing, \citet{chen2025utility} 
require similar prices for customers with similar baseline utilities, and characterise optimal fair pricing policies and the regret of demand learning. Although related in spirit, their constraint acts directly on continuous prices, whereas
ours acts on selection probabilities over discrete actions, inducing fundamentally different boundary
geometry and policy-learning challenges that preclude a direct application of their techniques.
Other pricing studies focus on disparities across customer groups and over time. \citet{cohen2025dynamic} consider price and demand fairness, 
while \citet{chen2026fairness} constrain relative price differences between groups under nonparametric demand. 
In online resource allocation, \citet{chen2021fairer} formulate the problem as a linear program and use the analytic centre of its offline optimal solution set to guide fair online allocations.

In this paper, we focus  on individual fairness across contexts, formalised through a global Lipschitz constraint on action selection probabilities in \Cref{def_lipschitz_fairness}, with more discussions in \Cref{sec_setup}.

\subsection{List of contributions}

We summarise the main contributions of the paper.
\begin{itemize}[noitemsep,nolistsep, leftmargin=*]
    \item To the best of our knowledge, this is the first paper to systematically study contextual bandits under a global covariate-based Lipschitz fairness constraint (\Cref{def_lipschitz_fairness}), extending the individual fairness concepts in \citet{dwork2012fairness} to bandit learning. Our framework  incorporates similarity-based fairness into contextual bandit problems, requiring the smoothness of the policy with respect to the underlying context space.

    \item  At the population level, we demonstrate that the Lipschitz fairness fundamentally reshapes the oracle policy by inducing a randomisation region around the decision boundary. The induced randomisation leads to a deviation from the unconstrained reward-maximising policy.  We use this deviation to explicitly quantify the population-level cost of fairness, i.e. the expected reward
loss relative to the unconstrained oracle (\Cref{thm:fairness_cost}).

    \item We establish a canonical representation of  optimal fair policies (\Cref{thm:canonical_extension}), showing how a policy's value on the decision boundary determines its reward-maximising Lipschitz extension over the entire context space. Moreover, we show that the resulting infinite-dimensional empirical boundary optimisation
reduces without relaxation gap to a finite-dimensional linear program. Guided by this structure, we develop a lazy adaptive greedy algorithm for continuous contexts that exploits the fairness induced exploration (\Cref{alg:lazy-adaptive-boundary-greedy,alg:lazy-one-dimensional-fair-greedy}). For more general context distributions, we develop a lazy polyhedral upper confidence bound (UCB) algorithm with efficient vertex-based optimisation (\Cref{alg:lazy-global-model-ucb,alg:lazy-one-dimensional-fair-ucb}).
    
    \item We show that both proposed algorithms satisfy the Lipschitz fairness constraint without slackness and attain a cumulative regret bound of order $\sqrt{T}+T^{(d-2)/(d-1)}\indc(d\geq 4)$, up to poly-logarithmic factors (\Cref{thm:fair_regret_control} and \ref{thm:main-regret}). The additional dimension-dependent term reflects the complexity of learning over the class of Lipschitz functions defined on the $(d-1)$-dimensional decision boundary, arising from the Lipschitz fairness constraints. We further establish a matching minimax lower bound (\Cref{thm:lower_bound}), demonstrating the rate optimality of both proposed algorithms and that, whenever $d \geq 4$, Lipschitz fairness incurs an unavoidable nonparametric   dimension-dependent cost that dominates  the intrinsic $\sqrt{T}$ parametric bandit complexity.

    \item Systematic numerical experiments on both simulated and real datasets are carried out, supporting our theoretical findings, and highlighting their practicality.
    
\end{itemize} The rest of the paper is organised as follows. Section \Cref{sec_formulation} introduces the problem formulation and the Lipschitz fairness constraint. Section \ref{sec_mix_cost} characterises the structural properties of optimal fair policies, including the cost of fairness and the canonical boundary representation. Section \ref{sec_learning_incomplete_information} develops fair learning algorithms and establishes their regret upper bounds under continuous and general context distributions. Section \ref{sec_lower} derives matching minimax lower bounds. Section \ref{sec_numerical}  presents numerical experiments on synthetic and real datasets. Section \ref{sec_conclusion} concludes. Additional algorithms, technical proofs and supplementary numerical results are provided in the Appendix.

\noindent\textbf{Notation.}
For any positive integer $m$, denote $[m]=\{1,\ldots,m\}$. For $a\in\mathbb R$, let $\lceil a\rceil$ and $\lfloor a\rfloor$ denote its ceiling and floor, and let $[a]_{[0,1]}=\min\{1,\max\{0,a\}\}$ denote the clipping of $a$ onto $[0,1]$. For $a,b\in\mathbb R$, write $a\vee b=\max\{a,b\}$ and $a\wedge b=\min\{a,b\}$, and let $\{a\}_+=a\vee0$. For $\nu \geq 1$ and any vector $v \in \mathbb{R}^p$, let $\|v\|_\nu$ denote its $\ell_\nu$-norm, and let $\|v\|_{\nu^*}$ denote its dual norm, where $1/\nu+1/\nu^*=1$, with the convention $1/\infty=0$. For a point $x$ and any sets $A$ and $B$, write $\dist_\nu(x,A)=\inf_{z\in A}\|x-z\|_\nu$, $\dist_\nu(A, B)=\inf_{y \in A, z\in B}\|y-z\|_\nu$, and $\operatorname{diam}_\nu(A)=\sup_{x,z\in A}\|x-z\|_\nu$. For any two symmetric matrices $A$ and $B$, write $A\succeq B$ when $A-B$ is positive semidefinite, and let $\lambda_{\min}(A)$ be the smallest eigenvalue of $A$. For any function $f : \mathcal{X} \rightarrow \mathbb{R}$ and a subset $A \subseteq \mathcal{X}$, denote by $f|_A : A \rightarrow \mathbb{R}$ the restriction of $f$ to~$A$. For two sequences of positive numbers $\{a_n\}$ and $\{b_n\}$, write $a_n\lesssim b_n$ if there exist a universal constant $C>0$ and some $n_0\in\mathbb N$ such that $a_n\le Cb_n$ for all $n\ge n_0$, and write $a_n\gtrsim b_n$ if $b_n\lesssim a_n$. We write $a_n\asymp b_n$ if both $a_n\lesssim b_n$ and $a_n\gtrsim b_n$. We further write $a_n\lesssim_{\log} b_n$ (or $a_n = O_{\log}(b_n)$) if $a_n\lesssim b_n$ holds up to poly-logarithmic factors. For random variables $\{X_n\}$ and positive numbers $\{a_n\}$, write $X_n=O_{\mathrm{p}}(a_n)$ if $\lim_{M\to\infty}\limsup_{n\to\infty}\mathbb P(|X_n|\ge Ma_n)=0$.  For any probability measure $P$, let $\|f\|_{L_2(P)}=\{\mathbb E_P(f^2)\}^{1/2}$. We write $\lambda_d(\cdot)$ for $d$-dimensional Lebesgue measure. 

\section{Problem formulation} \label{sec_formulation}
In this section, we formalise the problem setup for Lipschitz fair contextual bandits in \Cref{sec_setup}, followed by technical assumptions in \Cref{sec_assumption}. 

\subsection{Setup and \texorpdfstring{$(L,\nu)$}{(L,v)}-Lipschitz fairness} \label{sec_setup}

Consider a two-armed contextual bandit problem over a finite time horizon $T\in\mathbb{Z}_+$. The two-arm setting is widely studied and arises naturally in applications such as clinical trials, where the decision involves whether to administer a particular treatment. We focus on this setting for simplicity of exposition and to highlight the key ideas of our framework, with discussions on the extensions to multiple arms ($K\geq 3$) in \Cref{sec_conclusion}.
At each round $t\in[T]$, a context $X_t\in\mathcal X = [0,1]^d$ is observed before an action is chosen, with both arms sharing the same context.
After observing $X_t$, the learner chooses an action $A_t \in [2]$ and receives the realised reward $Y_t=R_{t,A_t}$ of the selected arm.

For $x \in [0, 1]^d$, write $\widetilde x=(1,x^\top)^\top \in \mathbb{R}^{d+1}$.  We assume the linear reward model 
$R_{t,k}=\widetilde X_t^\top\beta_k^*+ \xi_{t,k}$, $k \in [2]$, where $\beta_k^* \in\mathbb R^{d+1}$ is an unknown arm-specific parameter and $\xi_{t,k}$ is a mean-zero noise variable.  Write $\theta^* = \beta_1^*-\beta_2^* = (a^*, \gamma^{*\top})^{\top}$, where $a^* \in \mathbb{R}$ and $\gamma^* \in \mathbb{R}^d$ denote the intercept and slope differences, respectively.

\begin{remark}
    Our framework also extends to generalised linear rewards $R_{t,k}=\mu(\widetilde X_t^\top\beta_k^*)+\xi_{t,k}$, $k\in [2]$, with a known, strictly increasing inverse link function $\mu$, and mean-zero noise $\xi_{t,k}$ \citep[e.g.][]{li2017provably, filippi2010parametric}. Monotonicity of $\mu$ preserves the sign of the reward gap and consequently the linear decision boundary, allowing the structural characterisation of fair policies discussed in \Cref{sec_mix_cost} to extend under suitable regularity conditions. 
\end{remark}

A policy is modelled as a measurable function $p(\cdot):[0,1]^d \to[0,1]$, where $p(x)$ denotes the probability of pulling arm 1 at context $x$, i.e.~$p(x) = \mathbb{P}(A_t = 1|X=x)$. In this setting, at the population level, the standard reward-maximisation objective for bandits is to seek a policy $p$ that maximises the expected reward
\begin{align} \label{eq_reward}
R(p)=\mathbb E\left[\widetilde X^\top\beta_1^*p(X)+\widetilde X^\top\beta_2^*\{1-p(X)\}\right]= \mathbb E\left[\widetilde X^\top\beta_2^*+\Delta^* (X)p(X)\right],
\end{align}
where $\Delta^*(x)=\widetilde x^\top(\beta_1^*-\beta_2^*) = \widetilde x^\top \theta^*$ denotes the conditional mean reward gap between arms 1 and 2 at context $x$, and the expectation is taken with respect to the stochastic context $X$.

By \eqref{eq_reward}, the expected reward $R(p)$ depends on $p$ only through $\mathbb E\{\Delta^*(X)p(X)\}$. Without fairness constraints, the reward-maximising policy simply selects the arm with the larger conditional mean reward at each $x$, i.e.  by  comparing the sign of $\Delta^*(x)$ as
\begin{equation} \label{eq_unconstrained}
    p^*_{\rm un}(x)=\indc\{\Delta^*(x)>0\},
\end{equation}
where the policy values on the decision boundary $\mathcal{H}_{\theta^*} = \{x\in[0,1]^d:\Delta^*(x)=0\}$ do not affect $R(p)$ and ties may be broken uniformly at random.

To incorporate fairness into contextual bandits, we introduce the notion of $(L,\nu)$-Lipschitz fairness, following the individual fairness framework of \citet{dwork2012fairness}.

\begin{definition}[$(L,\nu)$-Lipschitz fairness]\label[definition]{def_lipschitz_fairness}
For $L\ge 0$ and $\nu\in[1,\infty]$, a policy $p(\cdot):[0,1]^d\to[0,1]$ is $(L,\nu)$-Lipschitz fair if
\begin{align*}
    |p(x)-p(x')|\le L\|x-x'\|_\nu,\qquad \text{for all }x,x'\in[0,1]^d.
\end{align*}
We denote the class of all $(L,\nu)$-Lipschitz fair policies by
\begin{align*}
    \mathcal P_{L,\nu}  = \left\{ p:[0,1]^d\to[0,1]:\, |p(x)-p(x')|\le L\|x-x'\|_\nu,\ \text{for all }x,x'\in[0,1]^d \right\}.
\end{align*}
\end{definition}

\Cref{def_lipschitz_fairness} ensures that individuals with similar contexts receive similar treatment probabilities through the $(L, \nu)$-Lipschitz constraints on $p$.  The parameter $L$ controls the strength of this constraint, with   smaller values imposing stronger smoothness constraints. When $L=\infty$, the fairness constraint is inactive, and recovers the unconstrained optimal policy $p^*_{\rm un}$ in~\eqref{eq_unconstrained}. When the constraint is active with  $L<\infty$, fairness  
prevents abrupt shifts in action probabilities across nearby contexts, which may require assigning positive
selection probability to the arm with the smaller conditional mean reward.

\begin{remark}
The fairness requirement in \Cref{def_lipschitz_fairness} is \emph{ex ante}.  Individuals with similar covariates receive similar treatment allocation probabilities, although their realised assignments may differ.   
\citet{wei2023fair} adopt a related perspective when studying group fairness in adaptive experiments.
\end{remark}

\begin{remark}[Choice of the similarity metric]
The similarity metric can be chosen according to the specific application \citep{dwork2012fairness}.  While we work with  the unweighted $\ell_\nu$-norm for concreteness, 
 our framework can readily accommodate alternative metrics, such as the weighted distance $d_W(x,x')=\|W(x-x')\|_\nu$, where $W$ is a suitable weight matrix reflecting task-specific~similarities.
\end{remark}

\begin{remark}
Under suitable regularity conditions, reward-maximising fair policies assign both arms marginal selection probabilities bounded away from zero (\Cref{lem:lipschitz-min-pull-proportion}).  This connects our Lipschitz fairness formulation  (\Cref{def_lipschitz_fairness}) to arm-level minimum pull-proportion fairness notions \citep[e.g.][]{chen2020fair,patil2021achieving}. 
\end{remark}

With \Cref{def_lipschitz_fairness}, for a fixed pair $(L,\nu)$, we define the optimal $(L,\nu)$-Lipschitz fair policy~by
\begin{align*}
   p^*_{L,\nu} \in \argmax_{p \in \mathcal P_{L,\nu}} \ R(p).
\end{align*}
When there is no ambiguity, write $p^*$ for $p^*_{L,\nu}$. In the online learning setting, our goal is to construct a sequence of estimated policies $\{\widehat{p}_t\}_{t=1}^T$ such that $\widehat{p}_{t} \in \mathcal P_{L,\nu}$ for all $t \in [T]$, while minimising the cumulative   regret defined as
\begin{align} \label{eq_accumulated_regret}
\mathfrak R_T(\{\widehat p_t\}_{t=1}^T) = \sum_{t=1}^T \left\{R(p^*)-R(\widehat p_t)\right\}.
\end{align}
Here, we emphasize that for a data-dependent policy $\widehat p_t$ measurable with respect to the history prior to round $t$,  $R(\widehat p_t)=\mathbb{E}\left[\widetilde X^\top\beta_2^*+\Delta^* (X)\widehat p_t|\widehat p_t\right]$ denotes its conditional population reward with respect to a new draw $X\sim P_X$.

\subsection{Assumptions} \label{sec_assumption}
In this subsection, we present the assumptions needed.

\begin{assumption}[Contexts]\label[assumption]{a_contexts}
The contexts $X_1,\ldots,X_T$ are independent and identically distributed (i.i.d.)~following an arbitrary unknown law~$P_X$ supported on $\cX=[0,1]^d$. Assume that $g(\cdot)$, the density of $
P_X$, exists and is   uniformly bounded above and away from zero, i.e.~$ c_g \leq \inf_{x \in[0,1]^d }g(x)  \leq \sup_{x \in[0,1]^d } g(x)  \leq C_g$, where $0 < c_g \leq C_g $ are absolute constants.
\end{assumption}

\Cref{a_contexts} is standard in the nonparametric and contextual bandit literature \citep[e.g.][]{perchet2013multi}. It ensures that the contexts are sufficiently spread across the covariate space, supporting our structural analysis and the theoretical guarantees of the proposed greedy algorithm in \Cref{alg:lazy-adaptive-boundary-greedy}. 
In \Cref{sec_alg_non_cts}, we relax it to \Cref{a_contexts_non_cts}, allowing discrete and mixed context distributions.

\begin{assumption}[Noise]\label[assumption]{a_noise}
For each $t \in [T]$ and $k\in[2]$, the noise variables $\xi_{t,k}$'s are independent across~$t$ and $k$, independent of the contexts $\{X_t\}_{t\ge 1}$ and the algorithmic randomness. Assume that $\xi_{t,k}$'s are mean-zero and sub-Gaussian, i.e.~there exists an absolute constant $\sigma >0$ such that $ \mathbb{E}[\exp(\lambda \xi_{t,k})]  \le \exp(\lambda^2\sigma^2/2)$ for all $\lambda\in\mathbb R$.
\end{assumption}

\begin{assumption}[Bounded parameters]\label[assumption]{a_beta}
There exists an absolute constant $0 <C_\beta<\infty$ such that, for $k\in[2]$, we have $ \|\beta_k^*\|_2\le C_\beta$.
\end{assumption}

\Cref{a_noise,a_beta} regulate the noise distribution and arm regression parameters, respectively. Both are commonly adopted in the linear contextual bandit literature \citep[e.g.][]{goldenshluger2013linear, wu2020stochastic, bastani2021mostly}.

\begin{assumption}[Boundary parameter sets]\label[assumption]{a_theta}
There exist absolute constants $A_0,B_0>0$, $0<c_0<C_0<B_0$, and $\rho_0\in(0,1/2)$ such that $\theta^*=\beta_1^*-\beta_2^*\in \Theta=\bigcup_{j=1}^d \Theta_j$, where, for each $j\in[d]$, $\gamma_{-j}$ denotes the vector obtained by removing the $j$-th coordinate of $\gamma$, and
\begin{align*}
& \Theta_j=\big\{\theta=(a, \gamma^{\top})^{\top}\in\mathbb R\times\mathbb R^d:\; |a|\le A_0,\ \|\gamma\|_2\le B_0,\ c_0\le \gamma_j\le C_0,\\
& \hspace{3cm} \text{and there exists }u_\theta\in[0,1]^{d-1}\text{ such that }
\rho_0 \le -(a+\gamma_{-j}^\top u_\theta)/\gamma_j \le 1- \rho_0
\big\},
\end{align*}
with the convention that $\gamma_{-1}^{\top}u_{\theta}=0$ if $d=1.$
\end{assumption}

\begin{figure}[!htbp]
    \centering
    \includegraphics[width=1\linewidth]{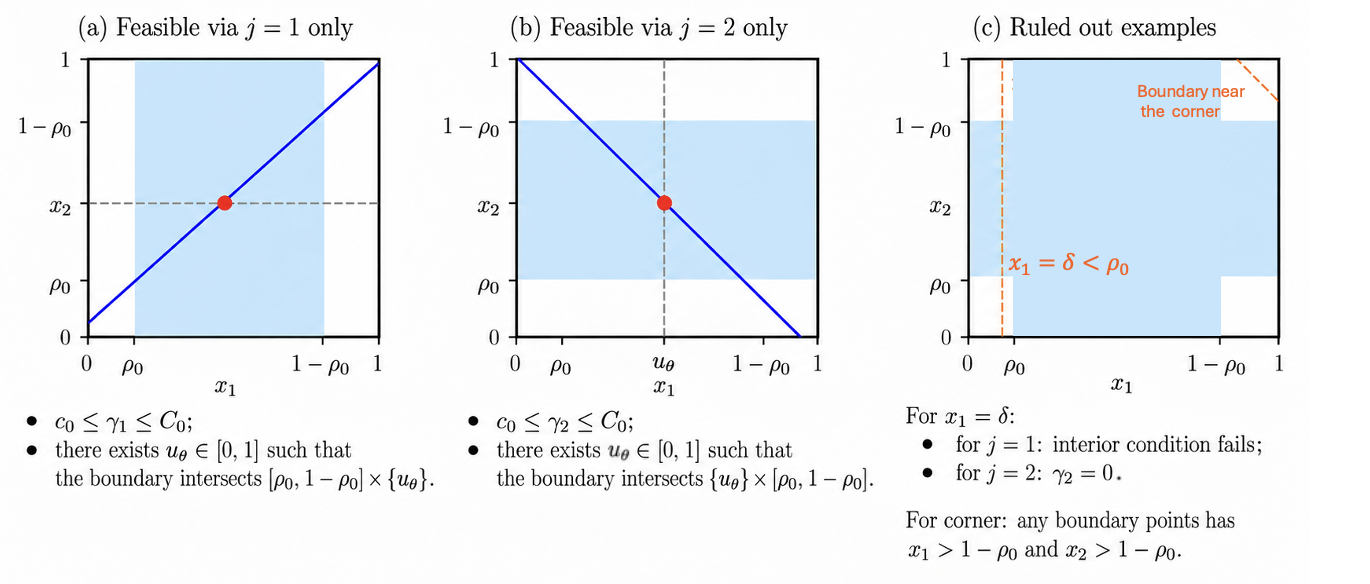}
    \caption{Illustration of \Cref{a_theta} when $d=2$. Plots (a) and (b) illustrate feasible decision boundaries satisfying the conditions for $j=1$ and $j=2$, respectively, with the red points marking the corresponding interior crossings. Plot (c) illustrates two examples of the ruled-out boundaries.}
    \label{fig_boundary_assumption}
\end{figure}

\Cref{a_theta} imposes geometric regularity on the decision boundary $\mathcal H_\theta=\{x\in[0,1]^d:a+\gamma^\top x=0\}$, with  $\theta=(a,\gamma^\top)^\top$.  Each $\Theta_j$, $j \in [d]$, requires the boundary to cross the interior interval $[\rho_0,1-\rho_0]$ along the $j$-th coordinate, with $\gamma_j\in[c_0,C_0]$ preventing the reward gap from vanishing
or degenerating along that coordinate,  which would otherwise reduce the effective dimension of the problem.    Taking the union over $j\in[d]$ allows this direction to vary with $\theta$. Together with the boundedness of $\|\gamma\|_2\leq B_0$, these conditions exclude boundaries confined to arbitrarily small regions near the edges or corners of the covariate space, where one arm dominates almost everywhere; see \Cref{fig_boundary_assumption} for an illustration.  Note that the constants in \Cref{a_theta} need not be known to the learning algorithms proposed in \Cref{sec_learning_incomplete_information}.  They enter only through the constants in the theoretical guarantees.

\section{Structural consequences of the Lipschitz fairness} \label{sec_mix_cost}

In this section, we characterise how the Lipschitz fairness constraint reshapes the optimal fair policy in \Cref{sec_boundary}, explicitly quantify the resulting cost of fairness in \Cref{sec_cost}, and further identify a canonical structure of the optimal fair policy in \Cref{sec_structural}, providing the key structural insight underlying the proposed online learning algorithm in \Cref{sec_learning_incomplete_information}.

\subsection{Boundary randomisation under Lipschitz fairness} \label{sec_boundary}

To understand how the Lipschitz fairness constraint characterises the optimal policy, we first introduce some notation on the geometry induced by a general reward-gap parameter $\theta \in \Theta$. For any $\theta \in \Theta$, the decision boundary $\mathcal H_\theta=\{x\in[0,1]^d:\widetilde x^\top\theta=0\}$ consists of contexts where the two arms have equal expected rewards. 
It separates the covariate space into the regions 
\[\mathcal X_\theta^+=\{x\in[0,1]^d:\widetilde x^\top\theta > 0\} \quad \text{and} \quad \mathcal X_\theta^-=\{x\in[0,1]^d:\widetilde x^\top\theta<0\},\]
where arms 1 and 2, respectively, have the higher expected reward under $\theta$. 

\begin{figure}[!htbp]
    \centering
    \begin{subfigure}[t]{0.48\textwidth}
        \centering
        \includegraphics[width=0.9\textwidth]{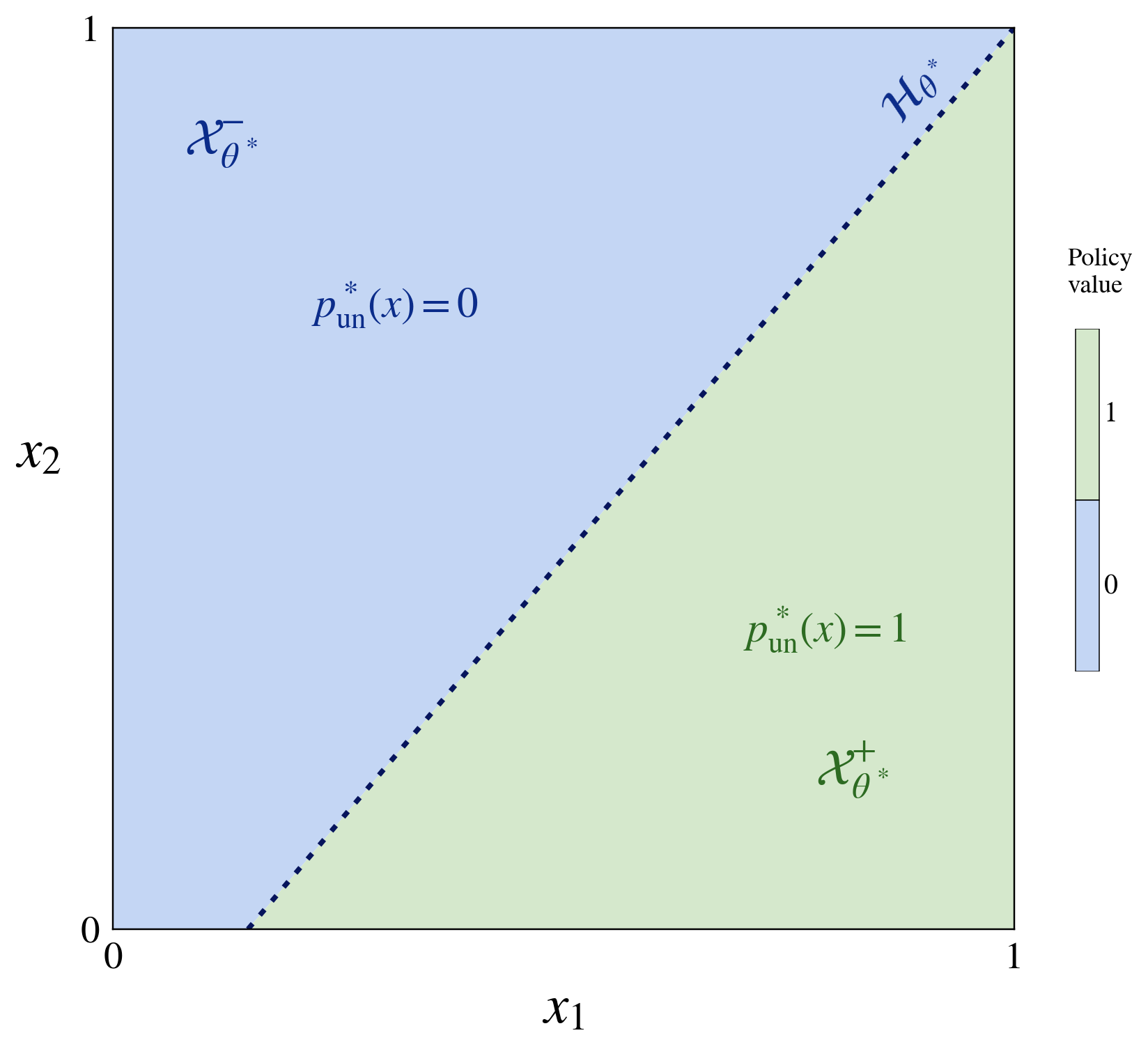}
        \caption{Unconstrained optimal policy.}
        \label{fig_unconstrained}
    \end{subfigure}
    \hfill
    \begin{subfigure}[t]{0.48\textwidth}
        \centering
        \includegraphics[width=0.9\textwidth]{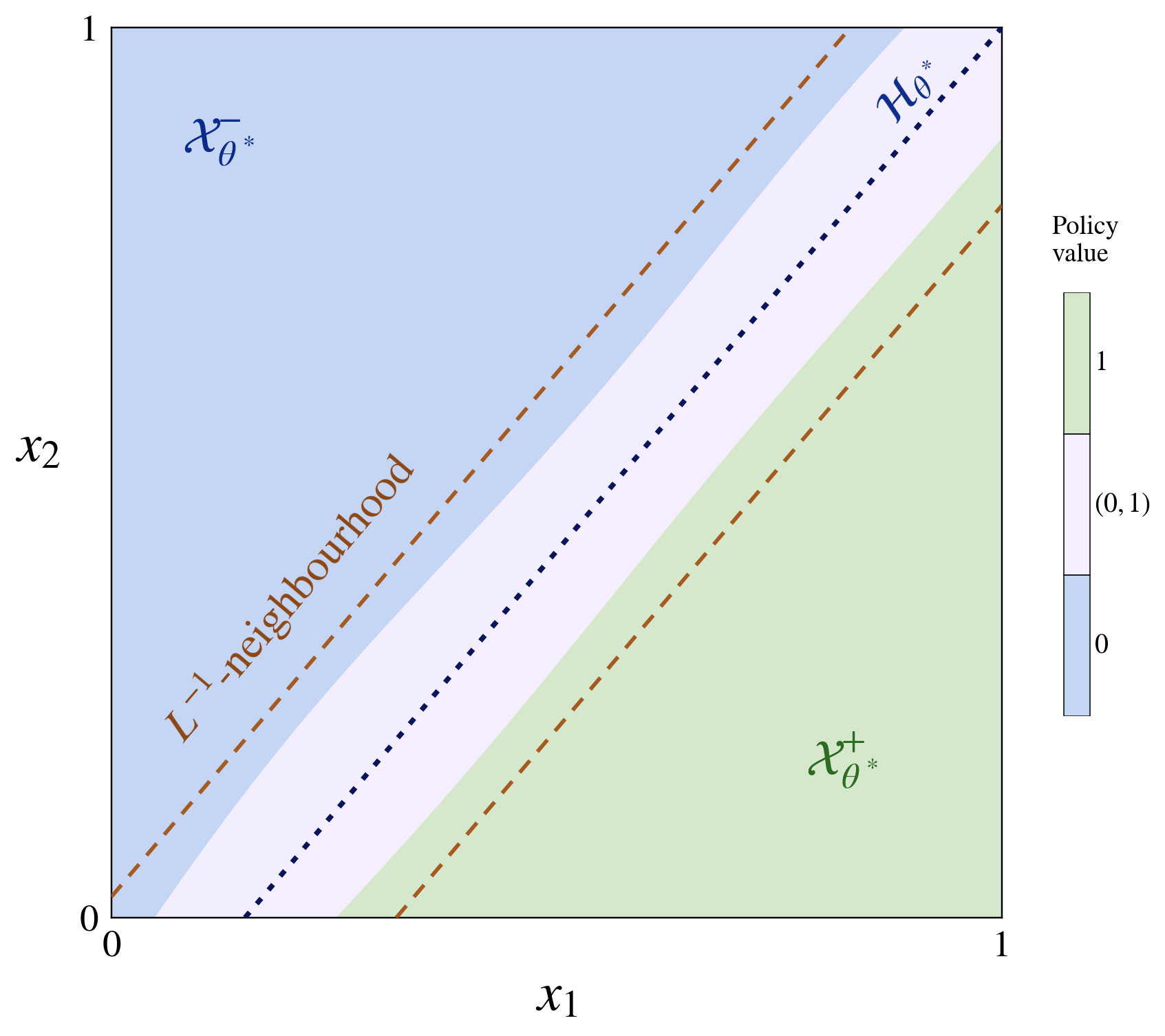}
        \caption{Optimal fair policy under the Lipschitz constraint.}
        \label{fig_mixing_region}
    \end{subfigure}
    \caption{Comparison between the unconstrained and fair optimal policies.}
    \label{fig_policy_comparison}
\end{figure}

For $\theta = \theta^*$, the unconstrained optimal policy $p^*_{\mathrm{un}}$ in \eqref{eq_unconstrained} switches abruptly from always selecting arm~1 on $\mathcal X_{\theta^*}^+$ to always selecting arm 2 on $\mathcal X_{\theta^*}^-$; see \Cref{fig_unconstrained}.
Lipschitz fairness prohibits this discontinuity. The following lemma characterises the resulting randomisation near the decision boundary.

\begin{lemma}\label[lemma]{lem:mixing}
Suppose Assumptions \ref{a_contexts} and \ref{a_theta} hold and fix an $L>0$. Define the randomisation region associated with the optimal policy $p^*$ by $\mathcal M(p^*)=\{x\in [0,1]^d:0<p^*(x)<1\}$, and denote $\theta^* = (a^*, \gamma^{*\top})^{\top}$, where $a^* \in \mathbb{R}$ and $\gamma^* \in \mathbb{R}^d$. It then holds that
\[
    \mathcal M(p^*)\subseteq\big\{x\in [0,1]^d:\dist_{\nu}(x,\mathcal H_{\theta^*})\le L^{-1}\big\} \subseteq \big\{x \in [0,1]^d: |\Delta^*(x)|/\|\gamma^*\|_{\nu^*} \le L^{-1}\big\}.
\]
Furthermore, if $\dist_{\nu}(x,\mathcal H_{\theta^*})>1/L$, then $p^*(x)\in\{0,1\}$, with $p^*(x)= \mathbbm{1}\{\Delta^*(x)>0\}$. 
\end{lemma}

The proof of \Cref{lem:mixing} is given in \Cref{sec_app_thm_boundary}. \Cref{lem:mixing} shows that, rather than switching directly between the two arms across the decision boundary, the optimal fair policy may interpolate between assignment probabilities $0$ and $1$. Such interpolation can only occur inside a $1/L$-neighbourhood of the decision boundary. The randomisation region is, therefore, contained in a slab around the decision boundary with width at most $2/L$ in $\ell_\nu$-distance. Outside this region, $p^*$ agrees with the unconstrained greedy rule $p^*_{\rm un}$, as illustrated in \Cref{fig_mixing_region}.

Under \Cref{a_beta}, $\|\gamma^*\|_{\nu^*} \lesssim 1$, so \Cref{lem:mixing} also implies that the randomisation only occurs where $|\Delta^*(x)|\lesssim 1/L$.  Thus, the optimal fair policy only randomises between arms when the reward gap is of order at most $1/L$, limiting the resulting reward loss.  This observation provides the key intuition behind the analysis of the cost of fairness in \Cref{thm:fairness_cost} in \Cref{sec_cost}.

\subsection{Cost of fairness} \label{sec_cost}
Let $R^*_{\rm un} = R(p^*_{\rm un})$ denote the reward for unconstrained optimal policy in \eqref{eq_unconstrained} and $ R_L^*=\sup_{p\in\mathcal P_{L,\nu}}R(p)$ the optimal reward under $(L,\nu)$-Lipschitz fairness constraint.  Building on the boundary randomisation result in \Cref{lem:mixing}, we quantify the cost of fairness $\operatorname{Cost}(L) =R_{\rm un}^*-R_L^*$ in the following theorem.

\begin{theorem}[Cost of fairness]\label[theorem]{thm:fairness_cost}
Suppose Assumptions \ref{a_contexts}, \ref{a_beta} and \ref{a_theta} hold. Let $d_\nu=d^{1/\nu}$, with the convention $d_\infty=1$.
\begin{enumerate}[noitemsep,nolistsep, leftmargin=*]
\item For every $L>0$, we have that $ \max\{0, \operatorname{Cost}(0)-Ld_\nu\mathbb E|\Delta^*(X)|\}\le \operatorname{Cost}(L)\le \min\{\operatorname{Cost}(0), 1/L^2\}$. 

\item Suppose additionally that $L\ge L_0$ for a constant $L_0<\infty$, depending only on the constants in Assumptions \ref{a_contexts} and \ref{a_theta}. Then, there exist constants $0<c_{\rm fair}\le C_{\rm fair}<\infty$ depending on $d$ and the absolute constants in Assumptions \ref{a_beta} and \ref{a_theta}, such that
\[
    c_{\rm fair}L^{-2}\le \operatorname{Cost}(L)\le C_{\rm fair}L^{-2}.
\]
\end{enumerate}
\end{theorem}
The proof of \Cref{thm:fairness_cost} is in \Cref{sec_app_cost}. Theorem~\ref{thm:fairness_cost} characterises the cost of fairness as the fairness parameter $L$ varies. 

When $L=0$, the fairness constraint forces the policy to be constant across the entire covariate space, so it is optimal to always select the arm with the larger expected reward.  The baseline cost is therefore
\begin{align*}
\operatorname{Cost}(0) = 2^{-1}\left(\mathbb E|\Delta^*(X)|-|\mathbb E\Delta^*(X)|\right) \le 2^{-1}\sqrt{\operatorname{Var}(\Delta^*(X))}
\le (4\sqrt{3})^{-1}\sqrt{C_g}B_0 =O(1),
\end{align*}
where the last inequality follows from \Cref{a_contexts,a_beta,a_theta}. Under these conditions, $\operatorname{Cost}(0)$ is strictly positive for each fixed $d$ and bounded above independent of $d$. Since increasing $L$ enlarges the admissible policy class, $\operatorname{Cost}(L)$ is non-increasing in $L$, with $\operatorname{Cost}(0)$ being its   maximum.

For small $L > 0$, \Cref{thm:fairness_cost} shows that $\operatorname{Cost}(0) - \operatorname{Cost}(L) = O(L)$, so the cost remains bounded away from zero for  sufficiently small $L$.
As $L$ increases further, the cost is controlled by the local behaviour of the policy near the decision boundary. By \Cref{lem:mixing}, randomisation is confined to a slab of width at most of order $O(L^{-1})$, within which the reward gap $|\Delta^*(X)|=O(L^{-1})$, yielding a total reward loss of order at most $O(L^{-2})$.  Thus, $\operatorname{Cost}(L)\to 0$ as $L\to\infty$.

For a constant $L_0$ depending on the boundary separation constant $\rho_0$ in \Cref{a_theta}, when $L \geq L_0$, the upper bound of order $O(L^{-2})$ is sharp.  By \Cref{a_theta}, the decision boundary $\mathcal{H}_{\theta^*}$ contains an interior region away from the edge of the covariate domain.  Across this boundary, the unconstrained oracle $p^*_{\rm un}$ switches between arms, whereas a Lipschitz fair policy requires a distance of at least order $O(L^{-1})$ to transition from probability 0 to 1.  Since the reward gap grows linearly with the distance from the boundary, this restriction creates an unavoidable loss of order $\int_0^{L^{-1}} s\;\mathrm ds\asymp L^{-2}$. The condition $L\ge L_0$ ensures that the boundary band of width $1/L$ remains inside the covariate domain, so that the local lower bound construction is not affected by boundary effects.

\begin{remark}
    In general, the optimal fair policy need not be unique or admit a simple closed-form expression.  In \Cref{sec_app_example}, we consider an example with explicit expressions for an optimal policy and the cost of fairness. Numerical results in \Cref{sec:numerical-simulation} further illustrate the theoretical findings in \Cref{thm:fairness_cost}.
\end{remark}

\subsection{Canonical structures} \label{sec_structural}

We next establish a canonical representation of the optimal fair policy, determined by its values on 
the decision boundary.  We formulate this representation for any 
$\theta\in\Theta$. Taking $\theta=\theta^*$
characterises optimal policies under the true reward model, whereas allowing $\theta$ to vary supports learning with unknown reward parameters in 
\Cref{sec_learning_incomplete_information}.

Motivated by the form of the unconstrained optimal policy in \eqref{eq_unconstrained} and the boundary randomisation phenomenon in \Cref{lem:mixing}, for any $\theta\in\Theta$, 
reward maximisation favours increasing $p$ towards $1$ on $\mathcal X_\theta^+$ and decreasing on $\mathcal X_\theta^-$.
Although policy values on 
$\mathcal H_\theta$ contribute no reward gap, they constrain the feasible policies elsewhere through the Lipschitz fairness constraint. 
This motivates a two-stage
construction: first specify the values on $\mathcal H_\theta$ and then extend  them off the boundary to the rest of the covariate space through an extremal Lipschitz extension.

For any parameter $\theta \in \Theta$, define $\mathcal Q_{\theta,L}$ the collection of $(L,\nu)$-Lipschitz fair policies over boundary,~i.e.
\begin{equation} \label{eq_set_voundary_policy}
    \mathcal Q_{\theta,L}  =
    \left\{ q:\mathcal H_\theta\to[0,1]:  |q(h)-q(h')|\le L\|h-h'\|_{\nu}\ \text{for all } h,h'\in\mathcal H_\theta  \right\}.
\end{equation}
For $q\in\mathcal Q_{\theta,L}$ and $x \in [0,1]^d$, define the upper and lower $(L,\nu)$-McShane envelopes \citep{mcshane1934extension} by 
\begin{equation} \label{eq:envelope-definitions}
    U_{\theta,q}(x)=\inf_{h\in\mathcal H_\theta}\{q(h)+L\|x-h\|_{\nu}\}
    \quad \mbox{and} \quad
    M_{\theta,q}(x)=\sup_{h\in\mathcal H_\theta}\{q(h)-L\|x-h\|_{\nu}\}.
\end{equation}
The boundary-canonical policy associated with $(\theta,q)$ takes the largest feasible value on $\mathcal{X}_{\theta}^+$ and the smallest on $\mathcal{X}_{\theta}^-$, i.e.
\begin{align} \label{eq_canonical}
    p_{\theta,q}(x)
    =
    \begin{cases}
    \min\{1,U_{\theta,q}(x)\}, & x\in\mathcal X_\theta^+,\\
    \max\{0,M_{\theta,q}(x)\}, & x\in\mathcal X_\theta^-.
    \end{cases}
\end{align}
Note that $U_{\theta,q}(x)=M_{\theta,q}(x)=q(x)$ for every $x\in \mathcal{H}_{\theta}$, so $p_{\theta,q}$ agrees continuously along $\mathcal{H}_{\theta}$. 
The following theorem establishes a canonical representation of optimal Lipschitz fair policies through this construction.

\begin{theorem}[Canonical representation]\label[theorem]{thm:canonical_extension}
Let $\theta\in\mathbb R^{d+1}$ satisfy $\mathcal H_\theta\neq\emptyset$. For every $q\in\mathcal Q_{\theta,L}$, it holds that $p_{\theta,q}$ in \eqref{eq_canonical} belongs to $\mathcal P_{L,\nu}$. Moreover, for any $p\in\mathcal P_{L,\nu}$, let $q=p|_{\mathcal H_{\theta}}$, then the canonical extension does not decrease the pointwise reward contribution, i.e.
\[
    (\widetilde x^\top\theta)p_{\theta,q}(x)\ge(\widetilde x^\top\theta)p(x),\qquad x\in\mathcal X.
\]
Consequently, we have that
\begin{equation}\label{eq:canonical-objective-equivalence}
    \sup_{p\in\mathcal P_{L,\nu}}\mathbb E\{(\widetilde X^\top\theta)p(X)\} =\sup_{q\in\mathcal Q_{\theta,L}}\mathbb E\{(\widetilde X^\top\theta)p_{\theta,q}(X)\},
\end{equation}
and both optimisation problems in \eqref{eq:canonical-objective-equivalence} admit at least one maximiser.
\end{theorem}

The proof of \Cref{thm:canonical_extension} is in \Cref{sec_app_thm_canonical}. For each fixed $\theta$, \Cref{thm:canonical_extension} shows that replacing a feasible policy by its canonical extension $p_{\theta,p|_{\mathcal H_\theta}}$ cannot decrease the objective associated with the reward gap $\widetilde x^\top\theta$. For a given candidate parameter, this reduces optimisation over all Lipschitz fair policies to optimisation over their boundary values. Since $\theta^*\in\Theta$, the canonical class $\{p_{\theta,q}:\theta\in\Theta,\ q\in\mathcal Q_{\theta,L}\}\subseteq\mathcal P_{L,\nu}$ contains at least one policy attaining the optimal population reward.

The key reason why the canonical representation cannot decrease the objective value is that, once the policy values $q$ on $\mathcal H_\theta$ are fixed, any $(L,\nu)$-Lipschitz extension $p$ satisfying $p|_{\mathcal H_\theta}=q$ must satisfy
\[
M_{\theta,q}(x)\le p(x)\le U_{\theta,q}(x),\qquad x\in\mathcal X.
\]
The canonical policy attains the upper bound on $\mathcal X_\theta^+$ and the lower bound on $\mathcal X_\theta^-$, maximising the pointwise reward while preserving $q$ and Lipschitz fairness.
When $\theta=\theta^*$, this pushes the policy towards the unconstrained optimal policy $p^*_{\rm un}$ in \eqref{eq_unconstrained} as far as the fairness permits.
This canonical representation guides the construction of \Cref{alg:lazy-adaptive-boundary-greedy} in \Cref{sec_learning_incomplete_information}.

Importantly,  \Cref{thm:canonical_extension} reduces functional optimisation over $d$-dimensional Lipschitz policies $\mathcal P_{L,\nu}$ to optimisation over $(d-1)$-dimensional boundary traces $\mathcal Q_{\theta,L}$. Since $\theta^*\in\Theta$ ensures $\mathcal H_{\theta^*}\neq\emptyset$, the canonical family $\{p_{\theta,q}:\theta\in\Theta,\ q\in\mathcal Q_{\theta,L}\}\subseteq\mathcal P_{L,\nu}$ contains a population-optimal fair policy for any context law $P_X$. This dimension reduction both motivates our algorithms and explains why the nonparametric learning complexity in \Cref{sec_learning_incomplete_information} is governed by the intrinsic boundary dimension $d-1$ instead of the ambient dimension $d$.

\section{Optimal fair policy with incomplete information} \label{sec_learning_incomplete_information} 

In this section, we study online learning of the optimal Lipschitz fair policy under incomplete information, where both the arm-specific reward parameters $\{\beta_1^*, \beta_2^*\}$ and the context distribution $P_X$ are unknown. A distinctive feature of our problem is that, whereas the unconstrained oracle $p^*_{\rm un}$ in \eqref{eq_unconstrained} depends solely on $\theta^* = \beta_1^* - \beta_2^*$, an optimal fair policy  $p^*$ also depends on $P_X$, and the optimal boundary policy values $q^*$ in \eqref{eq:canonical-objective-equivalence} balances reward gains on opposite sides of $\mathcal H_{\theta^*}$ weighted by the covariate distribution near that boundary.  Consequently, even if the reward parameters $(\beta_1^*, \beta_2^*)$ were known, identifying  an optimal fair policy  $p^*$ would still require functional information about the unknown context distribution $P_X$. This introduces an additional nonparametric component to the learning problem, despite the linear specification of the reward functions. 

In the rest of this section, we first introduce a lazy fair greedy algorithm for continuously distributed contexts in \Cref{sec_algorithm}, then establish its fairness and regret guarantees in \Cref{sec_upper}. In \Cref{sec_alg_non_cts}, we relax the density assumption to allow for general context distributions, including discrete and mixed contexts, and develop a UCB-based algorithm with matching theoretical guarantees.

\subsection{Policy learning over continuous feature spaces} \label{sec_algorithm}

Building upon the canonical representation in \Cref{thm:canonical_extension} and the fairness-induced exploration property established in \Cref{lem:lipschitz-min-pull-proportion}, we introduce in \Cref{alg:lazy-adaptive-boundary-greedy}, a lazy fair greedy algorithm for learning under incomplete information.  We focus on the case $d \geq 2$, where the policy is constructed from values on an estimated $(d-1)$-dimensional decision boundary.  The analogous algorithm for $d = 1$, where the boundary is a single point, is given in \Cref{alg:lazy-one-dimensional-fair-greedy} in Appendix.

\begin{algorithm}[!htbp]
\caption{Lazy fair greedy algorithm with boundary learning ($d \ge 2$)}
\label{alg:lazy-adaptive-boundary-greedy}
\begin{algorithmic}
\Require{Time horizon $T \in \mathbb{Z}_+$, fairness parameter $(L,\nu)$, ridge parameter $\lambda > 0$, exploration period $T_o \in \mathbb{Z}_+$}
\For{$k\in\{1,2\}$}
\State $V_k\gets \lambda I_{d+1}$, $b_k\gets (0, \ldots, 0)^{\top} \in \mathbb{R}^{d+1}$.
\EndFor
\State Set $s\gets T_o$.
\State \textbf{Exploration:}
\For{$t \in \{1,2,\ldots,2T_o\}$}
\State Observe $X_t$, set $\widetilde X_t=(1,X_t^\top)^\top$.
\State Sample $A_t = 1$ with probability $1/2$ and $A_t = 2$ otherwise. Observe $Y_t$.
\State Update $V_{A_t}\gets V_{A_t}+\widetilde X_t\widetilde X_t^\top$ and $b_{A_t}\gets b_{A_t}+\widetilde X_tY_t$.
\EndFor
\State \textbf{Adaptive lazy update:}
\For{$t \in \{2T_o+1,\ldots,T\}$}
\State Observe $X_t$ and set $\widetilde X_t=(1,X_t^\top)^\top$.
\If{$t\ge2s$}
\State Set $s\gets t$.
\For{$k\in\{1,2\}$}
\State Set $\widehat\beta_{k,s-1}\gets V_k^{-1}b_k$.
\EndFor
\State Set $\widehat\theta_s\gets\widehat\beta_{1,s-1}-\widehat\beta_{2,s-1}$.
\State Obtain $\widehat p_s$ from \Cref{alg:greedy-boundary-optimisation} and set $\widehat p\gets\widehat p_s$.
\EndIf
\State Sample $A_t = 1$ with probability $\widehat p(X_t)$ and $A_t = 2$ otherwise. Observe $Y_t$.
\State Update $V_{A_t}\gets V_{A_t}+\widetilde X_t\widetilde X_t^\top$ and $b_{A_t}\gets b_{A_t}+\widetilde X_tY_t$.
\EndFor
\end{algorithmic}
\end{algorithm}

\begin{algorithm}[!htbp]
\caption{Boundary optimisation at update time $s$}\label{alg:greedy-boundary-optimisation}
\begin{algorithmic}
\Require{$X_1,\ldots,X_{s-1}$, $\widehat\theta_s$, $L$, and $\nu$.}
\State Solve the finite-dimensional linear program
\begin{equation}\label{eq_alg_lp}
\begin{aligned}
\widehat z_s=(\widehat z_{s,1},\ldots,\widehat z_{s,s-1})\in\argmax_{z_1,\ldots,z_{s-1}}&\quad\frac1{s-1}\sum_{i=1}^{s-1}\widetilde X_i^\top\widehat\theta_s z_i\\
\mathrm{s.t.}&\quad 0\le z_i\le1,\qquad i=1,\ldots,s-1,\\
&\quad |z_i-z_j|\le L\|X_i-X_j\|_\nu,\qquad i,j=1,\ldots,s-1.
\end{aligned}
\end{equation}
\State Set $\widehat p_s^{\rm McExt}(x)\gets\min\{1,\min_{1\le i\le s-1}\{\widehat z_{s,i}+L\|x-X_i\|_\nu\}\}$ for $x\in[0,1]^d$.
\State Set $\widehat q_s\gets\widehat p_s^{\rm McExt}|_{\mathcal H_{\widehat\theta_s}}$ and $\widehat p_s\gets p_{\widehat\theta_s,\widehat q_s}$, with $p_{\widehat\theta_s,\widehat q_s}$ defined in \eqref{eq_canonical}.
\State \Return $\widehat p_s$.
\end{algorithmic}
\end{algorithm}

There are two main components in \Cref{alg:lazy-adaptive-boundary-greedy}: an initial exploration phase followed by an adaptive lazy policy-update phase. 

\textbf{Exploration.} During the first $2T_o$ rounds, the algorithm uses the constant policy $p_t(\cdot)=1/2$ and samples each arm with equal probability, independently across rounds and of the contexts and reward noise. This provides sufficient data to accurately estimate both arm-specific parameters $\beta_1^*$ and $\beta_2^*$ with high probability (\Cref{lem:initial-exploration-estimation}), and avoids projecting subsequent estimators $\widehat\theta_s$ onto $\Theta$, which would require prior knowledge of the constants in \Cref{a_theta}.  Note that  $p_t(\cdot)=1/2$ trivially belongs to $\mathcal P_{L,\nu}$, which is fair.

\textbf{Adaptive lazy updates.} After the exploration phase, the algorithm updates the ridge estimators and constructs a fair plug-in policy using the canonical representation in \Cref{thm:canonical_extension}. Under continuous contexts   in \Cref{a_contexts} and the boundary regularity conditions in \Cref{a_theta},  any canonical policy automatically induces a randomisation region with  positive Lebesgue measure and thus ensures sufficient covariate diversity for both arms  (\Cref{lem:algorithm-policy-diversity}). This property guarantees that both arms can continue to receive sufficiently informative observations, enabling  estimation of their reward parameters during the greedy exploitation  phase without additional forced exploration (\Cref{lem:lipschitz-min-pull-proportion}).

Policy updates are performed lazily using a dyadic schedule with the policy $\widehat{p}$ fixed between two consecutive update rounds. This ensures that the linear program in   \eqref{eq_alg_lp} is solved   at most $\lceil \log_2 (T)\rceil$ times, which  relieves the computational burden.  Similar ideas based on lazy updates have been commonly used in related problems to reduce the number of policy updates while preserving regret guarantees, including linear stochastic bandits \citep[e.g.][]{abbasi2011improved} and fair dynamic pricing \citep{chen2022privacy}. 

At each update time $s$, we fix $\widehat\theta_s$ and chooses the boundary values $q\in\mathcal Q_{\widehat\theta_s,L}$ that maximise the empirical plug-in reward gap,
\begin{equation}\label{eq_alg_optimisation}
\widehat q_s\in\argmax_{q\in\mathcal Q_{\widehat\theta_s,L}}\frac1{s-1}\sum_{i=1}^{s-1}\widetilde X_i^\top\widehat\theta_s p_{\widehat\theta_s,q}(X_i),
\end{equation} 
where $\mathcal Q_{\theta,L}$ and $p_{\theta,q}$ are defined in \eqref{eq_set_voundary_policy} and \eqref{eq_canonical}, respectively. Although \eqref{eq_alg_optimisation} is posed  over an infinite-dimensional  class of Lipschitz functions on $\mathcal H_{\widehat\theta_s}$, \Cref{alg:greedy-boundary-optimisation} computes an exact maximiser via  the finite-dimensional linear program \eqref{eq_alg_lp} without any relaxation gap (\Cref{lem:finite-lp-empirical-canonical-optimizer}).

The linear program in \eqref{eq_alg_lp} in \Cref{alg:greedy-boundary-optimisation} can be interpreted as first finding an optimal discretised fair policy on the observed contexts $\{X_i\}_{i\in[s-1]}$, where the variables $\{\widehat{z}_{s,i}\}_{i \in [s-1]}$ represent the probability of choosing arm $1$ at the sample point $X_i$, i.e.~$\widehat{p}(X_i)$. Once these values on the observed contexts have been determined, the McShane extension then interpolates $\{\widehat{z}_{s,i}\}_{i \in [s-1]}$  to  a globally defined $(L,\nu)$-Lipschitz policy $\widehat p^{\rm McExt}$ over the entire covariate space. Restricting $\widehat p^{\rm McExt}$ to the decision boundary produces a feasible and optimal boundary function $\widehat q_s\in\mathcal Q_{\widehat{\theta}_s,L}$ whose canonical extension $\widehat p_s = p_{\widehat\theta_s,\widehat q_s}$ solves the infinite-dimensional problem \eqref{eq_alg_optimisation}.

\begin{remark}[Equivalence between \eqref{eq_alg_lp} and \eqref{eq_alg_optimisation}]
The equivalence can be derived as follows.
\begin{itemize}[noitemsep,nolistsep, leftmargin=*]
    \item The optimal value of \eqref{eq_alg_lp} provides an upper bound on that of \eqref{eq_alg_optimisation}.  By \Cref{thm:canonical_extension}, every feasible boundary function $q\in\mathcal Q_{\widehat\theta_s,L}$ induces an $(L,\nu)$-Lipschitz canonical policy $p_{\widehat\theta_s,q}$, whose sample values $z_i=p_{\widehat\theta_s,q}(X_i)$, $i\in[s-1]$, satisfy the constraints in \eqref{eq_alg_lp}. Consequently, \eqref{eq_alg_lp} is a relaxation of \eqref{eq_alg_optimisation}, and provides an upper bound. 
    \item The optimal value of \eqref{eq_alg_optimisation} provides an upper bound on that of \eqref{eq_alg_lp}.  The empirical Lipschitz constraints in \eqref{eq_alg_lp} allow its linear programming solution $\{\widehat z_{s,i}\}_{i\in[s-1]}$ to be extended to a globally fair McShane extension which satisfies $\widehat p_s^{\rm McExt}(X_i)=\widehat{z}_{s,i}$ for all $i\in[s-1]$. Restricting $\widehat p_s^{\rm McExt}$ to the estimated decision boundary then yields $\widehat q_s=\widehat p_s^{\rm McExt}|_{\mathcal H_{\widehat\theta_s}}\in\mathcal Q_{\widehat\theta_s,L}$.  By \Cref{thm:canonical_extension}, $(\widetilde X_i^\top\widehat\theta_s)p_{\widehat\theta_s,\widehat q_s}(X_i) \ge (\widetilde X_i^\top\widehat\theta_s)\widehat p_s^{\rm McExt}(X_i) = (\widetilde X_i^\top\widehat\theta_s)\widehat z_{s,i}$, so \eqref{eq_alg_optimisation} attains at least the optimal value of \eqref{eq_alg_lp}.
\end{itemize}
\end{remark}

\begin{remark}[Computational complexity]
  Over the $T$ horizon, the computational complexity of \Cref{alg:lazy-adaptive-boundary-greedy} scales as $O(T^2+\log (T)\mathrm{Cost}_{\mathrm{LP}}(T, T^2) + T\mathrm{Cost}_{\mathrm{canonical}})$, where $T^2$ is the cost of computing $T^2$ constraints in \Cref{alg:greedy-boundary-optimisation}, $\mathrm{Cost}_{\mathrm{LP}}(T,T^2)$ is the cost of solving a linear program with $T$ variables and $T^2$ constraints, and $\mathrm{Cost}_{\mathrm{canonical}}$ is the cost of evaluating $\widehat p_s(X_t)$ via \eqref{eq_canonical} in \Cref{alg:greedy-boundary-optimisation},  which involves the global McShane extension, restriction to the boundary, and then canonical extension. Computational details under the $\ell_1$-metric are provided in Appendix \ref{sec:app-boundary-computation}. 
\end{remark}

\subsection{Fairness and regret control} \label{sec_upper}

In this subsection, we present the main theorem analysing the fairness and regret of Algorithms \ref{alg:lazy-adaptive-boundary-greedy}~and~\ref{alg:lazy-one-dimensional-fair-greedy}.  

\begin{theorem}[Fairness and regret control]
\label[theorem]{thm:fair_regret_control}
Suppose that Assumptions \ref{a_contexts}-\ref{a_theta} hold and $L$ is an absolute constant satisfying $0<L<L_0'$ for some $L'_0 >0$. For any $\eta \in (0,1/2)$, take $\lambda=C_\lambda$ and $T_o = C_1\lceil\sqrt{T\log(T/\eta)}\rceil$, where $C_1, C_{\lambda}>0$ are absolute constants. Then the following statements hold for Algorithms \ref{alg:lazy-adaptive-boundary-greedy} and \ref{alg:lazy-one-dimensional-fair-greedy}.
\begin{enumerate}[noitemsep,nolistsep, leftmargin=*]
    \item  \textbf{(Fairness.)} The collection of policies $\{\widehat{p}_t\}_{t=1}^T$ satisfies $(L,\nu)$-Lipschitz fairness, i.e.~$\widehat{p}_{t} \in \mathcal P_{L,\nu}$ for all $t \in [T]$.

    \item  \textbf{(Regret.)} There exist constants $C_2, C_3,C_4 >0 $, depending only on $d,L,c_g,C_g,\sigma$, $C_\beta,\lambda, \nu$ and the constants in \Cref{a_theta}, such that whenever $T\ge C_4\log(T/\eta)$, it holds with probability at least $1-\eta$ that
    \[
    \mathfrak R_T(\{\widehat p_t\}_{t=1}^T)
    \le C_2\sqrt{T\log(T/\eta)}
    +C_3\begin{cases}
    \sqrt T\log(eT^2/\eta), & d= 1,2,\\[2pt]
    \sqrt T\log^{\frac{3}{2}}(eT^2/\eta), & d=3,\\[2pt]
    T^{\frac{d-2}{d-1}}\log(eT^2/\eta), & d\ge4,
    \end{cases}
    \]
    where $\mathfrak R_T(\{\widehat p_t\}_{t=1}^T)$ is the accumulated regret defined in \eqref{eq_accumulated_regret}.
\end{enumerate}
\end{theorem}

We defer the detailed proof of \Cref{thm:fair_regret_control} to \Cref{sec_app_proof_regret}, with a sketch at the end of this subsection. \Cref{thm:fair_regret_control} demonstrates that both Algorithms \ref{alg:lazy-adaptive-boundary-greedy} and \ref{alg:lazy-one-dimensional-fair-greedy} satisfy the fairness constraints and achieve a regret sublinear in $T$.
\Cref{thm:fair_regret_control} considers a fixed constant fairness parameter $L \in (0, L'_0)$, for some $L'_0>0$, allowing us to study how regret scales with the time horizon under a fixed fairness requirement.  

Up to poly-logarithmic factors, the upper bound on regret in \Cref{thm:fair_regret_control} can be rewritten as
\begin{align*}
    \mathfrak R_T = O_{\mathrm{p}}\Big(\sqrt{T} + T^{\frac{d-2}{d-1}}\indc\{d\geq 4\}\Big),
\end{align*}
which at a high level, is of the form 
\[\mathfrak R_T = \text{ intrinsic bandit regret } + \text{ fairness-induced regret}.\]

This decomposition highlights the coexistence of a parametric bandit learning cost and a nonparametric fairness cost. The first term, of order $\sqrt T$, captures the classical cost of learning finite-dimensional linear reward parameters under bandit feedback, which represents the intrinsic minimax parametric bandit regret  \citep[e.g.][]{dani2008stochastic,chu2011contextual}. The second term, of order $T^{(d-2)/(d-1)}$ for $d\geq 4$, captures the additional statistical cost of enforcing Lipschitz fairness by learning the unknown $(d-1)$-dimensional decision boundary.  In particular, for $n$ historical contexts,   the uniform  empirical-process  deviation over the $(d-1)$-dimensional boundary-canonical  Lipschitz  class decays at the nonparametric rate $n^{-1/(d-1)}$ when $d-1 \ge 3$, which sums over $n\in[T]$ to $O(T^{(d-2)/(d-1)})$.  Notably, similar dimension-dependent rates also arise in other  nonparametric contextual bandit problems involving Lipschitz function classes, e.g.~\citet{cesa2017algorithmic} and \citet{mao2018contextual}, where the regret scales with the ambient dimension $d.$ By contrast, in our setting, \Cref{thm:canonical_extension} reduces the effective nonparametric dimension to $d-1$. Consequently, for $d\leq 3$, the fairness-induced term is absorbed into the parametric $\sqrt{T}$ rate  up to logarithmic terms \citep[see, e.g.,][Theorem~2.14.2]{van1996weak}. For $d\ge4$, it strictly dominates the overall regret, making learning the fair policy on the boundary the main statistical difficulty.

\textbf{Sketch of the proof.} We conclude this subsection with a sketch of the proof of the regret upper bound in \Cref{thm:fair_regret_control} for $d\geq2$, while the fairness guarantee is a direct consequence of \Cref{thm:canonical_extension}.

Under the boundedness assumptions in Assumptions~\ref{a_contexts} and \ref{a_beta}, the initial exploration incurs $O(T_o) = O_{\rm\log}(\sqrt{T})$ regret and localises $\widehat\theta_s$ within an enlarged regular neighbourhood $\widetilde\Theta \supset \Theta$ across all lazy updates with high probability (\Cref{lem:localization-and-estimation}). In each subsequent lazy update phase starting at $s$, the regret has two sources: (i) estimation error in $\theta^{\ast}$, which leads to a mismatch between the true and estimated decision boundaries $\mathcal H_{\theta^*}$ and $\mathcal H_{\widehat \theta}$, and (ii) empirical-process statistical error from replacing the population reward objective by its empirical plug-in counterpart in \eqref{eq_alg_optimisation}.  Since the lazy update batches grow geometrically, summing the per-update bounds over all batches inflates the regret only by logarithmic factors.

Specifically, at each update  time $s$, we compare the learned policy $\widehat p_s$ with the optimal fair policy~$p^*$ through an intermediate policy defined by the canonical extension over the estimated boundary. Let $q_s^\dagger=p^{*}|_{\mathcal H_{\widehat\theta_s}}$ be the restriction of the oracle policy to the estimated decision boundary, and define $p_s^\dagger=p_{\widehat\theta_s,q_s^\dagger}$. Since $p^*$ is $(L,\nu)$-Lipschitz, it holds that $q_s^\dagger\in\mathcal Q_{\widehat\theta_s,L}$, and $p_s^\dagger$ is a feasible competitor for the optimisation problem in \eqref{eq_alg_optimisation}. The instantaneous regret is then decomposed as
\begin{align*}
R(p^{*})-R(\widehat p_s) =& \{R(p^{*})-R(p_s^\dagger)\} + \{R(p_s^\dagger)-R(\widehat p_s)\}\\
=& \text{ boundary mismatch error } + \text{ empirical reward estimation error}.
\end{align*}
\begin{itemize}[noitemsep,nolistsep, leftmargin=*]
\item  Boundary mismatch error: To bound $R(p^{*})-R(p_s^\dagger)$, we partition $\mathcal X$ according to whether $x$ lies on the same side or opposite sides of $\mathcal H_{\theta^*}$ and $\mathcal H_{\widehat\theta_s}$. On the same-side event,  both $p^*(x)$ and $p_s^\dagger(x)$ use the   same envelope type, hence 
$|p^*(x) - p_s^\dagger(x)|$ can be upper bounded by $O(\|\widehat\theta_s - \theta^*\|_2)$ 
using the Lipschitz property (\textbf{Steps 2 and 3} in \Cref{lem:expected-transfer-no-common-chart}) and Hausdorff stability of $\mathcal{H}_{\widehat{\theta}_s}$ (\Cref{lem:hausdorff}). On the opposite-side event, since $|p^*(x) - p_s^\dagger(x)| \le 1$,  we use a crude pointwise bound of order $O(1)$.  This event has small probability due to the stability of the decision boundary.  \Cref{lem:hausdorff} controls the Hausdorff distance between the true and estimated decision boundaries by $O(\|\widehat{\theta}_s - \theta^*\|_2)$, confining the sign disagreement to a thin slab of the covariate around the true boundary under \Cref{a_contexts} (\Cref{lem:slab_probability}). 
Combining the Lipschitz control on the same-side event with the slab probability bound on the sign-disagreement event shows that the loss caused by replacing $\mathcal H_{\theta_s^*}$ with $\mathcal H_{\widehat\theta_s}$ is only of the same order as the parameter estimation error $\|\widehat\theta_s-\theta^*\|_2$ (\Cref{lem:expected-transfer-no-common-chart}).

\item 
Empirical reward estimation error: Using the empirical optimality of $\widehat{q}_s$ (and hence $\widehat{p}_s$) in \eqref{eq_alg_optimisation}, and that $|\widetilde{X}_i^{\top}(\widehat{\theta}_s-\theta^*)|\lesssim \|\widehat\theta_s - \theta^*\|_2$, then up to the estimation error of $O(\|\widehat\theta_s-\theta^*\|_2)$, the term     $R(p_s^\dagger)-R(\widehat p_s)$ can be bounded by the uniform empirical deviation over the canonical policy class, 
\begin{equation}\label{eq-uniform-deviation}
\sup_{\theta \in \Theta,\ q\in\mathcal Q_{\theta,L}}\Big|\frac1n\sum_{i=1}^n \widetilde X_i^\top \theta^* p_{\theta,q}(X_i)-\mathbb E\big \{\widetilde X^\top \theta^* p_{\theta,q}(X)\big\}\Big|.
\end{equation}
Since $\mathcal{Q}_{\widehat\theta_s,L}$ is parametrised by functions living only on the decision boundary, the complexity of $\mathcal{Q}_{\widehat\theta_s,L}$ scales with the intrinsic boundary dimension $d-1$, instead of the ambient dimension $d$ (\Cref{lem:l2-bracketing-policy}). This low-dimensional representation reduces the learning problem to an $L$-Lipschitz function class on the $(d-1)$-dimensional decision boundary, giving bracketing entropy of order $O_{\log}(\epsilon^{-(d-1)})$.
Letting~$n$ be the effective sample size, the bracketing maximal inequality \citep[Theorem~2.14.2]{van1996weak} leads to a high probability bound on $\eqref{eq-uniform-deviation}$, which is of order $O_{\log}(n^{-1/2})$ when $d \in \{2, 3\}$, and $O_{\log}(n^{-1/(d-1)})$ when $d \geq 4$. Summing over the dyadic schedule  gives the rate in \Cref{lem:control-delta-L-n} and explains the dimension-dependent term in \Cref{thm:fair_regret_control}.
\end{itemize}

\subsection{Extensions beyond continuous covariates} \label{sec_alg_non_cts}

Our analysis so far is established for continuously distributed contexts, detailed in \Cref{a_contexts}.  It 
allows us to control the probability mass near the decision boundary, crucial in quantifying the cost of fairness in \Cref{sec_cost}. It also guarantees sufficient covariate diversity for both arms under the canonical policies, enabling the greedy learning approach and its regret guarantees in \Cref{sec_algorithm,sec_upper}.  To allow
discrete and mixed discrete-continuous covariates, 
we relax \Cref{a_contexts} as follows.

\begin{assumptionprime}[General contexts] \label{a_contexts_non_cts}
The contexts $X_1,\ldots,X_T$ are i.i.d.~from an arbitrary, unknown Borel probability law $P_X$ supported on $\cX=[0,1]^d$. Let $\mathcal S=\supp(P_X)$ and $\Sigma_X=\mathbb E(\widetilde{X}_1\widetilde{X}_1^\top)$. We assume that there exists an absolute constant $\kappa_X>0$ such that $\lambda_{\min}(\Sigma_X)\ge\kappa_X$.
\end{assumptionprime}

\Cref{a_contexts_non_cts} permits continuous, discrete and mixed context distributions, and the support $\mathcal S$ can be unknown. We continue to define fair policies on the entire covariate space $\cX$ without loss of generality, since any Lipschitz fair policy defined on~$\mathcal S$ can be extended to $\cX$ while preserving its Lipschitz constant and its values on $\mathcal S$ (\Cref{lem:support-extension} in \Cref{sec_app_structural_general}). Conversely, any globally defined fair policy remains fair when restricted to $\mathcal S$. Since $X\in\mathcal S$ almost surely, the expected reward depends only on policy values on $\mathcal S$, and hence the two formulations yield the same optimal reward.  Even when the decision boundary $\mathcal H_\theta$ does not intersect the support $\mathcal S$, the canonical representation in \Cref{thm:canonical_extension} remains valid for arbitrary context distributions, since its construction relies on the Lipschitz geometry rather than the distribution of contexts. We defer structural results under general context distributions, including the boundary randomisation property restricted to $\mathcal S$, to \Cref{sec_app_structural_general}.

Although the structural characterisation remains applicable, the covariate diversity property underlying the greedy algorithm is no longer guaranteed under \Cref{a_contexts_non_cts}:  for instance, $\mathcal S$ may place zero probability mass inside the randomisation region around $\mathcal H_{\theta^*}$, causing a plug-in canonical policy to act deterministically on $\mathcal S$ and under-explore one arm.  To overcome this difficulty, we develop a UCB-based learning algorithm that explicitly accounts for uncertainty in the reward parameters, achieving fairness and regret guarantees for general context distributions. For $d \geq 2$, the algorithm and its optimistic update are presented in \Cref{alg:lazy-global-model-ucb,alg:vertex-oracle}.  The $d = 1$ analogue is deferred to \Cref{alg:lazy-one-dimensional-fair-ucb} in \Cref{sec:one-d-ucb}.

\begin{algorithm}[!htbp]
\caption{Lazy polyhedral UCB based algorithm $(d \geq 2)$}\label{alg:lazy-global-model-ucb}
\begin{algorithmic}
\Require{$T,T_o\in\mathbb N$, $2T_o<T$, $\lambda>0$, failure probability $\eta\in(0,1)$, $\sigma$, $C_\beta$, $L$, $\nu$.}
\State Set $s_0\gets2T_o+1$ and tolerance $\eta_{\rm c}\gets\eta/8$. Define
\begin{equation}\label{eq_ucb_alpha}
\alpha_T\gets\sqrt{2\sigma^2\left\{\log\left(\frac{2}{\eta_{\rm c}}\right)+\frac{d+1}{2}\log\left(1+\frac T\lambda\right)\right\}}+C_\beta\sqrt\lambda.
\end{equation}
\State Initialise $V_k\gets\lambda I_{d+1}$, $B_k\gets0$, and $\ell_k\gets\log\det(V_k)$ for $k\in\{1,2\}$.
\State For $t\ge1$ and $k\in\{1,2\}$, write
\begin{equation*}
V_{t,k}=\lambda I_{d+1}+\sum_{i=1}^{t-1}\widetilde X_i\widetilde X_i^\top\indc\{A_i=k\},\qquad \widehat\beta_{t,k}=V_{t,k}^{-1}\sum_{i=1}^{t-1}\widetilde X_iY_i\indc\{A_i=k\}.
\end{equation*}
\For{$t=1,\ldots,2T_o$}
\State Observe $X_t$ and set $\widetilde X_t=(1,X_t^\top)^\top$.
\State Set $\widehat p_t(\cdot)\gets 1/2$, sample $A_t=1$ with probability $1/2$ and $A_t=2$ otherwise, and observe $Y_t$.
\State Update
\begin{equation*}
V_{A_t}\gets V_{A_t}+\widetilde X_t\widetilde X_t^\top,\qquad B_{A_t}\gets B_{A_t}+\widetilde X_tY_t,\qquad \ell_{A_t}\gets\log\det(V_{A_t}).
\end{equation*}
\EndFor
\State Set $s\gets s_0$ and for $k\in\{1,2\}$,  set $V_{s,k}\gets V_k$ and $\ell_{s,k}\gets\ell_k$.
\For{$t=s_0,\ldots,T$}
\If{$t=s_0$ \textbf{or} $t\ge2s$ \textbf{or} $\exists k\in\{1,2\}:\ell_k-\ell_{s,k}\ge\log(2)$}
\State Set $s\gets t$ and and, for $k\in\{1,2\}$, set $V_{s,k}\gets V_k$, $\widehat\beta_{s,k}\gets V_k^{-1}B_k$, and $\ell_{s,k}\gets\ell_k$. Obtain $p_s$ from \Cref{alg:vertex-oracle} and set $\widehat p\gets p_s$.
\EndIf
\State Observe $X_t$ and set $\widetilde X_t=(1,X_t^\top)^\top$.
\State Sample $A_t=1$ with probability $\widehat p(X_t)$ and $A_t=2$ otherwise, and observe $Y_t$.
\State Update
\begin{equation*}
V_{A_t}\gets V_{A_t}+\widetilde X_t\widetilde X_t^\top,\qquad B_{A_t}\gets B_{A_t}+\widetilde X_tY_t,\qquad \ell_{A_t}\gets\log\det(V_{A_t}).
\end{equation*}
\EndFor
\end{algorithmic}
\end{algorithm}

\begin{algorithm}[!htbp]
\caption{Vertex optimisation at update time $s$}\label{alg:vertex-oracle}
\begin{algorithmic}
\Require{$X_1,\ldots,X_{s-1}$, $V_{s,k},\widehat\beta_{s,k}$ for $k\in\{1,2\}$, $\alpha_T$, $L$, and $\nu$.}
\State For $k\in\{1,2\}$, form the vertex set as
\begin{equation}\label{eq:confidence-vertices}
\mathcal V_{s,k}=\left\{\widehat\beta_{s,k}\pm\sqrt{d+1}\alpha_T V_{s,k}^{-1/2}e_r:r\in[d+1]\right\},
\end{equation}
where $e_r$ is the $r$th standard basis vector and $V_{s,k}^{-1/2}$ is the symmetric inverse square root.
\For{each pair $v=(v_1,v_2)\in\mathcal V_{s,1}\times\mathcal V_{s,2}$}
\State Set $\theta_v=v_1-v_2$ and solve (with $\mathcal Z_{s-1,L}$ defined in \eqref{eq:ZnL})
\begin{equation}\label{eq:vertex-lp}
z_v\in\argmax_{z\in\mathcal Z_{s-1,L}}
\frac1{s-1}\sum_{i=1}^{s-1}(\widetilde X_i^\top\theta_v)z_i.
\end{equation}
\State Set
\begin{equation*}
S_v=\frac1{s-1}\sum_{i=1}^{s-1}\widetilde X_i^\top v_2
+\frac1{s-1}\sum_{i=1}^{s-1}(\widetilde X_i^\top\theta_v)z_{v,i}.
\end{equation*}
\EndFor
\State Select
\begin{equation*}
\widehat v=(\widehat v_1,\widehat v_2)
\in\argmax_{v\in\mathcal V_{s,1}\times\mathcal V_{s,2}}S_v.
\end{equation*}
\State Set $\bar\beta_s=(\bar\beta_{s,1},\bar\beta_{s,2})=(\widehat v_1,\widehat v_2)$ and define $\theta_s=\bar\beta_{s,1}-\bar\beta_{s,2}$. Let $z_s=z_{\widehat v}$ be the optimiser of \eqref{eq:vertex-lp} corresponding to the selected vertex pair $\widehat v$, with $z_{s,i}=z_{\widehat v,i}$ for $i\in[s-1]$.
\State Define
\begin{equation*}
\overline p_s(x)=\min\left\{1,\min_{1\le i\le s-1}
\left[z_{s,i}+L\|x-X_i\|_\nu\right]\right\}.
\end{equation*}
\State Construct $q_s=\overline p_s|_{\mathcal H_{\theta_s}}$ and \Return $p_s=p_{\theta_s,q_s}$.
\end{algorithmic}
\end{algorithm}

Unlike the greedy approach in \Cref{alg:lazy-adaptive-boundary-greedy}, which constructs policies using point estimates of the reward parameters, the proposed algorithm adopts an optimistic learning strategy inspired by UCB techniques \citep[e.g.][]{abbasi2011improved}. By optimising over confidence sets of plausible reward parameters, the algorithm explicitly accounts for estimation uncertainty and encourages exploration when the reward parameters are insufficiently learned, without relying on the covariate diversity property established under \Cref{a_contexts}.

\noindent \textbf{Confidence set construction.}
At each update time $s$, instead of directly using the point estimators~$\widehat\beta_{s,k}$ as in \Cref{alg:lazy-adaptive-boundary-greedy}, we construct a confidence set $\mathcal C_{s,k}$ for each arm and optimise over all candidate reward models contained in these sets, which preserve statistical coverage while making the optimistic policy optimisation tractable. 
With $V_{s,k}$ and $\widehat\beta_{s,k}$ defined in \Cref{alg:lazy-global-model-ucb}, we form
\begin{equation}\label{eq:candidate-set}
\mathcal C_{s,k}=\left\{b:\left\|V_{s,k}^{1/2}(b-\widehat\beta_{s,k})\right\|_1\le\sqrt{d+1}\alpha_T\right\},
\end{equation}
which replaces the usual $\ell_2$-ball of $V_{s,k}^{1/2}(b-\widehat\beta_{s,k})$ with radius $\alpha_T$ in the standard UCB algorithms by an $\ell_1$-ball with enlarged radius $\sqrt{d+1}\alpha_T$. This modification preserves the confidence coverage while replacing the ellipsoidal region with a polytope of $2(d+1)$ explicit vertices, motivating the design of the computationally-efficient vertex optimisation in \Cref{alg:vertex-oracle}. 

\noindent \textbf{Vertex optimisation.}  
The model parameters and policy are estimated via the optimisation problem 
\begin{align}
& (\bar\beta_{s,1}^{\circ},\bar\beta_{s,2}^{\circ},q^\circ_s) \in\argmax_{\substack{b_k\in\mathcal C_{s,k},\ k\in\{1,2\}, \, \mathcal{H}_{b_1 - b_2} \neq \emptyset \\q\in\mathcal Q_{b_1-b_2,L}}} \frac1{s-1}\sum_{i=1}^{s-1}\left[\widetilde X_i^\top b_1p_{b_1-b_2,q}(X_i)+\widetilde X_i^\top b_2\{1-p_{b_1-b_2,q}(X_i)\}\right]  \nonumber \\
= & \argmax_{\substack{b_k\in\mathcal C_{s,k},\ k\in\{1,2\} \\ \mathcal{H}_{b_1 - b_2} \neq \emptyset}} \left\{\frac1{s-1}\sum_{i=1}^{s-1}\widetilde X_i^\top b_2+\max_{q\in\mathcal Q_{b_1 - b_2,L}}\frac1{s-1}\sum_{i=1}^{s-1}\{\widetilde X_i^\top(b_1 - b_2)\}p_{b_1 - b_2,q}(X_i)\right\},
\label{eq:ideal-update}
\end{align}
consisting of two optimisation problems.
\begin{itemize}[noitemsep,nolistsep, leftmargin=*]
    \item The inner optimisation over $\mathcal Q_{b_1 - b_2,L}$ has an exact linear program representation. Specifically, just as the greedy boundary problem \eqref{eq_alg_optimisation} in \Cref{alg:lazy-adaptive-boundary-greedy} reduces without relaxation gap to the finite-dimensional linear program \eqref{eq_alg_lp}, for any fixed pair $(b_1, b_2)$ with $\theta = b_1 - b_2$, the inner maximisation in \eqref{eq:ideal-update} reduces to 
\begin{equation}\label{eq:Hs}
H_s(\theta)=\max_{z\in\mathcal Z_{s-1,L}}\frac1{s-1}\sum_{i=1}^{s-1}(\widetilde X_i^\top\theta)z_i,
\end{equation}
where
\begin{equation}\label{eq:ZnL}
\mathcal Z_{s-1,L}=\{z\in[0,1]^{s-1}:|z_i-z_j|\le L\|X_i-X_j\|_\nu,\ i,j\in[s-1]\}.
\end{equation}
An optimiser $z$ of \eqref{eq:Hs} can be interpolated to a global Lipschitz policy by the McShane extension and then replaced by its boundary-canonical representative without decreasing the empirical objective.  A detailed formulation of the equivalence between the inner optimisation in \eqref{eq:ideal-update} and  \eqref{eq:Hs} is left in \Cref{lem:finite-lp-empirical-canonical-optimizer}.

\item The outer optimisation over $\mathcal{C}_{s, k}$'s in \eqref{eq:ideal-update}, with the notation in \eqref{eq:Hs}, can be written as
\begin{equation}\label{eq:Fs}
\max_{b_k\in\mathcal C_{s,k},\ k\in\{1,2\}}F_s(b_1,b_2), \quad \mbox{with }F_s(b_1,b_2)=\frac1{s-1}\sum_{i=1}^{s-1}\widetilde X_i^\top b_2+H_s(b_1-b_2),
\end{equation}
which, as shown in \Cref{lem:vertex-optimisation}, can be reduced to an equivalent optimisation over the vertex pairs of confidence sets. 
\begin{lemma}[Vertex optimisation]\label[lemma]{lem:vertex-optimisation}
For every update time $s$, the function $F_s$ in \eqref{eq:Fs} is jointly convex and continuous, and we have that 
\begin{equation}\label{eq:vertex-max}
\max_{b_k\in\mathcal C_{s,k},\ k\in\{1,2\}}F_s(b_1,b_2)=\max_{v_k\in \mathcal V_{s,k},\ k\in\{1,2\}}F_s(v_1,v_2),
\end{equation}
where $\mathcal C_{s,k}$ and $\mathcal V_{s,k}$ are defined in \eqref{eq:candidate-set} and \eqref{eq:confidence-vertices} respectively.
Moreover, if the selected model pair $(\bar\beta_{s,1},\bar\beta_{s,2})$ satisfies $\mathcal H_{\bar\beta_{s,1}-\bar\beta_{s,2}}\neq\emptyset$, then the output of \Cref{alg:vertex-oracle} solves the original optimistic update in \eqref{eq:ideal-update}, i.e.
\begin{equation}\label{eq:vertex-opt}
\frac1{s-1}\sum_{i=1}^{s-1}r_{\bar\beta_s,p_s}(X_i)=\sup_{b_k\in\mathcal C_{s,k},\ q\in\mathcal Q_{b_1-b_2,L}}\frac1{s-1}\sum_{i=1}^{s-1}r_{b,p_{b_1-b_2,q}}(X_i),
\end{equation}
where $r_{b,p}(x)=\widetilde x^\top b_1p(x)+\widetilde x^\top b_2\{1-p(x)\}$. 
\end{lemma}

The proof of \Cref{lem:vertex-optimisation} is in \Cref{sec_app_vertex_optimisation}. It follows from \Cref{lem:initial-localisation} that, a sufficiently long initial exploration phase ensures that all candidate reward-gap parameters remain within the regular parameter neighbourhood $\Theta$ with high probability.  Their decision boundaries are therefore nonempty, so \Cref{lem:vertex-optimisation} guarantees that \Cref{alg:vertex-oracle} provides an exact solution to the optimistic update in~\eqref{eq:ideal-update} through a finite number of linear programs on this event.

The equivalence of the vertex optimisation in \eqref{eq:vertex-max} follows directly from the construction of the polyhedral confidence sets in \eqref{eq:candidate-set}. Under the transformation $u=V_{s,k}^{1/2}(b-\widehat\beta_{s,k})$, the confidence set $\mathcal C_{s,k}$ becomes an $\ell_1$-ball with radius $\sqrt{d+1}\alpha_T$, which is the convex hull of its $2(d+1)$ vertices. Mapping these vertices back through $b=\widehat\beta_{s,k}+V_{s,k}^{-1/2}u$ gives the vertex set $\mathcal V_{s,k}$ in the original parameter space. Hence, every $b\in\mathcal C_{s,k}$ can be written as a convex combination of the vertices in $\mathcal V_{s,k}$ and maximising the jointly convex function $F_s$ over $\mathcal C_{s,1}\times\mathcal C_{s,2}$ is therefore equivalent to maximising it over the finite set of vertex pairs $\mathcal V_{s,1}\times\mathcal V_{s,2}$.

\end{itemize}

\begin{remark}[Computational complexity]
    The vertex formulation ensures that extending the method beyond continuously distributed covariates incurs only an additional computational cost of order $O(d^2)$, compared to \Cref{alg:lazy-adaptive-boundary-greedy}. In particular, each confidence polytope $\mathcal V_{s,k}$ has $2(d+1)$ vertices, so \Cref{alg:vertex-oracle} only needs to consider $4(d+1)^2$ vertex pairs and solves one linear program for each~pair. 
\end{remark}

\begin{theorem}[Regret control]\label[theorem]{thm:main-regret}
    Suppose Assumptions \ref{a_contexts_non_cts}, \ref{a_noise}-\ref{a_theta} hold and $L$ is an absolute constant satisfying $0<L<L_0'$ for some $L'_0 >0$. Implement \Cref{alg:lazy-global-model-ucb} ($d\ge 2$) or \Cref{alg:lazy-one-dimensional-fair-ucb} ($d=1$), with $\lambda = C_\lambda$, where $C_{\lambda} >0$ is an absolute constant, $T_o=\lceil\sqrt T\rceil$, and $\alpha_T$ in \eqref{eq_ucb_alpha}. Then for any $\eta \in (0,1/2)$, we have with probability at least $1-\eta$ that
    \begin{align}
    \mathfrak R_T
    \le\; C_1\alpha_T\sqrt{T\left\{\log\left(1+\frac{T}{\lambda}\right)+\log(8/\eta)\right\}} +C_2\begin{cases}
    \sqrt T\log(T/\eta), & d=1,2,\\
    \sqrt T\log^{3/2}(T/\eta), & d=3,\\
    T^{(d-2)/(d-1)}\log(T/\eta), & d\ge4.
    \end{cases}
    \notag
    \end{align}
    whenever $T\ge16\vee C_3(1+\sqrt{d+1})^4\alpha_T^4/\kappa_X^2$ and $C_1, C_2, C_3 >0$ are constants depending on $d,\nu,L,\sigma,C_\beta$ and the constants in \Cref{a_theta}.     
\end{theorem}

The regret upper bound in \Cref{thm:main-regret} matches the continuous context distribution rate in \Cref{thm:fair_regret_control}, up to poly-logarithmic factors, 
despite the absence of guaranteed covariate diversity for the greedy~approach.

The detailed proof of \Cref{thm:main-regret} is in \Cref{sec_app_ucb_regret}, with a sketch below.
The proof follows a similar strategy to that of \Cref{thm:fair_regret_control}, particularly in controlling the uniform deviation between empirical and population rewards through the complexity of the boundary-canonical policy class. The main difference lies in the treatment of parameter estimation uncertainty under the UCB framework. Specifically, let $\widehat R_s(b,p)=(s-1)^{-1}\sum_{i=1}^{s-1}r_{b,p}(X_i)$ and $R_b(p)=\mathbb E[r_{b,p}(X)]$, with $R(p)=R_{\beta^*}(p)$. Since the true reward parameters lie within the confidence sets with high probability, the selected model pair $(\bar\beta_s,p_s)$ in \eqref{eq:vertex-opt} satisfies $\widehat R_s(\bar\beta_s,p_s)\ge\widehat R_s(\beta^*,p^*)$ (\Cref{lem:empirical-optimism}). Consequently, the instantaneous regret can be decomposed and bounded by
\begin{align*}
R(p^*)-R(p_s)
\le{}& \{R(p^*)-\widehat R_s(\beta^*,p^*)\} +\{\widehat R_s(\bar\beta_s,p_s)-R_{\bar\beta_s}(p_s)\} +\{R_{\bar\beta_s}(p_s)-R(p_s)\},
\end{align*}
where the first two terms can be controlled using the uniform empirical-process bounds in \Cref{lem:ucb-control-delta}, and the last term captures the uncertainty in the estimated reward parameters. Unlike the proof of \Cref{thm:fair_regret_control}, which requires controlling the boundary mismatch error, the UCB analysis bounds the remaining parameter estimation error through arm-specific confidence widths. Writing $s=\tau(t)$ for the most recent update time before $t$ and $w_{s,k}(x)=\min\{1,\|\widetilde x\|_{V_{s,k}^{-1}}\}$, we obtain from the proof of \Cref{lem:one-step-regret}~that
\begin{align*}
R_{\bar\beta_s}(p_s)-R(p_s)
&\lesssim_{\log} \mathbb E\left[p_s(X_t)w_{s,1}(X_t)+\{1-p_s(X_t)\}w_{s,2}(X_t)\mid\mathcal F_{t-1}\right] = \mathbb E\left[w_{s,A_t}(X_t)\mid\mathcal F_{t-1}\right],
\end{align*}
where the  equality uses the fact that, conditional on $(\mathcal F_{t-1}, X_t)$, arm $1$ is pulled with probability $p_s(X_t)$. This  allows us to express the uncertainty under the randomised fair policy in terms of the confidence width of the arm actually selected, whose cumulative contribution can be bounded  using \Cref{lem:lazy-width-sum-lnu}.

\section{Optimality and minimax lower bounds} \label{sec_lower}
To establish the rate optimality of the algorithms developed in \Cref{sec_learning_incomplete_information}, and to quantify the inevitable learning regret under the fairness constraint, in this section, we derive minimax lower bounds.

While the upper bounds in \Cref{sec_learning_incomplete_information} hold for general $\ell_{\nu}$-Lipschitz fairness,  the minimax lower bound depends on the metric, $\nu$ in particular, through the geometry of the admissible policy classes. Here we focus on the $\ell_1$ formulation, which fixes a concrete metric geometry of the policy class and serves as a representative case to understand the unavoidable boundary-learning difficulty induced by Lipschitz fairness.

Let $\mathcal G_{\rm cts}$ denote the collection of data generating distributions satisfying Assumptions \ref{a_contexts}, \ref{a_noise}, \ref{a_beta} and~\ref{a_theta}, and $\mathcal G_{\rm gen}$ the collection satisfying Assumptions \ref{a_contexts_non_cts}, \ref{a_noise}, \ref{a_beta} and \ref{a_theta}. The continuous-context class is a subclass of the general-context class, i.e.~$\mathcal G_{\rm cts}\subseteq\mathcal G_{\rm gen}$, since under \Cref{a_contexts}, the context density is uniformly bounded away from zero on $[0,1]^d$, which implies that $\lambda_{\min}\{\mathbb E(\widetilde X\widetilde X^\top)\}$ is bounded away from zero and hence \Cref{a_contexts_non_cts} is satisfied. It therefore suffices to only establish the lower bound over $\mathcal G_{\rm cts}$.

Denote $\mathcal F_{t-1}=\sigma\{(X_s,A_s,Y_s):1\le s\le t-1\}$ as the observed history before round $t$. We define $\Pi_{L,1}$ as the class of non-anticipating $(L,1)$-Lipschitz fair learning algorithms $\pi=\{p_t\}_{t\ge1}$, where at each round $t$, it selects an $\mathcal F_{t-1}$-measurable random function $p_t:[0,1]^d\to[0,1]$ such that $p_t\in\mathcal P_{L,1}$ almost surely, i.e.
\[
    \Pi_{L,1}=\Bigl\{\{p_t\}_{t=1}^\infty: p_t \text{ is } \mathcal F_{t-1}\text{-measurable and } p_t\in \mathcal{P}_{L,1} \text{ a.s. for every } t\Bigr\}.
\]

For a policy sequence $\pi=\{p_t\}_{t\ge1}$, define its   expected regret under $G$ by
\[
    \mathfrak R_T(\pi;G)=\mathbb E_G^\pi\left[\sum_{t=1}^T\{R_G(p_G^*)-R_G(p_t)\}\right],
\]
where $R_G(p)=\mathbb E_G[p(X)\widetilde X^\top\beta_1^*+\{1-p(X)\}\widetilde X^\top\beta_2^*]$ and  $p_G^* \in \argmax_{p\in\mathcal P_{L,1}} R_G(p)$.

\begin{theorem}[Minimax lower bound]\label[theorem]{thm:lower_bound}
 For any fixed  $L \in(0,1)$, 
 there is a constant $c_{d,L}>0$ depending only on $d$ and $L$, such that
\[
    \inf_{\pi\in\Pi_{L,1}} \sup_{G\in\mathcal G_{\rm cts}} \mathfrak R_T(\pi;G)\ge c_{d,L}
    \begin{cases}
    \sqrt T, & d\in\{1,2,3\},\\[2pt]
    T^{\frac{d-2}{d-1}}, & d\ge4.
    \end{cases}
\]
\end{theorem}

The upper bound in \Cref{thm:fair_regret_control} and the lower bound in \Cref{thm:lower_bound} together show that Algorithm \ref{alg:lazy-adaptive-boundary-greedy} for $d\geq 2$ and Algorithm \ref{alg:lazy-one-dimensional-fair-greedy} for $d=1$, given in the appendix, achieve optimal regret up to poly-logarithmic factors. In particular, the extra fairness-induced regret of the order $T^{(d-2)/(d-1)}$ is tight and inevitable when $d \geq 4$. 

The detailed proof of \Cref{thm:lower_bound} is deferred to \Cref{sec_app_lower_bound}. 
The key idea is to construct covariate distributions that are statistically difficult to distinguish but induce different optimal fair policies. In the unconstrained problem, once the reward gap $\Delta^*(x)$ is fixed, the oracle policy is determined pointwise by the sign of $\Delta^*(x)$ and is therefore independent of the covariate density. Under Lipschitz fairness, however, the policy must transition smoothly across the decision boundary, and the optimal shape of this transition depends on how probability mass is distributed near the boundary. We exploit this feature by keeping the reward gap fixed, specifically $\Delta^*(x)=x_1-1/2,$ and perturbing only the covariate density locally around the boundary $\mathcal{H}_{\theta^*}=\{x\in[0,1]^d:x_1=1/2\}$. These local density perturbations shift probability mass across the boundary and hence alter the optimal fairness correction, while keeping the  KL divergence between environments  small.

For $d\in\{1,2\}$, a two-point density perturbation is sufficient to capture the parametric difficulty at rate $\sqrt{T}$, leading to a Le Cam argument \citep[e.g.][]{yu1997assouad}. The two environments favour different local shifts of the boundary value,  forcing any fair policy to incur $\Omega(\sqrt{T})$ regret under at least one of them.  For $d\geq3$, the $(d-1)$-dimensional decision boundary accommodates a packing of $m^{d-1}$ spatially  separated local perturbation  regions along the tangential coordinates $(x_2,\ldots,x_d)$. By placing independent signed density perturbations across these regions and applying Assouad's lemma \citep[e.g.][]{yu1997assouad}, any fair policy must simultaneously estimate all $m^{d-1}$ local boundaries. Hence,  by balancing local statistical detectability against the aggregate reward loss induced by the Lipschitz constraint, we choose $m \asymp T^{1/(d-1)}$ to obtain the sharp nonparametric lower bound of order $\Omega(T^{(d-2)/(d-1)})$, which matches $\sqrt{T}$ at $d=3$ and strictly dominates $\sqrt{T}$ for $d\ge 4$.

\section{Numerical experiments}\label{sec_numerical}

In \Cref{sec_learning_incomplete_information}, we proposed \Cref{alg:lazy-adaptive-boundary-greedy} for continuous contexts and \Cref{alg:lazy-global-model-ucb} for general contexts, with their one-dimensional counterparts in \Cref{alg:lazy-one-dimensional-fair-greedy,alg:lazy-one-dimensional-fair-ucb}, respectively.In this section, we  evaluate the empirical performance of our framework, focusing on the lazy fair  greedy algorithms (\Cref{alg:lazy-adaptive-boundary-greedy,alg:lazy-one-dimensional-fair-greedy})  owing to  their computational efficiency; an additional  mixed-context experiment in \Cref{sec:app-mixed-design} further illustrates their empirical performance beyond  purely  continuous contexts.  Simulation and real-data analysis are collected in \Cref{sec:numerical-simulation,sec:numerical-real-data}, respectively.   The code and datasets used in the experiments are available online at \url{https://github.com/chrisfanwang/fair-contextual-bandits}.

\subsection{Simulation studies}\label{sec:numerical-simulation}

We consider the two-arm linear reward model in \Cref{sec_setup}, with $d\in\{1,2,3,4\}$ and independent noise random variables $\xi_{t,k}\sim\mathcal N(0,0.3^2)$. Throughout, fairness is measured using the $\ell_1$-distance ($\nu=1$ in \Cref{def_lipschitz_fairness}). Contexts are sampled independently across $t$ from one of the following two distributions:
\begin{itemize}[noitemsep,nolistsep, leftmargin=*]
    \item \textbf{Uniform design} with $X_t$ uniformly drawn from $[0,1]^d$.
    \item \textbf{Truncated Gaussian design} where each coordinate of $X_t$ is independently drawn from $\mathcal N(0.5,0.18^2)$ and then truncated to $[0,1]$.
\end{itemize}
For each context distribution, we consider two reward regimes:
\begin{itemize}[noitemsep,nolistsep, leftmargin=*]
    \item \textbf{Homogeneous regime} where
    \[
        \beta_1^*  =\bigg(-\frac{\sqrt d}{4}, \frac{1}{2\sqrt d}\mathbf{1}_d^\top\bigg)^\top  \quad \mbox{and} \quad  \beta_2^*=-\beta_1^*.
    \]
    The resulting reward gap is $\Delta^*(x) = d^{-1/2}\sum_{j=1}^d(x_j-1/2)$,   so that all covariates contribute equally to the gap.
    \item \textbf{Heterogeneous regime} with the coefficients in \Cref{tab:numerical-heterogeneous}, varying in magnitude and sign across covariates and arms.
\end{itemize}
These choices yield four baseline settings in which each arm is optimal with probability $1/2$.  

To examine learning when one arm is optimal over a smaller region, we consider two additional $d=2$ settings under the uniform design:
\begin{itemize}[noitemsep,nolistsep, leftmargin=*]
    \item \textbf{Arm 1 small.} We set $\beta_1^*=(-0.41,0.135,0.5)^\top$ and $\beta_2^*=-\beta_1^*$, giving $\Delta^*(x)=0.27x_1+x_2-0.82$.  Arm 1 is optimal with probability $0.315$.
    \item \textbf{Arm 2 small.} We set $\beta_1^*=(-0.1,0.7,0.4)^\top$ and $\beta_2^*=(0.1,-0.3,-0.6)^\top$, giving $\Delta^*(x)=x_1+x_2-0.2$. Arm 2 is optimal with probability $0.02$.
\end{itemize}

\begin{table}[t]
\centering
\caption{Coefficients in the heterogeneous reward regime.  The first component denotes the intercept.}
\label{tab:numerical-heterogeneous}
\begin{tabular}{clll}
\toprule
$d$ & $\beta_1^*$ & $\beta_2^*$ & $\Delta^*(x)$\\
\midrule
1 & $(-0.3,0.8)^\top$ & $(0.3,-0.4)^\top$ & $1.2x_1-0.6$\\
2 & $(-0.5,0.7,0.4)^\top$ & $(0.5,-0.3,-0.6)^\top$ & $x_1+x_2-1$\\
3 & $(0,0.6,-0.7,0.15)^\top$ & $(0,-0.4,0.3,0.15)^\top$ & $x_1-x_2$\\
4 & $(0,0.6,0.35,-0.3,-0.6)^\top$ & $(0,-0.4,-0.65,0.7,0.4)^\top$ & $x_1+x_2-x_3-x_4$\\
\bottomrule
\end{tabular}
\end{table}

\medskip
\noindent
\textbf{Learning, oracle policies and evaluation.} We implement \Cref{alg:lazy-adaptive-boundary-greedy} and its one-dimensional counterpart in \Cref{alg:lazy-one-dimensional-fair-greedy}. We consider horizons $T\in\{1024, 2048, 4096, 8192, 16384, 32768\}$, with a separate experiment for each horizon because the exploration length depends on $T$. The initial exploration phase lasts $2T_o(T)$ rounds, where $T_o(T)=\lceil\sqrt{T\log(4T)}\rceil$, following the choice in \Cref{thm:fair_regret_control} with $C_1=1$. We use ridge regularisation with $\lambda=1$ throughout.

The baseline reward gaps can be written as $\Delta^*(x)=c\sum_{j=1}^d s_j(x_j-1/2)$, where $c>0$ and $s_j\in\{-1,0,1\}$. Under both the uniform and truncated Gaussian designs, an optimal fair policy is
\begin{equation}\label{eq-num-oracle}
 p_{L,1}^*(x)=\bigg[\frac12+L\sum_{j=1}^d s_j(x_j-1/2)\bigg]_{[0,1]},
\end{equation}
where $[x]_{[a, b]} = \max\{a, \, \min\{x, \, b\}\}$.  The Arm 1 small and Arm 2 small settings also admit explicit oracle policies; see \Cref{sec:app-simulation}.

We evaluate the population reward $R(p)$ in \eqref{eq_reward} by integrating over the covariate distribution. Because analytical integration becomes computationally demanding as the dimension or policy complexity grows, we use exact analytical integration for $d \in \{1,2\}$ (as well as for the oracle cost-of-fairness curves under the uniform design across all $d\in\{1,2,3,4\}$) and quasi-Monte Carlo integration for $d\in\{3,4\}$. For regret evaluation under either design, we apply the same integration rule to the oracle and learned policies. All quasi-Monte Carlo evaluations use $16384$ Halton points for $d=3$ and $32768$ for $d=4$, after discarding the first $128$ burin-in points. Under the truncated-Gaussian design, we transform these points coordinate-wise using the inverse cumulative distribution function of the truncated normal distribution. We compute cumulative regret $\mathfrak R_T$, defined in \eqref{eq_accumulated_regret}, by summing the population reward differences between the oracle and learned policies over all rounds, including initial exploration.

\begin{figure}[ht]
\centering
\includegraphics[width=.8\linewidth]{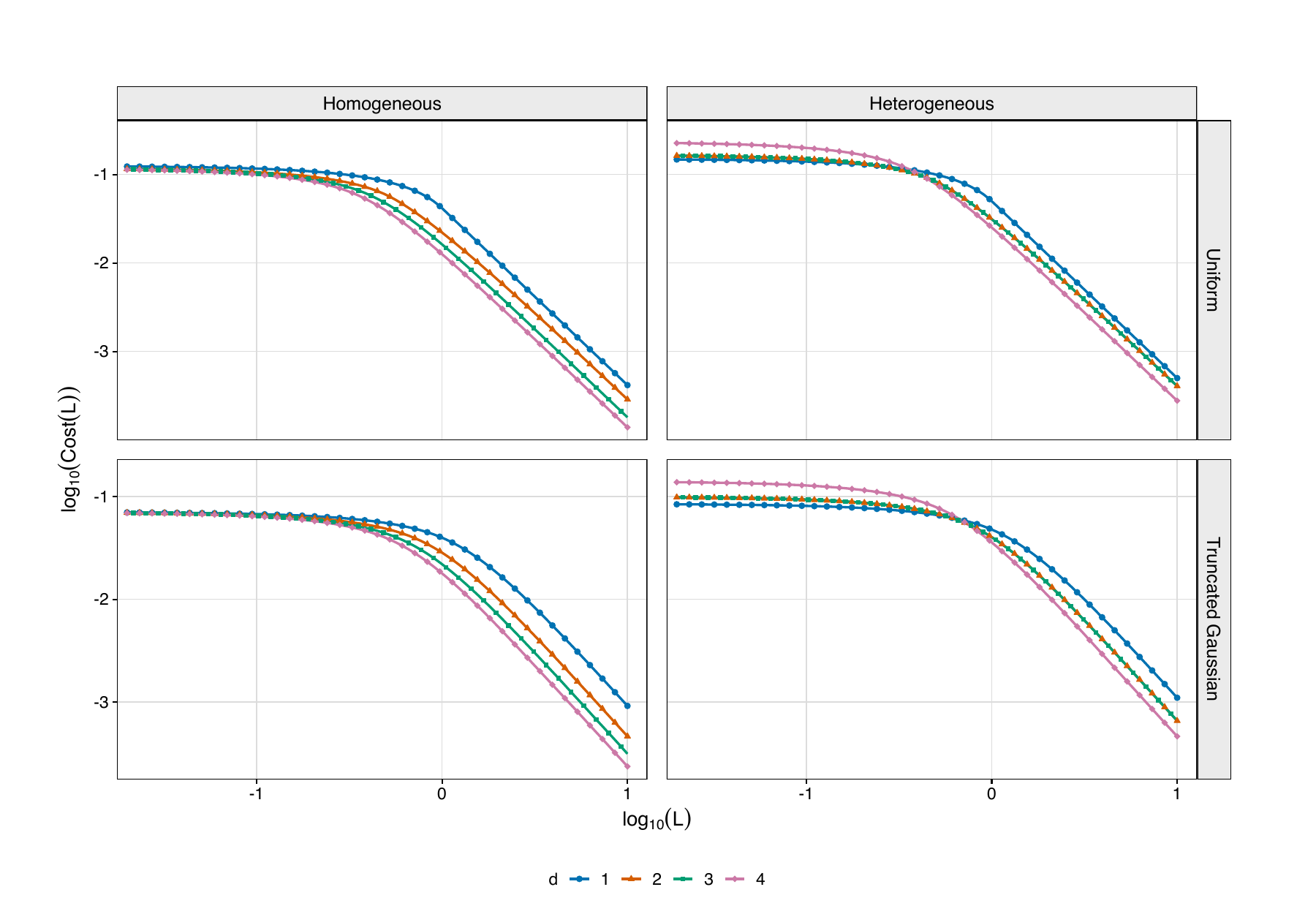}
\caption{Cost of fairness, $\mathrm{Cost}(L)$, on log--log axes. Columns show homogeneous (left) and heterogeneous (right) rewards; rows show uniform (top) and truncated-Gaussian (bottom) covariates.}
\label{fig:cost-fairness-loglog}
\end{figure}

\begin{figure}[ht]
\centering
\includegraphics[width=.8\linewidth]{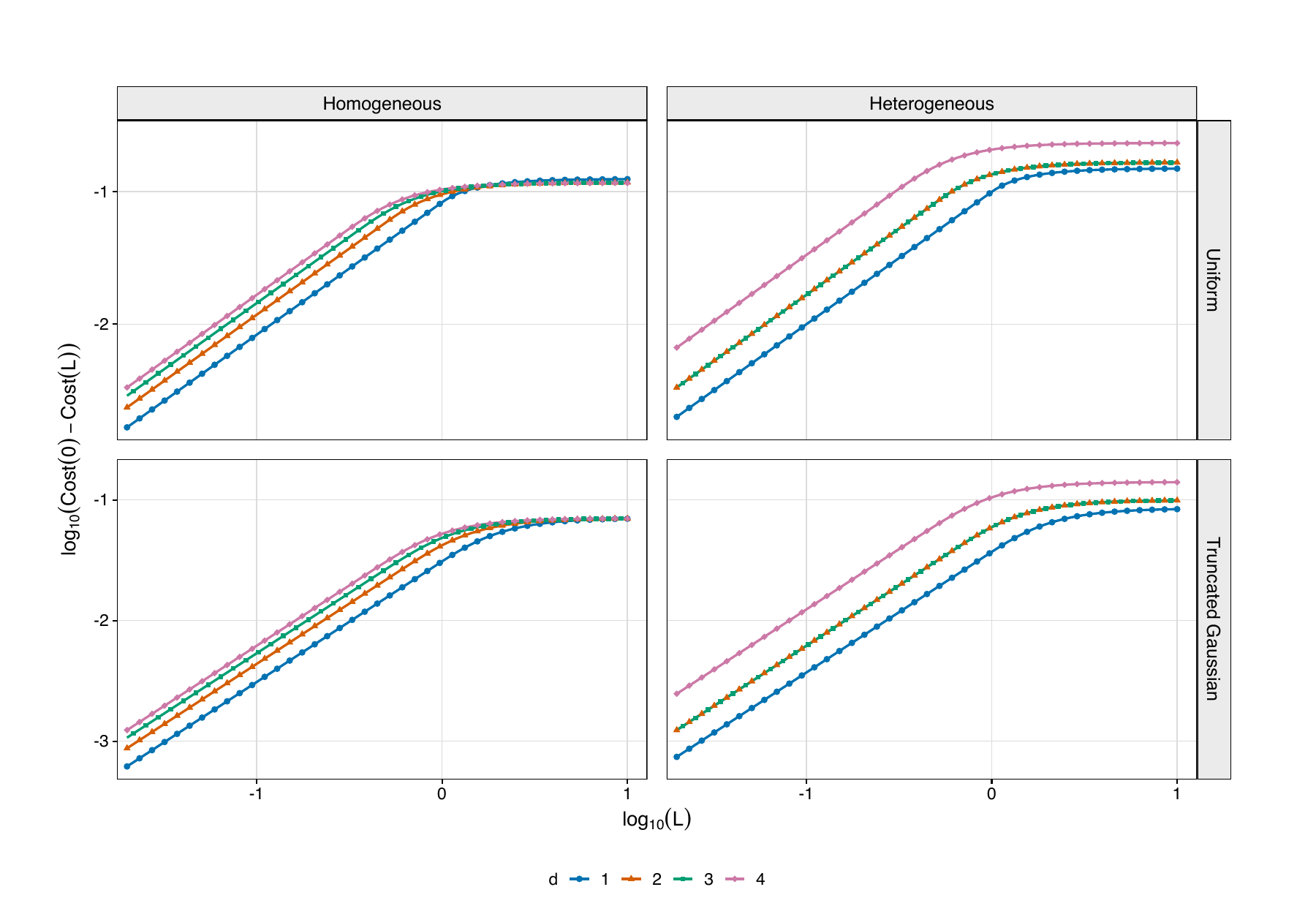}
\caption{Improvement over the best constant policy, $\mathrm{Cost}(0)-\mathrm{Cost}(L)$, on log--log axes. Columns show homogeneous (left) and heterogeneous (right) rewards; rows show uniform (top) and truncated-Gaussian (bottom) covariates.}
\label{fig:cost-improvement-loglog}
\end{figure}

\medskip
\noindent
\textbf{Cost of fairness.}
We first examine $\mathrm{Cost}(L)$, defined in \Cref{sec_cost} as the reward difference between the unconstrained and fair oracles.  \Cref{fig:cost-fairness-loglog} shows the results for the four baseline settings with $d \in \{1,2,3,4\}$, using $81$ logarithmically spaced values of $L$ from $0.02$ to $10$. The curves show a similar pattern across distributions, reward regimes and dimensions. For small $L$ ($L\leq0.1$), the cost remains close to $\mathrm{Cost}(0)$, the cost associated with the best constant policy. As $L$ increases, the fairness constraint becomes less restrictive and the cost decreases. Over $3\leq L\leq10$, the curves are approximately linear on log-log axes, with slopes close to $-2$, consistent with the $L^{-2}$ rate in \Cref{thm:fairness_cost}.

For the same settings, \Cref{fig:cost-improvement-loglog} reports $\mathrm{Cost}(0)-\mathrm{Cost}(L)=R_L^*-R_0^*$, the reward gain achieved by the optimal policy under the Lipschitz fairness constraint relative to the best constant policy. For small $L$, the gain increases approximately linearly with $L$ (slope near $1$ on log-log axes), matching the $O(L)$ upper bound in \Cref{thm:fairness_cost}. As $L$ increases further, the gain approaches $\mathrm{Cost}(0)$ and the curves flatten. In the special case of $d=1$ with uniform covariates and homogeneous rewards, the results in both figures also agree with the closed-form expressions derived  in \Cref{ex:cost_fairness_1d}.

\begin{figure}[ht]
\centering
\includegraphics[width=.8\linewidth]{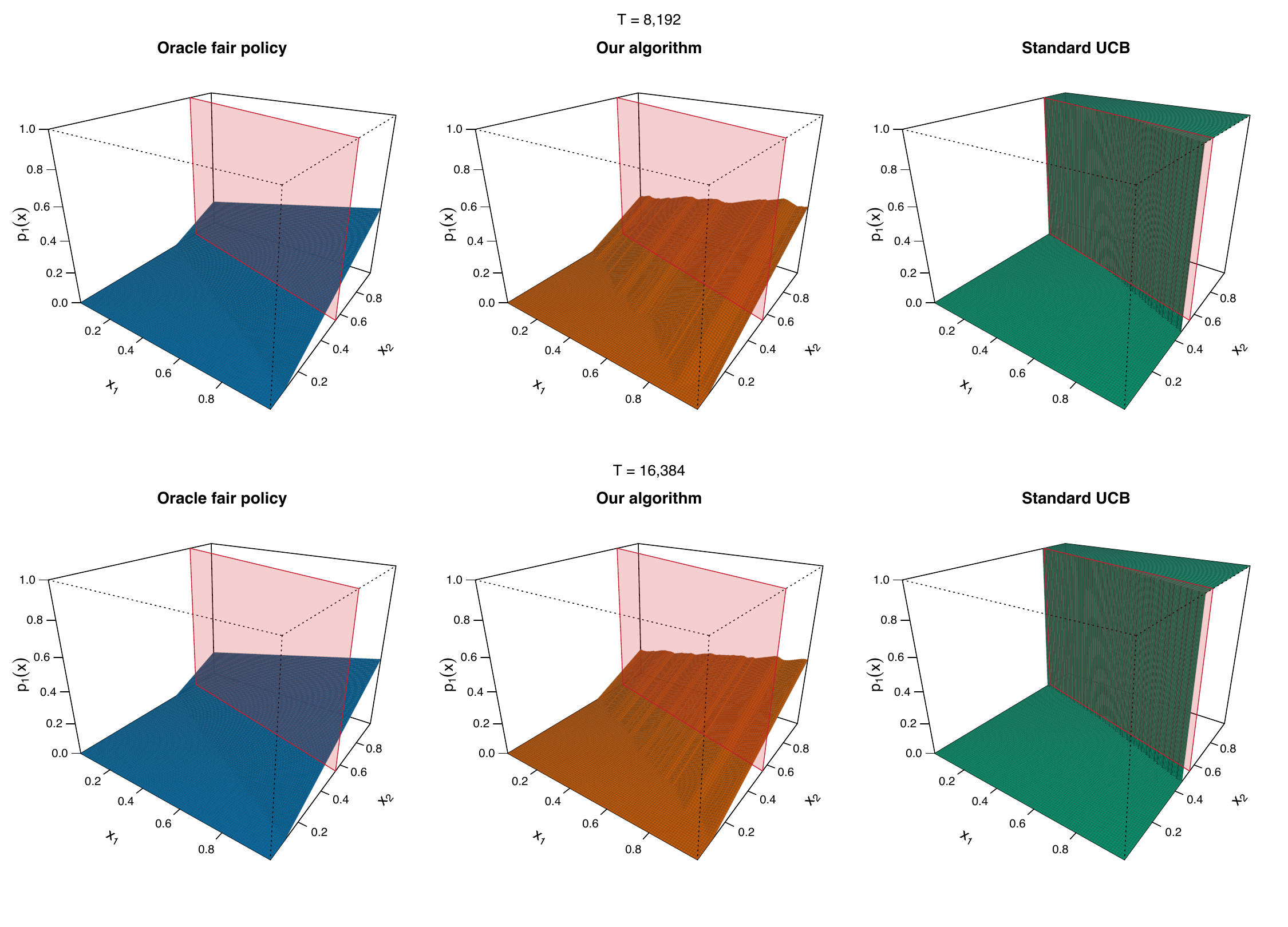}
\caption{Probability $p(x)$ of choosing arm 1 in the Arm 1 small setting with uniform covariates, $d=2$ and $L=0.5$. Columns show the fair oracle (left), our learned policy (middle) and unconstrained Standard UCB (right); rows show separate experiments at $T=8192$ (top) and $T=16384$ (bottom). The red plane marks the true decision boundary.}
\label{fig:num-policy}
\end{figure}

\medskip
\noindent
\textbf{Policy visualisation.}
We next compare the structure of three policies: the fair oracle $p^*_{L,1}$, our learned fair policy $\widehat{p}_T$ and unconstrained standard UCB. \Cref{fig:num-policy} shows action-selection probabilities in the Arm 1 small setting with uniform covariates, $d=2$ and $L=0.5$. For standard UCB, the probability of selecting arm~1 jumps between 0 and 1 across the estimated decision boundary. The fair oracle and our learned policy instead exhibit a smooth transition, illustrating the boundary randomisation and geometry structure described in \Cref{lem:mixing} and \Cref{thm:canonical_extension}. We also observe that the learned policy is closer to the oracle at $T=16384$ than at $T=8192$. The population reward gap relative to the oracle, $R(p_{L,1}^*)-R(\widehat p_T)$, decreases from $5.54\times10^{-4}$ at $T=8192$ to $2.94\times10^{-5}$ at $T=16384$, a reduction of approximately $94.7\%$, indicating improved performance at the longer horizon.

\begin{figure}[ht]
\centering
\includegraphics[width=.8\linewidth]{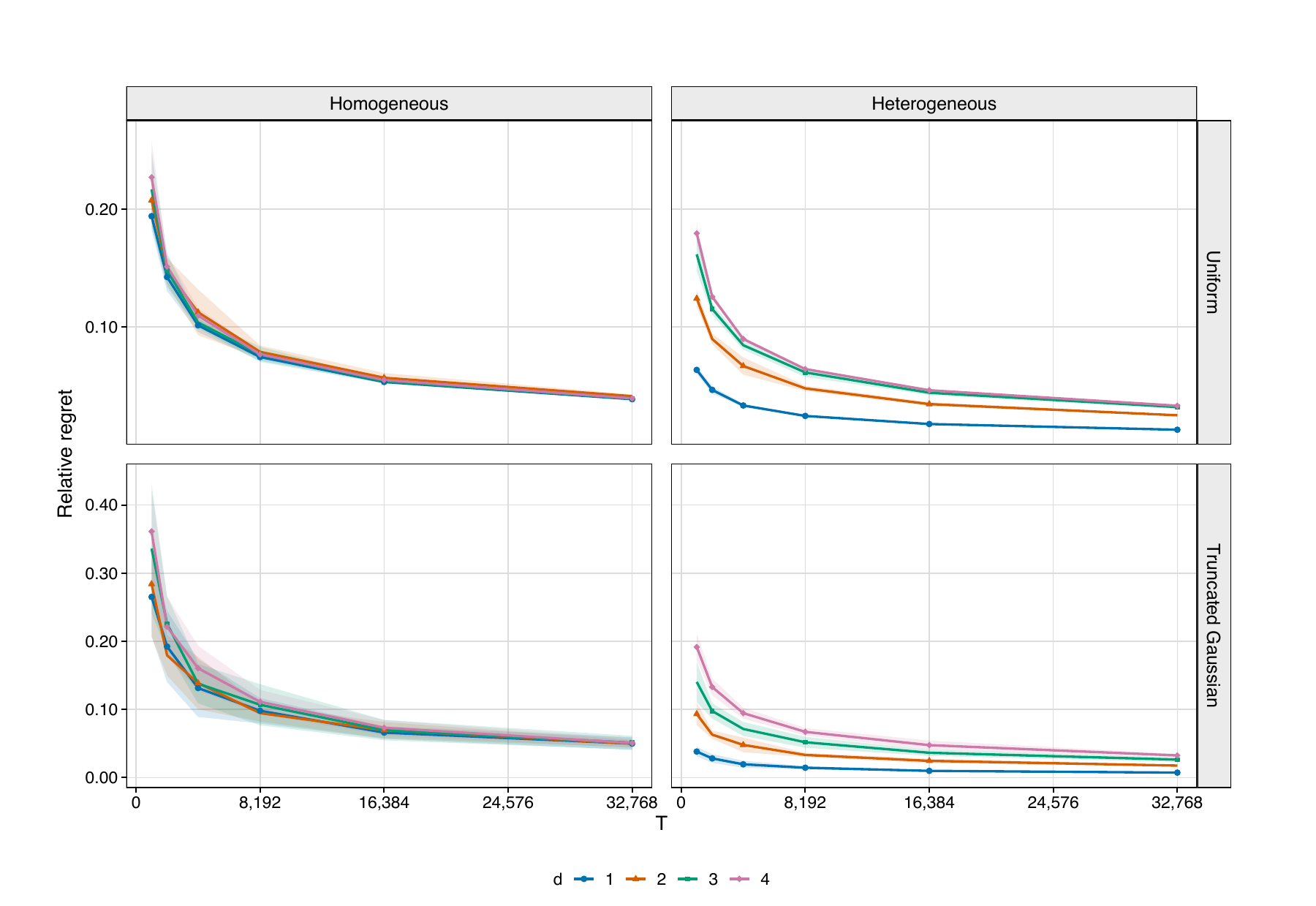}
\caption{Relative regret for $d=1,2,3,4$ at $L=0.5$. Columns show homogeneous (left) and heterogeneous (right) rewards; rows show uniform (top) and truncated-Gaussian (bottom) covariates. Lines show means over $50$ independent Monte Carlo trials; shaded bands show the mean $\pm$ one standard deviation.}
\label{fig:relative-regret-by-d}
\end{figure}

\medskip
\noindent
\textbf{Regret.} To evaluate learning performance, we first examine the relative regret, defined as
\[
\mathrm{RelReg}_T = \frac{\mathfrak R_T}{T R(p_{L,1}^*)},
\]
which is our algorithm's cumulative expected reward loss relative to the fair oracle, normalised by the oracle's cumulative expected reward. \Cref{fig:relative-regret-by-d} shows the results at $L=0.5$, based on $50$ Monte Carlo trials. Across  all four baseline settings and dimensions $d\in\{1,2,3,4\}$, the relative regret decreases steadily with $T$, consistent with the sublinear regret bounds in \Cref{thm:fair_regret_control}. Additional experiments in Appendix \ref{sec:app-simulation} further examines the sensitivity of regret to the fairness parameter $L\in\{0.25,0.5,1,2,4\}$ under uniform covariates and homogeneous rewards with $d=2$.

Appendix \ref{sec:app-simulation} further plots the unnormalised cumulative regret $\mathfrak R_T$ against the time horizon $T$ on log-log axes for all four baseline settings and $d\in\{1,2,3,4\}$. The log-log curves are approximately linear, with empirical slopes close to $1/2$ even for $d=4$, where the worst-case theoretical upper bound scales as $T^{2/3}$ up to logarithmic factors. We conjecture that the discrepancy in dimension dependence arises because the symmetric baseline designs yield the clipped affine policy in \eqref{eq-num-oracle}, whose boundary decision   $p^{*}_{L,1}|\mathcal{H}_{\theta^*}$ equals to the constant $1/2$ uniformly, whereas  the minimax rate $T^{(d-2)/(d-1)}$ reflects the difficulty of learning a spatially
varying $d-1$-dimensional boundary. To probe this nonparametric regime,  we construct a mixed discrete-continuous design with $d=4$ in Appendix \ref{sec:app-mixed-design}. It consists of a categorical covariate with six levels and another three continuous covariates, with conditional category probabilities differing across regions of the continuous covariate space. The fair oracle's arm selection probability is thus nonconstant along the decision boundary. In this setting, we find the cumulative regret curve is approximately linear on log-log axes, with an empirical slope of about $0.65$ with respect to $T$, numerically matching the theoretical exponent $2/3$ in \Cref{thm:fair_regret_control} and \Cref{thm:lower_bound}.

\subsection{Real data analysis} \label{sec:numerical-real-data}

In this section, we evaluate the proposed methods on two binary classification datasets: Adult \citep{becker1996adult} and Jobs \citep{tankha2023jobs}. After data preprocessing, the Adult dataset contains $43303$ census records with $12$ demographic and socioeconomic covariates (such as age, education, occupation and working hours), where the task is to predict whether annual income exceeds \$50000. The Jobs dataset contains $73462$ records with $11$ personal and professional covariates (such as age, education, programming experience, previous salary and computer skills) with the binary employment status as the response. Further details on sample selection, covariate definitions and preprocessing details are provided in Appendix \ref{sec:app-real-data}.

For each dataset,  we hold out a fixed yet randomly sampled  $10\%$ subset to fit a covariate preprocessing transformation based on  principal component analysis (PCA) \citep{jolliffe2016principal}, which retains the first three principal components and rescales and clips their scores to $[0,1]^3$. The fitted   fitted transformation is then applied to all remaining 90\% of records, which serve as the online sampling dataset, and 
Lipschitz fairness is enforced with respect to the $\ell_1$-distance on the transformed covariates $X_t\in [0,1]^3$.

Following \citet{ma2026locally}, we formulate each classification task as a two-arm bandit problem.  Selecting arm~1 assigns the predicted label $1$ to the current record, while selecting arm~2 assigns the predicted label $0$. For a record with true label $B_t\in\{0,1\}$, we define the arm reward as $R_{t,1}=B_t$ and $R_{t,2}=1-B_t$. At round $t$, the algorithm observes the transformed covariates $X_t$, selects arm $A_t$ and receives $Y_t=R_{t,A_t}$, which equals one for a correct prediction and zero otherwise. Since the population reward functions are unknown, regret relative to the fair oracle cannot be evaluated directly. We therefore report cumulative online accuracy,
\[
\mathrm{Acc}_T=\frac{1}{T}\sum_{t=1}^{T}Y_t,
\]
the proportion of correct predictions over all $T$ rounds, including initial exploration.

\begin{figure}[H]
\centering
\includegraphics[width=.8\linewidth]{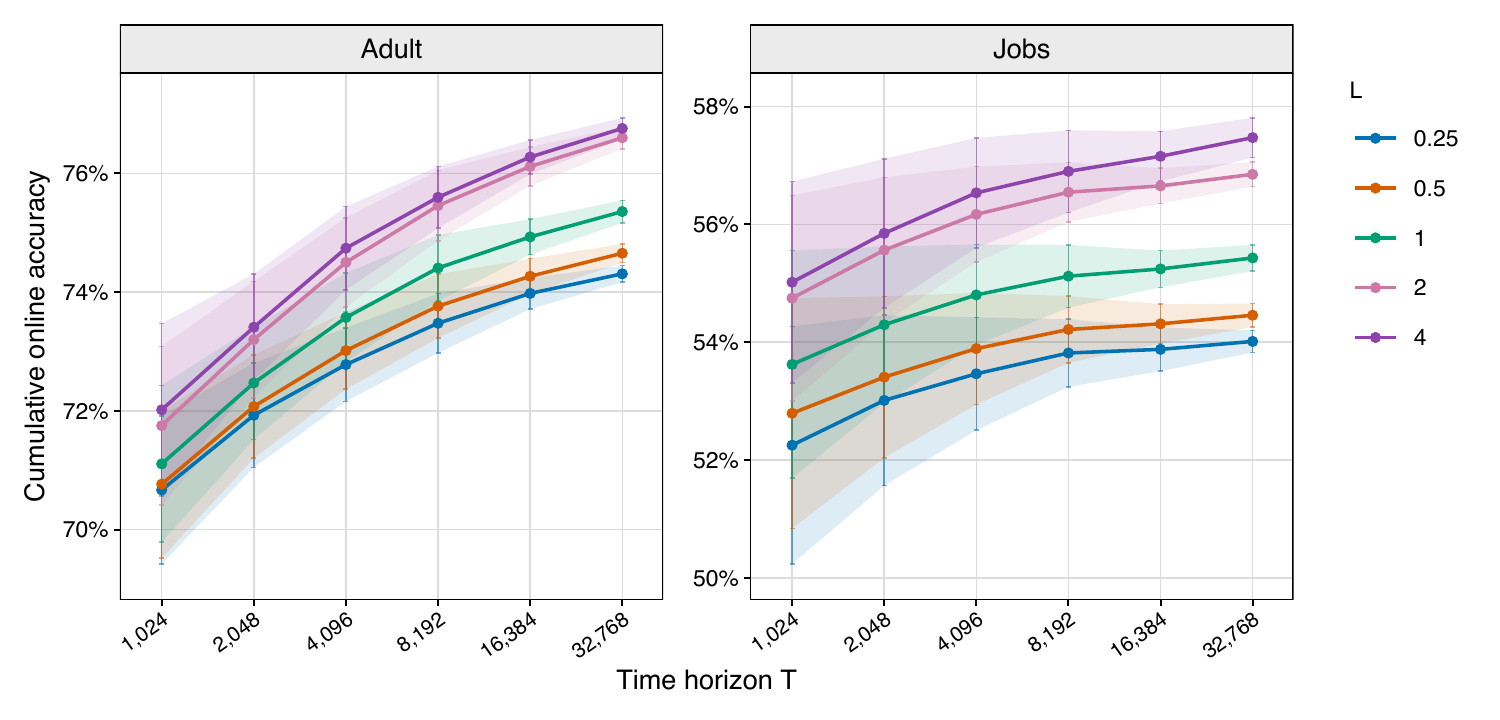}
\caption{Cumulative online accuracy against $T$ for $L\in\{0.25,0.5,1,2,4\}$. Panels show Adult (left) and Jobs (right). Lines show means over $50$ independent Monte Carlo trials; shaded bands show the mean $\pm$ one standard deviation.}
\label{fig:real-accuracy}
\end{figure}

We consider $T\in\{1024,2048,4096,8192,16384,32768\}$ and $L\in\{0.25,0.5,1,2,4\}$ with $50$ independent replications for each dataset and $(T,L)$ pair.   Within each replication, we randomly permute the online sampling pool and feed its first $T$ records to a separate learner for each $(T,L)$ pair sharing the same permutation across all values of $T$ and $L$, using the same exploration schedule $T_o(T)$ as in \Cref{sec:numerical-simulation}. \Cref{fig:real-accuracy} reports mean cumulative online accuracy against $T$. For both datasets, mean accuracy increases with $T$, indicating improved predictive performance at longer horizons. Accuracy is also consistently higher for larger $L$.  At $T=32768$, increasing $L$ from $0.25$ to $4$ raises mean accuracy from $74.31\%$ to $76.76\%$ for Adult and from $54.01\%$ to $57.47\%$ for Jobs. This reflects the theoretical trade-off between the strength of the individual fairness constraint and the rewards.

\section{Conclusions} \label{sec_conclusion}
In this paper, we study optimal learning in contextual bandits under a covariate-based Lipschitz fairness constraint.  Building upon the canonical boundary representations, we develop computationally efficient lazy greedy and UCB-type learning algorithms and establish matching upper and lower minimax regret bounds up to poly-logarithmic factors. Together, our results characterise the fundamental  cost of individual fairness at both the population and sequential learning levels, revealing a dimension phase transition at $d = 4$, beyond which nonparametric boundary estimation becomes more challenging than parametric reward learning.

We envisage several potential extensions for future research. First, one may generalise the current framework to multi-armed contextual bandits with $K \geq 3$. Our results in both Sections~\ref{sec_boundary} and \ref{sec_structural} reveal that the decision boundary is central to both the structure and learning difficulty of fair policies. In the setting when $K \geq 3$, however, the geometry may involve up to the order of $O(K^2)$ pairwise decision boundaries, which substantially complicates both the canonical representation and theoretical analysis. A potential remedy is to assume lower-dimensional structures of the fairness constraint, for example imposing Lipschitz continuity on expected rewards or utilities \citep[e.g.][]{chen2025utility}. Second, it would be interesting to move beyond fairness within each period and study temporal fairness constraints across time \citep[e.g.][]{gupta2021individual, cohen2025dynamic}. Under these settings, decisions are required to remain fair not only across covariates at a given time but also across nearby rounds over time. Deriving practical sequential learning algorithms and theoretical guarantees on the inevitable cost of fairness under these scenarios remains an intriguing area for further~investigation.

\section*{Acknowledgements}
Xue is supported by the EPSRC programme grant DASS (EP/Z531327/1). Jiang is supported in part  by the NSFC (Nos.~72522009, 12331009 and 72432002), and by the National Key R\&D Program of China (No.~2024YFA1015700). Yu is partially supported by the EPSRC programme grant DASS (EP/Z531327/1), UKRI-JST grant Statistical Safeguarding and Leverhulme Trust (Philip Leverhulme Prize).

\bibliographystyle{apalike}
\bibliography{reference}

\newpage
\appendix 
\section*{Appendices}
Additional algorithms and all technical details are collected in the Appendices. Proofs of the results in Sections~\ref{sec_mix_cost}, \ref{sec_learning_incomplete_information} and \ref{sec_lower} are provided in Appendices~\ref{sec_app_mix_cost}, \ref{sec_app_learning_incomplete_information} and \ref{sec_app_lower_results}, respectively. Appendix~\ref{sec_app_one_dimensional} collects the one-dimensional results, with a shared setup followed by separate presentations of the greedy algorithm for continuous covariates (analogue of \Cref{alg:lazy-adaptive-boundary-greedy}) and the UCB algorithm for general covariates (analogue of \Cref{alg:lazy-global-model-ucb}), their theoretical guarantees, and their proofs. Appendix~\ref{sec_app_structural_general} contains additional structural results for arbitrary context laws, including discrete and mixed distributions. Appendix~\ref{sec:app-boundary-computation} describes how to evaluate the boundary-canonical policy at arbitrary covariates under the $\ell_1$-metric. Appendix~\ref{sec:app-numerical} provides additional simulation results and preprocessing details for the real-data analysis.

Throughout the appendix, with a slight abuse of notation, unless specifically stated otherwise, let $C_1, C_2, \ldots >0$ denote constants whose values may vary from place to place and depend on $d$.

\section{Proofs of results in Section \ref{sec_mix_cost}} \label[appendix]{sec_app_mix_cost}

\subsection{Proof of Lemma \ref{lem:mixing}} \label[appendix]{sec_app_thm_boundary}
\begin{proof}[Proof of \Cref{lem:mixing}]
    Note that under \Cref{a_theta}, the decision boundary $\mathcal{H}_{\theta^*}$ is non-empty. Hence, \Cref{lem:mixing} is a direct consequence of \Cref{lem:mixing_confinement} and \Cref{cor:slab_width}.
\end{proof}

\begin{theorem}\label[theorem]{lem:mixing_confinement}
Suppose \Cref{a_contexts} holds and let $p^*\in\argmax_{p\in\mathcal P_{L,\nu}}R(p)$. Define the randomisation region~by 
\[\mathcal M(p^*)=\{x\in\mathcal X:0<p^*(x)<1\}\]
and write $\dist_{\nu}(x,A)=\inf_{z\in A}\|x-z\|_{\nu}$ for any $A\subseteq\mathcal X$.

If $\mathcal H_{\theta^*}=\emptyset$, then $\Delta^*$ has constant sign on $\mathcal X$ and the optimal fair policy is deterministic: $p^*(x)=1$ if $\Delta^*(x)>0$ for all $x\in\mathcal X$, and $p^*(x)=0$ if $\Delta^*(x)<0$ for all $x\in\mathcal X$.

If $\mathcal H_{\theta^*}\neq\emptyset$, then it holds that
\[
    \mathcal M(p^*)\subseteq\left\{x\in\mathcal X:\dist_{\nu}(x,\mathcal H_{\theta^*})\le \frac1L\right\}.
\]
Equivalently, if $\dist_{\nu}(x,\mathcal H_{\theta^*})>1/L$, then $p^*(x)\in\{0,1\}$, with $p^*(x)=1$ when $\Delta^*(x)>0$ and $p^*(x)=0$ when $\Delta^*(x)<0$.
\end{theorem}

\begin{proof}
For every policy $p$, we can write
\[
    R(p)=\mathbb E[\widetilde X^\top\beta_2^*]+\mathbb E[\Delta^*(X)p(X)].
\]
The first term is independent of $p$, so it is enough to consider only the second term.

\noindent \textbf{Case 1: $\mathcal H_{\theta^*}=\emptyset$.} In this case, by the continuity of $\Delta^*$, we have that $\Delta^*$ has constant sign on $\mathcal X$. If $\Delta^*(x)>0$ for all $x\in\mathcal X$, then the objective is pointwise increasing in $p(x)$, and hence $p^*(x)=1$ is optimal. If $\Delta^*(x)<0$ for all $x\in\mathcal X$, the same argument gives $p^*(x)=0$.

\noindent \textbf{Case 2: $\mathcal H_{\theta^*}\neq\emptyset$.} This case can happen, for example under the construction of \Cref{a_theta}. Fix $x_0\in\mathcal X$ such that $\dist_{\nu}(x_0,\mathcal H_{\theta^*})>1/L$.  We first consider the case $\Delta^*(x_0)>0$. Since $\Delta^*$ is continuous, there exists $r>0$ such that $A=\overline B_{\nu}(x_0,r)\cap\mathcal X$ satisfies $\Delta^*(x)>0$ for all $x\in A$ and $\dist_{\nu}(A,\mathcal H_{\theta^*})>1/L$. Define
\begin{equation} \label{lem:mixing_confinement_eq1}
    \phi_A(x)=\max\{0,1-L\dist_{\nu}(x,A)\}.
\end{equation}
Since $x\mapsto\dist_{\nu}(x,A)$ is 1-Lipschitz with respect to $\|\cdot\|_{\nu}$, we have that $\phi_A$ is $L$-Lipschitz. Moreover, by construction, $\phi_A\in[0,1]$ and $\phi_A(x)=1$ for all $x\in A$. We claim that $\Delta^*(x)>0$ whenever $\phi_A(x)>0$. Indeed, if $\phi_A(x)>0$, then $\dist_{\nu}(x,A)<1/L$. Suppose $\Delta^*(x)\le0$, then for some $a\in A$ with $\|x-a\|_{\nu}<1/L$, the line segment between $a$ and $x$ intersects $\mathcal H_{\theta^*}$ by convexity of $\mathcal X$ and linearity of $\Delta^*$. This gives $\dist_{\nu}(A,\mathcal H_{\theta^*})<1/L$, which is a contradiction.

Define 
\[\widetilde p(x)=\max\{p^*(x),\phi_A(x)\}.\]
Since both $p^*$ and $\phi_A$ are $L$-Lipschitz and the pointwise maximum of two $L$-Lipschitz functions is $L$-Lipschitz, we have that $\widetilde p\in\mathcal P_{L,\nu}$. Also, by construction, we have that $\widetilde p\ge p^*$ everywhere and $\widetilde p = 1$ on $A$. 

Suppose for contradiction that $p^*(x_0)<1$. Since $p^*$ is continuous, after shrinking $r$ if necessary, there exists an $\eta>0$ such that $p^*(x)\le1-\eta$ for all $x\in A$. Hence, we have $\widetilde p(x)-p^*(x)\ge\eta$ on $A$. Since $\widetilde p-p^*$ can be positive only when $\phi_A>0$, the claim above gives $\Delta^*(x)\{\widetilde p(x)-p^*(x)\}\ge0$ for all $x\in\mathcal X$. Therefore,
\begin{align*}
    R(\widetilde p)-R(p^*)
    &=\mathbb E[\Delta^*(X)\{\widetilde p(X)-p^*(X)\}]\\
    &\ge \int_A g(x)\Delta^*(x)\{\widetilde p(x)-p^*(x)\}\,dx\\
    &\ge \eta\int_A g(x)\Delta^*(x)\,dx>0,
\end{align*}
which contradicts the optimality of $p^*$. Hence $p^*(x_0)=1$.

The case $\Delta^*(x_0)<0$ is analogous. Choose $A=\overline B_{\nu}(x_0,r)\cap\mathcal X$ such that $\Delta^*(x)<0$ for all $x\in A$ and $\dist_{\nu}(A,\mathcal H_{\theta^*})>1/L$. With the same $\phi_A$ as in \eqref{lem:mixing_confinement_eq1}, it holds that $\Delta^*(x)<0$ whenever $\phi_A(x)>0$. Set 
\[\widetilde p(x)=\min\{p^*(x),1-\phi_A(x)\}.\] 
If $p^*(x_0)>0$, then after shrinking $r$ there exists $\eta>0$ such that $p^*(x)\ge\eta$ on $A$, so $\widetilde p(x)=0$ on $A$ and $p^*(x)-\widetilde p(x)\ge\eta$ on $A$. Since $\widetilde p-p^*$ can be negative only when $\phi_A>0$, we again get
\[
    R(\widetilde p)-R(p^*)=\mathbb E[\Delta^*(X)\{\widetilde p(X)-p^*(X)\}]>0,
\]
a contradiction. Therefore $p^*(x_0)=0$.

Thus, whenever $\dist_{\nu}(x,\mathcal H_{\theta^*})>1/L$, the optimal policy is deterministic with the correct sign. This proves the claim.
\end{proof}

\begin{corollary}\label[corollary]{cor:slab_width}
Denote $\theta^* = (a^*, \gamma^{*\top})^{\top}$, where $a^* \in \mathbb{R}$ and $\gamma^* \in \mathbb{R}^d$. Suppose the assumptions of \Cref{lem:mixing_confinement} hold and $\gamma^*\neq0$. Let $\nu^*$ be the dual exponent of $\nu$, i.e.~$1/\nu+1/\nu^*=1$, with the convention that $\nu^*=\infty$ when $\nu=1$ and $\nu^*=1$ when $\nu=\infty$. Define the full affine hyperplane $\overline{\mathcal H}_{\theta^*}=\{z\in\mathbb R^d:(1,z^\top)^\top\theta^*=0\}$. Then, for any $x\in\mathcal X$, we have that
\[
    \dist_{\nu}(x,\overline{\mathcal H}_{\theta^*})=\frac{|\Delta^*(x)|}{\|\gamma^*\|_{\nu^*}}.
\]
Consequently, the randomisation region satisfies
\[
    \mathcal M(p^*)\subseteq\left\{x\in\mathcal X:\frac{|\Delta^*(x)|}{\|\gamma^*\|_{\nu^*}}\le \frac1L\right\}.
\]
Equivalently, with $u_{\nu}^*=\gamma^*/\|\gamma^*\|_{\nu^*}$,
\[
    \mathcal M(p^*)\subseteq\left\{x\in\mathcal X:\left|\frac{a^*}{\|\gamma^*\|_{\nu^*}}+x^\top u_{\nu}^*\right|\le \frac1L\right\}.
\]
In particular, the randomisation region is contained in a slab around the decision boundary whose full thickness, measured in $\ell_{\nu}$-distance, is at most $2/L$.
\end{corollary}

\begin{proof}
For any $h\in\overline{\mathcal H}_{\theta^*}$, H\"{o}lder's inequality gives
\[
    |\Delta^*(x)|=|{\gamma^*}^\top(x-h)|\le \|\gamma^*\|_{\nu^*}\|x-h\|_{\nu}.
\]
Taking the infimum over $h\in\overline{\mathcal H}_{\theta^*}$ gives $\dist_{\nu}(x,\overline{\mathcal H}_{\theta^*})\ge|\Delta^*(x)|/\|\gamma^*\|_{\nu^*}$. Conversely, by duality, there exists $w\in\mathbb R^d$ such that $\|w\|_{\nu}=1$ and ${\gamma^*}^\top w=\|\gamma^*\|_{\nu^*}$. Setting $h=x-\{\Delta^*(x)/\|\gamma^*\|_{\nu^*}\}w$ gives $h\in\overline{\mathcal H}_{\theta^*}$ and $\|x-h\|_{\nu}=|\Delta^*(x)|/\|\gamma^*\|_{\nu^*}$. Therefore the distance identity holds.

Since $\overline{\mathcal H}_{\theta^*}$ contains $\mathcal H_{\theta^*}$, we have $\dist_{\nu}(x,\overline{\mathcal H}_{\theta^*})\le\dist_{\nu}(x,\mathcal H_{\theta^*})$. Combining this with \Cref{lem:mixing_confinement} gives
\[
    \mathcal M(p^*)\subseteq\left\{x\in\mathcal X:\dist_{\nu}(x,\overline{\mathcal H}_{\theta^*})\le \frac1L\right\}
    =
    \left\{x\in\mathcal X:\frac{|\Delta^*(x)|}{\|\gamma^*\|_{\nu^*}}\le \frac1L\right\}.
\]
The final display follows by writing $\Delta^*(x)=a^*+x^\top\gamma^*$.
\end{proof}

\subsection{Proof of Theorem \ref{thm:fairness_cost}} \label[appendix]{sec_app_cost}

\begin{proof}[Proof of \Cref{thm:fairness_cost}]
The proposition is a direct consequence of Lemmas \ref{lem:fairness_cost_upper}, \ref{lem:fairness_cost_lower} and \ref{lem:small_L_fairness_cost}. It remains only to justify $\mathbb  E|\Delta^*(X)| \lesssim 1$. Note that
$|\Delta^*(X)|\le \|\widetilde X\|_2\|\beta_1^*-\beta_2^*\|_2\le 2C_\beta\sqrt{d+1}$ by Assumptions \ref{a_contexts} and~\ref{a_beta}. Therefore, the Theorem holds.
\end{proof}

\begin{lemma}[Upper bound on the cost of fairness] \label[lemma]{lem:fairness_cost_upper}
Suppose Assumptions \ref{a_contexts} and \ref{a_theta} hold. Write $\theta^*=(a^*,{\gamma^*}^\top)^\top$ and let $\nu^*$ be the dual exponent of $\nu$. Then
\[
    0\le \operatorname{Cost}(L) \le  \mathbb E\left[|\Delta^*(X)|\mathbf 1\left\{\dist_{\nu}(X,\mathcal H_{\theta^*})\le \frac1L\right\}\right].
\]
Consequently, there exists a constant $C_{\rm fair}>0$, depending only on $C_g$ and the constants in \Cref{a_theta}, such that
\[
    \operatorname{Cost}(L)\le C_{\rm fair}\frac{\|\gamma^*\|_{\nu^*}^2}{L^2}.
\]
\end{lemma}

\begin{proof}
By definition, we have that 
\[
    \operatorname{Cost}(L) = \mathbb E[\{\Delta^*(X)\}_+-\Delta^*(X)p^*(X)].
\]
Note that the cost is non-negative. Indeed, if $\Delta^*(x)\ge0$, then we have that
\[
    \{\Delta^*(x)\}_+-\Delta^*(x)p^*(x)  = \Delta^*(x)\{1-p^*(x)\}\ge0.
\]
If $\Delta^*(x)<0$, then
\[
    \{\Delta^*(x)\}_+-\Delta^*(x)p^*(x) =  -\Delta^*(x)p^*(x)\ge0.
\]
Thus $\{\Delta^*(x)\}_+-\Delta^*(x)p^*(x)\ge0$ for all $x\in\mathcal X$.

Moreover, by \Cref{lem:mixing}, if $\dist_{\nu}(x,\mathcal H_{\theta^*})>1/L$, then $p^*(x)=1$ whenever $\Delta^*(x)>0$ and $p^*(x)=0$ when $\Delta^*(x)<0$. Hence, $\{\Delta^*(x)\}_+-\Delta^*(x)p^*(x) = 0$ outside the set $\{x\in\mathcal X:\dist_{\nu}(x,\mathcal H_{\theta^*})\le1/L\}$. Since $\{\Delta^*(x)\}_+-\Delta^*(x)p^*(x)$ is always bounded above by $|\Delta^*(x)|$, we obtain that
\[
    \operatorname{Cost}(L)
    \le
    \mathbb E\left[|\Delta^*(X)|\mathbf 1\left\{\dist_{\nu}(X,\mathcal H_{\theta^*})\le \frac1L\right\}\right].
\]

It remains to bound the expectation. Since $\overline{\mathcal H}_{\theta^*}$ contains $\mathcal H_{\theta^*}$, we have $\dist_{\nu}(x,\overline{\mathcal H}_{\theta^*})\le\dist_{\nu}(x,\mathcal H_{\theta^*})$. By \Cref{cor:slab_width}, we have that
\[
    \dist_{\nu}(x,\overline{\mathcal H}_{\theta^*})=\frac{|\Delta^*(x)|}{\|\gamma^*\|_{\nu^*}}.
\]
Thus, it holds that
\[
    \left\{x\in\mathcal X:\dist_{\nu}(x,\mathcal H_{\theta^*})\le \frac1L\right\}  \subseteq \left\{x\in\mathcal X:|\Delta^*(x)|\le \frac{\|\gamma^*\|_{\nu^*}}{L}\right\}.
\]

By \Cref{a_theta}, there exists $j^*\in\{1,\ldots,d\}$ such that $\theta^*\in\Theta_{j^*}$. We will write $x=(u,z)\in[0,1]^{d-1}\times[0,1]$, where $z=x_{j^*}$, and write $ \Delta^*(x)=a^*+{\gamma^*_{-j^*}}^\top u+\gamma^*_{j^*}z$. Set $r_L=\|\gamma^*\|_{\nu^*}/L$. 

Since $g(x)\le C_g$ by \Cref{a_contexts},
\begin{align*}
    \mathbb E\left[|\Delta^*(X)|\mathbf 1\{|\Delta^*(X)|\le r_L\}\right]
    &\le C_g\int_{[0,1]^d}|\Delta^*(x)|\mathbf 1\{|\Delta^*(x)|\le r_L\}\, \mathrm{d}x\\
    &= C_g\int_{[0,1]^{d-1}}\int_0^1 |a^*+{\gamma^*_{-j^*}}^\top u+\gamma^*_{j^*}z|\mathbf 1\{|a^*+{\gamma^*_{-j^*}}^\top u+\gamma^*_{j^*}z|\le r_L\}\, \mathrm{d}z\,\mathrm{d}u.
\end{align*}

By \Cref{a_theta}, there exists $j^*\in\{1,\ldots,d\}$ such that $|\gamma^*_{j^*}|\ge c_0$. Write $x=(u,z)\in[0,1]^{d-1}\times[0,1]$, where $z=x_{j^*}$ and $u=x_{-j^*}$, and set $A_u=a^*+{\gamma^*_{-j^*}}^\top u$. Then
\begin{align*}
    \operatorname{Cost}(L)
    &\le C_g\int_{[0,1]^{d-1}}\int_0^1 |A_u+\gamma^*_{j^*}z|\mathbf 1\{|A_u+\gamma^*_{j^*}z|\le r_L\}\,\mathrm{d}z\,\mathrm{d}u\\
    &= C_g\int_{[0,1]^{d-1}}\frac{1}{|\gamma^*_{j^*}|}\int_{I_u}|y|\mathbf 1\{|y|\le r_L\}\,\mathrm{d}y\,\mathrm{d}u\\
    &\le C_g\int_{[0,1]^{d-1}}\frac{1}{|\gamma^*_{j^*}|}\int_{-r_L}^{r_L}|y|\,\mathrm{d}y\,\mathrm{d}u\\
    &= C_g\int_{[0,1]^{d-1}}\frac{r_L^2}{|\gamma^*_{j^*}|}\,du\\
    &\leq \frac{C_g}{c_0}r_L^2  = \frac{C_g}{c_0}\frac{\|\gamma^*\|_{\nu^*}^2}{L^2},
\end{align*}
where $I_u$ is the interval between $A_u$ and $A_u+\gamma^*_{j^*}$, the first equality follows from the change of variable $y=A_u+\gamma^*_{j^*}z$ for fixed $u$, the first inequality follows because $|y|\mathbf 1\{|y|\le r_L\}$ is supported on $[-r_L,r_L]$, and the second inequality follows from \Cref{a_theta}.

The lemma holds with $C_{\rm fair}=C_g/c_0$.
\end{proof}

\begin{lemma}[Lower bound on the cost of fairness when L is large]\label[lemma]{lem:fairness_cost_lower}
Suppose Assumptions \ref{a_contexts} and \ref{a_theta} hold. Write $\theta^*=(a^*,{\gamma^*}^\top)^\top$. Then there exist constants $c_{\rm fair}>0$ and $L_0<\infty$, depending only on the constants in Assumptions \ref{a_contexts} and \ref{a_theta}, such that, for all $L\ge L_0$, we have that
\[
    \operatorname{Cost}(L)\ge \frac{c_{\rm fair}}{L^2}.
\]
\end{lemma}

\begin{proof}
Since $\theta^*\in\Theta$, there exists $j\in\{1,\ldots,d\}$ such that $\theta^*\in\Theta_j$. By relabelling coordinates inside this chart, write $x=(u,z)\in[0,1]^{d-1}\times[0,1]$ and write $\gamma^*=(b^*,c^*)$, where $c_0\le c^*\le C_0$. Define
\[
    \phi_{\theta^*}(u)=-\frac{a^*+{b^*}^\top u}{c^*}.
\]
By \textbf{Step 1} in the proof of \Cref{lem:uniform-interior-region}, it holds that $\phi_{\theta^*}$ is $K$-Lipschitz with respect to $\|\cdot\|_\nu$, where $K<\infty$ depends only on $d,\nu,B_0,c_0$. Moreover, there exist constants $m_0>0$ and $\rho>0$, depending only on the constants in \Cref{a_theta}, and a measurable set $S\subseteq[0,1]^{d-1}$ with $\lambda_{d-1}(S)\ge m_0$, such that
\begin{align} \label{lem:fairness_cost_lower_eq1}
    \rho\le \phi_{\theta^*}(u)\le 1-\rho,\qquad u\in S.
\end{align}
Thus the decision boundary contains a uniformly positive-measure interior patch.

Let $p\in\mathcal P_{L,\nu}$ be an arbitrary policy. For $u\in S$, write $h(u)=(u,\phi_{\theta^*}(u))$ and set $q(u)=p(h(u))$. For $s\ge0$, define
\[
    x_+(u,s)=(u,\phi_{\theta^*}(u)+s),
    \qquad
    x_-(u,s)=(u,\phi_{\theta^*}(u)-s).
\]
Note that, by \eqref{lem:fairness_cost_lower_eq1}, for $0\le s\le \min\{\rho,1/L\}$, both points belong to $\mathcal X$. Moreover, we have that
\[
    \Delta^*(x_+(u,s))=c^*s,
    \qquad
    \Delta^*(x_-(u,s))=-c^*s.
\]
Since $p$ is $L$-Lipschitz, we further have that
\[
    p(x_+(u,s)) \le p(h(u)) + L \|h(u) -x_+(u,s)\|_{\nu} = q(u)+Ls,
    \qquad
    p(x_-(u,s))\ge q(u)-Ls.
\]
Define the pointwise loss relative to the unconstrained oracle by $\ell_p(x)=\{\Delta^*(x)\}_+-\Delta^*(x)p(x)$. Then, we have that
\[
    \ell_p(x_+(u,s))=c^*s\{1-p(x_+(u,s))\}\ge c^*s(1-q(u)-Ls)_+,
\]
and
\[
    \ell_p(x_-(u,s))=c^*s p(x_-(u,s))\ge c^*s(q(u)-Ls)_+.
\]
For $L\ge1/\rho$, which gives $\rho \ge 1/L$, integrating over $s\in[0,1/L]$ gives
\begin{align*}
    \int_0^{1/L}\{\ell_p(x_+(u,s))+\ell_p(x_-(u,s))\}\,\mathrm{d}s
    &\ge c^*\int_0^{1/L}s(1-q(u)-Ls)_+\,\mathrm{d} s+c^*\int_0^{1/L}s(q(u)-Ls)_+\,\mathrm{d}s\\
    &=\frac{c^*}{6L^2}\{(1-q(u))^3+q(u)^3\}\\
    &\ge \frac{c_0}{24L^2},
\end{align*}
where the lasy inequality follows from the fact that $c^*\ge c_0$ in \Cref{a_theta} and $(1-q)^3+q^3\ge1/4$ for $q\in[0,1]$.

Using the density lower bound $g(x)\ge c_g$ from \Cref{a_contexts} and integrating over $u\in S$, we obtain
\[
    R_{\rm un}^*-R(p)
    =
    \mathbb E[\ell_p(X)]
    \ge
    c_g\lambda_{d-1}(S)\frac{c_0}{24L^2}
    \ge
    \frac{c_gm_0c_0}{24L^2}.
\]
Since this holds for every $p\in\mathcal P_{L,\nu}$, it also holds for the optimal fair policy. Therefore
\[
    \operatorname{Cost}(L)=R_{\rm un}^*-R_L^*\ge \frac{c_gm_0c_0}{24L^2}.
\]
Taking $c_{\rm fair}=c_gm_0c_0/24$ and $L_0=1/\rho$ proves the result.
\end{proof}

\begin{lemma}[Small-$L$ cost of fairness]\label[lemma]{lem:small_L_fairness_cost}
Let $R_L^*=\sup_{p\in \mathcal P_{L,\nu}}R(p)$ and define $\operatorname{Cost}(L)=R_{\rm un}^*-R_L^*$. Let $d_\nu=d^{1/\nu}$, with the convention $d_\infty=1$, and let $C_{d,\nu}=d^{\max\{1/\nu - 1/2,0\}}$. Then the following statements hold.
\begin{enumerate}
\item For $L=0$, we have that
\[
    \operatorname{Cost}(0)=\mathbb E[\{\Delta^*(X)\}_+]-\{\mathbb E[\Delta^*(X)]\}_+.
\]
\item Under \Cref{a_contexts}, for every $L\ge0$, we have that
\[
    \operatorname{Cost}(0)-Ld_\nu\mathbb E|\Delta^*(X)|\le \operatorname{Cost}(L)\le \operatorname{Cost}(0).
\]
Equivalently,
\[
    0\le \operatorname{Cost}(0)-\operatorname{Cost}(L)\le Ld_\nu\mathbb E|\Delta^*(X)|.
\]

\end{enumerate}
\end{lemma}

\begin{proof}
\noindent \textbf{Proof of (i).} When $L=0$, every policy in $\mathcal P_{0,\nu}$ is constant on $\mathcal X$. Hence any $p\in\mathcal P_{0,\nu}$ can be written as $p(x) = q$ for some $q\in[0,1]$. Since
\[
    R(p)=\mathbb E[\widetilde X^\top\beta_2^*]+\mathbb E[\Delta^*(X)p(X)],
\]
it holds that
\begin{equation} 
     R_0^*=\mathbb E[\widetilde X^\top\beta_2^*]+\sup_{q\in[0,1]}q\mathbb E[\Delta^*(X)]=\mathbb E[\widetilde X^\top\beta_2^*]+\{\mathbb E[\Delta^*(X)]\}_+.
\end{equation}
On the other hand, by the definition, we have that $ R_{\rm un}^*=\mathbb E[\widetilde X^\top\beta_2^*]+\mathbb E[\{\Delta^*(X)\}_+]$. Therefore, direct calculation gives
\[
    \operatorname{Cost}(0)=R_{\rm un}^*-R_0^*=\mathbb E[\{\Delta^*(X)\}_+]-\{\mathbb E[\Delta^*(X)]\}_+.
\]

Since $\mathcal P_{0,\nu}\subseteq\mathcal P_{L,\nu}$ for every $L\ge0$, we have $R_L^*\ge R_0^*$. Hence $\operatorname{Cost}(L)=R_{\rm un}^*-R_L^*\le R_{\rm un}^*-R_0^*=\operatorname{Cost}(0)$.

\noindent \textbf{Proof of (ii).} Fixing any $p\in\mathcal P_{L,\nu}$ and any $x_0\in\mathcal X$, we set $q=p(x_0)$. Since $p$ is $L$-Lipschitz with respect to $\|\cdot\|_\nu$, and since $\mathcal X\subseteq[0,1]^d$ under $\Cref{a_contexts}$, we have that $|p(x)-q|\le L\|x-x_0\|_\nu\le Ld_\nu$ for every $x\in\mathcal X$. Therefore, we have that
\begin{align*}
    R(p)
    &=\mathbb E[\widetilde X^\top\beta_2^*]+\mathbb E[\Delta^*(X)p(X)]\\
    &=\mathbb E[\widetilde X^\top\beta_2^*]+q\mathbb E[\Delta^*(X)]+\mathbb E[\Delta^*(X)\{p(X)-q\}]\\
    &\le \mathbb E[\widetilde X^\top\beta_2^*]+\{\mathbb E[\Delta^*(X)]\}_++\mathbb E[|\Delta^*(X)|\,|p(X)-q|]\\
    &\le \mathbb E[\widetilde X^\top\beta_2^*]+\{\mathbb E[\Delta^*(X)]\}_+ + Ld_\nu\mathbb E|\Delta^*(X)|\\
    &=R_0^*+Ld_\nu\mathbb E|\Delta^*(X)|,
\end{align*}
where the first inequality follows as $q <1$. Taking the supremum over $p\in\mathcal P_{L,\nu}$ gives $R_L^*\le R_0^*+Ld_\nu\mathbb E|\Delta^*(X)|$. Hence, it holds that
\begin{align*}
    \operatorname{Cost}(L)
    &=R_{\rm un}^*-R_L^* \ge R_{\rm un}^*-R_0^*-Ld_\nu\mathbb E|\Delta^*(X)| =\operatorname{Cost}(0)-Ld_\nu\mathbb E|\Delta^*(X)|.
\end{align*}

\end{proof}

\subsubsection{Examples for Section \ref{sec_cost}} \label[appendix]{sec_app_example}
In general, the optimal fair policy need not be unique and may not admit a simple closed-form expression. To provide a better understanding of the cost of fairness described in \Cref{sec_cost}, in this section, we consider a special example in which both the optimal policy and the cost of fairness can be computed explicitly. 
\begin{example}[One-dimensional example]\label[example]{ex:cost_fairness_1d}
Let $\mathcal X=[0,1]$, $X\sim\mathrm{Unif}[0,1]$, and $\Delta^*(x)=x-1/2$. Consider the class of policies $\mathcal P_{L,1}=\{p:[0,1]\to[0,1]: |p(x)-p(y)|\le L|x-y|,\ x,y\in[0,1]\}$.
\begin{itemize}
    \item When $L=0$, every feasible policy is constant, and hence $\operatorname{Cost}(0)=1/8$. 

    \item When $0<L\le1$, the optimal fair policy is not unique. For any $c\in[L/2,1-L/2]$, a policy of the form $p_c(x)=c+L(x-1/2)$ is optimal, and we have that $\operatorname{Cost}(L)=1/8 - L/12$. 

    \item When $L>1$, an optimal fair policy is $ p^*(x)= [ \frac12+L(x-\frac12)]_{[0,1]}$. In this case, we have that $\operatorname{Cost}(L)=1/(24L^2)$.
\end{itemize}
Thus, we have that 
\[\operatorname{Cost}(0) - \operatorname{Cost}(L) = \frac{L}{12}, \text{ when } 0<L\le1, \text{ and } \operatorname{Cost}(L)= \frac{1}{24L^2}, \text{ when } L>1.\]

\end{example}



\begin{proof}[Proof of \Cref{ex:cost_fairness_1d}]
    Let $\mathcal X=[0,1]$, $X\sim\mathrm{Unif}[0,1]$, and $\Delta^*(x)=x-1/2$. Consider the fair policy class
\[
    \mathcal P_L=\{p:[0,1]\to[0,1]: |p(x)-p(y)|\le L|x-y|,\ x,y\in[0,1]\}.
\]
We focus on the non-constant part of the reward, where
\[
    R(p)=\int_0^1 \left(x-\frac12\right)p(x)\,\mathrm{d}x.
\]
Note that the unconstrained oracle chooses arm $1$ when $x>1/2$ and arm $2$ when $x<1/2$, and hence
\[
    R_{\rm un}^*=\int_{1/2}^1\left(x-\frac12\right)\,\mathrm{d}x=\frac18.
\]

\noindent \textbf{When $L=0$.} In this case, every feasible policy must be constant on $[0,1]$. Thus any $p\in\mathcal P_0$ can be written as $p(x)= q$ for some fixed $q\in[0,1]$. Hence
\[
    R(p)=\int_0^1\left(x-\frac12\right)q\,\mathrm{d}x=q\int_0^1\left(x-\frac12\right)\,\mathrm{d}x=0.
\]
Therefore $R_0^*=0$, and consequently, $\operatorname{Cost}(0) = 1/8$.

\noindent \textbf{When $0<L\le1$.} We now compute the optimal fair value $R_L^*=\sup_{p\in\mathcal P_L}R(p)$. For any $p\in\mathcal P_L$, write $q=p(1/2)$. By the Lipschitz constraint, for $x\ge1/2$, $p(x)\le q+L(x-1/2)$, while for $x\le1/2$, $p(x)\ge q+L(x-1/2)$. Since $x-1/2$ is non-negative on $[1/2,1]$ and non-positive on $[0,1/2]$, both inequalities imply
\[
    \left(x-\frac12\right)p(x)\le \left(x-\frac12\right)\left\{q+L\left(x-\frac12\right)\right\},\qquad x\in[0,1].
\]
Therefore, we have that
\[
    R(p)\le q\int_0^1\left(x-\frac12\right)\,\mathrm{d}x+L\int_0^1\left(x-\frac12\right)^2\,\mathrm{d}x=\frac{L}{12}.
\]

In this case, upper bound of $L/12$ is attainable. Indeed, for any $c\in[L/2,1-L/2]$, the policy
\[
    p_c(x)=c+L\left(x-\frac12\right)
\]
belongs to $\mathcal P_L$, and we have
\[
    R(p_c)=\int_0^1\left(x-\frac12\right)\left\{c+L\left(x-\frac12\right)\right\}\,\mathrm{d}x=\frac{L}{12}.
\]
Hence, for $0<L\le1$, we have that $R_L^*=L/12$ and $\operatorname{Cost}(L) = 1/8 - L/12$.

\noindent \textbf{When $L>1$.} In this case, by the canonical structure in \Cref{thm:canonical_extension}, it is enough to optimise over the boundary value $q=p(1/2)$. The corresponding canonical policy is
\[
    p_q(x)=
    \begin{cases}
    \max\{0,q+L(x-1/2)\}, & x\le1/2,\\
    \min\{1,q+L(x-1/2)\}, & x\ge1/2.
    \end{cases}
\]
We compute the loss relative to the unconstrained oracle.

On the right side of the boundary, i.e.~$x \geq 1/2$, write $x=1/2+t$ with $t\in[0,1/2]$. Then $\Delta^*(x)=t>0$, so the oracle chooses arm $1$ and obtains contribution $t$, while the fair policy obtains contribution $t p_q(1/2+t)$. Hence the pointwise loss is
\[
    t\{1-p_q(1/2+t)\}=t\{1-\min(1,q+Lt)\}=t(1-q-Lt)_+.
\]
On the left side of the boundary, i.e.~$x \leq 1/2$, write $x=1/2-t$ with $t\in[0,1/2]$. Then $\Delta^*(x)=-t<0$, so the oracle chooses arm $2$ and obtains contribution $0$, while the fair policy obtains contribution $-t p_q(1/2-t)$. Hence the pointwise loss is
\[
    0-\{-t p_q(1/2-t)\}=t p_q(1/2-t)=t\max\{0,q-Lt\}=t(q-Lt)_+.
\]
Consequently, combining the two sides, the loss relative to the unconstrained oracle is
\[
    R_{\rm un}^*-R(p_q)=\int_0^{1/2}t(1-q-Lt)_+\,\mathrm{d}t+\int_0^{1/2}t(q-Lt)_+\,\mathrm{d}t.
\]
Note that the right-hand side is convex in $q$ and symmetric around $q=1/2$, so it is minimised at $q=1/2$. Therefore, when $L>1$, the optimal policy is
\[
    p^*(x)=
    \begin{cases}
    0, & x\le \frac12-\frac{1}{2L},\\
    \frac12+L\left(x-\frac12\right), & \left|x-\frac12\right|\le \frac{1}{2L},\\
    1, & x\ge \frac12+\frac{1}{2L}.
    \end{cases}
\]
For this policy,
\begin{align*}
    R_{\rm un}^*-R_L^*
    &=2\int_0^{1/(2L)}t\left(\frac12-Lt\right)\,\mathrm{d}t\\
    &=\frac{1}{24L^2}.
\end{align*}

Combining all three cases above, the cost of fairness is
\[
    \operatorname{Cost}(L)=R_{\rm un}^*-R_L^*
    =
    \begin{cases}
    \frac{1}{8}, & L=0,\\
    \frac18-\frac{L}{12}, & 0<L\le1,\\
    \frac{1}{24L^2}, & L>1.
    \end{cases}
\]

\end{proof}


\subsection{Proof of Theorem \ref{thm:canonical_extension}} \label[appendix]{sec_app_thm_canonical}\label[appendix]{sec_app_structural}
\begin{proof}[Proof of \Cref{thm:canonical_extension}]


Fix $\theta\in\mathbb R^{d+1}$ such that $\mathcal H_\theta\neq\emptyset$ and $q\in\mathcal Q_{\theta,L}$.

\noindent\textbf{Step 1: the envelopes agree with the trace on the boundary.} 
Let $x\in\mathcal H_\theta$. Since $x$ itself is admissible in the infimum defining $U_{\theta,q}$,
\[
    U_{\theta,q}(x) \le q(x)+L\|x-x\|_{\nu} = q(x).
\]
On the other hand, for every $h\in\mathcal H_\theta$, the $L$-Lipschitz property of $q$ gives $q(h)\ge q(x)-L\|h-x\|_{\nu}$, and hence
\[
    q(h)+L\|x-h\|_{\nu}\ge q(x).
\]
Taking the infimum over $h$ gives $U_{\theta,q}(x)\ge q(x)$. Therefore $U_{\theta,q}(x)=q(x)$. Similarly, using $h=x$ gives $M_{\theta,q}(x)\ge q(x)$, while $q(h)\le q(x)+L\|h-x\|_{\nu}$ gives
\[
    q(h)-L\|x-h\|_{\nu}\le q(x)
\]
for every $h\in\mathcal H_\theta$. Taking the supremum over $h$ gives $M_{\theta,q}(x)\le q(x)$. Thus $M_{\theta,q}(x)=q(x)$.

\medskip
\noindent\textbf{Step 2: the envelopes are $L$-Lipschitz.}
For any $x,y\in\mathcal X$ and every $h\in\mathcal H_\theta$, the triangle inequality gives
\[
    q(h)+L\|x-h\|_{\nu} \le q(h)+L\|y-h\|_{\nu}+L\|x-y\|_{\nu}.
\]
Taking the infimum over $h$ yields $U_{\theta,q}(x)\le U_{\theta,q}(y)+L\|x-y\|_{\nu}$. Interchanging $x$ and $y$ gives
\[
    |U_{\theta,q}(x)-U_{\theta,q}(y)| \le  L\|x-y\|_{\nu}.
\]
Similarly, for every $h\in\mathcal H_\theta$,
\[
    q(h)-L\|x-h\|_{\nu} \le  q(h)-L\|y-h\|_{\nu}+L\|x-y\|_{\nu}.
\]
Taking the supremum over $h$ and then interchanging $x$ and $y$ gives
\[
    |M_{\theta,q}(x)-M_{\theta,q}(y)|  \le L\|x-y\|_{\nu}.
\]

\medskip
\noindent\textbf{Step 3: the unclipped piecewise envelope is globally $L$-Lipschitz.}
Define
\[
    F_{\theta,q}(x)
    =
    \begin{cases}
    U_{\theta,q}(x), & x\in\mathcal X_\theta^+,\\
    M_{\theta,q}(x), & x\in\mathcal X_\theta^-.
    \end{cases}
\]
If $x$ and $y$ lie on the same side of $\mathcal H_\theta$, the claim follows from Step 2. It remains to consider the cross-boundary case. Let $x\in\mathcal X_\theta^+$ and $y\in\mathcal X_\theta^-$. If $\widetilde x^\top\theta=0$, take $h=x$. Otherwise, set
\[
    a_0=\frac{\widetilde x^\top\theta}{\widetilde x^\top\theta-\widetilde y^\top\theta}
    \qquad\text{and}\qquad
    h=(1-a_0)x+a_0y .
\]
Since $\widetilde x^\top\theta>0$ and $\widetilde y^\top\theta<0$, we have $a_0\in(0,1)$. By convexity of $\mathcal X=[0,1]^d$, $h\in\mathcal X$. Moreover,
\[
    \widetilde h^\top\theta = (1-a_0)\widetilde x^\top\theta+a_0\widetilde y^\top\theta = 0,
\]
so $h\in\mathcal H_\theta$. Since $h$ lies on the line segment between $x$ and $y$,
\[
    \|x-h\|_{\nu}+\|y-h\|_{\nu}=\|x-y\|_{\nu}.
\]
Using this particular boundary point $h$ in the definitions of $U_{\theta,q}$ and $M_{\theta,q}$ gives
\[
    U_{\theta,q}(x)-M_{\theta,q}(y) \le  q(h)+L\|x-h\|_{\nu}-\{q(h)-L\|y-h\|_{\nu}\} =  L\|x-y\|_{\nu}.
\]
For the reverse direction, for every $g\in\mathcal H_\theta$,
\[
    q(g)-L\|y-g\|_{\nu}  \le q(h)+L\|g-h\|_{\nu}-L\|y-g\|_{\nu} \le q(h)+L\|y-h\|_{\nu},
\]
where the last inequality follows from $\|g-h\|_{\nu}\le\|g-y\|_{\nu}+\|y-h\|_{\nu}$. Taking the supremum over $g$ gives
\[
    M_{\theta,q}(y)\le q(h)+L\|y-h\|_{\nu}.
\]
Similarly, for every $g\in\mathcal H_\theta$,
\[
    q(g)+L\|x-g\|_{\nu} \ge q(h)-L\|g-h\|_{\nu}+L\|x-g\|_{\nu} \ge q(h)-L\|x-h\|_{\nu}.
\]
Taking the infimum over $g$ gives
\[
    U_{\theta,q}(x)\ge q(h)-L\|x-h\|_{\nu}.
\]
Therefore
\[
    M_{\theta,q}(y)-U_{\theta,q}(x) \le L\|y-h\|_{\nu}+L\|x-h\|_{\nu} = L\|x-y\|_{\nu}.
\]
Thus $F_{\theta,q}$ is $L$-Lipschitz on $\mathcal X$.

\medskip
\noindent\textbf{Step 4: clipping preserves the Lipschitz constant.}
Since $q$ takes values in $[0,1]$, we have $U_{\theta,q}(x)\ge0$ and $M_{\theta,q}(x)\le1$ for every $x\in\mathcal X$. Hence $ p_{\theta,q}(x)=[F_{\theta,q}(x)]_{[0,1]}$.

The clipping map $u\mapsto[u]_{[0,1]}$ is 1-Lipschitz on $\mathbb R$, so for all $x,y\in\mathcal X$,
\[
    |p_{\theta,q}(x)-p_{\theta,q}(y)|  \le |F_{\theta,q}(x)-F_{\theta,q}(y)|  \le L\|x-y\|_{\nu}.
\]
Therefore $p_{\theta,q}\in\mathcal P_{L, \nu}$.

\medskip
\noindent\textbf{Step 5: domination by the canonical representative.}
Fix $\theta$ with $\mathcal H_\theta\neq\emptyset$ and $p\in\mathcal P_{L,\nu}$, and set $q=p|_{\mathcal H_\theta}$. If $x\in\mathcal X_{\theta}^+$, then for every $h\in\mathcal H_{\theta}$,
\[
    p(x)\le q(h)+L\|x-h\|_{\nu}.
\]
Taking the infimum over $h$ gives $p(x)\le U_{\theta,q}(x)$. Since $p(x)\le1$, we get
\[
    p(x)\le \min\{1,U_{\theta,q}(x)\}=p_{\theta,q}(x).
\]
If $x\in\mathcal X_{\theta}^-$, then for every $h\in\mathcal H_{\theta}$,
\[
    p(x)\ge q(h)-L\|x-h\|_{\nu}.
\]
Taking the supremum over $h$ gives $p(x)\ge M_{\theta,q}(x)$. Since $p(x)\ge0$, we get
\[
    p(x)\ge \max\{0,M_{\theta,q}(x)\}=p_{\theta,q}(x).
\]
By definition, $(\widetilde x^\top\theta)\ge0$ on $\mathcal X_{\theta}^+$ and $(\widetilde x^\top\theta)<0$ on $\mathcal X_{\theta}^-$. Therefore the two inequalities above imply
\[
    (\widetilde x^\top\theta)p_{\theta,q}(x) \ge (\widetilde x^\top\theta)p(x) \qquad\text{for all }x\in\mathcal X.
\]
\medskip

\noindent\textbf{Step 6: reduction of the optimiser.}
Fix any distribution of $X$ on $\mathcal X$ and retain the parameter $\theta$ with $\mathcal H_\theta\neq\emptyset$. Taking expectations in the pointwise inequality in \textbf{Step 5} shows that replacing any $p\in\mathcal P_{L,\nu}$ by $p_{\theta,p|_{\mathcal H_\theta}}$ cannot decrease the objective $\mathbb E\{(\widetilde X^\top\theta)p(X)\}$.

The class $\mathcal P_{L,\nu}$ is uniformly bounded and equicontinuous on the compact set $\mathcal X$. It is also closed under uniform limits: if $p_n\in\mathcal P_{L,\nu}$ and $p_n\to p$ uniformly, then $0\le p\le1$ and
\[
    |p(x)-p(y)|=\lim_{n\to\infty}|p_n(x)-p_n(y)|\le L\|x-y\|_\nu,\qquad x,y\in\mathcal X.
\]
Thus $\mathcal P_{L,\nu}$ is compact in the sup norm by Arzela--Ascoli \citep[e.g.~Theorem 7.25 in][]{rudin1976principles}. The objective is continuous in this norm because, for $p,r\in\mathcal P_{L,\nu}$,
\[
    \left|\mathbb E\left[(\widetilde X^\top\theta)\{p(X)-r(X)\}\right]\right| \le\sqrt{d+1}\|\theta\|_2\|p-r\|_\infty.
\]
An optimiser $p^\circ\in\mathcal P_{L,\nu}$ therefore exists. Set $q_\theta=p^\circ|_{\mathcal H_\theta}$. By \textbf{Step 5},
\[
    \mathbb E\{(\widetilde X^\top\theta)p_{\theta,q_\theta}(X)\} \ge\mathbb E\{(\widetilde X^\top\theta)p^\circ(X)\}.
\]
The reverse inequality follows from optimality of $p^\circ$ and feasibility of $p_{\theta,q_\theta}$. Hence $q_\theta$ attains the same optimal value. 
\end{proof}

\section{Proofs of results in Section \ref{sec_learning_incomplete_information}} \label[appendix]{sec_app_learning_incomplete_information}

\subsection{Proof of Lemma \ref{lem:finite-lp-empirical-canonical-optimizer}} \label[appendix]{sec_app_proof_empirical_optimisation}

\begin{lemma}
\label[lemma]{lem:finite-lp-empirical-canonical-optimizer}
Fix $\theta\in\mathbb R^{d+1}$ with $\mathcal H_\theta\neq\emptyset$. For any given contexts $X_1,\ldots,X_n\in[0,1]^d$, consider the finite-dimensional linear program
\begin{align} 
\widehat z = (\widehat{z}_1, \ldots, \widehat{z}_n) \in \argmax_{z_1,\ldots,z_n} &\quad \frac1n\sum_{i=1}^n \theta^\top \widetilde X_i z_i,  \label{eq_empirical_optimisation} \\
\mathrm{s.t.} & \quad 0\le z_i\le 1,\qquad i=1,\ldots,n, \notag\\
&\quad |z_i-z_j|\le L\|X_i-X_j\|_\nu,\qquad i,j=1,\ldots,n. \notag
\end{align}
Define the McShane extension of $\widehat z$ by $\widehat p(x)=\min\{1,\min_{1\le j\le n}\{\widehat z_j+L\|x-X_j\|_\nu\}\}$, for all $x\in[0,1]^d$. Let $\widehat q= \widehat p|_{\mathcal H_\theta}$. Then, we have that $\widehat q\in\mathcal Q_{\theta,L}$ and
\[
\widehat q \in \argmax_{q\in\mathcal Q_{\theta,L}}\frac1n\sum_{i=1}^n \widetilde X_i^\top \theta  p_{\theta,q}(X_i).
\]
\end{lemma}

\begin{proof}
We prove the result for every $\theta\in\mathbb R^{d+1}$ such that $\mathcal H_\theta\neq\emptyset$, which includes both $\Theta$ and $\widetilde\Theta$. 

Define
\[
\mathrm{OPT}_{\theta,n}=\sup_{q\in\mathcal Q_{\theta,L}}\frac1n\sum_{i=1}^n \theta^\top \widetilde X_i p_{\theta,q}(X_i)
\]
and let $\mathrm{OPT}_{\mathrm{LP}}$ denote the optimal value of \eqref{eq_empirical_optimisation}. The feasible region of \eqref{eq_empirical_optimisation} is a non-empty compact polytope and hence an optimiser $\widehat z$ exists.

\medskip
\noindent
\textbf{Step 1: Reduction from canonical policies to the finite LP.}
Fix any $q\in\mathcal Q_{\theta,L}$. By \Cref{thm:canonical_extension}, we have $p_{\theta,q}\in\mathcal P_{L,\nu}$. Therefore, the sample values
\[
z_i=p_{\theta,q}(X_i),\quad i=1,\ldots,n,
\]
satisfy $0\le z_i\le1$ for $i=1,\ldots,n$, and
\[
|z_i-z_j|=|p_{\theta,q}(X_i)-p_{\theta,q}(X_j)|\le L\|X_i-X_j\|_\nu,\quad i,j=1,\ldots,n.
\]
Thus $(z_1,\ldots,z_n)$ is a feasible solution to \eqref{eq_empirical_optimisation}. It then follows that
\[
\frac1n\sum_{i=1}^n \theta^\top \widetilde X_i p_{\theta,q}(X_i)\le \mathrm{OPT}_{\mathrm{LP}}.
\]
Taking the supremum over $q\in\mathcal Q_{\theta,L}$ gives
\begin{equation}
\label{eq:finite-lp-opt-upper}
\mathrm{OPT}_{\theta,n}\le \mathrm{OPT}_{\mathrm{LP}}.
\end{equation}

\noindent\textbf{Step 2: McShane interpolation of the linear program (LP) optimiser.}
For each $j \in \{1, \ldots, n\}$, the map $x\mapsto \widehat z_j+L\|x-X_j\|_\nu$ is $L$-Lipschitz. Hence its pointwise minimum is also $L$-Lipschitz. Thus, we have $\widehat p\in\mathcal P_{L,\nu}$.

We next verify that $\widehat p$ interpolates the LP solution. Fix $i\in\{1,\ldots,n\}$. Note that 
\[
\widehat p(X_i) =  \min_{1\le j\le n}\{\widehat z_j+L\|X_i-X_j\|_\nu\}\le \widehat z_i+L\|X_i-X_i\|_\nu=\widehat z_i.
\]
On the other hand, by feasibility of $\widehat z$, for every $j=1,\ldots,n$, we have that
\[
\widehat z_j+L\|X_i-X_j\|_\nu\ge \widehat z_i.
\]
Taking the minimum over $j$ and using $0\le \widehat z_i\le1$ gives
\[
\widehat p(X_i)\ge \widehat z_i.
\]
Therefore, we have that $\widehat p(X_i)=\widehat z_i$ for all $i=1,\ldots,n$. Consequently,
\begin{equation}
\label{eq:finite-lp-extension-attains-lp}
\frac1n\sum_{i=1}^n \theta^\top \widetilde X_i\widehat p(X_i)=\frac1n\sum_{i=1}^n \theta^\top \widetilde X_i \widehat z_i=\mathrm{OPT}_{\mathrm{LP}}.
\end{equation}

\medskip

\noindent\textbf{Step 3: Canonicalisation of the interpolated LP optimiser.} Let $\widehat q=\widehat p|_{\mathcal H_\theta}$. Since $\widehat p\in\mathcal P_{L,\nu}$, we have $\widehat q\in\mathcal Q_{\theta,L}$. Applying \Cref{thm:canonical_extension} gives
\begin{equation}\label{eq:finite-lp-canonical-domination}
    (\widetilde x^\top\theta)p_{\theta,\widehat q}(x)\ge(\widetilde x^\top\theta)\widehat p(x),\qquad x\in[0,1]^d.
\end{equation}

Applying \eqref{eq:finite-lp-canonical-domination} at the sample points and using \eqref{eq:finite-lp-extension-attains-lp}, we obtain
\[
\frac1n\sum_{i=1}^n \theta^\top \widetilde X_i  p_{\theta,\widehat q}(X_i)\ge \frac1n\sum_{i=1}^n \theta^\top \widetilde X_i  \widehat p(X_i)=\mathrm{OPT}_{\mathrm{LP}}.
\]
Since $\widehat q\in\mathcal Q_{\theta,L}$, the left-hand side is bounded above by $\mathrm{OPT}_{\theta,n}$. Hence
\begin{equation}
\label{eq:finite-lp-opt-lower}
\mathrm{OPT}_{\theta,n}\ge \mathrm{OPT}_{\mathrm{LP}}.
\end{equation}
Combining \eqref{eq:finite-lp-opt-upper} and \eqref{eq:finite-lp-opt-lower} gives
\[
\mathrm{OPT}_{\theta,n}=\mathrm{OPT}_{\mathrm{LP}}.
\]
Moreover, the preceding display shows that the feasible trace $\widehat q=\widehat p|_{\mathcal H_\theta}$ attains this common value. Therefore, we have that
\[
\widehat q \in \argmax_{q\in\mathcal Q_{\theta,L}}\frac1n\sum_{i=1}^n \theta^\top \widetilde X_i  p_{\theta,q}(X_i).
\]
This completes the proof.
\end{proof}

\subsection{Proof of Theorem \ref{thm:fair_regret_control}} \label[appendix]{sec_app_proof_regret}
\begin{proof}[Proof of \Cref{thm:fair_regret_control}]
    Note that the fairness guarantee follows directly from the canonical representation in \Cref{thm:canonical_extension}, while regret control is a direct consequence of Theorems~\ref{thm:lazy-regret-upper-bound} and \ref{thm:one-d-regret}.
\end{proof}

\subsubsection{Enlarged parameter space}
\begin{lemma} \label[lemma]{lem:enlarged_theta}
    For constants $A_0, B_0, C_0, c_0, \rho_0$ in \Cref{a_theta}, denote $\widetilde A_0 = 3A_0/2$, $\widetilde B_0 =3B_0/2$, $\widetilde c_0 = c_0/2$, $\widetilde C_0 = 3C_0/2$, and $\widetilde\rho_0 = \rho_0/2$. Define the enlarged parameter space $\widetilde\Theta=\bigcup_{j=1}^d \widetilde\Theta_j$, where, for each $j\in\{1,\ldots,d\}$,
    \begin{align*}
    \widetilde\Theta_j=\big\{(a, \gamma^{\top})^{\top}\in\mathbb R\times\mathbb R^d:\;& |a|\le \widetilde A_0,\ \|\gamma\|_2\le \widetilde B_0,\ \widetilde c_0\le \gamma_j\le \widetilde C_0,\\
    &\text{and there exists }u_\theta\in[0,1]^{d-1}\text{ such that }
    -(1-\widetilde\rho_0)\gamma_j\le a+\gamma_{-j}^\top u_\theta\le -\widetilde\rho_0\gamma_j
    \big\}.
    \end{align*}
    Then it holds that $\theta^* \in \widetilde\Theta$. Moreover, there exists a constant $r_{\rm loc}>0$, depending only on $\theta^*$ and the constants in \Cref{a_theta}, such that $ \{\theta:\|\theta-\theta^*\|_2\le 2r_{\rm loc}\}\subseteq \widetilde\Theta_{j^*}\subseteq\widetilde\Theta$.
\end{lemma}

\begin{proof}
Write $\theta^*=(a^*, \gamma^{*\top})^{\top}$. By \Cref{a_theta}, there exist an index $j^*\in\{1,\ldots,d\}$ and a point $u^*\in[0,1]^{d-1}$ such that $|a^*|\le A_0$, $\|\gamma^*\|_2\le B_0$, $c_0\le \gamma^*_{j^*}\le C_0$, and
\[
-(1-\rho_0)\gamma^*_{j^*}\le a^*+{\gamma^*_{-j^*}}^\top u^*\le -\rho_0\gamma^*_{j^*}.
\]
Since $\widetilde A_0=3A_0/2$, $\widetilde B_0=3B_0/2$, $\widetilde c_0=c_0/2$, $\widetilde C_0=3C_0/2$, and $\widetilde\rho_0=\rho_0/2$, the bounds $|a^*|\le \widetilde A_0$, $\|\gamma^*\|_2\le \widetilde B_0$, and $\widetilde c_0\le \gamma^*_{j^*}\le \widetilde C_0$ trivially holds. Moreover, since $\gamma^*_{j^*}\ge c_0>0$ and $\widetilde\rho_0<\rho_0$, we also have
\[
-(1-\widetilde\rho_0)\gamma^*_{j^*}\le a^*+{\gamma^*_{-j^*}}^\top u^*\le -\widetilde\rho_0\gamma^*_{j^*}.
\]
Thus $\theta^*\in\widetilde\Theta_{j^*}\subseteq\widetilde\Theta$.

It remains to show that a small neighbourhood of $\theta^*$ is contained in $\widetilde\Theta_{j^*}$. Note that the functions $\theta\mapsto |a|$, $\theta\mapsto \|\gamma\|_2$, $\theta\mapsto \gamma_{j^*}$, and $\theta\mapsto a+\gamma_{-j^*}^\top u^*$ are continuous. Since all the inequalities defining $\widetilde\Theta_{j^*}$ hold with slack at $\theta^*$, there exists $r_{\rm loc}>0$ such that, whenever $\|\theta-\theta^*\|_2\le 2r_{\rm loc}$, the same inequalities remain valid with $j=j^*$. Therefore, the lemma holds.
\end{proof}

Throughout the rest of the proof, we work over the enlarged parameter space $\widetilde\Theta$ instead of $\Theta$. This enlargement ensures that the estimator $\widehat\theta$ belongs to the parameter space under consideration, thereby avoiding the need to project it back onto $\Theta$

For $\theta=(a, \gamma^{\top})^{\top}\in\widetilde\Theta_j$, write $s_\theta(x)=\widetilde x^\top\theta=a+\gamma^\top x$ and $x=(u^{(j)},z_j)$, where $z_j=x_j$ and $u^{(j)}$ is obtained from $x$ by deleting its $j$-th coordinate. The corresponding chart representation is
\[
    \phi_{\theta,j}(u)=-\frac{a+\gamma_{-j}^\top u}{\gamma_j}.
\]

Let $\{\mathcal F_t\}_{t\ge0}$ be the filtration generated by the history up to the end of round $t$, namely
\[
    \mathcal F_t=\sigma\{X_i,A_i,Y_i:1\le i\le t\},
    \qquad
    \mathcal F_0=\{\emptyset,\Omega\}.
\]
For each $t>2T_o$, let $\tau(t)$ denote the most recent lazy update time before the action at round $t$ is chosen:
\[
    \tau(t)=\max\left\{s\le t:\text{\Cref{alg:lazy-adaptive-boundary-greedy} recomputes }(\widehat\theta_s,\widehat q_s)\text{ at the beginning of round }s\right\}.
\]

\subsubsection{Upper bound on regret}
For any $\mathcal F_{t-1}$-measurable policy $p$, define  
\[
R_t(p)=\mathbb E\{\widetilde X^\top\beta_1^*p(X)+\widetilde X^\top\beta_2^*(1-p(X))\mid\mathcal F_{t-1}\}
=\mathbb E\{\widetilde X^\top\beta_2^*\}+\mathbb E\{\widetilde X^\top\theta^*p(X)\mid\mathcal F_{t-1}\},
\]
where $X$ is an independent copy of the context and $\widetilde X=(1,X^\top)^\top$.

Let
\[q^*\in\argmax _{q\in\mathcal Q_{\theta^*,L}}R(p_{\theta^*,q}),\quad p^*=p_{\theta^*,q^*}.\]
By \Cref{thm:canonical_extension}, it holds that $\sup_{p\in\mathcal P_{L, \nu}}R(p)  = R(p_{\theta^*,q^*})$.

For $t>2T_o$, let $\tau(t)$ be the most recent lazy update time before the action at round $t$ is chosen. 

Denote the cumulative regret
\[
\mathfrak R_T = \sum_{t=1}^{2T_o}\{R_t(p^*)-R_t(\widehat p_t)\} + \sum_{t=2T_o+1}^T\{R_t(p^*)-R_t(\widehat p_{\tau(t)})\}.
\]

\begin{theorem}[Regret upper bound for the lazy adaptive boundary greedy algorithm]
\label[theorem]{thm:lazy-regret-upper-bound}
Suppose Assumptions \ref{a_contexts}-\ref{a_theta} hold. Take $\lambda=C_\lambda$ and $T_o = C_1\lceil\sqrt{T\log(T/\eta)}\rceil$, for any $\eta \in (0,1/2)$. Assume that $T\ge C_0\log(T/\eta)$, where $C_0>0$ is a sufficiently large constant depending only on the constants listed below. Then, with probability at least $1-\eta$, there exists constants $C_1, C_2 >0 $, depending only on $d,L,c_g,C_g,\sigma$, $C_\beta,\lambda, \nu$, the constants in Assumption~\ref{a_contexts} and the constants in \Cref{a_theta}, such that
\[
\mathfrak R_T
\le C_1\sqrt{T\log(T/\eta)}
+C_2\begin{cases}
\sqrt T\log(eT^2/\eta), & d=2,\\
\sqrt T\log^{3/2}(eT^2/\eta), & d=3,\\
T^{1-\frac1{d-1}}\log(eT^2/\eta), & d\ge4.
\end{cases}
\]

\end{theorem}

\begin{proof}
Let $\mathcal S_T$ denote the set of lazy update times up to time $T$. Since the algorithm updates only when $t\ge2s$, the update times grow geometrically. Hence, for every $\alpha>0$, we have that $\sum_{s\in\mathcal S_T}s^\alpha\le C_\alpha T^\alpha$, where $C_{\alpha} = (1-2^{-\alpha})^{-1}$.

Write $\ell_T=\log(T/\eta)$. Since $\ell_T\ge\log 2$, there is a constant $K>0$, depending only on the fixed model parameters, such that
\[
C_T\log\{2(d+1)T/\eta\}\vee
\frac{C_{\rm loc}\log(2T/\eta)}{r_{\rm loc}^2}
\le K\ell_T,
\]
where $r_{\rm loc}$ is the constant defined in \Cref{lem:enlarged_theta}. Choose $C_0\ge \max\{4,16C_1^2,(K/C_1)^2\}$. The condition $T\ge C_0\ell_T$ then gives
\[
T_o=C_1\lceil\sqrt{T\ell_T}\rceil
\ge C_1\sqrt{C_0}\,\ell_T
\ge K\ell_T.
\]
Thus the exploration lower bounds in \Cref{prop:lazy-regret-decomposition} and \Cref{lem:localization-and-estimation} hold, including with confidence parameter $\eta/2$. Moreover, the same horizon condition ensures $T\ge2$, so $\lceil\sqrt{T\ell_T}\rceil\le2\sqrt{T\ell_T}$ and
\[
2T_o\le4C_1\sqrt{T\ell_T}
\le \frac{4C_1}{\sqrt{C_0}}T
\le T.
\]
Hence the exploration phase fits within the time horizon. If $2T_o=T$, the claimed regret bound follows directly from the bounded one-round regret and $T_o\lesssim\sqrt{T\ell_T}$. We therefore consider $2T_o<T$ below.

By \Cref{lem:localization-and-estimation}, with probability at least $1-\eta/2$, for every $s\in\mathcal S_T$, we have that $\widehat\theta_s\in\widetilde\Theta_{j^*}$ and 
\begin{equation} \label{thm:lazy-regret-upper-bound_eq2}
    \|\widehat\theta_s-\theta^*\|_2 \le C_2\sqrt{\frac{\log(T/\eta)}{s-1}}.
\end{equation}
By \Cref{lem:control-delta-L-n}, together with a union bound over $n=1,\ldots,T$, with probability at least $1-\eta/2$, for every $n=1,\ldots,T$, $\Delta_n^L\le C\rho_{n,d-1}(\eta/T)$, where $\Delta_n^L$ and $\rho_n$ are defined in \Cref{lem:control-delta-L-n}. Equivalently, for every $s\in\mathcal S_T$, we have that 
\begin{align} \label{thm:lazy-regret-upper-bound_eq1}
    \Delta_{s-1}^L
\le C_3
\begin{cases}
(s-1)^{-1/2}\log(eT^2/\eta), & d=2,\\
(s-1)^{-1/2}\log^{3/2}(eT^2/\eta), & d=3,\\
(s-1)^{-\frac1{d-1}}\log(eT^2/\eta), & d\ge4.
\end{cases}
\end{align}
By a union bound argument, the rest of the proof is constructed conditioning on the intersection of the high-probability events in \eqref{thm:lazy-regret-upper-bound_eq2} and \eqref{thm:lazy-regret-upper-bound_eq1}.

For any post-initialisation round $t>2T_o$, write $s=\tau(t)$ for the lazy update time at which the policy used at round $t$ was computed. By \Cref{prop:lazy-regret-decomposition}, we have that
\[
R_s(p^*)-R_s(\widehat p_s)
\le C_{\theta,L}\|\widehat\theta_s-\theta^*\|_2+2\Delta_{s-1}^L.
\]
Therefore,
\begin{equation} \label{thm:lazy-regret-upper-bound_eq3}
    R_s(p^*)-R_s(\widehat p_s)
\le C_2\sqrt{\frac{\log(T/\eta)}{s-1}}+2\Delta_{s-1}^L.
\end{equation}
Note that the regret from the first $2T_o$ forced-exploration rounds is bounded by $C_4T_o$ under Assumptions \ref{a_contexts} and \ref{a_beta}. Moreover, by the lazy update mechanism, if $s\in\mathcal S_T$ is an update time, then the policy $\widehat p_s$ is kept fixed until the next update. Since the next update occurs no later than time $2s$, we have that $\#\{t\le T:\tau(t)=s\}\le s$. Therefore, applying the one-step bound to every round with $\tau(t)=s$ in \eqref{thm:lazy-regret-upper-bound_eq3} gives
\begin{align*}
    \mathfrak R_T
&\le C_4T_o+\sum_{s\in\mathcal S_T}\#\{t\le T:\tau(t)=s\}\left\{C_2\sqrt{\frac{\log(T/\eta)}{s-1}}+2\Delta_{s-1}^L\right\}\\
&\le C_4T_o+C_2\sum_{s\in\mathcal S_T}s\sqrt{\frac{\log(T/\eta)}{s-1}}+2\sum_{s\in\mathcal S_T}s\Delta_{s-1}^L.
\end{align*}

Since $s-1\asymp s$ for $s\in\mathcal S_T$, the estimation term satisfies
\[
\sum_{s\in\mathcal S_T}s\sqrt{\frac{\log(T/\eta)}{s-1}}
\le C_5\sqrt{\log(T/\eta)}\sum_{s\in\mathcal S_T}s^{1/2}
\le C_6\sqrt{T\log(T/\eta)}.
\]

It remains to control the empirical-process term. If $d=2$, then $m=d-1=1$, so
\[
\sum_{s\in\mathcal S_T}s\Delta_{s-1}^L
\le C_3\log(eT^2/\eta)\sum_{s\in\mathcal S_T}s^{1/2}
\le C_7\sqrt T\log(eT^2/\eta).
\]
If $d=3$, then $m=d-1=2$, so  
\[
\sum_{s\in\mathcal S_T}s\Delta_{s-1}^L
\le C_3\log^{3/2}(eT^2/\eta)\sum_{s\in\mathcal S_T}s^{1/2}
\le C_8\sqrt T\log^{3/2}(eT^2/\eta).
\]
If $d\ge4$, then $m=d-1\ge3$, so
\[
\sum_{s\in\mathcal S_T}s\Delta_{s-1}^L
\le C_3\log(eT^2/\eta)\sum_{s\in\mathcal S_T}s^{1-\frac1{d-1}}
\le C_9 T^{1-\frac1{d-1}}\log(eT^2/\eta).
\]
Combining the preceding bounds yields  
\[
\mathfrak R_T
\lesssim T_o+\sqrt{T\log(T/\eta)}
+\begin{cases}
\sqrt T\log(eT^2/\eta), & d=2,\\
\sqrt T\log^{3/2}(eT^2/\eta), & d=3,\\
T^{1-\frac1{d-1}}\log(eT^2/\eta), & d\ge4.
\end{cases}
\]
This completes the proof.
\end{proof}

\begin{proposition}[Regret decomposition]
\label[proposition]{prop:lazy-regret-decomposition}
Suppose Assumptions \ref{a_contexts}-\ref{a_theta} hold and for any $\eta \in (0,1/2)$, 
\[T_o\ge C_T\log\{(d+1)T/\eta\}\vee \frac{C_{\rm loc}\log(T/\eta)}{r_{\rm loc}^2}.\]
Then, for every lazy update time $s\in\mathcal S_T$, we have with probability at least $1-\eta$ that
\[
R_s(p^*)-R_s(\widehat p_s)\le C_{\theta,L}\|\widehat\theta_s-\theta^*\|_2+2\Delta_{s-1}^{L},
\]
where $\Delta_{s-1}^{L}$ is defined in \Cref{lem:control-delta-L-n}, $C_{\theta,L}$ is a constant depending only on $L,d,\nu$ and the constants in Assumptions~\ref{a_contexts} and \ref{a_theta}.
\end{proposition}

\begin{proof}

Fix $t>2T_o$ and write $s=\tau(t)$. Let $\delta_s=\|\widehat\theta_s-\theta^*\|_2$. Conditional on $\mathcal F_{s-1}$, the quantities $\widehat\theta_s$, $\widehat q_s$ and $\widehat p_s$ are fixed. 

Define
\[
\mathcal E=\left\{\widehat\theta_s\in\widetilde\Theta_{j^*}\text{ for every lazy update time }s\le T\right\}.
\]
Then by \Cref{lem:localization-and-estimation}, it holds that $\mathbb{P}(\mathcal E)\geq 1-\eta$. The rest of the proof is constructed conditioning on $\mathcal{E}$ happening.

On $\mathcal E$, by construction of the algorithm, it holds that $\widehat\theta_s\in\widetilde\Theta_{j^*}$, $\widehat q_s\in\mathcal Q_{\widehat\theta_s,L}$ and $\widehat p_s=p_{\widehat\theta_s,\widehat q_s}$.

Under the same construction as \Cref{lem:expected-transfer-no-common-chart}, for $q^*\in\arg\max_{q\in\mathcal Q_{\theta^*,L}}R(p_{\theta^*,q})$ and $p^*=p_{\theta^*,q^*}$, define the global McShane extension of $q^*$ by
\[
\bar q^*(x)=\left[\inf_{h\in\mathcal H_{\theta^*}}\{q^*(h)+L\|x-h\|_{\nu}\}\right]_{[0,1]},\qquad x\in[0,1]^d.
\]
Define the trace on the estimated boundary by
\[
q_s^\dagger(h)=\bar q^*(h),\qquad h\in\mathcal H_{\widehat\theta_s},
\]
and set $p_s^\dagger=p_{\widehat\theta_s,q_s^\dagger}$.

\noindent \textbf{Decomposition of regret.} 
For any update time $s \in \mathcal{S}_T$, the policy played is $\widehat p_s$, and we can upper bounded the regret by
\begin{align*}
R_s(p^*)-R_s(\widehat p_s)&=\mathbb E\left[\widetilde X^\top\theta^*\{p^*(X)-\widehat p_s(X)\}\mid\mathcal F_{s-1}\right]\\
&=\mathbb E\left[\widetilde X^\top\theta^*\{p^*(X)-p_s^\dagger(X)\}\mid\mathcal F_{s-1}\right]+\mathbb E\left[\widetilde X^\top\theta^*\{p_s^\dagger(X)-\widehat p_s(X)\}\mid\mathcal F_{s-1}\right],
\end{align*}
where the first term is further bounded by
\[
\mathbb E\left[\widetilde X^\top\theta^*\{p^*(X)-p_s^\dagger(X)\}\mid\mathcal F_{s-1}\right]\le C_{\rm tr,L}\delta_s
\]
as shown in \Cref{lem:expected-transfer-no-common-chart}.

For the second term, by adding and subtracting the empirical averages, we have that
\begin{align*}
&\mathbb E\left[\widetilde X^\top\theta^*\{p_s^\dagger(X)-\widehat p_s(X)\}\mid\mathcal F_{s-1}\right]\\
&=\left[\mathbb E\{\widetilde X^\top\theta^*p_s^\dagger(X)\mid\mathcal F_{s-1}\}-\frac1{s-1}\sum_{i=1}^{s-1}\widetilde X_i^\top\theta^*p_s^\dagger(X_i)\right]\\
&\quad+\frac1{s-1}\sum_{i=1}^{s-1}\widetilde X_i^\top\theta^*\{p_s^\dagger(X_i)-\widehat p_s(X_i)\}\\
&\quad+\left[\frac1{s-1}\sum_{i=1}^{s-1}\widetilde X_i^\top\theta^*\widehat p_s(X_i)-\mathbb E\{\widetilde X^\top\theta^*\widehat p_s(X)\mid\mathcal F_{s-1}\}\right].
\end{align*}
Note that under $\mathcal{E}$, we have that $\widehat\theta_s \in \widetilde \Theta$. Moreover, since $q_s^\dagger\in\mathcal Q_{\widehat\theta_s,L}$ and $\widehat q_s\in\mathcal Q_{\widehat\theta_s,L}$, both $p_s^\dagger$ and $\widehat p_s$ belong to the class $\mathcal F_L$ defining $\Delta_{s-1}^L$ in \Cref{lem:control-delta-L-n}. Therefore the two bracketed terms are both bounded above by $\Delta_{s-1}^L$, and hence
\[
\mathbb E\left[\widetilde X^\top\theta^*\{p_s^\dagger(X)-\widehat p_s(X)\}\mid\mathcal F_{s-1}\right]\le2\Delta_{s-1}^L+\frac1{s-1}\sum_{i=1}^{s-1}\widetilde X_i^\top\theta^*\{p_s^\dagger(X_i)-\widehat p_s(X_i)\}.
\]
It remains to control the empirical term. Note that we can decompose it as
\begin{align*}
\frac1{s-1}\sum_{i=1}^{s-1}\widetilde X_i^\top\theta^*\{p_s^\dagger(X_i)-\widehat p_s(X_i)\}&=\frac1{s-1}\sum_{i=1}^{s-1}\widetilde X_i^\top\widehat\theta_s\{p_s^\dagger(X_i)-\widehat p_s(X_i)\}\\
&\quad+\frac1{s-1}\sum_{i=1}^{s-1}\widetilde X_i^\top(\theta^*-\widehat\theta_s)\{p_s^\dagger(X_i)-\widehat p_s(X_i)\}.
\end{align*}
The first term is nonpositive by empirical optimality of $\widehat q_s$. Indeed, $q_s^\dagger\in\mathcal Q_{\widehat\theta_s,L}$ is a feasible competitor in the empirical optimization over $\mathcal Q_{\widehat\theta_s,L}$, so
\[
\frac1{s-1}\sum_{i=1}^{s-1}\widetilde X_i^\top\widehat\theta_s\widehat p_s(X_i)\ge \frac1{s-1}\sum_{i=1}^{s-1}\widetilde X_i^\top\widehat\theta_s p_s^\dagger(X_i).
\]
Thus, it holds that
\[
\frac1{s-1}\sum_{i=1}^{s-1}\widetilde X_i^\top\widehat\theta_s\{p_s^\dagger(X_i)-\widehat p_s(X_i)\}\le0.
\]
For the second term, since $0\le p_s^\dagger,\widehat p_s\le1$ and $\widetilde X_i=(1,X_i^\top)^\top$ with $X_i\in[0,1]^d$, we have that
\begin{align*}
\frac1{s-1}\sum_{i=1}^{s-1}\widetilde X_i^\top(\theta^*-\widehat\theta_s)\{p_s^\dagger(X_i)-\widehat p_s(X_i)\}&\le \frac1{s-1}\sum_{i=1}^{s-1}|\widetilde X_i^\top(\theta^*-\widehat\theta_s)|\\
&\le \sqrt{d+1}\|\widehat\theta_s-\theta^*\|_2=\sqrt{d+1}\delta_s.
\end{align*}
Combining the preceding bounds thus gives
\[
R_s(p^*)-R_s(\widehat p_s)\le C_{\rm tr,L}\delta_s+2\Delta_{s-1}^L+\sqrt{d+1}\delta_s.
\]
Therefore,
\[
R_s(p^*)-R_s(\widehat p_s)\le C_{\theta,L}\|\widehat\theta_s-\theta^*\|_2+2\Delta_{s-1}^L,
\]
with $C_{\theta,L}=C_{\rm tr,L}+\sqrt{d+1}$. This completes the proof. 
\end{proof}

\begin{lemma}
\label[lemma]{lem:expected-transfer-no-common-chart}
Suppose Assumptions \ref{a_contexts}-\ref{a_theta} hold. For any $\theta'\in\widetilde\Theta_{j^*}$, define
\[
\bar q^*(x)=\left[\inf_{h\in\mathcal H_{\theta^*}}\{q^*(h)+L\|x-h\|_{\nu}\}\right]_{[0,1]},\qquad x\in[0,1]^d,
\]
and define the restricted trace on $\mathcal H_{\theta'}$ by setting $q^{\dagger}_{\theta'}(h')=\bar q^*(h')$, $\ h'\in\mathcal H_{\theta'}$.

Set $p^\dagger_{\theta'}=p_{\theta',q^\dagger_{\theta'}}$. Then it holds that $q^\dagger_{\theta'}\in\mathcal Q_{\theta',L}$, and there exists a constant $C_{\rm tr,L}<\infty$, depending only on $d,\nu,L$, the constants in Assumption~\ref{a_contexts} and the constants in \Cref{a_theta}, such that
\[
\mathbb E\left[|\widetilde X^\top\theta^*|\,|p^*(X)-p^\dagger_{\theta'}(X)|\right]\le C_{\rm tr,L}\|\theta'-\theta^*\|_2.
\]
\end{lemma}

\begin{proof}
Let $\delta=\|\theta'-\theta^*\|_2$. Since $\theta^*,\theta'\in\widetilde\Theta_{j^*}$, in the proof below, we worked with the fixed $j^*$-th chart $\widetilde\Theta_{j^*}$.

\noindent \textbf{Step 1: $q^\dagger_{\theta'}$ is a feasible $L$-Lipschitz trace.} The follow of this step follows from a similar argument as the one in the proof of \Cref{thm:canonical_extension}. Define $G(x)=\inf_{h\in\mathcal H_{\theta^*}}\{q^*(h)+L\|x-h\|_{\nu}\}$, where $q^*\in\argmax _{q\in\mathcal Q_{\theta^*,L}}R(p_{\theta^*,q})$. For any $x,y\in[0,1]^d$ and any $h\in\mathcal H_{\theta^*}$, the triangle inequality gives
\[
q^*(h)+L\|x-h\|_{\nu}\le q^*(h)+L\|y-h\|_{\nu}+L\|x-y\|_{\nu}.
\]
Taking the infimum over $h\in\mathcal H_{\theta^*}$ gives $G(x)\le G(y)+L\|x-y\|_{\nu}$. Reversing the roles of $x$ and $y$ gives $|G(x)-G(y)|\le L\|x-y\|_{\nu}$. Since the clipping map $z\mapsto [z]_{[0,1]}$ is nonexpansive, $\bar q^*$ is also $L$-Lipschitz on $[0,1]^d$ and takes values in $[0,1]$. Therefore, for any $h'_1,h'_2\in\mathcal H_{\theta'}$, it holds that
\[
|q^\dagger_{\theta'}(h'_1)-q^\dagger_{\theta'}(h'_2)|=|\bar q^*(h'_1)-\bar q^*(h'_2)|\le L\|h'_1-h'_2\|_{\nu}.
\]
Thus $q^\dagger_{\theta'}\in\mathcal Q_{\theta',L}$. 

\noindent \textbf{Step 2: Upper bound on $\|U_{\theta^*,q^*}-U_{\theta',q^\dagger_{\theta'}}\|_{\infty}$.} For any $x \in [0,1]^d$, define
\[
U_{\theta^*,q^*}(x)=\inf_{h\in\mathcal H_{\theta^*}}\{q^*(h)+L\|x-h\|_{\nu}\}, \quad U_{\theta',q^\dagger_{\theta'}}(x) = \inf_{h\in\mathcal H_{\theta'}}\{q^\dagger_{\theta'}(h)+L\|x-h\|_{\nu}\}.
\]
By \Cref{lem:hausdorff}, for every $h\in\mathcal H_{\theta^*}$, there exists $h'\in\mathcal H_{\theta'}$ such that $\|h-h'\|_{\nu}\le C_H\delta$. Note that the constant $C_H$ already contains the norm-comparison factor $C_{d,\nu}$ from \Cref{lem:norm-comparison}, as shown explicitly in the proof of \Cref{lem:hausdorff}. Hence, for such a pair, it holds that
\begin{align*}
U_{\theta',q^\dagger_{\theta'}}(x)&\le q^\dagger_{\theta'}(h')+L\|x-h'\|_{\nu}\\
&=\bar q^*(h')+L\|x-h'\|_{\nu} \\
&\le q^*(h)+L\|h'-h\|_{\nu}+L\|x-h'\|_{\nu}\\
& \leq  q^*(h)+L\|h'-h\|_{\nu}+L\|h-h'\|_{\nu} +L\|x-h\|_{\nu} \\
&\le q^*(h)+L\|x-h\|_{\nu}+2LC_H\delta.
\end{align*}
Taking the infimum over $h\in\mathcal H_{\theta^*}$ thus gives
\[
U_{\theta',q^\dagger_{\theta'}}(x)\le U_{\theta^*,q^*}(x)+2LC_H\delta.
\]

Conversely, for every $h'\in\mathcal H_{\theta'}$ and every $h\in\mathcal H_{\theta^*}$, the triangle inequality gives
\begin{align*}
q^*(h)+L\|h'-h\|_{\nu}+L\|x-h'\|_{\nu}\ge q^*(h)+L\|x-h\|_{\nu}.
\end{align*}
Taking the infimum over $h\in\mathcal H_{\theta^*}$ gives
\begin{equation} \label{lem:expected-transfer-no-common-chart_eq1}
    G(h')+L\|x-h'\|_{\nu}\ge U_{\theta^*,q^*}(x).
\end{equation}
Since $q^*$ takes values in $[0,1]$, we have $G(h')\ge 0$. Moreover, by \Cref{lem:hausdorff}, for every $h'\in\mathcal H_{\theta'}$ there exists $h\in\mathcal H_{\theta^*}$ such that $\|h-h'\|_{\nu}\le C_H\delta$, and hence
\[
G(h')\le q^*(h)+L\|h'-h\|_{\nu}\le 1+LC_H\delta.
\]
Therefore, since $\bar q^*(h')$ is obtained by truncating $G(h')$ to $[0,1]$, we have that
\[
0\le G(h')-\bar q^*(h')\le LC_H\delta.
\]
Combining this with \eqref{lem:expected-transfer-no-common-chart_eq1} yields
\[
\bar q^*(h')+L\|x-h'\|_{\nu}\ge U_{\theta^*,q^*}(x)-LC_H\delta.
\]
Taking the infimum over $h'\in\mathcal H_{\theta'}$ gives
\[
U_{\theta',q^\dagger_{\theta'}}(x)\ge U_{\theta^*,q^*}(x)-LC_H\delta.
\]

Therefore, we have that, for any $x \in  [0,1]^d$,
\[
|U_{\theta^*,q^*}(x)-U_{\theta',q^\dagger_{\theta'}}(x)|\le 2LC_H\delta.
\]

\noindent \textbf{Step 3: Upper bound on $\|M_{\theta^*,q^*}-M_{\theta',q^\dagger_{\theta'}}\|_{\infty}$.} For any $x \in [0,1]^d$, define
\[M_{\theta^*,q^*}(x)=\sup_{h\in\mathcal H_{\theta^*}}\{q^*(h)-L\|x-h\|_{\nu}\},\quad M_{\theta',q^\dagger_{\theta'}}(x) = \sup_{h\in\mathcal H_{\theta'}}\{q^\dagger_{\theta'}(h)-L\|x-h\|_{\nu}\}.\]

First, by \Cref{lem:hausdorff}, for every $h'\in\mathcal H_{\theta'}$, we can choose $h\in\mathcal H_{\theta^*}$ such that $\|h-h'\|_{\nu}\le C_H\delta$. Since $\bar q^*(h')\le q^*(h)+L\|h'-h\|_{\nu}$, we have
\begin{align*}
\bar q^*(h')-L\|x-h'\|_{\nu}&\le q^*(h)+L\|h'-h\|_{\nu}-L\|x-h'\|_{\nu}\\
&\le q^*(h)-L\|x-h\|_{\nu}+2LC_H\delta\\
&\le M_{\theta^*,q^*}(x)+2LC_H\delta,
\end{align*}
where the second inequality follows from the triangle inequality. Taking the supremum over $h'\in\mathcal H_{\theta'}$ gives
\[
M_{\theta',q^\dagger_{\theta'}}(x)\le M_{\theta^*,q^*}(x)+2LC_H\delta.
\]
For the reverse inequality, by \Cref{lem:hausdorff}, for any fix $h\in\mathcal H_{\theta^*}$, we can choose $h'\in\mathcal H_{\theta'}$ such that $\|h-h'\|_{\nu}\le C_H\delta$. For any $\widetilde h\in\mathcal H_{\theta^*}$, since $q^*$ is $L$-Lipschitz on $\mathcal H_{\theta^*}$, we have that
\[
q^*(\widetilde h)+L\|h'-\widetilde h\|_{\nu}\ge q^*(h)-L\|h-\widetilde h\|_{\nu}+L\|h'-\widetilde h\|_{\nu}\ge q^*(h)-L\|h-h'\|_{\nu}.
\]
Taking the infimum over $\widetilde h\in\mathcal H_{\theta^*}$ gives $\bar q^*(h')\ge q^*(h)-LC_H\delta$. Hence, we have that
\begin{align*}
M_{\theta',q^\dagger_{\theta'}}(x)&\ge \bar q^*(h')-L\|x-h'\|_{\nu}\\
&\ge q^*(h)-LC_H\delta-L\|x-h'\|_{\nu}\\
&\ge q^*(h)-L\|x-h\|_{\nu}-2LC_H\delta.
\end{align*}
Taking the supremum over $h\in\mathcal H_{\theta^*}$ gives
\[
M_{\theta',q^\dagger_{\theta'}}(x)\ge M_{\theta^*,q^*}(x)-2LC_H\delta.
\]
Therefore,
\[
|M_{\theta^*,q^*}(x)-M_{\theta',q^\dagger_{\theta'}}(x)|\le 2LC_H\delta.
\]

\noindent \textbf{Step 4: Upper bound on $\mathbb E[|\widetilde X^\top\theta^*|\,|p^*(X)-p^\dagger_{\theta'}(X)|]$.} For $\theta=(a, \gamma^{\top})^{\top}\in\widetilde\Theta$, define $s_\theta(x)=\widetilde x^\top\theta=a+\gamma^\top x$. Also, define the sign-disagreement set 
\[
S(\theta^*,\theta')=\{x\in[0,1]^d:s_{\theta^*}(x)\ge 0,\ s_{\theta'}(x)<0\}\cup\{x\in[0,1]^d:s_{\theta^*}(x)<0,\ s_{\theta'}(x)\ge 0\}.
\]
On $S(\theta^*,\theta')^c$, both policies use the same type of envelope, either both upper envelopes or both lower envelopes. Since $y\mapsto \min\{1,y\}$ and $y\mapsto \max\{0,y\}$ are nonexpansive, results in \textbf{Steps 2 and 3} imply that
\[
|p^*(x)-p^\dagger_{\theta'}(x)|\le 2LC_H\delta, \text{ for all } x\in S(\theta^*,\theta')^c.
\]
On $S(\theta^*,\theta')$, we will use the crude bound $|p^*(x)-p^\dagger_{\theta'}(x)|\le 1$.

Hence we have that 

\begin{align*}
&\mathbb E\left[|\widetilde X^\top\theta^*|\,|p^*(X)-p^\dagger_{\theta'}(X)|\right]\\
&\quad=\mathbb E\left[|\widetilde X^\top\theta^*|\,|p^*(X)-p^\dagger_{\theta'}(X)|\mathbf 1\{X\in S(\theta^*,\theta')^c\}\right]\\
&\qquad+\mathbb E\left[|\widetilde X^\top\theta^*|\,|p^*(X)-p^\dagger_{\theta'}(X)|\mathbf 1\{X\in S(\theta^*,\theta')\}\right]\\
&\quad\le 2C_{\beta}\sqrt{d+1}\left\{2LC_H\delta\,P\{X\in S(\theta^*,\theta')^c\}+P\{X\in S(\theta^*,\theta')\}\right\}\\
&\quad\le 2C_{\beta}\sqrt{d+1}\left\{2LC_H\delta+P\{X\in S(\theta^*,\theta')\}\right\}\\
&\quad\le 2C_{\beta}\sqrt{d+1}\left\{2LC_H\delta+C_{\rm slab}\delta\right\}\\
&\quad=2C_{\beta}\sqrt{d+1}(2LC_H+C_{\rm slab})\delta,
\end{align*}
where the first inequality follows from the fact that $|\widetilde X^\top\theta^*|\le \sqrt{d+1}\|\theta^*\|_2 \le 2C_{\beta}\sqrt{d+1}$ and the third inequality follows from \Cref{lem:slab_probability}. Thus, the claim is proved.
\end{proof}

\begin{lemma}[Hausdorff stability]
\label[lemma]{lem:hausdorff}
Suppose \Cref{a_theta} holds. Then there exists a constant $C_H<\infty$, depending only on $d,\nu$ and constants $\widetilde A_0,\widetilde B_0,\widetilde c_0,\widetilde C_0, \widetilde\rho_0$ in \Cref{lem:enlarged_theta}, such that, for every $j\in\{1,\ldots,d\}$ and all $\theta,\theta'\in\widetilde\Theta_j$,
\begin{equation*}
d_{H,\nu}(\mathcal H_\theta,\mathcal H_{\theta'})
\le C_H\|\theta-\theta'\|_2.
\end{equation*}
\end{lemma}

\begin{proof}
Fix $j\in\{1,\ldots,d\}$. By relabelling coordinates, without loss of generality, we assume that $j=d$. Write $\theta=(a, \gamma^{\top})^{\top}$ and $\theta'=(a', \gamma^{'\top})^{\top}$, where $a,a' \in \mathbb{R}$ and $\gamma, \gamma' \in \mathbb{R}^{d}$. Then write $\gamma=(b,c)$ and $\gamma'=(b',c')$ in this chart, where $b,b' \in \mathbb{R}^{d-1}$ and $c,c' \in \mathbb{R}$. Note that, for every $u\in[0,1]^{d-1}$, it holds that
\begin{align} \notag
|\phi_\theta(u)-\phi_{\theta'}(u)|
&= \left| \frac{a+b^\top u}{c} -\frac{a'+b'^\top u}{c'}\right| \\ \notag
&\le \frac{|a-a'|+|(b-b')^\top u|}{\widetilde c_0}+\frac{|a'+b'^\top u|}{\widetilde c_0^2}|c-c'| \\ \label{lem:hausdorff_eq1}
&\le C_\phi\|\theta-\theta'\|_2,
\end{align}
where
\[
C_\phi=\frac{1+\sqrt{d-1}}{\widetilde c_0}+\frac{\widetilde A_0+\widetilde B_0\sqrt{d-1}}{\widetilde c_0^2}.
\]
Set $\eta=C_\phi\|\theta-\theta'\|_2$ and let $C_{d,\nu}=d^{\max\{1/\nu-1/2,0\}}$ be the constant in \Cref{lem:norm-comparison}. We show that every point of $\mathcal H_\theta$ is within $\ell_\nu$-distance $C\eta$ of $\mathcal H_{\theta'}$. The reverse direction then follows by symmetry.

Take any $h \in\mathcal H_\theta$ and, by \Cref{a_theta}, we can always find $u \in [0,1]^{d-1}$ such that it can be decomposed as $h = (u,\phi_\theta(u))$. By the definition of $\mathcal H_\theta$, $0\le \phi_\theta(u)\le1$. Therefore, by \eqref{lem:hausdorff_eq1} and the triangle inequality, we have that $-\eta\le \phi_{\theta'}(u)\le1+\eta$.

\noindent \textbf{Case 1:} If $0\le \phi_{\theta'}(u)\le1$, then $h'=(u,\phi_{\theta'}(u))\in\mathcal H_{\theta'}$, and hence
\[
    \|h-h'\|_\nu=|\phi_\theta(u)-\phi_{\theta'}(u)|\le\eta.
\]

\noindent \textbf{Case 2:} Next suppose $\phi_{\theta'}(u)<0$. By the definition of $\widetilde\Theta_j$, there exists $u_{\theta'}\in[0,1]^{d-1}$ such that $\widetilde\rho_0\le \phi_{\theta'}(u_{\theta'})\le1-\widetilde\rho_0$. For $\lambda\in[0,1]$, define $u_\lambda=(1-\lambda)u+\lambda u_{\theta'}$. Since $\phi_{\theta'}$ is affine, the map $\lambda\mapsto \phi_{\theta'}(u_\lambda)$ is affine. Since $\phi_{\theta'}(u)<0$ and $\phi_{\theta'}(u_{\theta'})\ge\widetilde\rho_0$, there exists $\lambda_0\in[0,1]$ such that $\phi_{\theta'}(u_{\lambda_0})=0$. Moreover,
\begin{equation*}
\lambda_0=\frac{-\phi_{\theta'}(u)}{\phi_{\theta'}(u_{\theta'})-\phi_{\theta'}(u)}\le\frac{\eta}{\widetilde\rho_0+\eta}\le\frac{\eta}{\widetilde\rho_0}.
\end{equation*}
Set $v=u_{\lambda_0}$ and $h'=(v,0)$. Then $h'\in\mathcal H_{\theta'}$. We first bound the distance in the Euclidean norm and then convert it to the fairness norm. Since $u,u_{\theta'}\in[0,1]^{d-1}$,
\[
    \|v-u\|_2
    \le
    \lambda_0\|u_{\theta'}-u\|_2
    \le
    \frac{\sqrt{d-1}}{\widetilde\rho_0}\eta.
\]
By \Cref{lem:norm-comparison},
\[
    \|v-u\|_\nu
    \le
    C_{d,\nu}\|v-u\|_2
    \le
    \frac{C_{d,\nu}\sqrt{d-1}}{\widetilde\rho_0}\eta.
\]
Since $\phi_{\theta'}(u)<0$ and $|\phi_\theta(u)-\phi_{\theta'}(u)|\le\eta$, we also have $\phi_\theta(u)\le\eta$. Hence, using the triangle inequality for $\|\cdot\|_\nu$,
\begin{align*}
\|h-h'\|_\nu
&\le \|u-v\|_\nu+|\phi_\theta(u)| \\
&\le \left(1+\frac{C_{d,\nu}\sqrt{d-1}}{\widetilde\rho_0}\right)\eta.
\end{align*}

\noindent \textbf{Case 3:} Finally suppose $\phi_{\theta'}(u)>1$. Again, by the definition of $\widetilde\Theta_j$, there exists $u_{\theta'}\in[0,1]^{d-1}$ such that $\widetilde\rho_0\le \phi_{\theta'}(u_{\theta'})\le1-\widetilde\rho_0$. Along the same segment $u_\lambda=(1-\lambda)u+\lambda u_{\theta'}$, there exists $\lambda_0\in[0,1]$ such that $\phi_{\theta'}(u_{\lambda_0})=1$. Moreover,
\begin{equation*}
\lambda_0=\frac{\phi_{\theta'}(u)-1}{\phi_{\theta'}(u)-\phi_{\theta'}(u_{\theta'})}\le\frac{\eta}{\widetilde\rho_0+\eta}\le\frac{\eta}{\widetilde\rho_0}.
\end{equation*}
Set $v=u_{\lambda_0}$ and $h'=(v,1)$. Then $h'\in\mathcal H_{\theta'}$. As in \textbf{Case 2},
\[
    \|v-u\|_\nu
    \le
    C_{d,\nu}\|v-u\|_2
    \le
    \frac{C_{d,\nu}\sqrt{d-1}}{\widetilde\rho_0}\eta,
\]
where the first inequality follows from \Cref{lem:norm-comparison}. Since $\phi_{\theta'}(u)>1$ and $|\phi_\theta(u)-\phi_{\theta'}(u)|\le\eta$, we have $\phi_\theta(u)\ge1-\eta$. Hence
\begin{align*}
\|h-h'\|_\nu
&\le
\|u-v\|_\nu+|\phi_\theta(u)-1| \\
&\le
\left(1+\frac{C_{d,\nu}\sqrt{d-1}}{\widetilde\rho_0}\right)\eta.
\end{align*}

Combining the three cases gives
\begin{equation*}
\sup_{h\in\mathcal H_\theta}\operatorname{dist}_{\nu}(h,\mathcal H_{\theta'})=\sup_{h\in\mathcal H_\theta}\inf_{\widetilde h\in \mathcal H_{\theta'}}\|h-\widetilde h\|_{\nu}\le C C_\phi\|\theta-\theta'\|_2,
\end{equation*}
where $C<\infty$ depends only on $d$, $\nu$ and $\widetilde\rho_0$. By symmetry, we also have
\begin{equation*}
\sup_{h'\in\mathcal H_{\theta'}}\operatorname{dist}_{\nu}(h',\mathcal H_\theta)
\le C C_\phi\|\theta-\theta'\|_2.
\end{equation*}
Therefore, the lemma follows.
\end{proof}

\subsubsection{Control of complexity of functional class} \label[appendix]{sec_app_cts_complexity}

\begin{lemma}[Probability bound for disagreement]
\label[lemma]{lem:slab_probability}
Suppose Assumptions \ref{a_contexts}-\ref{a_theta} hold. For $\theta=(a, \gamma^{\top})^{\top}\in\widetilde\Theta$, define $s_\theta(x)=\widetilde x^\top\theta=a+\gamma^\top x$. Then there exists a constant $C_{\rm slab}<\infty$, depending only on $d,\widetilde c_0$ and $C_g$, such that, for all $\theta,\theta'\in\widetilde\Theta$,
\[
\mathbb{P}\left[\{s_\theta(X)\ge 0,\ s_{\theta'}(X)<0\}\cup\{s_\theta(X)<0,\ s_{\theta'}(X)\ge 0\}\right]\le C_{\rm slab}\|\theta-\theta'\|_2.
\]
\end{lemma}

\begin{proof}
Fix $\theta=(a, \gamma^{\top})^{\top}$ and $\theta'=(a', \gamma^{'\top})^{\top}$ in $\widetilde\Theta$, and write $\delta=\|\theta-\theta'\|_2$. Define the sign-disagreement set
\[
S(\theta,\theta')=\{x\in[0,1]^d:s_\theta(x)\ge 0,\ s_{\theta'}(x)<0\}\cup\{x\in[0,1]^d:s_\theta(x)<0,\ s_{\theta'}(x)\ge 0\}.
\]

We first show that $S(\theta,\theta')$ is contained in a slab around $\mathcal H_\theta$. Note that if $x\in S(\theta,\theta')$, then it is the case that $s_\theta(x)$ and $s_{\theta'}(x)$ have different signs. If $s_\theta(x)\ge0$ and $s_{\theta'}(x)<0$, then $|s_\theta(x)|=s_\theta(x)\le s_\theta(x)-s_{\theta'}(x)$. If $s_\theta(x)<0$ and $s_{\theta'}(x)\ge0$, then $|s_\theta(x)|=-s_\theta(x)\le s_{\theta'}(x)-s_\theta(x)$. Therefore, in both cases, we have that
\[
|s_\theta(x)|\le |s_\theta(x)-s_{\theta'}(x)|.
\]

Since $\|\widetilde x\|_2\le \sqrt{d+1}$ on $[0,1]^d$, the Cauchy--Schwarz inequality gives
\[
|s_\theta(x)-s_{\theta'}(x)|\le \sqrt{d+1}\|\theta-\theta'\|_2=\sqrt{d+1}\delta.
\]
Therefore, it holds that
\begin{align} \label{lem:slab_probability_eq1}
    S(\theta,\theta')\subseteq\{x\in[0,1]^d:|s_\theta(x)|\le \sqrt{d+1}\delta\}.
\end{align}
Next we bound the probability of this slab. Since $\theta\in\widetilde\Theta$, there exists at least one index $j$ such that $\theta\in\widetilde\Theta_j$. By relabelling coordinates, write $x=(u,z)\in[0,1]^{d-1}\times[0,1]$ and $s_\theta(x)=a+b^\top u+cz$, where $c=\gamma_j\ge\widetilde c_0$. Let $r>0$. For fixed $u\in[0,1]^{d-1}$, the set of $z\in[0,1]$ such that $|a+b^\top u+cz|\le r$ is an interval, possibly empty, of length at most $2r/c$, and hence at most $2r/\widetilde c_0$. Therefore, by the density upper bound in \Cref{a_contexts} and the Fubini's theorem, we have that
\begin{align}  \notag
\mathbb{P}\{|s_\theta(X)|\le r\} &\le C_g\lambda_d\{x\in[0,1]^d:|s_\theta(x)|\le r\} \\\label{lem:slab_probability_eq2}
&\le C_g\int_{[0,1]^{d-1}}\frac{2r}{\widetilde c_0}\,\mathrm{d}u =\frac{2C_g}{\widetilde c_0}r.
\end{align}
Taking $r=\sqrt{d+1}\delta$ gives
\[
\mathbb{P}\{X\in S(\theta,\theta')\}\le  \mathbb{P} \{|s_\theta(X)|\le \sqrt{d+1}\delta\} \le \frac{2C_g\sqrt{d+1}}{\widetilde c_0}\delta.
\]
Thus the claim holds with $C_{\rm slab}=2C_g\sqrt{d+1}/\widetilde c_0$.
\end{proof}

\begin{lemma}[Transfer in $L_2(P)$]
\label[lemma]{lem:l2_transfer_slab}
Suppose Assumptions \ref{a_contexts}-\ref{a_theta} hold. Then there exists a constant $C_{\rm tr}<\infty$ such that, for every $j\in\{1,\ldots,d\}$, every $\theta,\theta'\in\widetilde\Theta_j$ and every $q\in \mathcal Q_{\theta,L}$, there exists $q'\in\mathcal Q_{\theta',L}$ satisfying
\[
\|p_{\theta,q}-p_{\theta',q'}\|_{L_2(P)}\le C_{\rm tr}\|\theta-\theta'\|_2^{1/2}.
\]
\end{lemma}

\begin{proof}
The proof follows from a similar argument as the proof for \Cref{lem:expected-transfer-no-common-chart}. We only include the difference here.

Fix $j\in\{1,\ldots,d\}$, $\theta,\theta'\in\widetilde\Theta_j$ and $q\in \mathcal Q_{\theta,L}$. Write $\delta=\|\theta-\theta'\|_2$. Define the global McShane extension of $q$ from $\mathcal H_\theta$ to $[0,1]^d$ by
\begin{equation*}
G(x)=\inf_{h\in \mathcal H_\theta}\{q(h)+L\|x-h\|_{\nu}\}, \qquad \bar q(x)=[G(x)]_{[0,1]}, \qquad x\in[0,1]^d,
\end{equation*}
where $[y]_{[0,1]}=\min\{1,\max\{0,y\}\}$. Define the transferred trace on $\mathcal H_{\theta'}$ by $q'(h')=\bar q(h')$ for all $h'\in \mathcal H_{\theta'}$.

\noindent\textbf{Step 1: Feasibility of the transferred trace.} By the same argument as in \textbf{Step 1} of the proof of \Cref{lem:expected-transfer-no-common-chart}, with $(\theta^*,q^*)$ replaced by $(\theta,q)$, the function $G$ is $L$-Lipschitz on $[0,1]^d$. Since the clipping map $y\mapsto [y]_{[0,1]}$ is nonexpansive, $\bar q$ is also $L$-Lipschitz and takes values in $[0,1]$. Therefore, its restriction $q'=\bar q|_{\mathcal H_{\theta'}}$ belongs to $\mathcal Q_{\theta',L}$.

\noindent\textbf{Step 2: Upper-envelope comparison.} Define
\begin{equation*}
U_{\theta,q}(x)=\inf_{h\in \mathcal H_\theta}\{q(h)+L\|x-h\|_{\nu}\}, \qquad U_{\theta',q'}(x)=\inf_{h'\in \mathcal H_{\theta'}}\{q'(h')+L\|x-h'\|_{\nu}\}.
\end{equation*}
By the same argument as in \textbf{Step 2} of the proof of \Cref{lem:expected-transfer-no-common-chart}, with $(\theta^*,q^*)$ replaced by $(\theta,q)$. Therefore, for every $x\in[0,1]^d$, we can show that
\begin{equation*}
|U_{\theta,q}(x)-U_{\theta',q'}(x)|\le 2LC_H\delta.
\end{equation*}

\noindent\textbf{Step 3: Lower-envelope comparison.} Define
\begin{equation*}
M_{\theta,q}(x)=\sup_{h\in \mathcal H_\theta}\{q(h)-L\|x-h\|_{\nu}\}, \qquad M_{\theta',q'}(x)=\sup_{h'\in \mathcal H_{\theta'}}\{q'(h')-L\|x-h'\|_{\nu}\}.
\end{equation*}
By the same argument as in \textbf{Step 3} of the proof of \Cref{lem:expected-transfer-no-common-chart}, again with $(\theta^*,q^*)$ replaced by $(\theta,q)$, we obtain, for every $x\in[0,1]^d$,
\begin{equation*}
|M_{\theta,q}(x)-M_{\theta',q'}(x)|\le 2LC_H\delta.
\end{equation*}

\noindent\textbf{Step 4: $L_2(P)$ bound.} Under the same notation as in \Cref{lem:slab_probability}, define the sign-disagreement set
\begin{equation*}
S(\theta,\theta')=\{x\in[0,1]^d:s_\theta(x)\ge 0,\ s_{\theta'}(x)<0\}\cup\{x\in[0,1]^d:s_\theta(x)<0,\ s_{\theta'}(x)\ge 0\}.
\end{equation*}
On $S(\theta,\theta')^c$, the two policies use the same type of clipped envelope; either both use the upper envelopes or both use the lower envelopes. Since the maps $y\mapsto \min\{1,y\}$ and $y\mapsto \max\{0,y\}$ are nonexpansive, \textbf{Steps 2 and 3} imply that
\begin{equation*}
|p_{\theta,q}(x)-p_{\theta',q'}(x)|\le 2LC_H\delta, \qquad x\in S(\theta,\theta')^c.
\end{equation*}
On $S(\theta,\theta')$, we use the crude bound $|p_{\theta,q}(x)-p_{\theta',q'}(x)|\le 1$. Hence, by \Cref{lem:slab_probability}, we have
\begin{align*}
\|p_{\theta,q}-p_{\theta',q'}\|_{L_2(P)}^2&=\mathbb E\left[\left|p_{\theta,q}(X)-p_{\theta',q'}(X)\right|^2\right]\\
&=\mathbb E\left[\left|p_{\theta,q}(X)-p_{\theta',q'}(X)\right|^2\mathds{1}\{X\in S(\theta,\theta')^c\}\right]\\
&\quad +\mathbb E\left[\left|p_{\theta,q}(X)-p_{\theta',q'}(X)\right|^2\mathds{1}\{X\in S(\theta,\theta')\}\right]\\
&\le (2LC_H)^2\delta^2+P\{X\in S(\theta,\theta')\}\\
&\le (2LC_H)^2\delta^2+C_{\rm slab}\delta.
\end{align*}
Moreover, by the definition of $\widetilde\Theta$, it holds that $\delta\le \operatorname{diam}(\widetilde\Theta)<\infty$. Therefore, we have that
\begin{align*}
(2LC_H)^2\delta^2+C_{\rm slab}\delta&\le \{(2LC_H)^2\operatorname{diam}(\widetilde\Theta)+C_{\rm slab}\}\delta := C_{\rm tr}\delta.
\end{align*}
Taking square roots gives
\begin{equation*}
\|p_{\theta,q}-p_{\theta',q'}\|_{L_2(P)}\le C_{\rm tr}\delta^{1/2}=C_{\rm tr}\|\theta-\theta'\|_2^{1/2}.
\end{equation*}
\end{proof}

\begin{lemma}[Uniform trace entropy]
\label[lemma]{lem:trace_entropy_weak}
Suppose \Cref{a_theta} holds. Then there exists a constant $C_Q<\infty$, depending only on $d,\nu,L,\widetilde A_0,\widetilde B_0,\widetilde c_0,\widetilde C_0$ and $\widetilde\rho_0$, such that, for every fixed $\theta\in\widetilde\Theta$ and every $\varepsilon\in(0,1)$,
\[
\log N(\varepsilon,\mathcal Q_{\theta,L},\|\cdot\|_\infty)\le C_Q\varepsilon^{-(d-1)}\log(C_Q/\varepsilon).
\]
\end{lemma}

\begin{proof}

Fix $\theta\in\widetilde\Theta$. If $L=0$, then $\mathcal Q_{\theta,0}$ consists of constant functions with values in $[0,1]$, and the result follows by discretising this interval. Suppose henceforth that $L>0$. Since $\mathcal H_\theta\subseteq[0,1]^d$ is contained in the affine hyperplane $\{x:s_\theta(x)=0\}$, it is a bounded subset of a $(d-1)$-dimensional affine space. By the Euclidean covering bound for bounded subsets of finite-dimensional affine spaces \citep[e.g.~Proposition 4.2.12 in][]{vershynin2018high}, for every $r\in(0,1]$, there exists an $(r/C_{d,\nu})$-net $\mathcal N_r=\{h_1,\ldots,h_N\}\subseteq\mathcal H_\theta$ with respect to $\|\cdot\|_2$ such that $ N\le C_Q r^{-(d-1)}$, where $C_Q<\infty$ depends only on $d$ and $\nu$. By \Cref{lem:norm-comparison}, we have that $\|x-y\|_\nu\le C_{d,\nu}\|x-y\|_2$ for all $x,y\in\mathbb R^d$. Hence $\mathcal N_r$ is an $r$-net of $\mathcal H_\theta$ with respect to $\|\cdot\|_\nu$. For $\varepsilon\in(0,1)$, we choose $r=1\wedge\{\varepsilon/(8L)\}$.  

We next discretise $[0,1]$ using a grid $\mathcal G_\varepsilon\subset[0,1]$.  Let $\mathcal G_\varepsilon\subset[0,1]$ be a finite grid with grid size at most $\varepsilon/4$, meaning that for every $y\in[0,1]$, there exists $g\in\mathcal G_\varepsilon$ such that $|y-g|\le \frac{\varepsilon}{4}$. Moreover, we have that $ |\mathcal G_\varepsilon| \leq C/\varepsilon$, where $C<\infty$ is a universal constant.

For each $q\in \mathcal Q_{\theta,L}$, choose $v_j\in\mathcal G_\varepsilon$ such that $|v_j-q(h_j)|\le \varepsilon/4$ for $j=1,\ldots,N$. Define
\[
\widetilde q(h)=\left[\min_{1\le j\le N}\{v_j+L\|h-h_j\|_{\nu}\}\right]_{[0,1]},\qquad h\in\mathcal H_\theta.
\]
Note that $\widetilde q$ is $L$-Lipschitz and takes values in $[0,1]$.

We now show that $\widetilde q$ approximates $q$ uniformly. For any $h\in\mathcal H_\theta$, choose $h_j\in\mathcal N_r$ such that $\|h-h_j\|_{\nu}\le r$. Then we have that
\[
\widetilde q(h)\le v_j+L\|h-h_j\|_{\nu}\le q(h_j)+\frac{\varepsilon}{4}+Lr \le q(h)+\frac{\varepsilon}{2},
\]
where the last inequality follows from $2Lr\le\varepsilon/4$. For the lower bound, for every $j$,
\[
v_j+L\|h-h_j\|_{\nu}\ge q(h_j)-\frac{\varepsilon}{4}+L\|h-h_j\|_{\nu}\ge q(h)-\frac{\varepsilon}{4},
\]
where the last inequality follows since $q$ is $L$-Lipschitz over $\mathcal H_\theta$. Taking the infimum over $j$ and then clipping to $[0,1]$ cannot decrease the value below $q(h)-\varepsilon/4$, since $q(h)\in[0,1]$. Thus we have that
\[
\widetilde q(h)\ge q(h)-\frac{\varepsilon}{4}.
\]
Combining the results, it then holds that $\|\widetilde q-q\|_\infty\le \varepsilon$.

Note that the number of possible choices of $(v_1,\ldots,v_N)$ is at most $(C/\varepsilon)^N$. Since each coordinate $v_j$ is chosen independently from the grid $\mathcal G_\varepsilon$, and $|\mathcal G_\varepsilon|\le C/\varepsilon$, the number of possible vectors $(v_1,\ldots,v_N)$ is at most $|\mathcal G_\varepsilon|^N\le (C/\varepsilon)^N$. Consequently, it holds that 
\[
\log N(\varepsilon, \mathcal Q_{\theta,L},\|\cdot\|_\infty)\le N\log(C/\varepsilon)\le C_Q\varepsilon^{-(d-1)}\log(C_Q/\varepsilon).
\]
\end{proof}

\begin{lemma}[$L_2(P)$ bracketing entropy]
\label[lemma]{lem:l2-bracketing-policy}
Suppose Assumptions \ref{a_contexts}-\ref{a_theta} hold and $d \ge 2$. Let
\[
\mathcal F_L=\{x\mapsto s_{\theta^*}(x)p_{\theta,q}(x):\theta\in\widetilde\Theta,\ q\in\mathcal Q_{\theta,L}\},
\]
where $s_{\theta^*}(x)=\widetilde x^\top\theta^*$. Then there exists a constant $C_{\mathcal F}<\infty$ such that, for every $\varepsilon\in(0,1)$,
\[
\log N_{[ ]}(\varepsilon,\mathcal F_L,L_2(P))\le C_{\mathcal F}\varepsilon^{-(d-1)}\log(C_{\mathcal F}/\varepsilon),
\]
where $N_{[]}(\varepsilon,\mathcal F_L,L_2(P))$ denotes the $L_2(P)$ bracketing number.
\end{lemma}

\begin{proof}
Write $m=d-1$ and define $B_*=\sup_{x\in[0,1]^d}|s_{\theta^*}(x)|$. By \Cref{a_theta}, it holds that $B_*<\infty$. If $B_*=0$, then every function in $\mathcal F_L$ is identically zero and the claim is trivial. Suppose $B_*>0$ and set $\bar B=1+B_*$.

For each $r\in\{1,\ldots,d\}$, define
\[
\mathcal F_{L,r}=\{x\mapsto s_{\theta^*}(x)p_{\theta,q}(x):\theta\in\widetilde\Theta_r,\ q\in\mathcal Q_{\theta,L}\}.
\]
Then $\mathcal F_L=\bigcup_{r=1}^d\mathcal F_{L,r}$. We first bound the bracketing number of each fixed chart $\mathcal F_{L,r}$.

Fix $r\in\{1,\ldots,d\}$ and $\varepsilon\in(0,1)$. Let $\delta_\theta=c_1\varepsilon^2$, where $c_1>0$ will be chosen sufficiently small. Let $\{\theta_1,\ldots,\theta_M\}$ be a $\delta_\theta$-net of $\widetilde\Theta_r$ with respect to $\|\cdot\|_2$. Since $\widetilde\Theta_r\subset\mathbb R^{d+1}$ is bounded by the definition of $\widetilde\Theta_r$, we have that
\[
\log M\le (d+1)\log(C'/\delta_\theta)\le C\log(C/\varepsilon),
\]
where $C$ is a constant depending on $d$.

For each fixed representative boundary $\mathcal H_{\theta_j}$, \Cref{lem:trace_entropy_weak} gives a finite $c_2\varepsilon/\bar B$-cover of $Q_{\theta_j,L}$ in the sup-norm, where $c_2>0$ is chosen sufficiently small. Thus there exists a finite set $\mathcal Q_j\subset\mathcal Q_{\theta_j,L}$ such that, for every $q^\circ\in\mathcal Q_{\theta_j,L}$, there exists $\widetilde q\in\mathcal Q_j$ satisfying $\|q^\circ-\widetilde q\|_\infty\le c_2\varepsilon/\bar B$. Moreover, it holds that $\log|\mathcal Q_j|\le C\varepsilon^{-m}\log(C/\varepsilon)$.

Now fix any $f\in\mathcal F_{L,r}$, and we can write $f(x)=s_{\theta^*}(x)p_{\theta,q}(x)$, for some $\theta\in\widetilde\Theta_r$ and $q\in\mathcal Q_{\theta,L}$. Choose~$\theta_j$~such that $\|\theta-\theta_j\|_2\le\delta_\theta$. Since $\theta$ and $\theta_j$ lie in the same chart $\widetilde\Theta_r$, by the construction in \textbf{Steps 2,3 and 4} of the proof of \Cref{lem:l2_transfer_slab}, there exists $q^\circ\in\mathcal Q_{\theta_j,L}$ such that, for all $x \in [0,1]^d$,
\begin{align}\label{lem:l2-bracketing-policy_eq1}
    |p_{\theta,q}(x)-p_{\theta_j,q^\circ}(x)| \le 2LC_H\delta_\theta + \indc_{S(\theta, \theta_j)}(x) \leq 2LC_H\delta_\theta + \indc_{A_j}(x),
\end{align}
where
\[
S(\theta,\theta_j)=\{x:s_\theta(x)\ge0,\ s_{\theta_j}(x)<0\}\cup\{x:s_\theta(x)<0,\ s_{\theta_j}(x)\ge0\}
\]
is the disagreement set,
\[
A_j = \{x: |s_{\theta_j}(x)|\le \sqrt{d+1}\delta_{\theta}\},
\]
and the last inequality follows as $S(\theta_j, \theta) \subseteq A_j$ by  \eqref{lem:slab_probability_eq1} in the proof of \Cref{lem:slab_probability}.

Choose $\widetilde q\in\mathcal Q_j$ such that $\|q^\circ-\widetilde q\|_\infty\le c_2\varepsilon/\bar B$. Since the envelope map is nonexpansive in the boundary trace, we have that
\begin{align}\label{lem:l2-bracketing-policy_eq2}
    \|p_{\theta_j,q^\circ}-p_{\theta_j,\widetilde q}\|_\infty\le c_2\varepsilon/\bar B.
\end{align}
Therefore, by \eqref{lem:l2-bracketing-policy_eq1}, \eqref{lem:l2-bracketing-policy_eq2} and the triangle inequality, for every $x\in[0,1]^d$,
\[
|p_{\theta,q}(x)-p_{\theta_j,\widetilde q}(x)|\le C\delta_\theta+c_2\varepsilon/\bar B+\indc_{A_j}(x).
\]
Since $|s_{\theta^*}(x)|\le B_*\le\bar B$, this implies
\[
|s_{\theta^*}(x)\{p_{\theta,q}(x)-p_{\theta_j,\widetilde q}(x)\}|\le \bar B\{C\delta_\theta+c_2\varepsilon/\bar B+\indc_{A_j}(x)\}.
\]

For each pair $(j,\widetilde q)$, define $f_{j,\widetilde q}(x)=s_{\theta^*}(x)p_{\theta_j,\widetilde q}(x)$.

Define the bracket
\[
\ell_{j,\widetilde q}(x)=f_{j,\widetilde q}(x)-\bar B\{C\delta_\theta+c_2\varepsilon/\bar B+\indc_{A_j}(x)\},
\]
and
\[
u_{j,\widetilde q}(x)=f_{j,\widetilde q}(x)+\bar B\{C\delta_\theta+c_2\varepsilon/\bar B+\indc_{A_j}(x)\}.
\]
We have that
\[
\ell_{j,\widetilde q}(x)\le f(x)\le u_{j,\widetilde q}(x)
\]
for every $x\in[0,1]^d$. Moreover, the $L_2(P)$ width of this bracket satisfies
\[
\|u_{j,\widetilde q}-\ell_{j,\widetilde q}\|_{L_2(P)}\le C\bar B\{\delta_\theta+\varepsilon/\bar B+\mathbb{P}(A_j)^{1/2}\}\leq C\bar B\{\delta_\theta+\varepsilon/\bar B + C^{1/2}\delta_\theta^{1/2}\},
\]
where the last inequality follows as
$\mathbb{P}(A_j) \leq C\delta_\theta$ by the same calculations as \eqref{lem:slab_probability_eq2} in the proof of \Cref{lem:slab_probability}.

Since $\delta_\theta=c_1\varepsilon^2$, choosing $c_1,c_2>0$ sufficiently small gives $\|u_{j,\widetilde q}-\ell_{j,\widetilde q}\|_{L_2(P)}\le\varepsilon$.

Thus the brackets indexed by $(j,\widetilde q)$ form an $L_2(P)$ $\varepsilon$-bracketing cover of $\mathcal F_{L,r}$, where the number of brackets is at most $M\max_{1\le j\le M}|\mathcal Q_j|$. Therefore, for every fixed $r\in\{1,\ldots,d\}$, we have that
\[
\log N_{[]}(\varepsilon,\mathcal F_{L,r},L_2(P))\le C\log(C/\varepsilon)+C\varepsilon^{-m}\log(C/\varepsilon)\le C_{\mathcal F}\varepsilon^{-m}\log(C_{\mathcal F}/\varepsilon).
\]
Since $\mathcal F_L=\bigcup_{r=1}^d\mathcal F_{L,r}$, the union of the chart-wise bracketing covers is an $\varepsilon$-bracketing cover of $\mathcal F_L$. Hence,
\[
N_{[]}(\varepsilon,\mathcal F_L,L_2(P))\le\sum_{r=1}^d N_{[]}(\varepsilon,\mathcal F_{L,r},L_2(P)).
\]
Taking logarithms and absorbing $\log(d)$ into the constant prove the claim.
\end{proof}

\begin{lemma}[Control of $\Delta_n^L$]
\label[lemma]{lem:control-delta-L-n}
Suppose Assumptions \ref{a_contexts}-\ref{a_theta} hold. Let
\[
\mathcal F_L=\{x\mapsto s_{\theta^*}(x)p_{\theta,q}(x):\theta\in\widetilde\Theta,\ q\in\mathcal Q_{\theta,L}\},
\]
where $s_{\theta^*}(x)=\widetilde x^\top\theta^*$. Define
\[
\Delta_n^L=\sup_{\theta\in\widetilde\Theta,\ q\in\mathcal Q_{\theta,L}}\left|\frac1n\sum_{i=1}^n s_{\theta^*}(X_i)p_{\theta,q}(X_i)-\mathbb E\{s_{\theta^*}(X)p_{\theta,q}(X)\}\right|.
\]
Then there exists a constant $C<\infty$, depending only on $d,L,\widetilde A_0,\widetilde B_0,\widetilde c_0, \widetilde C_0, \widetilde \rho_0,c_g,C_g$, such that, for every $n\ge2$ and every $\eta\in(0,1)$, with probability at least $1-\eta$,
\[
\Delta_n^L\le C\rho_{n,d-1}(\eta),
\]
where, for $m=d-1$,
\[
\rho_{n,m}(\eta)=
\begin{cases}
n^{-1/2}\log(en/\eta), & m=1,\\
n^{-1/2}\log^{3/2}(en/\eta), & m=2,\\
n^{-1/m}\log(en/\eta), & m\ge3.
\end{cases}
\]
\end{lemma}

\begin{proof}
Write $m=d-1$ and $\theta^*=(a^*, \gamma^{*\top})^{\top}$. For every $f\in\mathcal F_L$, we can rewrite it as $f(x)=s_{\theta^*}(x)p_{\theta,q}(x)$ for some $\theta\in\widetilde\Theta$ and $q\in\mathcal Q_{\theta,L}$. Since $0\le p_{\theta,q}\le1$ and $\theta^*\in\Theta$, \Cref{a_theta} gives
\[
|s_{\theta^*}(x)|\le |a^*|+\|\gamma^*\|_2\|x\|_2\le A_0+B_0\sqrt d=:B.
\]
Therefore, we have that  $\|f\|_\infty \le \|s_{\theta^*}\|_{\infty} \leq  B$ for every $f\in\mathcal F_L$.

By the bracketing maximal inequality for empirical processes \citep[e.g.~Theorem 2.14.2 in][]{van1996weak}, we have that
\[
\mathbb E\Delta_n^L\le C_1\inf_{\alpha\in(0,1)}\left\{\alpha+\frac1{\sqrt n}\int_\alpha^1\sqrt{1+\log N_{[]}(\varepsilon,\mathcal F_L,L_2(P))}\,\mathrm{d}\varepsilon\right\}.
\]
Moreover, by \Cref{lem:l2-bracketing-policy}, we have that
\[
\log N_{[]}(\varepsilon,\mathcal F_L,L_2(P))\le C_2\varepsilon^{-m}\log(C_2/\varepsilon).
\]
Therefore, we have that
\[
\mathbb E\Delta_n^L\le C_3\inf_{\alpha\in(0,1)}\left\{\alpha+\frac1{\sqrt n}\int_\alpha^1\varepsilon^{-m/2}\{\log(C_2/\varepsilon)\}^{1/2}\,\mathrm{d}\varepsilon\right\}.
\]

\noindent \textbf{Case 1: $m=1$.} In this case, 
\begin{align*}
\int_\alpha^1\varepsilon^{-1/2}\{\log(C_2/\varepsilon)\}^{1/2}\,\mathrm{d}\varepsilon
&\le \int_0^1\varepsilon^{-1/2}\{\log(C_2/\varepsilon)\}^{1/2}\,\mathrm{d}\varepsilon \\
&\le \sqrt{\log C_2}\int_0^1\varepsilon^{-1/2}\,\mathrm{d}\varepsilon+\int_0^1\varepsilon^{-1/2}\{\log(1/\varepsilon)\}^{1/2}\,\mathrm{d}\varepsilon \\
&=2\sqrt{\log C_2}+\int_0^\infty u^{1/2}e^{-u/2}\,\mathrm{d}u =:C_4
\end{align*}
Therefore, taking $\alpha=n^{-1/2}$ gives $\mathbb E\Delta_n^L\le C_5n^{-1/2}$.

\noindent \textbf{Case 2: $m=2$.} In this case, we have that
\[
\int_\alpha^1\varepsilon^{-1}\{\log(C_2/\varepsilon)\}^{1/2}\,d\varepsilon\le C_6\{\log(1/\alpha)\}^{3/2}.
\]
Taking $\alpha=n^{-1/2}$ gives $\mathbb E\Delta_n^L\le C_7n^{-1/2}(\log n)^{3/2}$.

\noindent \textbf{Case 3: $m\geq 3$.} In this case, we have that
\[
\int_\alpha^1\varepsilon^{-m/2}\{\log(C_2/\varepsilon)\}^{1/2}\,d\varepsilon\le C_8\alpha^{1-m/2}\{\log(1/\alpha)\}^{1/2}.
\]
Taking $\alpha=n^{-1/m}(\log n)^{1/m}$ gives $\mathbb E\Delta_n^L\le C_9n^{-1/m}(\log n)^{1/m}$.

Thus, for all $m\ge1$, we have that
\[
\mathbb E\Delta_n^L\le C_{10}
\begin{cases}
n^{-1/2}, & m=1,\\
n^{-1/2}(\log n)^{3/2}, & m=2,\\
n^{-1/m}(\log n)^{1/m}, & m\ge3.
\end{cases}
\]

Since $\|f\|_\infty\le B$ for all $f\in\mathcal F_L$, the Bennett concentration inequality for bounded empirical processes \citep[e.g.~Theorem 2.3 in][]{bousquet2002bennett} implies that, with probability at least $1-\eta$,
\[
\Delta_n^L\le 2\mathbb E\Delta_n^L+C B\sqrt{\frac{\log(1/\eta)}{n}}+C B\frac{\log(1/\eta)}{n}.
\]
Since $\log(en/\eta)\ge1$ and taking $m=d-1$, the lemma follows.
\end{proof}

\subsubsection{Performance of estimators}
\begin{lemma}[Initial estimation error after exploration]\label[lemma]{lem:initial-exploration-estimation}
Suppose Assumptions \ref{a_contexts}-\ref{a_theta} hold. During the first $2T_o$ rounds, use the constant policy $p_t(\cdot)=1/2$, with independent arm assignments across rounds and independent of the contexts and reward noise. For $\eta\in(0,1/2)$, define $\widehat\beta_{k,2T_o}=V_{k,2T_o}^{-1}b_{k,2T_o}$. Suppose there exist constants $C_{\rm init},C_{\rm T}<\infty$, depending only on $d,c_g,\sigma,\lambda$ and $C_\beta$, such that $T_o\ge C_{\rm T}\log\{(d+1)/\eta\}$. Then, taking $\lambda=C_\lambda$, where $C_\lambda>0$ is an absolute constant, it holds with probability at least $1-\eta$ that
\[
\|\widehat\beta_{k,2T_o}-\beta_k^*\|_2\le C_{\rm init}\sqrt{\frac{\log(T/\eta)}{T_o}},\qquad k\in\{1,2\}.
\]
Consequently, if $s_0$ is the first lazy update time, then we have with probability at least $1-\eta$ that
\[
\|\widehat\theta_{s_0}-\theta^*\|_2\le 2C_{\rm init}\sqrt{\frac{\log(T/\eta)}{T_o}}.
\]
\end{lemma}

\begin{proof}
For $k\in\{1,2\}$, let $I_k^0=\{i\le2T_o:A_i=k\}$ and $N_{0,k}=|I_k^0|$. Let $\lambda_X=c_g/(128d^2)$. By \Cref{l_covariate_diverse} with $A=[0,1]^d$, we have $\lambda_{\min}\{\mathbb E(\widetilde X\widetilde X^\top)\}\ge\lambda_X$. Define
\begin{align*}
\mathcal E_{\rm alloc}&=\{N_{0,1}\ge T_o/2,\ N_{0,2}\ge T_o/2\},\\
\mathcal E_{\rm init}&=\{\lambda_{\min}(V_{k,2T_o})\ge\lambda+\lambda_XT_o/4,\ \text{for every }k\in\{1,2\}\}.
\end{align*}
Under $p_t(\cdot)=1/2$, the exploration actions are independent of the contexts and reward noise, and $N_{0,k}\sim\operatorname{Bin}(2T_o,1/2)$ for each $k\in\{1,2\}$. By Hoeffding's inequality \citep[e.g.~Theorem~2.2.6 in][]{vershynin2018high} and a union bound over $k\in\{1,2\}$, we have that
\[
\mathbb P(\mathcal E_{\rm alloc}^c)\le2\exp\left(-\frac{T_o}{4}\right).
\]

Fix $k\in\{1,2\}$ and condition on the exploration actions $A_1,\ldots,A_{2T_o}$. Since these actions are sampled independently of the contexts, the matrices $\{\widetilde X_i\widetilde X_i^\top:i\in I_k^0\}$ remain independent and each has the same distribution as $\widetilde X\widetilde X^\top$. Moreover, $\widetilde X_i\widetilde X_i^\top\succeq0$ and $\lambda_{\max}(\widetilde X_i\widetilde X_i^\top)=\|\widetilde X_i\|_2^2\le d+1$. Therefore, on $\mathcal E_{\rm alloc}$, we have that
\[
\mu_{\min}:=\lambda_{\min}\left\{\sum_{i\in I_k^0}\mathbb E(\widetilde X_i\widetilde X_i^\top\mid A_1,\ldots,A_{2T_o})\right\}\ge N_{0,k}\lambda_X\ge\frac{\lambda_XT_o}{2}.
\]
By the matrix Chernoff inequality \citep[e.g.~Theorem~5.1.1 in][]{tropp2015introduction}, applied with matrix dimension $d+1$, upper eigenvalue bound $d+1$ and $\varepsilon=1/2$, it holds that
\begin{align*}
&\mathbb P\left\{\lambda_{\min}\left(\sum_{i\in I_k^0}\widetilde X_i\widetilde X_i^\top\right)\le\frac{\mu_{\min}}{2}\ \middle|\ A_1,\ldots,A_{2T_o}\right\}\\
&\qquad\le(d+1)\left(\frac{e^{-1/2}}{(1/2)^{1/2}}\right)^{\mu_{\min}/(d+1)}\le(d+1)\exp\left\{-\frac{\mu_{\min}}{8(d+1)}\right\}.
\end{align*}
Since $\mu_{\min}\ge\lambda_XT_o/2$ on $\mathcal E_{\rm alloc}$ and $V_{k,2T_o}=\lambda I_{d+1}+\sum_{i\in I_k^0}\widetilde X_i\widetilde X_i^\top$, a union bound over $k\in\{1,2\}$ and the tower property give
\[
\mathbb P(\mathcal E_{\rm init}^c\cap\mathcal E_{\rm alloc})\le2(d+1)\exp\left\{-\frac{\lambda_XT_o}{16(d+1)}\right\}.
\]
Combining the preceding two bounds, we obtain
\[
\mathbb P(\mathcal E_{\rm init}^c)\le2\exp\left(-\frac{T_o}{4}\right)+2(d+1)\exp\left\{-\frac{\lambda_XT_o}{16(d+1)}\right\}\le\frac{\eta}{2},
\]
where the last inequality follows whenever $T_o\ge4\log(8/\eta)$ and $T_o\ge\{16(d+1)/\lambda_X\}\log\{8(d+1)/\eta\}$. Both inequalities hold under $T_o\ge C_{\rm T}\log\{(d+1)/\eta\}$ for $C_{\rm T}$ sufficiently large. Consequently, with probability at least $1-\eta/2$, we have that
\begin{equation}\label{lem:conditional-min-eigen-V_eq1}
\lambda_{\min}(V_{k,2T_o})\ge\lambda+\frac{\lambda_XT_o}{4},\qquad k\in\{1,2\}.
\end{equation}

Next, for $k\in\{1,2\}$, write $S_k^0=\sum_{i\in I_k^0}\widetilde X_i\xi_{i,k}$. For $j=1,\ldots,d+1$, let $S_{k,j}^0=\sum_{i\in I_k^0}\widetilde X_{i,j}\xi_{i,k}$ denote the $j$th coordinate of $S_k^0$. Conditional on the exploration actions and contexts, the noise variables remain independent, mean-zero and $\sigma$-sub-Gaussian by \Cref{a_noise}. Since $|\widetilde X_{i,j}|\le1$ and $N_{0,k}\le2T_o$, the general Hoeffding inequality for sums of independent sub-Gaussian random variables \citep[e.g.~Theorem~2.6.3 in][]{vershynin2018high} gives, for every $u>0$,
\[
\mathbb P\left\{|S_{k,j}^0|\ge u\ \middle|\ A_1,\ldots,A_{2T_o},X_1,\ldots,X_{2T_o}\right\}\le2\exp\left\{-\frac{cu^2}{\sigma^2T_o}\right\},
\]
where $c>0$ is an absolute constant. Taking expectation over the actions and contexts and a union bound over $k\in\{1,2\}$ and $j=1,\ldots,d+1$, we obtain
\[
\mathbb P\left\{\max_{k\in\{1,2\}}\max_{1\le j\le d+1}|S_{k,j}^0|\ge u\right\}\le4(d+1)\exp\left\{-\frac{cu^2}{\sigma^2T_o}\right\}.
\]
Choosing $u=C_1\sqrt{T_o\log(T/\eta)}$ with $C_1$ sufficiently large, we have with probability at least $1-\eta/2$ that
\[
\|S_k^0\|_2\le\sqrt{d+1}\max_{1\le j\le d+1}|S_{k,j}^0|\le C_2\sqrt{T_o\log(T/\eta)},\qquad k\in\{1,2\}.
\]
Combining this noise bound with \eqref{lem:conditional-min-eigen-V_eq1}, with probability at least $1-\eta$, both bounds hold. We work under these high-probability events in the rest of the proof.

Since $Y_i=\widetilde X_i^\top\beta_k^*+\xi_{i,k}$ for $i\in I_k^0$, the ridge estimator satisfies
\[
\widehat\beta_{k,2T_o}-\beta_k^*=V_{k,2T_o}^{-1}S_k^0-\lambda V_{k,2T_o}^{-1}\beta_k^*.
\]
Therefore, for any $k\in\{1,2\}$, we have that
\begin{align*}
\|\widehat\beta_{k,2T_o}-\beta_k^*\|_2&\le\left\|V_{k,2T_o}^{-1}S_k^0\right\|_2+\lambda\left\|V_{k,2T_o}^{-1}\beta_k^*\right\|_2\\
&\le\left\|V_{k,2T_o}^{-1}\right\|_{\op}\left(\|S_k^0\|_2+\lambda\|\beta_k^*\|_2\right)\\
&\le\frac{4C_2\sqrt{T_o\log(T/\eta)}}{\lambda_XT_o}+\frac{4C_\lambda C_\beta}{\lambda_XT_o}\\
&\le C_{\rm init}\sqrt{\frac{\log(T/\eta)}{T_o}},
\end{align*}
where the last inequality follows by taking $T_o\ge C_{\rm T}\log\{(d+1)/\eta\}$ and enlarging $C_{\rm init}$ if necessary.

Finally, at the first lazy update time $s_0$, the algorithm uses the estimator based on the first $2T_o$ rounds. Hence, on the same event, we have that
\[
\|\widehat\theta_{s_0}-\theta^*\|_2\le\|\widehat\beta_{1,2T_o}-\beta_1^*\|_2+\|\widehat\beta_{2,2T_o}-\beta_2^*\|_2\le2C_{\rm init}\sqrt{\frac{\log(T/\eta)}{T_o}}.
\]
This proves the claim.
\end{proof}

\begin{lemma}[Estimation at lazy updates]\label[lemma]{lem:localization-and-estimation}
Suppose Assumptions \ref{a_contexts}-\ref{a_theta} hold.  Let $\mathcal S_T$ be the set of lazy update times up to time $T$, and let $s_0$ be the first lazy update time. Suppose that there exists a constant $C_{\rm loc}<\infty$, depending only on $d,L,c_g,\sigma,\lambda,C_\beta$ and the constants in \Cref{a_theta}, such that 
\[T_o\ge C_T\log\{(d+1)/\eta\}\vee \frac{C_{\rm loc}\log(T/\eta)}{r_{\rm loc}^2}.\]
Then, it holds with probability at least $1-\eta$ that
\[
\widehat\theta_s\in\widetilde\Theta_{j^*}\qquad\text{for every }s\in\mathcal S_T.
\]
Moreover, with probability at least $1-\eta$, for every $s\in\mathcal S_T$ and $k\in\{1,2\}$,
\[
\|\widehat\beta_{k,s-1}-\beta_k^*\|_2\le C_{\rm est}\sqrt{\frac{\log(T/\eta)}{s-1}},
\]
where $C_{\rm est}<\infty$ depends only on $d,L,c_g,\sigma,\lambda,C_\beta$ and the constants in \Cref{a_theta}. Moreover, we have that 
\[
\|\widehat\theta_s-\theta^*\|_2 \le 2C_{\rm est}\sqrt{\frac{\log(T/\eta)}{s-1}}.
\]
\end{lemma}

\begin{proof}
By \Cref{lem:initial-exploration-estimation}, if
\[
T_o\ge C_T\log\{(d+1)/\eta\}\vee \frac{C_{\rm loc}\log(T/\eta)}{r_{\rm loc}^2},
\]
then, with probability at least $1-\eta/4$, we have that
\[
\|\widehat\theta_{s_0}-\theta^*\|_2\le r_{\rm loc}.
\]

Moreover, by the same allocation and conditional matrix Chernoff arguments as in the proof of \Cref{lem:initial-exploration-estimation}, it also holds with probability at least $1-\eta/4$ that
\[
\lambda_{\min}(V_{k,2T_o})\ge\lambda+\frac{\lambda_XT_o}{4},\qquad k\in\{1,2\}.
\]

Now consider the post-initialization rounds. For $i>2T_o$, define $Z_{i,k}=\widetilde X_i\widetilde X_i^\top\mathbf 1\{A_i=k\}$. Then, we have that $Z_{i,k}\succeq0$ and $\lambda_{\max}(Z_{i,k})\le d+1$. Fix $t\ge4T_o$ and $k\in\{1,2\}$. Suppose that $\widehat\theta_r\in\widetilde\Theta$ for every $r\in\mathcal S_T$ with $r\le t$. For each $i=2T_o+1,\ldots,t$, let $\tau(i)$ be the lazy update time at which the policy used at round $i$ is computed. Then $\tau(i)\in\mathcal S_T$ and $\tau(i)\le i\le t$. Hence by the assumption that $\widehat\theta_r\in\widetilde\Theta$ for every $r\in\mathcal S_T$ with $r\le t$, it holds that $\widehat\theta_{\tau(i)}\in\widetilde\Theta$. Therefore, the policy used at round $i$ is generated by a parameter in $\widetilde\Theta$, and by \Cref{lem:algorithm-policy-diversity}, we have that
\[
\mathbb{E}[Z_{i,k}\mid\mathcal F_{i-1}]\succeq\lambda_{\rm div}I_{d+1}.
\]
Therefore, we have
\[
\lambda_{\min}\left(\sum_{i=2T_o+1}^t\mathbb{E}[Z_{i,k}\mid\mathcal F_{i-1}]\right)\ge\lambda_{\rm div}(t-2T_o).
\]
Applying the matrix concentration inequality \citep[e.g.~Lemma~9 in][with $R=d+1$, $\gamma=1/2$ and $\mu=\lambda_{\rm div}(t-2T_o)$]{bastani2021mostly}  gives
\begin{align*}
&\mathbb{P}\left\{
\lambda_{\min}\left(\sum_{i=2T_o+1}^tZ_{i,k}\right)\le\frac{\lambda_{\rm div}}{2}(t-2T_o),\
\lambda_{\min}\left(\sum_{i=2T_o+1}^t\mathbb{E}[Z_{i,k}\mid\mathcal F_{i-1}]\right)\ge\lambda_{\rm div}(t-2T_o)
\right\}\\
&\qquad\le(d+1)\left(\frac{e^{-1/2}}{(1/2)^{1/2}}\right)^{\lambda_{\rm div}(t-2T_o)/(d+1)}\le C_1(d+1)\exp\left\{-c\frac{\lambda_{\rm div}(t-2T_o)}{d+1}\right\}.
\end{align*}

Taking a union bound over $t=4T_o,\ldots,T$ and $k\in\{1,2\}$, then it holds with probability at least $1-\eta/4$ that if $\widehat\theta_r\in\widetilde\Theta$ for every $r\in\mathcal S_T$ with $r\le t$,
\[
\lambda_{\min}\left(\sum_{i=2T_o+1}^tZ_{i,k}\right)\ge\frac{\lambda_{\rm div}}{2}(t-2T_o),
\]
whenever $T_o\ge C_2\log\{(d+1)T/\eta\}$. Combining the preceding forced-exploration and post-initialization bounds, with probability at least $1-\eta/2$, it holds for every $t\in\{2T_o,\ldots,T\}$ and $k\in\{1,2\}$ that, if $\widehat\theta_r\in\widetilde\Theta$ for every $r\in\mathcal S_T$ with $r\le t$, then
\[
\lambda_{\min}(V_{k,t})\ge c_Vt.
\]

Next we control the noise term uniformly over time. This step requires a martingale argument because, after the forced-exploration phase, the action $A_i$ is adaptive. For $k\in\{1,2\}$ and $j=1,\ldots,d+1$, define
\[
M_{t,k,j}=\sum_{i=1}^t \widetilde X_{i,j}\xi_{i,k}\mathbf 1\{A_i=k\}.
\]
Let $\mathcal G_i$ denote the sigma-field generated by the history before round $i$, the current covariate $X_i$, the randomisation used by the policy at round $i$, and the chosen action $A_i$. In particular, conditional on $\mathcal G_i$, the quantity $\widetilde X_{i,j}\mathbf 1\{A_i=k\}$ is fixed and satisfies $|\widetilde X_{i,j}\mathbf 1\{A_i=k\}|\le1$. Since the action is chosen before observing the reward noise, and since $\xi_{i,k}$ is conditionally mean-zero and $\sigma$-sub-Gaussian, for every $u\in\mathbb R$, we have that
\[
\mathbb E\left[\exp\left\{u\widetilde X_{i,j}\xi_{i,k}\mathbf 1\{A_i=k\}\right\}\mid \mathcal G_i\right]\le \exp\left\{\frac{\sigma^2u^2}{2}\right\}.
\]
Therefore, by the tower property, for every fixed $t$, $k$ and $j$, the sum $M_{t,k,j}$ is sub-Gaussian with variance proxy $\sigma^2t$. 
Hence, for every $x>0$, we have that
\[
\mathbb P\{|M_{t,k,j}|\ge x\}\le 2\exp\left\{-\frac{x^2}{2\sigma^2t}\right\}.
\]
Taking
\[
x_t=\sigma\sqrt{2t\log\left(\frac{16(d+1)T}{\eta}\right)},
\]
and applying a union bound over $t=1,\ldots,T$, $k\in\{1,2\}$ and $j=1,\ldots,d+1$, we obtain
\[
\mathbb P\left\{\max_{k\in\{1,2\}}\max_{1\le j\le d+1}|M_{t,k,j}|\leq \sigma\sqrt{2t\log\left(\frac{16(d+1)T}{\eta}\right)}, \text{ for all  } t\in [T]\right\}\ge 1- \frac{\eta}{4}.
\]
Consequently, with probability at least $1-\eta/4$, for every $t=1,\ldots,T$ and $k\in\{1,2\}$,
\[
\left\|\sum_{i=1}^t \widetilde X_i\xi_{i,k}\mathbf 1\{A_i=k\}\right\|_2
\le \sqrt{d+1}\max_{1\le j\le d+1}|M_{t,k,j}|
\le \sigma\sqrt{2(d+1)t\log\left(\frac{16(d+1)T}{\eta}\right)}.
\]
This implies that,  with probability at least $1-\eta/4$, for every $t=1,\ldots,T$ and $k\in\{1,2\}$, we have that
\[
\left\|\sum_{i=1}^t \widetilde X_i\xi_{i,k}\mathbf 1\{A_i=k\}\right\|_2
\le C_\xi\sqrt{t\log(T/\eta)}.
\]
By a union bound argument, we work under these high-probability events in the rest of the proof.

\textbf{We prove by induction over the lazy update times that $\widehat\theta_s\in\widetilde\Theta_{j^*}$.} For the first lazy update time $s_0$, we have $\|\widehat\theta_{s_0}-\theta^*\|_2\le r_{\rm loc}$. By \Cref{lem:enlarged_theta}, we have that $\{\theta:\|\theta-\theta^*\|_2\le2r_{\rm loc}\}\subseteq\widetilde\Theta_{j^*}$. Hence $\widehat\theta_{s_0}\in\widetilde\Theta_{j^*}$. 

Now fix a lazy update time $s>s_0$ and suppose that $\widehat\theta_r\in\widetilde\Theta_{j^*}$ for every previous lazy update time $r<s$. Since $\widetilde\Theta_{j^*}\subseteq\widetilde\Theta$, the induction hypothesis gives that $\widehat\theta_r\in\widetilde\Theta$ for every $r\in\mathcal S_T$ with $r\le s-1$. Therefore, applying the eigenvalue bound with $t=s-1$ gives
\[
\lambda_{\min}(V_{k,s-1})\ge c_V(s-1),\qquad k\in\{1,2\}.
\]
For $k\in\{1,2\}$, write $S_{k,s-1}=\sum_{i=1}^{s-1}\widetilde X_i\xi_{i,k}\mathbf 1\{A_i=k\}$. As in the proof of \Cref{lem:initial-exploration-estimation}, the ridge estimator satisfies
\[
\widehat\beta_{k,s-1}-\beta_k^*=V_{k,s-1}^{-1}S_{k,s-1}-\lambda V_{k,s-1}^{-1}\beta_k^*.
\]
Therefore,
\begin{align*}
\|\widehat\beta_{k,s-1}-\beta_k^*\|_2
&\le \left\|V_{k,s-1}^{-1}S_{k,s-1}\right\|_2+\lambda\left\|V_{k,s-1}^{-1}\beta_k^*\right\|_2\\
&\le \left\|V_{k,s-1}^{-1}\right\|_{\op}\|S_{k,s-1}\|_2+\lambda\left\|V_{k,s-1}^{-1}\right\|_{\op}\|\beta_k^*\|_2\\
&\le \frac{C_\xi\sqrt{(s-1)\log(T/\eta)}}{c_V(s-1)}+\frac{C_\lambda C_\beta}{c_V(s-1)}\\
&\le C_{\rm est}\sqrt{\frac{\log(T/\eta)}{s-1}}.
\end{align*}
Hence,
\[
\|\widehat\theta_s-\theta^*\|_2\le \|\widehat\beta_{1,s-1}-\beta_1^*\|_2+\|\widehat\beta_{2,s-1}-\beta_2^*\|_2\le 2C_{\rm est}\sqrt{\frac{\log(T/\eta)}{s-1}}.
\]
Since $s-1\ge2T_o$, we have
\[
\|\widehat\theta_s-\theta^*\|_2\le 2C_{\rm est}\sqrt{\frac{\log(T/\eta)}{2T_o}}.
\]
By the lower bound on $T_o$, the final term is at most $r_{\rm loc}$. Therefore $\|\widehat\theta_s-\theta^*\|_2\le r_{\rm loc}$. Again using $\{\theta:\|\theta-\theta^*\|_2\le2r_{\rm loc}\}\subseteq\widetilde\Theta_{j^*}$, this implies $\widehat\theta_s\in\widetilde\Theta_{j^*}$.

Therefore, by induction, $\widehat\theta_s\in\widetilde\Theta_{j^*}$ for every $s\in\mathcal S_T$. The displayed estimation error bound was also proved at each lazy update time. This completes the proof.
\end{proof}

\subsubsection{Auxiliary results for covariate diversity} \label[appendix]{sec_app_covariate_diversity}
In this section, we provide proofs of covariate diversity used in the proof of \Cref{lem:localization-and-estimation}. By the construction of the algorithm, the policy played at round $t$ can be written as $\widehat p_t =p_{\widehat\theta_{\tau(t)},\widehat q_{\tau(t)}}$. 

\textbf{We would like to remark that throughout this subsection, all proofs are constructed under the condition that $\widehat\theta_{\tau(t)}\in\widetilde\Theta$.} These conditional results will then be used in the induction argument in the proof of \Cref{lem:localization-and-estimation}, where this event is shown to hold over all lazy update times on suitable large-probability events.

\begin{lemma}
\label[lemma]{lem:algorithm-policy-diversity}
Suppose Assumptions \ref{a_contexts}-\ref{a_theta} hold. There exists a constant $\lambda_{\rm div}=c_g\varepsilon_0v_0^3/(128d^2)$ depending only on $d,\nu,L,c_g$ and the constants in \Cref{a_theta}, such that, for every $t>2T_o$ satisfying $\widehat\theta_{\tau(t)}\in\widetilde\Theta$, we have that
\[
    \mathbb{E}\left[\widetilde X_t\widetilde X_t^\top\mathbf 1\{A_t=1\}\mid\mathcal F_{t-1}\right]\succeq \lambda_{\rm div}I_{d+1},
\]
and
\[
    \mathbb{E}\left[\widetilde X_t\widetilde X_t^\top\mathbf 1\{A_t=2\}\mid\mathcal F_{t-1}\right]\succeq \lambda_{\rm div}I_{d+1}.
\]
Equivalently, it holds that
\[
    \mathbb{E}\left[\widetilde X_t\widetilde X_t^\top\widehat p_t(X_t)\mid\mathcal F_{t-1}\right]\succeq \lambda_{\rm div}I_{d+1},
    \qquad
    \mathbb{E}\left[\widetilde X_t\widetilde X_t^\top\{1-\widehat p_t(X_t)\}\mid\mathcal F_{t-1}\right]\succeq \lambda_{\rm div}I_{d+1}.
\]
\end{lemma}

\begin{proof}
Fix $t>2T_o$ and write $\widetilde X_t=(1,X_t^\top)^\top$. By the lazy-update structure of \Cref{alg:lazy-adaptive-boundary-greedy}, the policy used at round $t$ is computed at the most recent lazy update time $\tau(t)$. Hence $\widehat p_t(\cdot)=p_{\widehat\theta_{\tau(t)},\widehat q_{\tau(t)}}(\cdot)$.

Since $\widehat\theta_{\tau(t)}\in\widetilde\Theta$ and $\widehat q_{\tau(t)}\in\mathcal Q_{\widehat\theta_{\tau(t)},L}$, \Cref{lem:uniform-interior-region} can be applied to the played policy. Define
\[
    B_t=\left\{x\in[0,1]^d:\varepsilon_0\le \widehat p_t(x)\le1-\varepsilon_0\right\}.
\]
Then $\lambda_d(B_t)\ge v_0$. Moreover, note that $\widehat p_t$ is $\mathcal F_{t-1}$-measurable, because at round $t$ the policy is either inherited from an earlier lazy update or recomputed using only observations up to time~$t-1$.

Since $X_t$ is independent of $\mathcal F_{t-1}$ and has density $g$ satisfying $g(x)\ge c_g$ on $[0,1]^d$, \Cref{l_covariate_diverse} gives
\[
    \mathbb{E}\left[\widetilde X_t\widetilde X_t^\top\mathbf 1\{X_t\in B_t\}\mid\mathcal F_{t-1}\right]\succeq c_g\frac{v_0^3}{128d^2}I_{d+1}.
\]
On $B_t$, we have $\widehat p_t(x)\ge\varepsilon_0$ and $1-\widehat p_t(x)\ge\varepsilon_0$. Therefore, it holds that
\begin{align*}
    \mathbb{E}\left[\widetilde X_t\widetilde X_t^\top\widehat p_t(X_t)\mid\mathcal F_{t-1}\right]
    &\succeq \varepsilon_0\mathbb{E}\left[\widetilde X_t\widetilde X_t^\top\mathbf 1\{X_t\in B_t\}\mid\mathcal F_{t-1}\right]\\
    &\succeq c_g\varepsilon_0\frac{v_0^3}{128d^2}I_{d+1}.
\end{align*}
The same argument gives
\[
    \mathbb{E}\left[\widetilde X_t\widetilde X_t^\top\{1-\widehat p_t(X_t)\}\mid\mathcal F_{t-1}\right]\succeq c_g\varepsilon_0\frac{v_0^3}{128d^2}I_{d+1}.
\]
Finally, by the action-sampling rule in \Cref{alg:lazy-adaptive-boundary-greedy},
\[
    \mathbb{P}(A_t=1\mid X_t,\mathcal F_{t-1})=\widehat p_t(X_t),
    \qquad
    \mathbb{P}(A_t=2\mid X_t,\mathcal F_{t-1})=1-\widehat p_t(X_t).
\]
Thus,
\[
    \mathbb{E}\left[\widetilde X_t\widetilde X_t^\top\mathbf 1\{A_t=1\}\mid\mathcal F_{t-1}\right]
    =
    \mathbb{E}\left[\widetilde X_t\widetilde X_t^\top\widehat p_t(X_t)\mid\mathcal F_{t-1}\right],
\]
and similarly for arm $2$. This proves the claim.
\end{proof}

\begin{lemma}\label[lemma]{lem:uniform-interior-region}
Suppose \Cref{a_theta} holds. Then there exist constants $\varepsilon_0>0$ and $v_0>0$, depending only on $d,\nu,L$ and the constants in \Cref{a_theta}, such that, for every $\theta\in\widetilde\Theta$ and every $q\in\mathcal Q_{\theta,L}$,
\begin{align*}
    \lambda_d\left(\left\{x\in[0,1]^d:\varepsilon_0\le p_{\theta,q}(x)\le 1-\varepsilon_0\right\}\right)\ge v_0,
\end{align*}
where $\lambda_d(\cdot)$ denotes the Lebesgue measure. In particular, the same conclusion holds for any reward-maximizing trace $q_\theta^\star\in\argmax_{q\in\mathcal Q_{\theta,L}}R(p_{\theta,q})$.
\end{lemma}

\begin{proof}
Fix $\theta=(a, \gamma^{\top})^{\top}\in\widetilde\Theta$. Then $\theta\in\widetilde\Theta_j$ for at least one $j\in\{1,\ldots,d\}$. By relabelling coordinates inside this chart, write $x=(u,z)\in[0,1]^{d-1}\times[0,1]$ and $\gamma=(b,c)$. Since $c>0$, define $\phi_\theta(u)=-(a+b^\top u)/c$. Denote
\[
    \mathcal H_\theta=\{(u,z)\in[0,1]^{d-1}\times[0,1]:a+b^\top u+cz=0\}=\{(u,\phi_\theta(u)):u\in[0,1]^{d-1},\ 0\le \phi_\theta(u)\le 1\}.
\]

\noindent\textbf{Step 1: Bounded slope plus one crossing gives a positive-measure interior patch.} Since $\phi_\theta(u)=-(a+b^\top u)/c$, norm equivalence on $\mathbb R^{d-1}$ implies that there is a constant $\kappa_{d,\nu}<\infty$, depending only on $d$ and $\nu$, such that
\[
    |\phi_\theta(u)-\phi_\theta(v)|\le \frac{\|b\|_2}{\widetilde c_0}\|u-v\|_2\le K\|u-v\|_{\nu}.
\]
Here $K=\widetilde B_0\kappa_{d,\nu}/\widetilde c_0$. Thus uniformly over $\theta$, we have $|\phi_\theta(u)-\phi_\theta(v)|\le K\|u-v\|_{\nu}$.

By the definition of $\widetilde\Theta_j$, there is some $u_\theta\in[0,1]^{d-1}$ such that $\widetilde\rho_0\le \phi_\theta(u_\theta)\le 1-\widetilde\rho_0$. Let $D_{\nu}=\sup_{v\in[-1,1]^{d-1}}\|v\|_{\nu}$. Define
\begin{equation} \label{eq:r_0}
     r_0=
    \begin{cases}
    1, & K=0,\\[4pt]
    \displaystyle \min\left\{1,\frac{\widetilde\rho_0}{2KD_{\nu}}\right\}, & K>0.
    \end{cases}
\end{equation}
Let $S_\theta=[0,1]^{d-1}\cap\{u:\|u-u_\theta\|_\infty\le r_0\}$. Then standard calculation gives $\lambda_{d-1}(S_\theta)\ge r_0^{d-1}$. Moreover, for every $u\in S_\theta$, it holds that $\|u-u_\theta\|_{\nu}\le D_{\nu}r_0$. Therefore, it holds that
\[
    |\phi_\theta(u)-\phi_\theta(u_\theta)|\le K\|u-u_{\theta}\|_{\nu}\le KD_{\nu}r_0\le \frac{\widetilde\rho_0}{2}.
\]
Hence, we have that
\[
    \frac{\widetilde\rho_0}{2}\le \phi_\theta(u)\le 1-\frac{\widetilde\rho_0}{2},\qquad \text{for all }u\in S_\theta.
\]
In particular, for every $u\in S_\theta$, the point $h(u)=(u,\phi_\theta(u))$ belongs to $\mathcal H_\theta$. We define $\rho=\widetilde\rho_0/2$.

\medskip

\noindent\textbf{Step 2: Choose a uniform interior probability level.} Let $C_{\rm gr}=1+K$, and we will use this constant to control how vertical distance compares to the $\ell_{\nu}$ distance from the full affine graph of $\phi_\theta$. Since $\mathcal H_\theta$ is contained in this full graph, for $x=(u,\phi_\theta(u)+s)$, for any $h=(v,\phi_\theta(v))\in\mathcal H_\theta$, we have that
\begin{align} \notag
s&=\phi_\theta(u)+s-\phi_\theta(u)\\ \notag
&\le |\phi_\theta(u)+s-\phi_\theta(v)|+|\phi_\theta(v)-\phi_\theta(u)|\\ \notag
&\le |\phi_\theta(u)+s-\phi_\theta(v)|+K\|u-v\|_{\nu}\\ \notag
&\le (1+K)\|x-h\|_{\nu}\\  \label{eq:C_gr}
&=C_{\rm gr}\|x-h\|_{\nu},
\end{align}
where the second inequality uses the $K$-Lipschitz property of $\phi_\theta$, and the third inequality follows because both $|\phi_\theta(u)+s-\phi_\theta(v)|$ and $\|u-v\|_{\nu}$ are bounded by $\|x-h\|_{\nu}$. Since this holds for every $h\in\mathcal H_\theta$, taking the infimum over $h\in\mathcal H_\theta$ gives
\begin{equation} \label{dist_x_H}
       \operatorname{dist}_{\nu}(x,\mathcal H_\theta)\ge \frac{s}{C_{\rm gr}},
\end{equation}
 
The case $x=(u,\phi_\theta(u)-s)$ is identical.

Now define $\varepsilon_0=\min\{1,L\rho\}/\{4(C_{\rm gr}+4)\}$. This ensures
\begin{align}
(C_{\rm gr}+4)\varepsilon_0&\le 1,\label{eq:uniform-interior-eps-upper-one}\\
\frac{(C_{\rm gr}+1)\varepsilon_0}{L}&\le \rho.\label{eq:uniform-interior-eps-vertical-room}
\end{align}
Note that \eqref{eq:uniform-interior-eps-vertical-room} ensures that all the small vertical intervals in \textbf{Step 4}, e.g.~\eqref{eq_vertical_inerval}, remain inside the cube.

\medskip

\noindent\textbf{Step 3: Partition the boundary trace.} For $u\in S_\theta$, write $h(u)=(u,\phi_\theta(u))\in\mathcal H_\theta$, and define $q_\theta(u)=q(h(u))$. Partition $S_\theta$ into
\[
    A_0=\{u\in S_\theta:q_\theta(u)\le 2\varepsilon_0\},\quad A_m=\{u\in S_\theta:2\varepsilon_0<q_\theta(u)<1-2\varepsilon_0\},\quad A_1=\{u\in S_\theta:q_\theta(u)\ge 1-2\varepsilon_0\}.
\]
These three sets cover $S_\theta$. We will show that for every $u\in S_\theta$, there is a vertical interval of length at least $\varepsilon_0/L$ on which $\varepsilon_0\le p_{\theta,q}(u,z)\le 1-\varepsilon_0$.

\medskip

\noindent\textbf{Step 4: Low trace values give an interior region above the boundary.} 
Take $u\in A_0$. Then $q_\theta(u)\le 2\varepsilon_0$. Consider points above the boundary:
\begin{align} \label{eq_vertical_inerval}
    x=(u,\phi_\theta(u)+s),
    \qquad
    s\in\left[\frac{C_{\rm gr}\varepsilon_0}{L},\frac{(C_{\rm gr}+1)\varepsilon_0}{L}\right].
\end{align}
This interval has length $\varepsilon_0/L$. Since $u\in S_\theta$, the boundary height satisfies $\phi_\theta(u)\le 1-\rho$. By \eqref{eq:uniform-interior-eps-vertical-room}, we have that $(C_{\rm gr}+1)\varepsilon_0/L\le \rho$, therefore, these points remain inside $[0,1]^d$. Moreover, it holds that $x\in\mathcal X_\theta^+$, as by construction, $a+\gamma_{-j}^\top u+\gamma_j\{\phi_\theta(u)+s\}
=\gamma_j s>0$. 

Since $x$ lies above the boundary, we have that $p_{\theta,q}(x)=\min\{1,U_{\theta,q}(x)\}$. Since $q\ge0$, it holds that
\[
    U_{\theta,q}(x)=\inf_{h\in\mathcal H_\theta}\{q(h)+L\|x-h\|_{\nu}\}\ge L\operatorname{dist}_{\nu}(x,\mathcal H_\theta).
\]
Therefore, it holds that $U_{\theta,q}(x)\ge Ls/C_{\rm gr}\ge \varepsilon_0$ by \eqref{dist_x_H}. On the other hand, using the specific boundary point $h(u)$,
\[
    U_{\theta,q}(x)\le q(h(u))+L\|x-h(u)\|_{\nu}.
\]
Since $\|x-h(u)\|_{\nu}=s$, we have
\[
    U_{\theta,q}(x)\le q_\theta(u)+Ls\le 2\varepsilon_0+(C_{\rm gr}+1)\varepsilon_0=(C_{\rm gr}+3)\varepsilon_0 \leq 1-\varepsilon_0,
\]
where the last inequality follows from \eqref{eq:uniform-interior-eps-upper-one}. Hence $\varepsilon_0\le p_{\theta,q}(x)\le 1-\varepsilon_0$. Thus, over every $u\in A_0$, there is a vertical interval of length $\varepsilon_0/L$ above the boundary where the policy is interior.

\medskip

\noindent\textbf{Step 5: High trace values give an interior region below the boundary.} Take $u\in A_1$. Then $q_\theta(u)\ge 1-2\varepsilon_0$. Consider points below the boundary:
\[
    x=(u,\phi_\theta(u)-s),
    \qquad
    s\in\left[\frac{C_{\rm gr}\varepsilon_0}{L},\frac{(C_{\rm gr}+1)\varepsilon_0}{L}\right].
\]
Again this interval has length $\varepsilon_0/L$. Since $\phi_\theta(u)\ge \rho$ and $(C_{\rm gr}+1)\varepsilon_0/L\le \rho$ by \eqref{eq:uniform-interior-eps-vertical-room}, these points remain inside the cube.

Since $x$ is below the boundary, $p_{\theta,q}(x)=\max\{0,M_{\theta,q}(x)\}$. Using the boundary point $h(u)$, we have that
\[
    M_{\theta,q}(x)\ge q(h(u))-L\|x-h(u)\|_{\nu}=q_\theta(u)-Ls.
\]
Thus,
\[
    M_{\theta,q}(x)\ge 1-2\varepsilon_0-(C_{\rm gr}+1)\varepsilon_0=1-(C_{\rm gr}+3)\varepsilon_0 \ge \varepsilon_0,
\]
where the last inequality follows from \eqref{eq:uniform-interior-eps-upper-one}. Also, since $q\le1$, we have that
\[
    M_{\theta,q}(x)=\sup_{h\in\mathcal H_\theta}\{q(h)-L\|x-h\|_{\nu}\}\le 1-L\operatorname{dist}_{\nu}(x,\mathcal H_\theta).
\]
Therefore, we have that $M_{\theta,q}(x)\le 1-Ls/C_{\rm gr}\le 1-\varepsilon_0$, and hence, $\varepsilon_0\le p_{\theta,q}(x)\le 1-\varepsilon_0$. So over every $u\in A_1$, there is a vertical interval of length $\varepsilon_0/L$ below the boundary where the policy is interior.

\medskip

\noindent\textbf{Step 6: Middle trace values give an interior region near the boundary.} Take $u\in A_m$. Then $2\varepsilon_0<q_\theta(u)<1-2\varepsilon_0$. Consider points just above the boundary:
\[
    x=(u,\phi_\theta(u)+s),
    \qquad
    s\in\left[0,\frac{\varepsilon_0}{L}\right].
\]
These points remain inside the cube since $\varepsilon_0/L\le \rho$, which follows from \eqref{eq:uniform-interior-eps-vertical-room}.

Since $x$ is above the boundary, $p_{\theta,q}(x)=\min\{1,U_{\theta,q}(x)\}$. Using $h(u)$ gives
\[
    U_{\theta,q}(x)\le q_\theta(u)+Ls\le q_\theta(u)+\varepsilon_0<1-\varepsilon_0.
\]
For the lower bound, for every $h\in\mathcal H_\theta$, we have that
\[
    q(h)\ge q(h(u))-L\|h-h(u)\|_{\nu}.
\]
Hence, it holds that
\[
    q(h)+L\|x-h\|_{\nu}\ge q(h(u))-L\|h-h(u)\|_{\nu}+L\|x-h\|_{\nu}.
\]
Using the triangle inequality, $\|h-h(u)\|_{\nu}\le \|h-x\|_{\nu}+\|x-h(u)\|_{\nu}$, so
\[
    q(h)+L\|x-h\|_{\nu}\ge q(h(u))-L\|x-h(u)\|_{\nu}.
\]
Taking the infimum over $h\in\mathcal H_\theta$ gives
\[
    U_{\theta,q}(x)\ge q_\theta(u)-Ls\ge q_\theta(u)-\varepsilon_0>\varepsilon_0.
\]
Therefore $\varepsilon_0\le p_{\theta,q}(x)\le 1-\varepsilon_0$. Thus over every $u\in A_m$, there is also a vertical interval of length $\varepsilon_0/L$ on which the policy is interior.

\medskip

\noindent\textbf{Step 7: Integrate over the base patch.} For every $u\in S_\theta$, one of the three cases $u\in A_0$, $u\in A_m$, or $u\in A_1$ holds. In each case, the vertical fiber $\{z\in[0,1]:(u,z)\in[0,1]^d\}$ contains a subinterval of length at least $\varepsilon_0/L$ on which $\varepsilon_0\le p_{\theta,q}(u,z)\le 1-\varepsilon_0$. Therefore, by Fubini's theorem,
\[
    \lambda_d\left(\left\{x\in[0,1]^d:\varepsilon_0\le p_{\theta,q}(x)\le 1-\varepsilon_0\right\}\right)\ge \lambda_{d-1}(S_\theta)\frac{\varepsilon_0}{L}.
\]
Since $\lambda_{d-1}(S_\theta)\ge r_0^{d-1}$, we get
\[
    \lambda_d\left(\left\{x\in[0,1]^d:\varepsilon_0\le p_{\theta,q}(x)\le 1-\varepsilon_0\right\}\right)\ge r_0^{d-1}\frac{\varepsilon_0}{L}.
\]
Taking $v_0=r_0^{d-1}\varepsilon_0/L$ proves the claim.
\end{proof}

\begin{lemma}[Covariate diversity on positive-measure sets]\label[lemma]{l_covariate_diverse}
Suppose \Cref{a_contexts} holds. Let $A\subseteq[0,1]^d$ be measurable and suppose $\lambda_d(A)\ge v_0>0$. Define $\widetilde X=(1,X^\top)^\top$. Then, we have that
\[
    \lambda_{\min}\left(\mathbb{E}[\widetilde X\widetilde X^\top\mathbf 1\{X\in A\}]\right)\ge c_g\frac{v_0^3}{128d^2}.
\]
\end{lemma}

\begin{proof}[Proof of \Cref{l_covariate_diverse}]
Let $M_A=\mathbb{E}[\widetilde X\widetilde X^\top\mathbf 1\{X\in A\}]$. Fix any unit vector $w=(w_0,\bar w)\in\mathbb R^{d+1}$, where $\bar w\in\mathbb R^d$. We first show that there exists a constant $C_d=4d$ such that, for every $t>0$,
\[
    \lambda_d\{x\in[0,1]^d:|w_0+\bar w^\top x|\le t\}\le C_dt.
\]
If $\|\bar w\|_2\ge 1/(2\sqrt d)$, choose $j\in\{1,\ldots,d\}$ such that $|\bar w_j|\ge \|\bar w\|_2/\sqrt d\ge 1/(2d)$. For fixed $x_{-j}$, the set of $x_j\in[0,1]$ satisfying $|w_0+\bar w^\top x|\le t$ is an interval of length at most $2t/|\bar w_j|\le 4dt$. Hence, by Fubini's theorem, $\lambda_d\{x\in[0,1]^d:|w_0+\bar w^\top x|\le t\}\le 4dt$.

If $\|\bar w\|_2<1/(2\sqrt d)$, then $|w_0|=(1-\|\bar w\|_2^2)^{1/2}$ and $|\bar w^\top x|\le \|\bar w\|_2\|x\|_2<1/2$ for all $x\in[0,1]^d$. Hence, it holds that  $|w_0+\bar w^\top x|>1/4$ from the triangle inequality for all $x\in[0,1]^d$. Therefore the slab is empty when $t\le1/4$, while for $t>1/4$ its Lebesgue measure is at most $1\le4t\le4dt$. Thus, in all cases, we have that
\[
    \lambda_d\{x\in[0,1]^d:|w_0+\bar w^\top x|\le t\}\le 4dt.
\]

Take $t_0=v_0/(8d)$. Then, we have that
\[
    \lambda_d\{x\in[0,1]^d:|w_0+\bar w^\top x|\le t_0\}\le \frac{v_0}{2}.
\]
Since $\lambda_d(A)\ge v_0$, we have
\[
    \lambda_d\left(A\cap\{x:|w_0+\bar w^\top x|>t_0\}\right)\ge \frac{v_0}{2}.
\]
Therefore, using the density lower bound $g(x)\ge c_g$ from \Cref{a_contexts},
\begin{align*}
    w^\top M_Aw&=\mathbb{E}[(w^\top\widetilde X)^2\mathbf 1\{X\in A\}]\\
    &=\int_A (w_0+\bar w^\top x)^2g(x)\,\mathrm{d}x\\
    &\ge c_g\int_A (w_0+\bar w^\top x)^2\,\mathrm{d}x\\
    &\ge c_g\int_{A\cap\{|w_0+\bar w^\top x|>t_0\}}(w_0+\bar w^\top x)^2\,\mathrm{d}x\\
    &\ge c_gt_0^2\frac{v_0}{2}
    =c_g\frac{v_0^3}{128d^2}.
\end{align*}
Taking the infimum over all unit vectors $w\in\mathbb R^{d+1}$ gives $\lambda_{\min}(M_A)\ge c_gv_0^3/(128d^2)$.
\end{proof}

\begin{lemma}[Exploration property]\label[lemma]{lem:lipschitz-min-pull-proportion}
Suppose Assumptions \ref{a_contexts} and \ref{a_theta} hold. Let $\theta\in \Theta$ and $q\in\mathcal Q_{\theta,L}$, with $\mathcal Q_{\theta,L}$ defined in \eqref{eq_set_voundary_policy}. Define $\varepsilon_L=\min\{1,L\rho_0/4\}$.

Then there exists a constant $c_{\rm pull}>0$, depending only on $d,\nu,c_g$ and the constants in \Cref{a_theta}, such that
\begin{equation*}
    \mathbb E\{p_{\theta,q}(X)\}\wedge \mathbb E\{1-p_{\theta,q}(X)\}\ge c_{\rm pull}\frac{\varepsilon_L^2}{L}.
\end{equation*}
In particular, if $L\ge 4/\rho_0$, then $\mathbb E\{p_{\theta,q}(X)\}\wedge \mathbb E\{1-p_{\theta,q}(X)\}\ge c_{\rm pull}/L$, whereas if $L<4/\rho_0$, then $ \mathbb E\{p_{\theta,q}(X)\}\wedge \mathbb E\{1-p_{\theta,q}(X)\}\ge c'_{\rm pull}L$, where $c'_{\rm pull} = c_{\rm pull} \rho_0^2/4$.
\end{lemma}
\begin{proof}[Proof of \Cref{lem:lipschitz-min-pull-proportion}]
Fix $\theta\in\widetilde\Theta$ and $q\in\mathcal Q_{\theta,L}$. By \Cref{lem:uniform-interior-region}, there exists a measurable set $B_{\theta,q}\subseteq[0,1]^d$ such that
\begin{equation*}
    \frac{\varepsilon_L}{4(C_{\rm gr}+4)}\le p_{\theta,q}(x)\le 1-\frac{\varepsilon_L}{4(C_{\rm gr}+4)},\qquad x\in B_{\theta,q},
\end{equation*}
where $C_{\rm gr}>0$ is the constant defined in \Cref{eq:C_gr}. Moreover, it holds that
\begin{equation*}
    \lambda_d(B_{\theta,q})\ge r_0^{d-1}\frac{\varepsilon_L}{4(C_{\rm gr}+4)L},
\end{equation*}
where $r_0>0$ is defined in \eqref{eq:r_0} and depends only on $d,\nu$ and the constants in \Cref{a_theta}. Since $p_{\theta,q}(x)\ge \varepsilon_L/\{4(C_{\rm gr}+4)\}$ on $B_{\theta,q}$, we have that
\begin{align} \label{lem:lipschitz-min-pull-proportion_eq1}
\mathbb E\{p_{\theta,q}(X)\}&\ge \mathbb E\{p_{\theta,q}(X)\indc(X\in B_{\theta,q})\} \ge \frac{\varepsilon_L}{4(C_{\rm gr}+4)}\mathbb P(X\in B_{\theta,q}).
\end{align}
By \Cref{a_contexts}, the density of $X$ is bounded below by $c_g$ on $[0,1]^d$. Hence, we have that
\begin{align} \label{lem:lipschitz-min-pull-proportion_eq2}
\mathbb P(X\in B_{\theta,q})&=\int_{B_{\theta,q}}g(x)\,dx \ge c_g\lambda_d(B_{\theta,q}) \ge c_g r_0^{d-1}\frac{\varepsilon_L}{4(C_{\rm gr}+4)L}.
\end{align}
Substituting \eqref{lem:lipschitz-min-pull-proportion_eq2} into \eqref{lem:lipschitz-min-pull-proportion_eq1} gives
\begin{equation*}
    \mathbb E\{p_{\theta,q}(X)\}\ge \frac{c_g r_0^{d-1}}{16(C_{\rm gr}+4)^2}\frac{\varepsilon_L^2}{L}.
\end{equation*}
Similarly, since $1-p_{\theta,q}(x)\ge \varepsilon_L/\{4(C_{\rm gr}+4)\}$ on $B_{\theta,q}$, we have that
\begin{align*}
\mathbb E\{1-p_{\theta,q}(X)\}&\ge \mathbb E[\{1-p_{\theta,q}(X)\}\indc(X\in B_{\theta,q})] \ge \frac{\varepsilon_L}{4(C_{\rm gr}+4)}\mathbb P(X\in B_{\theta,q}) \ge \frac{c_g r_0^{d-1}}{16(C_{\rm gr}+4)^2}\frac{\varepsilon_L^2}{L}.
\end{align*}
Therefore, we have that
\begin{equation*}
    \mathbb E\{p_{\theta,q}(X)\}\wedge \mathbb E\{1-p_{\theta,q}(X)\}\ge \frac{c_g r_0^{d-1}}{16(C_{\rm gr}+4)^2}\frac{\varepsilon_L^2}{L}.
\end{equation*}
Taking $c_{\rm pull}=c_g r_0^{d-1}/\{16(C_{\rm gr}+4)^2\}$ proves the first claim.

In the case when $L\ge2/\widetilde\rho_0$, we have that $\varepsilon_L=1$. Hence the previous display gives
\begin{equation*}
    \mathbb E\{p_{\theta,q}(X)\}\wedge \mathbb E\{1-p_{\theta,q}(X)\}\ge \frac{c_{\rm pull}}{L}.
\end{equation*}

On the other hand, when $L<2/\widetilde\rho_0$, we have that $\varepsilon_L=L\widetilde\rho_0/2$. Hence the previous display gives
\begin{align*}
    \mathbb E\{p_{\theta,q}(X)\}\wedge \mathbb E\{1-p_{\theta,q}(X)\}
    &\ge c_{\rm pull}\frac{(L\widetilde\rho_0/2)^2}{L} =\frac{c_{\rm pull}\widetilde\rho_0^2}{4}L =c'_{\rm pull}L,
\end{align*}
where $c'_{\rm pull}=c_{\rm pull}\widetilde\rho_0^2/4$. This completes the proof.
\end{proof}

\begin{lemma}[Norm comparison]
\label[lemma]{lem:norm-comparison}
Fix $\nu\in[1,\infty]$. Let $ C_{d,\nu}=d^{\max\{1/\nu-1/2,0\}}$, where $1/\infty=0$. Then, for every $z\in\mathbb R^d$,
\[
    \|z\|_\nu\le C_{d,\nu}\|z\|_2.
\]
The same bound also holds for vectors in $\mathbb R^{d-1}$, after viewing them as vectors in $\mathbb R^d$ by adding one zero coordinate.
\end{lemma}

\begin{proof}
If $\nu\ge2$, then the monotonicity of $\ell_p$ norms gives $\|z\|_\nu\le\|z\|_2$, so the claim holds because $C_{d,\nu}=1$. If $1\le\nu<2$, Holder's inequality gives
\[
    \|z\|_\nu = \left(\sum_{j=1}^d |z_j|^\nu\right)^{1/\nu}
    \le  d^{1/\nu-1/2}\left(\sum_{j=1}^d z_j^2\right)^{1/2}= C_{d,\nu}\|z\|_2.
\]
This proves the claim.
\end{proof}

\subsection{Proof of Lemma \ref{lem:vertex-optimisation}}\label[appendix]{sec_app_vertex_optimisation}
\begin{proof}[Proof of \Cref{lem:vertex-optimisation}]
By \eqref{eq:Hs} and \eqref{eq:Fs}, we have that
\begin{equation*}
F_s(b_1,b_2)=\max_{z\in\mathcal Z_{s-1,L}}\frac1{s-1}\sum_{i=1}^{s-1}\left[\widetilde X_i^\top b_1z_i+\widetilde X_i^\top b_2(1-z_i)\right].
\end{equation*}
For any fixed $z\in\mathcal Z_{s-1,L}$, the expression inside the maximum is affine in $(b_1,b_2)$, meaning that it depends linearly on $b_1$ and $b_2$. In particular, for any two model pairs $(b_1,b_2)$ and $(b_1',b_2')$ and any $\lambda\in[0,1]$, its value at the convex combination ${\lambda b_1+(1-\lambda)b_1',\lambda b_2+(1-\lambda)b_2'}$ is exactly
\begin{align*}
&\frac1{s-1}\sum_{i=1}^{s-1}\left[\widetilde X_i^\top{\lambda b_1+(1-\lambda)b_1'}z_i+\widetilde X_i^\top{\lambda b_2+(1-\lambda)b_2'}(1-z_i)\right]\\
=&\frac{\lambda}{s-1}\sum_{i=1}^{s-1}\left[\widetilde X_i^\top b_1z_i+\widetilde X_i^\top b_2(1-z_i)\right]+\frac{1-\lambda}{s-1}\sum_{i=1}^{s-1}\left[\widetilde X_i^\top b_1'z_i+\widetilde X_i^\top b_2'(1-z_i)\right].
\end{align*}
Since the feasible set $\mathcal Z_{s-1,L}$ is the same for every choice of $(b_1,b_2)$, taking the maximum over $z$ on the left-hand side cannot exceed taking the maximum separately in the two terms on the right-hand side. Therefore, it holds that  $F_s$ is jointly convex. Continuity follows from the compactness of $\mathcal Z_{s-1,L}$ and the continuity of the objective.

By \Cref{lem:polyhedral-confidence}, any feasible pair can be written as
\begin{equation*}
b_1=\sum_{v_1\in\mathcal V_{s,1}}a_{v_1}v_1,\qquad b_2=\sum_{v_2\in\mathcal V_{s,2}}c_{v_2}v_2,
\end{equation*}
where the coefficients are nonnegative and each collection sums to one. Combining these two representations gives
\begin{equation*}
(b_1,b_2)=\sum_{v_1\in\mathcal V_{s,1}}\sum_{v_2\in\mathcal V_{s,2}}a_{v_1}c_{v_2}(v_1,v_2).
\end{equation*}
The coefficients $a_{v_1}c_{v_2}$ are nonnegative and satisfy
\begin{equation*}
\sum_{v_1\in\mathcal V_{s,1}}\sum_{v_2\in\mathcal V_{s,2}}a_{v_1}c_{v_2}=\left(\sum_{v_1\in\mathcal V_{s,1}}a_{v_1}\right)\left(\sum_{v_2\in\mathcal V_{s,2}}c_{v_2}\right)=1.
\end{equation*}
Thus, every feasible pair $(b_1,b_2)$ is a convex combination of the vertex pairs $(v_1,v_2)$. Joint convexity therefore gives
\begin{equation*}
F_s(b_1,b_2)\le\sum_{v_1\in\mathcal V_{s,1}}\sum_{v_2\in\mathcal V_{s,2}}a_{v_1}c_{v_2}F_s(v_1,v_2)\le\max_{v_k\in\mathcal V_{s,k},\ k\in\{1,2\}}F_s(v_1,v_2),
\end{equation*}
where the second inequality follows because the middle expression is a convex combination of the values $F_s(v_1,v_2)$ and therefore cannot exceed their largest value. Since every vertex pair is feasible, this proves \eqref{eq:vertex-max}.

Suppose the selected model pair satisfies $\mathcal H_{\bar\beta_{s,1}-\bar\beta_{s,2}}\neq\emptyset$. By \Cref{thm:canonical_extension}, every feasible canonical policy $p_{b_1-b_2,q}$ belongs to $\mathcal P_{L,\nu}$, so its sample values belong to $\mathcal Z_{s-1,L}$. Hence, its empirical objective is bounded above by $F_s(b_1,b_2)$, and the supremum on the right-hand side of \eqref{eq:vertex-opt} is bounded above by the left-hand side of \eqref{eq:vertex-max}. For every vertex pair $v$, the optimality of $z_v$ gives
\begin{equation*}
S_v=\frac1{s-1}\sum_{i=1}^{s-1}\widetilde X_i^\top v_2+H_s(v_1-v_2)=F_s(v_1,v_2).
\end{equation*}
If $v^\circ$ maximises $F_s$ over the vertex pairs, the selected pair satisfies $S_{\widehat v}=S_{v^\circ}=F_s(v_1^\circ,v_2^\circ)$. The McShane extension of $z_s$ belongs to $\mathcal P_{L,\nu}$ and agrees with $z_s$ at the sample points. By \Cref{thm:canonical_extension}, its boundary-canonical representative does not decrease the empirical objective. Hence, it holds that
\begin{equation*}
\frac1{s-1}\sum_{i=1}^{s-1}r_{\bar\beta_s,p_s}(X_i)\ge S_{\widehat v}=\max_{b_k\in\mathcal C_{s,k},\ k\in\{1,2\}}F_s(b_1,b_2).
\end{equation*}
Since the returned model-policy pair is feasible in \eqref{eq:ideal-update}, the reverse inequality also holds. This proves \eqref{eq:vertex-opt}.
\end{proof}

\subsubsection{Polyhedral confidence sets}
\begin{lemma}[Polyhedral confidence sets]\label[lemma]{lem:polyhedral-confidence}
For every update time $s$ and $k\in\{1,2\}$, we have that
\begin{equation*}
\left\{b:\|b-\widehat\beta_{s,k}\|_{V_{s,k}}\le\alpha_T\right\}\subseteq\mathcal C_{s,k}=\operatorname{conv}(\mathcal V_{s,k}),
\end{equation*}
where $\mathcal{C}_{s,k}$ is defined in \eqref{eq:candidate-set},  $\mathcal{V}_{s,k}$ is defined in \eqref{eq:confidence-vertices}, and $\operatorname{conv}(\cdot)$ denotes the convex hull.

Suppose $\|\widehat\beta_{s,k}-\beta_k^*\|_{V_{s,k}}\le\alpha_T$, then it holds that $\beta_k^*\in\mathcal C_{s,k}$ and, for every $b_k\in\mathcal C_{s,k}$,
\begin{equation}
\|b_k-\beta_k^*\|_{V_{s,k}}\le(1+\sqrt{d+1})\alpha_T.
\end{equation}
\end{lemma}

\begin{proof}
For every $u\in\mathbb R^{d+1}$, we have $\|u\|_2\le\|u\|_1\le\sqrt{d+1}\|u\|_2$. Note that $\mathcal{C}_{s,k}=\{b:\|V_{s,k}^{1/2}(b-\widehat\beta_{s,k})\|_1\leq \sqrt{d+1}\alpha_T\}$, applying the inequality to $u=V_{s,k}^{1/2}(b-\widehat\beta_{s,k})$ proves the first claim.

We next show that $\mathcal C_{s,k}=\operatorname{conv}(\mathcal V_{s,k})$. Under the transformation $u=V_{s,k}^{1/2}(b-\widehat\beta_{s,k})$, the set $\mathcal C_{s,k}$ corresponds to the $\ell_1$-ball $B_1(\sqrt{d+1}\alpha_T) = \{u\in\mathbb R^{d+1}:\|u\|_1\le\sqrt{d+1}\alpha_T\}$, i.e.~$\mathcal C_{s,k}=\widehat\beta_{s,k}+V_{s,k}^{-1/2}B_1(\sqrt{d+1}\alpha_T)$. Note that this $\ell_1$-ball is the convex hull of its $2(d+1)$ vertices $\{\pm\sqrt{d+1}\alpha_Te_r:r\in[d+1]\}$. To see this, for any $u$ in the ball, the contribution of its $r$th coordinate can be written as 
$$\frac{|u_r|}{\sqrt{d+1}\alpha_T}\{\operatorname{sign}(u_r)\sqrt{d+1}\alpha_Te_r\}=u_re_r.$$ 
Summing over $r$ therefore recovers $u$, while the corresponding nonnegative coefficients have total weight $\|u\|_1/(\sqrt{d+1}\alpha_T)\le1$. If this total weight is strictly smaller than one, the remaining weight can be split equally between any pair of opposite vertices, e.g.~$\sqrt{d+1}\alpha_Te_1$ and $-\sqrt{d+1}\alpha_Te_1$. Since these two vertices sum to zero, assigning equal additional weight to them leaves $u$ unchanged while making the total weight equal to one. Hence every point in the $\ell_1$-ball is a convex combination of these $2(d+1)$ vertices. Conversely, all these vertices belong to the $\ell_1$-ball, which is convex, so the $\ell_1$-ball is exactly their convex hull.

Since $u=V_{s,k}^{1/2}(b-\widehat\beta_{s,k})$ by construction, we equivalently have $b=\widehat\beta_{s,k}+V_{s,k}^{-1/2}u$. Therefore, the vertices $\pm\sqrt{d+1}\alpha_Te_r$ are mapped to $\widehat\beta_{s,k}\pm\sqrt{d+1}\alpha_TV_{s,k}^{-1/2}e_r$. These transformed vertices form the set
\begin{equation*}
\mathcal V_{s,k}=\left\{\widehat\beta_{s,k}\pm\sqrt{d+1}\alpha_TV_{s,k}^{-1/2}e_r:r\in[d+1]\right\}.
\end{equation*}
Since affine transformations preserve convex combinations, this gives $\mathcal C_{s,k}=\operatorname{conv}(\mathcal V_{s,k})$.

Suppose $\|\widehat\beta_{s,k}-\beta_k^*\|_{V_{s,k}}\le\alpha_T$. The containment proved above gives $\beta_k^*\in\mathcal C_{s,k}$. Finally, for every $b_k\in\mathcal C_{s,k}$, we have $\|b_k-\widehat\beta_{s,k}\|_{V_{s,k}}\le\sqrt{d+1}\alpha_T$, and the triangle inequality gives
\begin{align*}
\|b_k-\beta_k^*\|_{V_{s,k}}
&\le\|b_k-\widehat\beta_{s,k}\|_{V_{s,k}}+\|\widehat\beta_{s,k}-\beta_k^*\|_{V_{s,k}}\\
&\le\sqrt{d+1}\alpha_T+\alpha_T=(1+\sqrt{d+1})\alpha_T.
\end{align*}
This completes the proof.
\end{proof}

\begin{lemma}[Computational complexity]\label[lemma]{lem:vertex-complexity}
It holds that each update in \Cref{alg:vertex-oracle} requires solving $4(d+1)^2$ LPs with $s-1$ variables and at most $s(s-1)$ linear inequalities. With $s_0=2T_o+1\le T$, the number $K_T$ of updates in \Cref{alg:lazy-global-model-ucb} satisfies
\begin{equation}\label{eq:vertex-update-count}
K_T\le1+\frac{\log(T/s_0)}{\log2}+\frac{2(d+1)}{\log2}\log\left(1+\frac T\lambda\right).
\end{equation}
Consequently, the total number of policy LPs is at most $4(d+1)^2K_T$. 
\end{lemma}

\begin{proof}
Each confidence polytope has $2(d+1)$ vertices, giving $4(d+1)^2$ pairs. For $n=s-1$, the constraints in \eqref{eq:ZnL} consist of $2n$ box inequalities and $n(n-1)$ pairwise Lipschitz inequalities, for a total of $n(n+1)=s(s-1)$. The two inverse square roots and the pairwise distances are computed once per update. 

To count the total number of updates, write $s_1=s_0<s_2<\cdots<s_{K_T}\le T$. We will update when index doubles or the determinant doubles. If $K_0$ updates are assigned to the time trigger, then 
\begin{equation*}
K_0\log2\le\sum_{j=1}^{K_T-1}\log(s_{j+1}/s_j)\le\log(T/s_0).
\end{equation*}
If $K_k$ updates are assigned to arm $k$, monotonicity of $V_{s,k}$ gives
\begin{equation*}
K_k\log2\le\sum_{j=1}^{K_T-1}\log\frac{\det(V_{s_{j+1},k})}{\det(V_{s_j,k})}\le\log\frac{\det(V_{T,k})}{\lambda^{d+1}}\le(d+1)\log\left(1+\frac T\lambda\right),
\end{equation*}
where the last inequality follows from the trace bound $\operatorname{tr}(V_{T,k})\le(d+1)(\lambda+T)$ and the arithmetic--geometric mean inequality for eigenvalues, similar as in \Cref{lem:lazy-width-sum-lnu}. Since $K_T=1+K_0+K_1+K_2$, this proves \eqref{eq:vertex-update-count}.
\end{proof}

\subsection{Proof of Theorem \ref{thm:main-regret}}\label[appendix]{sec_app_ucb_regret}
For this subsection, let $\tau(t)=\max\{s\le t:s\text{ is an update time of \Cref{alg:lazy-global-model-ucb}}\}$ denote the most recent update time before $t$, and let $\mathcal F_t=\sigma\{X_i,A_i,Y_i:1\le i\le t\}$.

\begin{proof}[Proof of \Cref{thm:main-regret}]
Let $B=\sqrt{d+1}(C_\beta+1)$ and fix $r_{\rm loc}\in(0,1]$ as in \Cref{lem:enlarged_theta}. By \Cref{thm:canonical_extension}, there exists a population-optimal canonical policy $p^*=p_{\theta^*,q^*}$. During the first $2T_o$ rounds, the algorithm uses $p_t= 1/2$. Since the conditional mean reward of each arm has magnitude at most $B$, the regret in each exploration round is at most $2B$. Hence, since $T_o=\lceil\sqrt T\rceil$, the total regret during exploration is at most $4BT_o\le C_1\sqrt T$, where $C_1=8B$.

\medskip
\noindent \textbf{High-probability events.}
We first construct a high-probability event $\mathcal E=\mathcal E_{\rm conf}\cap\mathcal E_{\rm emp}\cap\mathcal E_{\rm mart}\cap\mathcal E_{\rm init}$, on which the remaining analysis is carried out. 

Define
\begin{equation}\label{eq:main-regret-confidence-event}
\mathcal E_{\rm conf}=\left\{\|\widehat\beta_{t,k}-\beta_k^*\|_{V_{t,k}}\le\alpha_T,\ \text{for every }t\in[T]\text{ and }k\in\{1,2\}\right\}.
\end{equation}

Applying \Cref{cor:self-normalised-deterministic-lnu} with failure probability $\eta_{\rm c}=\eta/8$, it then holds with probability at least $1-\eta/8$ that for every $t\in[T]$ and $k\in\{1,2\}$,
\begin{equation*}
\|S_{t,k}\|_{V_{t,k}^{-1}}\le\sqrt{2\sigma^2\left\{\log\left(\frac{2}{\eta_{\rm c}}\right)+\frac{d+1}{2}\log\left(1+\frac{T}{\lambda}\right)\right\}}.
\end{equation*}
Combining this with \eqref{l_ucb_lnu_eq1} and $\|\beta_k^*\|_2\le C_\beta$ from \Cref{a_beta}, the definition of $\alpha_T$ in \eqref{eq_ucb_alpha} gives
\begin{equation*}
\|\widehat\beta_{t,k}-\beta_k^*\|_{V_{t,k}}\le\alpha_T
\end{equation*}
for every $t\in[T]$ and $k\in\{1,2\}$. Hence $P(\mathcal E_{\rm conf})\ge1-\eta/8$. On $\mathcal E_{\rm conf}$, \Cref{lem:polyhedral-confidence} gives $\beta_k^*\in\mathcal C_{s,k}$ for every update time $s$ and $k\in\{1,2\}$.

Next, with $\Delta_n^{\mathrm{UCB}}$ as in \eqref{eq_delta_ucb}, define
\begin{equation}\label{eq:main-regret-empirical-event}
\mathcal E_{\rm emp}=\left\{\Delta_n^{\mathrm{UCB}}\le C_2\rho_{n,d-1}\left(\frac{\eta}{8T}\right),\ \text{for every }n\in\{2,\ldots,T\}\right\}.
\end{equation}
By \Cref{lem:ucb-control-delta}, for each $n\in\{2,\ldots,T\}$ the displayed inequality fails with probability at most $\eta/(8T)$. A union bound over these $T-1$ values of $n$ therefore gives $P(\mathcal E_{\rm emp})\ge1-(T-1)\eta/(8T)\ge1-\eta/8$.

For $s_0\le t\le T$, write $W_t=\min\{1,\|\widetilde X_t\|_{V_{\tau(t),A_t}^{-1}}\}$, and define
\begin{equation}\label{eq:main-regret-martingale-event}
\mathcal E_{\rm mart}=\left\{\sum_{t=s_0}^T\mathbb E[W_t\mid\mathcal F_{t-1}]\le\sum_{t=s_0}^TW_t+\sqrt{2T\log(8/\eta)}\right\}.
\end{equation}
Since $W_t$ is $\mathcal F_t$-measurable and $0\le W_t\le1$ by construction, the variables $Z_t=\mathbb E[W_t\mid\mathcal F_{t-1}]-W_t$ form a martingale difference sequence with $|Z_t|\le1$. The Azuma--Hoeffding inequality \citep[e.g.~Corollary 2.20 in][]{wainwright2019high}, applied to at most $T$ increments, gives $P(\mathcal E_{\rm mart})\ge1-\eta/8$.

Finally, let $\mathcal E_{\rm init}$ be the event associated with the exploration phase in \eqref{eq:initial-design-event}. Since the first coordinate of $\widetilde X$ is one, we have that $\kappa_X\le1$. By construction of $\alpha_T$ in \eqref{eq_ucb_alpha}, we have that $\alpha_T^2\ge2\sigma^2\log(16/\eta)$. Choose $C_3$ sufficiently large so that the condition
\begin{equation*}
T\ge\frac{C_3(1+\sqrt{d+1})^4\alpha_T^4}{\kappa_X^2}
\end{equation*}
implies
\begin{equation*}
\sqrt T\ge8\log(16/\eta),\qquad \sqrt T\ge\frac{16(d+1)}{\kappa_X}\log\left(\frac{16(d+1)}{\eta}\right).
\end{equation*}
Since $T_o\ge\sqrt T$, \eqref{eq:initial-allocation-prob} and \eqref{eq:initial-design-prob} then give $P(\mathcal E_{\rm init}^c)\le\eta/4$. Set $\mathcal E=\mathcal E_{\rm conf}\cap\mathcal E_{\rm emp}\cap\mathcal E_{\rm mart}\cap\mathcal E_{\rm init}$. By a union bound,
\begin{equation}
P(\mathcal E)\ge1-\frac{\eta}{8}-\frac{\eta}{8}-\frac{\eta}{8}-\frac{\eta}{4}=1-\frac{5\eta}{8}\ge1-\eta.
\end{equation}

\medskip
\noindent \textbf{Analysis on $\mathcal{E}$.}  The rest of the proof is carried out on $\mathcal E$.  It remains to verify that, at every update time $s$, every difference $\theta=b_1-b_2$ with $b_k\in\mathcal C_{s,k}$ belongs to the fixed regular neighbourhood $\widetilde\Theta_{j^*}$. Together with the bound $\|\bar\beta_{s,k}\|_2\le C_\beta+1$ from \Cref{lem:initial-localisation} and the construction of $p_s$, this ensures that $r_{\bar\beta_s,p_s}(x) = \widetilde x^\top \bar\beta_{s,1}p_s(x)+\widetilde x^\top \bar\beta_{s,2}\{1-p_s(x)\}$  belongs to the deterministic class $\mathcal F_{\mathrm{UCB}}$ defined in \Cref{lem:ucb-bracketing}, so that the empirical-process bound in \Cref{lem:ucb-control-delta} applies. On $\mathcal E_{\rm conf}\cap\mathcal E_{\rm init}$, \eqref{eq:initial-euclidean-radius} gives, for every update time $s$ and every $b_k\in\mathcal C_{s,k}$,
\begin{equation*}
\|(b_1-b_2)-\theta^*\|_2\le\frac{2(1+\sqrt{d+1})\alpha_T}{\sqrt{\lambda+\kappa_XT_o/4}}\le\frac{4(1+\sqrt{d+1})\alpha_T}{\sqrt{\kappa_X}T^{1/4}}.
\end{equation*}
The assumed lower bound on $T$ yields
\begin{equation*}
\frac{4(1+\sqrt{d+1})\alpha_T}{\sqrt{\kappa_X}T^{1/4}}\le\frac{4}{C_3^{1/4}}\le2r_{\rm loc}.
\end{equation*}
Thus, by \Cref{lem:enlarged_theta}, $b_1-b_2\in\widetilde\Theta_{j^*}$ for every update time. Moreover, \eqref{eq:initial-euclidean-radius} and the same lower bound on $T$ give $\|b_k-\beta_k^*\|_2\le r_{\rm loc}\le1$, and hence $\|b_k\|_2\le C_\beta+1$ for every candidate.

Applying \Cref{lem:one-step-regret} with $s=\tau(t)$ and summing from $t=s_0$ to $T$ gives
\begin{align*}
\sum_{t=s_0}^T\{R(p^*)-R(p_{\tau(t)})\}\le{}&2\sum_{t=s_0}^T\Delta_{\tau(t)-1}^{\mathrm{UCB}} + \{(1+\sqrt{d+1})\alpha_T+2B\}\sum_{t=s_0}^T\mathbb E\left[\min\{1,\|\widetilde X_t\|_{V_{\tau(t),A_t}^{-1}}\}\mid\mathcal F_{t-1}\right].
\end{align*}
By the definition of $\mathcal E_{\rm emp}$ in \eqref{eq:main-regret-empirical-event} and \Cref{lem:error-sums}, with $u=\eta/(8T)$, it holds that
\begin{align*}
2\sum_{t=s_0}^T\Delta_{\tau(t)-1}^{\mathrm{UCB}}\le C_5\begin{cases}\sqrt T\log(eT/u),&d=2,\\ \sqrt T\log^{3/2}(eT/u),&d=3,\\ T^{(d-2)/(d-1)}\log(eT/u),&d\ge4.\end{cases}
\end{align*}
Since $u=\eta/(8T)$, we have $\log(eT/u)=\log(8eT^2/\eta)\le C_6\log(T/\eta)$ for $\eta<1/2$, where $C_6=6$. Consequently,
\begin{equation*}
2\sum_{t=s_0}^T\Delta_{\tau(t)-1}^{\mathrm{UCB}}\le C_7\begin{cases}\sqrt T\log(T/\eta),&d=2,\\ \sqrt T\log^{3/2}(T/\eta),&d=3,\\ T^{(d-2)/(d-1)}\log(T/\eta),&d\ge4.\end{cases}
\end{equation*}
For the remaining term, the event $\mathcal E_{\rm mart}$ in \eqref{eq:main-regret-martingale-event} and \Cref{lem:lazy-width-sum-lnu} give
\begin{align*}
\sum_{t=s_0}^T\mathbb E\left[\min\{1,\|\widetilde X_t\|_{V_{\tau(t),A_t}^{-1}}\}\mid\mathcal F_{t-1}\right]\le{}&C_8\sqrt{T(d+1)\log\left(1+\frac{T}{\lambda}\right)} + C_9\sqrt{T\log(8/\eta)}.
\end{align*}
Since $\alpha_T\ge C_\beta\sqrt\lambda$, while $\lambda$ is fixed, we have $(1+\sqrt{d+1})\alpha_T+2B\le C_{10}(1+\sqrt{d+1})\alpha_T$, where $C_{10}=1+2(C_\beta+1)/(C_\beta\sqrt\lambda)$. Therefore, using $\sqrt a+\sqrt b\le\sqrt{2(a+b)}$ and setting $C_{11}=\sqrt2\,C_{10}\max\{C_8,C_9\}$,
\begin{align*}
&\{(1+\sqrt{d+1})\alpha_T+2B\}\sum_{t=s_0}^T\mathbb E\left[\min\{1,\|\widetilde X_t\|_{V_{\tau(t),A_t}^{-1}}\}\mid\mathcal F_{t-1}\right]\\
\le & C_{11}(1+\sqrt{d+1})\alpha_T\sqrt{T\left\{(d+1)\log\left(1+\frac{T}{\lambda}\right)+\log(8/\eta)\right\}}.
\end{align*}
Combining the preceding bounds with the exploration regret proves the bound in \Cref{thm:main-regret}.

\end{proof}

\begin{corollary}\label[corollary]{cor:ucb-regret-rate}
Under the conditions of \Cref{thm:main-regret}, for fixed model constants, it holds that
\[
    \mathfrak R_T=\begin{cases}\widetilde O_{\mathrm p}(\sqrt T),&d=2,3,\\\widetilde O_{\mathrm p}(T^{(d-2)/(d-1)}),&d\ge4.\end{cases}
\]
\end{corollary}
\begin{proof}
By \eqref{eq_ucb_alpha}, $\alpha_T=O(\sqrt{\log(T/\eta)})$ for fixed $d,\lambda,\sigma,C_\beta$. Substituting this bound into \Cref{thm:main-regret} gives the stated rates.
\end{proof}

\subsubsection{Initial exploration and localisation}
\begin{lemma}[Exploration phase]\label[lemma]{lem:initial-localisation}
Suppose \Cref{a_contexts_non_cts,a_noise,a_beta,a_theta} hold. Suppose we play $p_t= 1/2$ in the exploration period when $t\le2T_o$ and select the confidence radius $\alpha_T$ as in \eqref{eq_ucb_alpha}. For $\eta\in(0,1/2)$, if
\begin{equation*}
T_o\ge\frac{16(d+1)}{\kappa_X}\log\left(\frac{16(d+1)}{\eta}\right),
\end{equation*}
then, with probability at least $1-\eta$, for every update time $s$, every $k\in\{1,2\}$, and every $b_k\in\mathcal C_{s,k}$, we have that
\begin{equation}\label{eq:initial-euclidean-radius}
\|b_k-\beta_k^*\|_2\le\frac{(1+\sqrt{d+1})\alpha_T}{\sqrt{\lambda+\kappa_XT_o/4}}.
\end{equation}
If we choose $T_o=\lceil\sqrt T\rceil$, take $r_{\rm loc}\le1$, and suppose $T\ge16\vee C_1(1+\sqrt{d+1})^4\alpha_T^4/(\kappa_X^2r_{\rm loc}^4)$, then, on the same event, every candidate satisfies $\|b_k\|_2\le C_\beta+1$ and every difference $b_1-b_2$ with $b_k\in\mathcal C_{s,k}$ belongs to the fixed regular neighbourhood $\widetilde\Theta_{j^*}$ in \Cref{lem:enlarged_theta}, for every update time $s$.
\end{lemma}

\begin{proof}
For $k\in\{1,2\}$, let $I_k^0=\{i\le2T_o:A_i=k\}$ and $N_{0,k}=|I_k^0|=\sum_{t=1}^{2T_o}\indc\{A_t=k\}$. Define
\begin{align}
\mathcal E_{\rm alloc}&=\{N_{0,1}\ge T_o/2,\ N_{0,2}\ge T_o/2\},\\
\mathcal E_{\rm init}&=\{\lambda_{\min}(V_{s_0,k})\ge\lambda+\kappa_XT_o/4,\ \text{for every }k\in\{1,2\}\}.\label{eq:initial-design-event}
\end{align}
Under $p_t= 1/2$, the exploration actions are independent of the contexts and $N_{0,1}\sim\operatorname{Bin}(2T_o,1/2)$, with $N_{0,2}=2T_o-N_{0,1}$. By Hoeffding's inequality \citep[e.g.~Theorem~2.2.6 in][]{vershynin2018high}, we have~that
\begin{equation*}
P(N_{0,k}<T_o/2)\le\exp\left(-\frac{T_o}{4}\right)\le\exp\left(-\frac{T_o}{8}\right),\qquad k\in\{1,2\}.
\end{equation*}
A union bound over $k\in\{1,2\}$ gives
\begin{equation}\label{eq:initial-allocation-prob}
P(\mathcal E_{\rm alloc}^c)\le2\exp\left(-\frac{T_o}{8}\right).
\end{equation}

Fix $k\in\{1,2\}$ and condition on the exploration actions $A_1,\ldots,A_{2T_o}$. Since the exploration actions are sampled independently of the contexts, the matrices $\{\widetilde X_i\widetilde X_i^\top:i\in I_k^0\}$ remain independent and each has the same distribution as $\widetilde X\widetilde X^\top$. Moreover, $\widetilde X_i\widetilde X_i^\top\succeq0$, $\lambda_{\max}(\widetilde X_i\widetilde X_i^\top)=\|\widetilde X_i\|_2^2\le d+1$, and $\mathbb E(\widetilde X_i\widetilde X_i^\top)=\Sigma_X\succeq\kappa_XI_{d+1}$ by \Cref{a_contexts_non_cts}. Therefore, on $\mathcal E_{\rm alloc}$, we have that
\begin{equation*}
\mu_{\min}:=\lambda_{\min}\left\{\sum_{i\in I_k^0}\mathbb E(\widetilde X_i\widetilde X_i^\top\mid A_1,\ldots,A_{2T_o})\right\}\ge N_{0,k}\kappa_X\ge\frac{\kappa_XT_o}{2}.
\end{equation*}

By the matrix Chernoff inequality \citep[e.g.~Theorem~5.1.1 in][]{tropp2015introduction}, applied with matrix dimension $d+1$, upper eigenvalue bound $d+1$, and $\varepsilon=1/2$, it holds that
\begin{align*}
&P\left\{\lambda_{\min}\left(\sum_{i\in I_k^0}\widetilde X_i\widetilde X_i^\top\right)\le\frac{\mu_{\min}}{2}\ \middle|\ A_1,\ldots,A_{2T_o}\right\}\\
&\qquad\le(d+1)\left(\frac{e^{-1/2}}{(1/2)^{1/2}}\right)^{\mu_{\min}/(d+1)}
\le(d+1)\exp\left\{-\frac{\mu_{\min}}{8(d+1)}\right\}.
\end{align*}
Since $\mu_{\min}\ge\kappa_XT_o/2$ on $\mathcal E_{\rm alloc}$, we have
\begin{equation*}
P\left\{\lambda_{\min}\left(\sum_{i\in I_k^0}\widetilde X_i\widetilde X_i^\top\right)\le\frac{\kappa_XT_o}{4}\ \middle|\ A_1,\ldots,A_{2T_o}\right\}\le(d+1)\exp\left\{-\frac{\kappa_XT_o}{16(d+1)}\right\}.
\end{equation*}
Recall from \Cref{alg:lazy-global-model-ucb} that $V_{s_0,k}=\lambda I_{d+1}+\sum_{i\in I_k^0}\widetilde X_i\widetilde X_i^\top$. Therefore, a union bound over $k\in\{1,2\}$ and tower property give
\begin{equation}\label{eq:initial-design-prob}
P(\mathcal E_{\rm init}^c\cap\mathcal E_{\rm alloc})\le2(d+1)\exp\left\{-\frac{\kappa_XT_o}{16(d+1)}\right\}.
\end{equation}
Combining \eqref{eq:initial-allocation-prob} and \eqref{eq:initial-design-prob}, we obtain
\begin{equation*}
P(\mathcal E_{\rm init}^c)\le2\exp\left(-\frac{T_o}{8}\right)+2(d+1)\exp\left\{-\frac{\kappa_XT_o}{16(d+1)}\right\}\le\frac{\eta}{4},
\end{equation*}
where the last inequality follows whenever
\begin{equation*}
T_o\ge8\log\left(\frac{16}{\eta}\right)\qquad\text{and}\qquad T_o\ge\frac{16(d+1)}{\kappa_X}\log\left(\frac{16(d+1)}{\eta}\right).
\end{equation*}
Since $e_1^\top\Sigma_Xe_1=1$, we have $\kappa_X\le1$, and therefore the assumed lower bound on $T_o$ implies both inequalities.

Recall the confidence event $\mathcal E_{\rm conf}$ in \eqref{eq:main-regret-confidence-event}. By the same concentration argument following that display, $P(\mathcal E_{\rm conf})\ge1-\eta/8$. A union bound therefore gives
\begin{equation*}
P(\mathcal E_{\rm init}\cap\mathcal E_{\rm conf})\ge1-\frac{\eta}{4}-\frac{\eta}{8}=1-\frac{3\eta}{8}\ge1-\eta.
\end{equation*}
The rest of the proof is carried out conditioning on $\mathcal E_{\rm init}\cap\mathcal E_{\rm conf}$ happening. By \Cref{lem:polyhedral-confidence}, $\beta_k^*\in\mathcal C_{s,k}$ and $\|b_k-\beta_k^*\|_{V_{s,k}}\le(1+\sqrt{d+1})\alpha_T$ for every update time $s$ and every $b_k\in\mathcal C_{s,k}$. Since $V_{s,k}\succeq V_{s_0,k}$ for every $s\ge s_0$, we also have $V_{s,k}\succeq(\lambda+\kappa_XT_o/4)I_{d+1}$. Therefore,
\begin{equation*}
\sqrt{\lambda+\kappa_XT_o/4}\,\|b_k-\beta_k^*\|_2\le\|b_k-\beta_k^*\|_{V_{s,k}}\le(1+\sqrt{d+1})\alpha_T,
\end{equation*}
which proves \eqref{eq:initial-euclidean-radius} for every update time $s$.

Finally, for $\theta=b_1-b_2$, the triangle inequality and \eqref{eq:initial-euclidean-radius} give
\begin{equation*}
\|\theta-\theta^*\|_2\le\frac{2(1+\sqrt{d+1})\alpha_T}{\sqrt{\lambda+\kappa_XT_o/4}}.
\end{equation*}
If $T_o=\lceil\sqrt T\rceil$, then we have that
\begin{equation*}
\frac{2(1+\sqrt{d+1})\alpha_T}{\sqrt{\lambda+\kappa_XT_o/4}}\le\frac{4(1+\sqrt{d+1})\alpha_T}{\sqrt{\kappa_X}T^{1/4}}.
\end{equation*}
Hence, whenever
\begin{equation*}
T\ge\frac{16(1+\sqrt{d+1})^4\alpha_T^4}{\kappa_X^2r_{\rm loc}^4},
\end{equation*}
we have $\|\theta-\theta^*\|_2\le2r_{\rm loc}$. It therefore follows from \Cref{lem:enlarged_theta} that $\theta\in\widetilde\Theta_{j^*}$, for every update time $s$ on the same event. Under the same lower bound on $T$, \eqref{eq:initial-euclidean-radius} gives $\|b_k-\beta_k^*\|_2\le r_{\rm loc}\le1$, and hence $\|b_k\|_2\le C_\beta+1$.

\end{proof}

\subsubsection{Complexity under arbitrary context laws}\label[appendix]{sec_app_ucb_complexity}
Throughout this subsubsection, we let $s_\theta(x)=\widetilde x^\top\theta$.
The geometric results in \Cref{lem:hausdorff} and the trace-entropy bound in \Cref{lem:trace_entropy_weak} depend only on the parameter and Lipschitz assumptions, so they apply to arbitrary context laws. However, the slab-probability bound in \Cref{lem:slab_probability} uses a bounded density and cannot be used here. We instead control policy differences pointwise, including on the sign-disagreement set, and obtain entropy bounds uniformly over all probability measures. This uniformity allows us to use the empirical measure $P_n$ in the conditional Rademacher argument below. 

\begin{lemma}[Localisation of sign disagreement]
\label[lemma]{lem:disagreement-localisation}
Suppose \Cref{a_theta} holds. Then there exists a constant $C_{\rm loc}<\infty$, depending only on $d,\nu$ and the constants in \Cref{lem:enlarged_theta}, such that, for every $j\in\{1,\ldots,d\}$ and all $\theta,\theta'\in\widetilde\Theta_j$,
\[
\sup_{x\in S(\theta,\theta')}\left\{\dist_\nu(x,\mathcal H_\theta)\vee\dist_\nu(x,\mathcal H_{\theta'})\right\}\le C_{\rm loc}\|\theta-\theta'\|_2,
\]
where
\[
S(\theta,\theta')=\{x\in[0,1]^d:s_\theta(x)\ge0,\ s_{\theta'}(x)<0\}\cup\{x\in[0,1]^d:s_\theta(x)<0,\ s_{\theta'}(x)\ge0\}.
\]
\end{lemma}

\begin{proof}
Fix $j\in\{1,\ldots,d\}$. By relabelling coordinates, without loss of generality, we assume that $j=d$. Write $\theta=(a,\gamma)$ and $\theta'=(a',\gamma')$, where $a,a'\in\mathbb R$ and $\gamma,\gamma'\in\mathbb R^{d}$. Then write $\gamma=(b,c)$ and $\gamma'=(b',c')$, where $b,b'\in\mathbb R^{d-1}$ and $c,c'\in\mathbb R$. For $u\in[0,1]^{d-1}$, define
\[
\phi_\theta(u)=-\frac{a+b^\top u}{c},\qquad \phi_{\theta'}(u)=-\frac{a'+{b'}^\top u}{c'}.
\]
Note that the calculation in \eqref{lem:hausdorff_eq1} gives
\begin{equation}\label{lem:disagreement-localisation_eq1}
\|\phi_\theta-\phi_{\theta'}\|_\infty\le C_\phi\|\theta-\theta'\|_2.
\end{equation}
Write $\delta=\|\theta-\theta'\|_2$ and set $\eta=C_\phi\delta$. Let $C_{d,\nu}=d^{\max\{1/\nu-1/2,0\}}$ be the constant in \Cref{lem:norm-comparison}.

Take any $x=(u,z)\in S(\theta,\theta')$. Since $s_\theta(u,z)=c\{z-\phi_\theta(u)\}$ and $s_{\theta'}(u,z)=c'\{z-\phi_{\theta'}(u)\}$ with $c,c'>0$, the value $z$ lies between $\phi_\theta(u)$ and $\phi_{\theta'}(u)$. Therefore, by \eqref{lem:disagreement-localisation_eq1},
\begin{equation}\label{lem:disagreement-localisation_eq2}
|z-\phi_\theta(u)|\vee|z-\phi_{\theta'}(u)|\le\eta.
\end{equation}
In the rest of the proof, we show that $x$ is within $\ell_\nu$-distance $C\eta$ of $\mathcal H_\theta$. The bound for $\mathcal H_{\theta'}$ then follows by symmetry.

\noindent \textbf{Case 1:} If $0\le\phi_\theta(u)\le1$, then $h=(u,\phi_\theta(u))\in\mathcal H_\theta$, and hence
\[
\|x-h\|_\nu=|z-\phi_\theta(u)|\le\eta.
\]

\noindent \textbf{Case 2:} Next suppose $\phi_\theta(u)<0$. By the definition of $\widetilde\Theta_j$, there exists $u_\theta\in[0,1]^{d-1}$ such that $\widetilde\rho_0\le\phi_\theta(u_\theta)\le1-\widetilde\rho_0$. For $\lambda\in[0,1]$, define $u_\lambda=(1-\lambda)u+\lambda u_\theta$. Since $\phi_\theta$ is affine, the map $\lambda\mapsto\phi_\theta(u_\lambda)$ is affine. Since $\phi_\theta(u)<0$ and $\phi_\theta(u_\theta)\ge\widetilde\rho_0$, there exists $\lambda_0\in[0,1]$ such that $\phi_\theta(u_{\lambda_0})=0$. Moreover, by \eqref{lem:disagreement-localisation_eq2} and $z\ge0$, we have that
\[
\lambda_0=\frac{-\phi_\theta(u)}{\phi_\theta(u_\theta)-\phi_\theta(u)}\le\frac{\eta}{\widetilde\rho_0+\eta}\le\frac{\eta}{\widetilde\rho_0},\qquad z\le\eta.
\]
Set $v=u_{\lambda_0}$ and $h=(v,0)$. Then $h\in\mathcal H_\theta$. We first bound the distance in the Euclidean norm and then convert it to the fairness norm. Since $u,u_\theta\in[0,1]^{d-1}$,
\[
\|v-u\|_2\le\lambda_0\|u_\theta-u\|_2\le\frac{\sqrt{d-1}}{\widetilde\rho_0}\eta.
\]
By \Cref{lem:norm-comparison},
\[
\|v-u\|_\nu\le C_{d,\nu}\|v-u\|_2\le\frac{C_{d,\nu}\sqrt{d-1}}{\widetilde\rho_0}\eta.
\]
Hence, by the triangle inequality,
\[
\|x-h\|_\nu\le\|u-v\|_\nu+z\le\left(1+\frac{C_{d,\nu}\sqrt{d-1}}{\widetilde\rho_0}\right)\eta.
\]

\noindent \textbf{Case 3:} Finally suppose $\phi_\theta(u)>1$. Again, by the definition of $\widetilde\Theta_j$, there exists $u_\theta\in[0,1]^{d-1}$ such that $\widetilde\rho_0\le\phi_\theta(u_\theta)\le1-\widetilde\rho_0$. Along the same segment $u_\lambda=(1-\lambda)u+\lambda u_\theta$, there exists $\lambda_0\in[0,1]$ such that $\phi_\theta(u_{\lambda_0})=1$. Moreover, by \eqref{lem:disagreement-localisation_eq2} and $z\le1$,
\[
\lambda_0=\frac{\phi_\theta(u)-1}{\phi_\theta(u)-\phi_\theta(u_\theta)}\le\frac{\eta}{\widetilde\rho_0+\eta}\le\frac{\eta}{\widetilde\rho_0},\qquad 1-z\le\eta.
\]
Set $v=u_{\lambda_0}$ and $h=(v,1)$. Then $h\in\mathcal H_\theta$. As in \textbf{Case 2},
\[
\|v-u\|_\nu\le C_{d,\nu}\|v-u\|_2\le\frac{C_{d,\nu}\sqrt{d-1}}{\widetilde\rho_0}\eta,
\]
where the first inequality follows from \Cref{lem:norm-comparison}. Hence
\[
\|x-h\|_\nu\le\|u-v\|_\nu+|z-1|\le\left(1+\frac{C_{d,\nu}\sqrt{d-1}}{\widetilde\rho_0}\right)\eta.
\]

Combining \textbf{Cases 1--3} gives $\dist_\nu(x,\mathcal H_\theta)\le C\eta$. Applying the same argument with $\theta'$ in place of $\theta$ gives $\dist_\nu(x,\mathcal H_{\theta'})\le C\eta$. Since $\eta=C_\phi\|\theta-\theta'\|_2$, the claim follows.
\end{proof}

\begin{lemma}
\label[lemma]{lem:mcshane-envelope-stability}
Suppose \Cref{a_theta} holds. Fix $j\in\{1,\ldots,d\}$, $\theta,\theta'\in\widetilde\Theta_j$ and $q\in\mathcal Q_{\theta,L}$. Define the global McShane extension of $q$ from $\mathcal H_\theta$ to $[0,1]^d$ by
\[
G(x)=\inf_{h\in\mathcal H_\theta}\{q(h)+L\|x-h\|_\nu\},\qquad \bar q(x)=[G(x)]_{[0,1]},\qquad x\in[0,1]^d,
\]
where $[y]_{[0,1]}=\min\{1,\max\{0,y\}\}$. Define the transferred trace on $\mathcal H_{\theta'}$ by $q'(h')=\bar q(h')$ for all $h'\in\mathcal H_{\theta'}$. Then $q'\in\mathcal Q_{\theta',L}$ and with $U_{\theta,q}$
and $M_{\theta,q}$ defined in \eqref{eq:envelope-definitions},
\[
\|U_{\theta,q}-U_{\theta',q'}\|_\infty\vee\|M_{\theta,q}-M_{\theta',q'}\|_\infty\le2LC_H\|\theta-\theta'\|_2,
\]
\end{lemma}

\begin{proof}
The construction is the same as in \textbf{Steps 1--3} of the proof of \Cref{lem:expected-transfer-no-common-chart}, with $(\theta^*,q^*)$ replaced by $(\theta,q)$. Those steps use only the Lipschitz property of the trace, nonexpansiveness of clipping, and the Hausdorff bound in \Cref{lem:hausdorff}, and they do not rely on the context density or optimality of the trace. Thus $q'\in\mathcal Q_{\theta',L}$ and both envelope differences are bounded by $2LC_H\|\theta-\theta'\|_2$, as claimed.
\end{proof}

\begin{lemma}[Gap between the clipped canonical envelopes]
\label[lemma]{lem:clipped-envelope-gap}
For every $\theta\in\widetilde\Theta$, $q\in\mathcal Q_{\theta,L}$ and $x\in[0,1]^d$, define $\overline U_{\theta,q}(x)=1\wedge U_{\theta,q}(x)$ and $\overline M_{\theta,q}(x)=0\vee M_{\theta,q}(x)$. Then, we have that
\[
0\le\overline U_{\theta,q}(x)-\overline M_{\theta,q}(x)\le2L\dist_\nu(x,\mathcal H_\theta).
\]
\end{lemma}

\begin{proof}
For any $h,g\in\mathcal H_\theta$, we have that
\[
q(h)-L\|x-h\|_\nu\le q(g)+L\|h-g\|_\nu-L\|x-h\|_\nu\le q(g)+L\|x-g\|_\nu.
\]
Taking the supremum over $h$ followed by the infimum over $g$ gives $M_{\theta,q}(x)\le U_{\theta,q}(x)$. Moreover, since $U_{\theta,q}\ge0$ and $M_{\theta,q}\le1$, clipping preserves the ordering and we have that
\[
0\le\overline U_{\theta,q}(x)-\overline M_{\theta,q}(x)\le U_{\theta,q}(x)-M_{\theta,q}(x).
\]

Moreover, since set $\mathcal H_\theta$ is nonempty and compact, let $h_x\in\mathcal H_\theta$ attain $\dist_\nu(x,\mathcal H_\theta)$. Then
\[
U_{\theta,q}(x)\le q(h_x)+L\dist_\nu(x,\mathcal H_\theta),\qquad M_{\theta,q}(x)\ge q(h_x)-L\dist_\nu(x,\mathcal H_\theta).
\]
Combining the preceding bounds proves the result.
\end{proof}

\begin{lemma}[Uniform transfer in the sup-norm]
\label[lemma]{lem:uniform-transfer-sup}
Suppose \Cref{a_theta} holds. Then there exists a constant $C_{\rm tr}<\infty$, depending only on $d,\nu,L$ and the constants in \Cref{lem:enlarged_theta}, such that, for every $j\in\{1,\ldots,d\}$, every $\theta,\theta'\in\widetilde\Theta_j$ and every $q\in\mathcal Q_{\theta,L}$, there exists $q'\in\mathcal Q_{\theta',L}$ satisfying
\[
\|p_{\theta,q}-p_{\theta',q'}\|_\infty\le C_{\rm tr}\|\theta-\theta'\|_2.
\]
\end{lemma}

\begin{proof}
Fix $j\in\{1,\ldots,d\}$, $\theta,\theta'\in\widetilde\Theta_j$ and $q\in\mathcal Q_{\theta,L}$. Write $\delta=\|\theta-\theta'\|_2$, and let $q'$ be the transferred trace constructed in \Cref{lem:mcshane-envelope-stability}. Define the sign-disagreement set
\[
S(\theta,\theta')=\{x\in[0,1]^d:s_\theta(x)\ge0,\ s_{\theta'}(x)<0\}\cup\{x\in[0,1]^d:s_\theta(x)<0,\ s_{\theta'}(x)\ge0\}.
\]

\noindent \textbf{Step 1: Control on $S(\theta,\theta')^c$.} On $S(\theta,\theta')^c$, the two policies use the same type of clipped envelope; either both use the upper envelopes or both use the lower envelopes. Since the maps $y\mapsto\min\{1,y\}$ and $y\mapsto\max\{0,y\}$ are nonexpansive, \Cref{lem:mcshane-envelope-stability} implies that
\begin{equation}\label{lem:uniform-transfer-sup_eq1}
|p_{\theta,q}(x)-p_{\theta',q'}(x)|\le2LC_H\delta,\qquad x\in S(\theta,\theta')^c.
\end{equation}

\noindent \textbf{Step 2: Control on $S(\theta,\theta')$.}

\noindent \textbf{Case 1:} Suppose $s_\theta(x)\ge0$ and $s_{\theta'}(x)<0$. Then $p_{\theta,q}(x)=\overline U_{\theta,q}(x)$ and $p_{\theta',q'}(x)=\overline M_{\theta',q'}(x)$. By Lemmas \ref{lem:disagreement-localisation}, \ref{lem:mcshane-envelope-stability}, and \ref{lem:clipped-envelope-gap}, we have that
\begin{align*}
|p_{\theta,q}(x)-p_{\theta',q'}(x)|&=|\overline U_{\theta,q}(x)-\overline M_{\theta',q'}(x)|\\
&\le|\overline U_{\theta,q}(x)-\overline M_{\theta,q}(x)|+|\overline M_{\theta,q}(x)-\overline M_{\theta',q'}(x)|\\
&\le2L\dist_\nu(x,\mathcal H_\theta)+2LC_H\delta\\
&\le2L(C_{\rm loc}+C_H)\delta.
\end{align*}

\noindent \textbf{Case 2:} Suppose $s_\theta(x)<0$ and $s_{\theta'}(x)\ge0$. Then $p_{\theta,q}(x)=\overline M_{\theta,q}(x)$ and $p_{\theta',q'}(x)=\overline U_{\theta',q'}(x)$. By the same argument as \textbf{Case 1}, we have that
\begin{align*}
|p_{\theta,q}(x)-p_{\theta',q'}(x)|&=|\overline M_{\theta,q}(x)-\overline U_{\theta',q'}(x)|\\
&\le|\overline M_{\theta,q}(x)-\overline U_{\theta,q}(x)|+|\overline U_{\theta,q}(x)-\overline U_{\theta',q'}(x)|\\
&\le2L\dist_\nu(x,\mathcal H_\theta)+2LC_H\delta\\
&\le2L(C_{\rm loc}+C_H)\delta.
\end{align*}
Combining \textbf{Cases 1 and 2} with \eqref{lem:uniform-transfer-sup_eq1} proves the claim with $C_{\rm tr}=2L(C_{\rm loc}+C_H)$.
\end{proof}

\begin{lemma}[Sup-norm entropy of the canonical policy class]
\label[lemma]{lem:sup-entropy-policy}
Suppose \Cref{a_theta} holds and $d\ge2$. Let $\mathcal P_L=\{p_{\theta,q}:\theta\in\widetilde\Theta,\ q\in\mathcal Q_{\theta,L}\}$. Then there exists a constant $C_{\mathcal P}<\infty$, depending only on $d,\nu,L$ and the constants in \Cref{lem:enlarged_theta}, such that, for every $\varepsilon\in(0,1)$,
\[
\log N(\varepsilon,\mathcal P_L,\|\cdot\|_\infty)\le C_{\mathcal P}\varepsilon^{-(d-1)}\log(C_{\mathcal P}/\varepsilon).
\]
\end{lemma}

\begin{proof}
Write $m=d-1$. If $L=0$, then every $q\in\mathcal Q_{\theta,0}$ is constant on $\mathcal H_\theta$, the upper and lower canonical envelopes equal that constant on $[0,1]^d$, and hence $\mathcal P_0$ is the class of constant functions with values in $[0,1]$. Discretising $[0,1]$ proves the result. Suppose henceforth that $L>0$.

For each $r\in\{1,\ldots,d\}$, define
\[
\mathcal P_{L,r}=\{p_{\theta,q}:\theta\in\widetilde\Theta_r,\ q\in\mathcal Q_{\theta,L}\}.
\]
Then $\mathcal P_L=\bigcup_{r=1}^d\mathcal P_{L,r}$. We first bound the covering number of each fixed chart $\mathcal P_{L,r}$.

Fix $r\in\{1,\ldots,d\}$ and $\varepsilon\in(0,1)$. Let $\delta_\theta=\varepsilon/(2C_{\rm tr})$. Let $\{\theta_1,\ldots,\theta_M\}$ be a $\delta_\theta$-net of $\widetilde\Theta_r$ with respect to $\|\cdot\|_2$. Since $\widetilde\Theta_r\subset\mathbb R^{d+1}$ is bounded by the definition of $\widetilde\Theta_r$, we have that
\[
\log M\le(d+1)\log(C'/\delta_\theta)\le C\log(C/\varepsilon),
\]
where $C$ is a constant depending on $d$.

For each fixed representative boundary $\mathcal H_{\theta_j}$, \Cref{lem:trace_entropy_weak} gives a finite $\varepsilon/2$-cover of $\mathcal Q_{\theta_j,L}$ in the sup-norm. Thus there exists a finite set $\mathcal Q_j\subset\mathcal Q_{\theta_j,L}$ such that, for every $q^\circ\in\mathcal Q_{\theta_j,L}$, there exists $\widetilde q\in\mathcal Q_j$ satisfying $\|q^\circ-\widetilde q\|_\infty\le\varepsilon/2$. Moreover, it holds that $\log|\mathcal Q_j|\le C\varepsilon^{-m}\log(C/\varepsilon)$.

Now fix any $p\in\mathcal P_{L,r}$, and we can write $p=p_{\theta,q}$ for some $\theta\in\widetilde\Theta_r$ and $q\in\mathcal Q_{\theta,L}$. Choose $\theta_j$ such that $\|\theta-\theta_j\|_2\le\delta_\theta$. Since $\theta$ and $\theta_j$ lie in the same chart $\widetilde\Theta_r$, by \Cref{lem:uniform-transfer-sup}, there exists $q^\circ\in\mathcal Q_{\theta_j,L}$ such that
\begin{equation}\label{lem:sup-entropy-policy_eq1}
\|p_{\theta,q}-p_{\theta_j,q^\circ}\|_\infty\le C_{\rm tr}\delta_\theta=\varepsilon/2.
\end{equation}
Choose $\widetilde q\in\mathcal Q_j$ such that $\|q^\circ-\widetilde q\|_\infty\le\varepsilon/2$. Since the envelope map is nonexpansive in the boundary trace, we have that
\begin{equation}\label{lem:sup-entropy-policy_eq2}
\|p_{\theta_j,q^\circ}-p_{\theta_j,\widetilde q}\|_\infty\le\varepsilon/2.
\end{equation}
Therefore, by \eqref{lem:sup-entropy-policy_eq1}, \eqref{lem:sup-entropy-policy_eq2} and the triangle inequality, the policies $p_{\theta_j,\widetilde q}$ form an $\varepsilon$-cover of $\mathcal P_{L,r}$.

Thus, for every fixed $r\in\{1,\ldots,d\}$, we have that
\[
\log N(\varepsilon,\mathcal P_{L,r},\|\cdot\|_\infty)\le C\log(C/\varepsilon)+C\varepsilon^{-m}\log(C/\varepsilon)\le C_{\mathcal P}\varepsilon^{-m}\log(C_{\mathcal P}/\varepsilon).
\]
Since $\mathcal P_L=\bigcup_{r=1}^d\mathcal P_{L,r}$, the union of the chart-wise covers is an $\varepsilon$-cover of $\mathcal P_L$. Hence,
\[
N(\varepsilon,\mathcal P_L,\|\cdot\|_\infty)\le\sum_{r=1}^dN(\varepsilon,\mathcal P_{L,r},\|\cdot\|_\infty).
\]
Taking logarithms and absorbing $\log(d)$ into the constant prove the claim.
\end{proof}

\begin{lemma}[$L_2(P)$ bracketing entropy]
\label[lemma]{lem:ucb-bracketing}
Suppose Assumptions \ref{a_beta} and \ref{a_theta} hold and $d\ge2$. For $b=(b_1,b_2)$ and a policy $p:[0,1]^d\to[0,1]$, define
\[
r_{b,p}(x)=\widetilde x^\top b_1p(x)+\widetilde x^\top b_2\{1-p(x)\}.
\]
Let
\[
\mathcal F_{\mathrm{UCB}}=\left\{x\mapsto r_{b,p_{\theta,q}}(x):\|b_k\|_2\le C_\beta+1,\ k\in\{1,2\},\ \theta=b_1-b_2\in\widetilde\Theta,\ q\in\mathcal Q_{\theta,L}\right\}.
\]
Let $P$ be an arbitrary Borel probability measure supported on $[0,1]^d$. Then there exists a constant $C_{\mathcal F}<\infty$, depending only on $d,\nu,L,C_\beta$ and the constants in \Cref{lem:enlarged_theta}, and independent of $P$, such that, for every $\varepsilon\in(0,1)$,
\[
\log N_{[]}(\varepsilon,\mathcal F_{\mathrm{UCB}},L_2(P))\le C_{\mathcal F}\varepsilon^{-(d-1)}\log(C_{\mathcal F}/\varepsilon),
\]
where $N_{[]}(\varepsilon,\mathcal F_{\mathrm{UCB}},L_2(P))$ denotes the $L_2(P)$ bracketing number.
\end{lemma}

\begin{proof}
Write $m=d-1$, $R_X=\sup_{x\in[0,1]^d}\|\widetilde x\|_2=\sqrt{d+1}$ and $B_\Theta=1\vee\sup_{\theta\in\widetilde\Theta}\|\theta\|_2$. By \Cref{lem:enlarged_theta}, it holds that $B_\Theta<\infty$. Let
\[
\delta_b=\frac{\varepsilon}{6R_X},\qquad \delta_\theta=\frac{\varepsilon}{6R_X},\qquad \delta_p=\frac{\varepsilon}{6R_XB_\Theta}.
\]

Let $\mathcal B_\beta=\{b\in\mathbb R^{d+1}:\|b\|_2\le C_\beta + 1\}$. Since $\mathcal B_\beta$ is a bounded subset of $\mathbb R^{d+1}$, there exists a $\delta_b$-net $\{v_1,\ldots,v_{N_b}\}$ of $\mathcal B_\beta$ under $\|\cdot\|_2$ satisfying
\[
\log N_b\le(d+1)\log(C/\delta_b)\le C\log(C/\varepsilon).
\]
Similarly, since $\widetilde\Theta$ is a bounded subset of $\mathbb R^{d+1}$, there exists a $\delta_\theta$-net $\{\vartheta_1,\ldots,\vartheta_{N_\theta}\}$ of $\widetilde\Theta$ under $\|\cdot\|_2$ satisfying
\[
\log N_\theta\le(d+1)\log(C/\delta_\theta)\le C\log(C/\varepsilon).
\]
By \Cref{lem:sup-entropy-policy}, there exists a $\delta_p$-cover $\{\pi_1,\ldots,\pi_{N_p}\}$ of $\mathcal P_L$ under the sup-norm with
\[
\log N_p\le C\delta_p^{-m}\log(C/\delta_p)\le C\varepsilon^{-m}\log(C/\varepsilon).
\]

Now fix any $f\in\mathcal F_{\mathrm{UCB}}$. We can write
\[
f(x)=\widetilde x^\top b_2+\widetilde x^\top\theta p_{\theta,q}(x)
\]
for some $b_1,b_2\in\mathcal B_\beta$, $\theta=b_1-b_2\in\widetilde\Theta$ and $q\in\mathcal Q_{\theta,L}$. Choose $v_\ell$, $\vartheta_j$ and $\pi_k$ such that
\[
\|b_2-v_\ell\|_2\le\delta_b,\qquad \|\theta-\vartheta_j\|_2\le\delta_\theta,\qquad \|p_{\theta,q}-\pi_k\|_\infty\le\delta_p.
\]
Define
\[
f_{\ell,j,k}(x)=\widetilde x^\top v_\ell+\widetilde x^\top\vartheta_j\pi_k(x).
\]
Since $0\le p_{\theta,q}\le1$, for every $x\in[0,1]^d$ we have that
\begin{align*}
|f(x)-f_{\ell,j,k}(x)|
&\le|\widetilde x^\top(b_2-v_\ell)|+|\widetilde x^\top(\theta-\vartheta_j)|p_{\theta,q}(x)+|\widetilde x^\top\vartheta_j|\,|p_{\theta,q}(x)-\pi_k(x)|\\
&\le R_X\delta_b+R_X\delta_\theta+R_XB_\Theta\delta_p =\frac{\varepsilon}{2}.
\end{align*}

For each triple $(\ell,j,k)$, define the bracket
\[
\ell_{\ell,j,k}(x)=f_{\ell,j,k}(x)-\varepsilon/2
\]
and
\[
u_{\ell,j,k}(x)=f_{\ell,j,k}(x)+\varepsilon/2.
\]
For every $f\in\mathcal F_{\mathrm{UCB}}$, there exists a triple $(\ell,j,k)$ such that
\[
\ell_{\ell,j,k}(x)\le f(x)\le u_{\ell,j,k}(x)
\]
for every $x\in[0,1]^d$. Moreover, the $L_2(P)$ width of this bracket satisfies
\[
\|u_{\ell,j,k}-\ell_{\ell,j,k}\|_{L_2(P)}=\varepsilon.
\]
Thus the brackets indexed by $(\ell,j,k)$ form an $L_2(P)$ $\varepsilon$-bracketing cover of $\mathcal F_{\mathrm{UCB}}$. Since the number of brackets is at most $N_bN_\theta N_p$, we have that
\begin{align*}
\log N_{[]}(\varepsilon,\mathcal F_{\mathrm{UCB}},L_2(P))
&\le\log N_b+\log N_\theta+\log N_p\\
&\le C\log(C/\varepsilon)+C\varepsilon^{-m}\log(C/\varepsilon)\\
&\le C_{\mathcal F}\varepsilon^{-m}\log(C_{\mathcal F}/\varepsilon).
\end{align*}
Since the construction is independent of $P$, the claim follows.
\end{proof}

\begin{lemma}[Conditional Rademacher complexity]
\label[lemma]{lem:rademacher-policy}
Suppose the conditions of \Cref{lem:ucb-bracketing} hold. For arbitrary deterministic points $x_1,\ldots,x_n\in[0,1]^d$, let $P_n=n^{-1}\sum_{i=1}^n\delta_{x_i}$ be the empirical measure and define
\[
\widehat{\mathfrak R}_n(\mathcal F_{\mathrm{UCB}})=\mathbb E_\sigma\left[\sup_{f\in\mathcal F_{\mathrm{UCB}}}\left|\frac1n\sum_{i=1}^n\sigma_i f(x_i)\right|\right],
\]
where $\sigma_1,\ldots,\sigma_n$ are independent Rademacher random variables. Then there exists a constant $C<\infty$, depending only on $d,\nu,L,C_\beta$ and the constants in \Cref{lem:enlarged_theta}, and independent of $x_1,\ldots,x_n$, such that, with $m=d-1$,
\[
\widehat{\mathfrak R}_n(\mathcal F_{\mathrm{UCB}})\le C\begin{cases}n^{-1/2},&m=1,\\n^{-1/2}\log^{3/2}(en),&m=2,\\\{\log(en)/n\}^{1/m},&m\ge3.\end{cases}
\]
\end{lemma}

\begin{proof}
For $n=1$, the bound follows from $\sup_{f\in\mathcal F_{\mathrm{UCB}}}\|f\|_\infty\le\sqrt{d+1}C_\beta$, so assume $n\ge2$. By \Cref{lem:l2-bracketing-policy}, selecting one function from every nonempty bracket gives, for every probability measure $Q$ supported on $[0,1]^d$,
\begin{equation}\label{lem:rademacher-policy_eq1}
\log N(\varepsilon,\mathcal F_{\mathrm{UCB}},L_2(Q))\le C_1\varepsilon^{-m}\log(C_1/\varepsilon),\qquad\varepsilon\in(0,1).
\end{equation}
In particular, this holds for $Q=P_n$. Let $\mathcal F_{\mathrm{UCB}}^\pm=\mathcal F_{\mathrm{UCB}}\cup(-\mathcal F_{\mathrm{UCB}})\cup\{0\}$ and consider the Rademacher process $Z_f=n^{-1}\sum_{i=1}^n\sigma_i f(x_i)$ indexed by $\mathcal F_{\mathrm{UCB}}^\pm$. Conditionally on $x_1,\ldots,x_n$, its increments satisfy
\begin{align*}
\mathbb E_\sigma\exp\{\lambda(Z_f-Z_g)\}
&=\mathbb E_\sigma\exp\left\{\frac{\lambda}{n}\sum_{i=1}^n\sigma_i\{f(x_i)-g(x_i)\}\right\}\\
&=\prod_{i=1}^n\mathbb E_{\sigma_i}\exp\left\{\frac{\lambda}{n}\sigma_i\{f(x_i)-g(x_i)\}\right\}\\
&\le\prod_{i=1}^n\exp\left\{\frac{\lambda^2\{f(x_i)-g(x_i)\}^2}{2n^2}\right\}\\
&=\exp\left\{\frac{\lambda^2}{2n}\|f-g\|_{L_2(P_n)}^2\right\},
\end{align*}
where the second equality follows from the independence of $\sigma_1,\ldots,\sigma_n$. Thus, the process is sub-Gaussian under the pseudometric $d_n(f,g)=\|f-g\|_{L_2(P_n)}/\sqrt n$.

Since $\mathcal F_{\mathrm{UCB}}^\pm$ is symmetric and contains the zero function, we have
\[
\mathbb E_\sigma\left[\sup_{f,g\in\mathcal F_{\mathrm{UCB}}^\pm}(Z_f-Z_g)\right]=2\mathbb E_\sigma\left[\sup_{f\in\mathcal F_{\mathrm{UCB}}^\pm}Z_f\right] \geq 2\widehat{\mathfrak R}_n(\mathcal F_{\mathrm{UCB}}).
\]
Let $D_n=\operatorname{diam}(\mathcal F_{\mathrm{UCB}}^\pm,d_n)$. If $d_n(f,g)\le\delta$, the Cauchy--Schwarz inequality gives
\[
|Z_f-Z_g|\le\|f-g\|_{L_2(P_n)}=\sqrt n\,d_n(f,g)\le\sqrt n\,\delta.
\]
Therefore, the truncated Dudley entropy bound \citep[e.g.~Theorem~5.22 in][]{wainwright2019high} implies
\begin{align*}
2\widehat{\mathfrak R}_n(\mathcal F_{\mathrm{UCB}})
&\le2\sqrt n\,\delta+32\int_{\delta/4}^{D_n}\sqrt{\log N(u,\mathcal F_{\mathrm{UCB}}^\pm,d_n)}\,\mathrm du\\
&=2\sqrt n\,\delta+\frac{32}{\sqrt n}\int_{\sqrt n\delta/4}^{\sqrt nD_n}\sqrt{\log N(\varepsilon,\mathcal F_{\mathrm{UCB}}^\pm,L_2(P_n))}\,\mathrm d\varepsilon,
\end{align*}
where the equality follows from $N(u,\mathcal F_{\mathrm{UCB}}^\pm,d_n)=N(\sqrt n\,u,\mathcal F_{\mathrm{UCB}}^\pm,L_2(P_n))$, and the change of variables $\varepsilon=\sqrt n\,u$. Moreover, since
\[
N(\varepsilon,\mathcal F_{\mathrm{UCB}}^\pm,L_2(P_n))\le2N(\varepsilon,\mathcal F_{\mathrm{UCB}},L_2(P_n))+1,
\]
and every $f\in\mathcal F_{\mathrm{UCB}}$ satisfies
\[
\|f\|_\infty\le\sup_{x\in[0,1]^d}\|\widetilde x\|_2C_\beta=\sqrt{d+1}C_\beta,
\]
we have that $\sqrt nD_n=\operatorname{diam}(\mathcal F_{\mathrm{UCB}}^\pm,L_2(P_n))$ is bounded by a fixed constant. Therefore, setting $\alpha=\sqrt n\delta/4$, rescaling the fixed upper endpoint to one, and absorbing universal constants, \eqref{lem:rademacher-policy_eq1} yields
\begin{equation}\label{lem:rademacher-policy_eq2}
\widehat{\mathfrak R}_n(\mathcal F_{\mathrm{UCB}})\le C_2\inf_{\alpha\in(0,1)}\left\{\alpha+\frac1{\sqrt n}\int_\alpha^1\varepsilon^{-m/2}\{\log(C_1/\varepsilon)\}^{1/2}\,\mathrm d\varepsilon\right\}.
\end{equation}

\noindent \textbf{Case 1: $m=1$.} In this case,
\begin{align*}
\int_\alpha^1\varepsilon^{-1/2}\{\log(C_1/\varepsilon)\}^{1/2}\,\mathrm d\varepsilon&\le\int_0^1\varepsilon^{-1/2}\{\log(C_1/\varepsilon)\}^{1/2}\,\mathrm d\varepsilon\\
&\le\sqrt{\log C_1}\int_0^1\varepsilon^{-1/2}\,\mathrm d\varepsilon+\int_0^1\varepsilon^{-1/2}\{\log(1/\varepsilon)\}^{1/2}\,\mathrm d\varepsilon\\
&=2\sqrt{\log C_1}+\int_0^\infty u^{1/2}e^{-u/2}\,\mathrm du=:C_3.
\end{align*}
Therefore, taking $\alpha=n^{-1/2}$ gives $\widehat{\mathfrak R}_n(\mathcal F_{\mathrm{UCB}})\le C_4n^{-1/2}$.

\noindent \textbf{Case 2: $m=2$.} In this case, we have that
\[
\int_\alpha^1\varepsilon^{-1}\{\log(C_1/\varepsilon)\}^{1/2}\,\mathrm d\varepsilon\le C_5\{\log(C_1/\alpha)\}^{3/2}.
\]
Taking $\alpha=n^{-1/2}$ gives $\widehat{\mathfrak R}_n(\mathcal F_{\mathrm{UCB}})\le C_6n^{-1/2}\{\log(en)\}^{3/2}$.

\noindent \textbf{Case 3: $m\ge3$.} In this case, we have that
\[
\int_\alpha^1\varepsilon^{-m/2}\{\log(C_1/\varepsilon)\}^{1/2}\,\mathrm d\varepsilon\le C_7\alpha^{1-m/2}\{\log(1/\alpha)\}^{1/2}.
\]
Taking $\alpha=n^{-1/m}\{\log(en)\}^{1/m}$ gives $\widehat{\mathfrak R}_n(\mathcal F_{\mathrm{UCB}})\le C_8\{\log(en)/n\}^{1/m}$.

Thus, for all $m\ge1$, we have that
\[
\widehat{\mathfrak R}_n(\mathcal F_{\mathrm{UCB}})\le C\begin{cases}n^{-1/2},&m=1,\\n^{-1/2}\log^{3/2}(en),&m=2,\\\{\log(en)/n\}^{1/m},&m\ge3.\end{cases}
\]
This completes the proof.
\end{proof}

\begin{lemma}[Control of $\Delta_n^{\mathrm{UCB}}$]
\label[lemma]{lem:ucb-control-delta}
Suppose Assumptions \ref{a_beta} and \ref{a_theta} hold and $d\ge2$. Let $\mathcal F_{\mathrm{UCB}}$ be the function class defined in \Cref{lem:ucb-bracketing}. Let $X_1,\ldots,X_n$ be i.i.d. from an arbitrary Borel probability measure $P_X$ supported on $[0,1]^d$. Define
\begin{equation}\label{eq_delta_ucb}
    \Delta_n^{\mathrm{UCB}}=\sup_{f\in\mathcal F_{\mathrm{UCB}}}\left|\frac1n\sum_{i=1}^nf(X_i)-\mathbb E\{f(X)\}\right|.
\end{equation}
Then there exists a constant $C<\infty$, depending only on $d,\nu,L,C_\beta$ and the constants in \Cref{lem:enlarged_theta}, and independent of $P_X$, such that, for every $n\ge2$ and every $\eta\in(0,1)$, with probability at least $1-\eta$,
\[
\Delta_n^{\mathrm{UCB}}\le C\rho_{n,d-1}(\eta),
\]
where $\rho_{n,m}$ is defined in \Cref{lem:control-delta-L-n}.
\end{lemma}

\begin{proof}
Write $m=d-1$. Since every $f\in\mathcal F_{\mathrm{UCB}}$
is a convex combination of the two linear reward functions
at each context, we have
\[
\|f\|_\infty\le\sqrt{d+1}(C_\beta+1)=:B.
\]

Note that the finite brackets of pointwise width $\varepsilon$ constructed
in \Cref{lem:ucb-bracketing} imply that
$\mathcal F_{\mathrm{UCB}}$ is separable in the sup-norm.
Consequently, the empirical and Rademacher suprema are measurable.
By the standard symmetrisation inequality,
\[
\mathbb E_X\Delta_n^{\mathrm{UCB}}
\le 2\mathbb E_X\widehat{\mathfrak R}_n(\mathcal F_{\mathrm{UCB}})
=2\mathfrak R_n(\mathcal F_{\mathrm{UCB}}).
\]

Replacing one observation changes $\Delta_n^{\mathrm{UCB}}$ by at most $2B/n$. The bounded-differences inequality therefore gives, with probability at least $1-\eta$,
\[
\Delta_n^{\mathrm{UCB}}\le\mathbb E_X\Delta_n^{\mathrm{UCB}}+B\sqrt{\frac{2\log(1/\eta)}{n}}.
\]
Combining this with symmetrisation gives the usual uniform convergence bound based on Rademacher complexity \citep[e.g.~Theorem~4.10 in][]{wainwright2019high}:
\[
\Delta_n^{\mathrm{UCB}}\le2\mathfrak R_n(\mathcal F_{\mathrm{UCB}})+B\sqrt{\frac{2\log(1/\eta)}{n}},
\]
where $\mathfrak R_n(\mathcal F_{\mathrm{UCB}})=\mathbb E_X\{\widehat{\mathfrak R}_n(\mathcal F_{\mathrm{UCB}})\}$. Note that the bound in \Cref{lem:rademacher-policy} holds for every deterministic sample, therefore by the tower property, it also holds after taking the expectation over $X_1,\ldots,X_n$. Combining the two bounds gives
\[
\Delta_n^{\mathrm{UCB}}\le C\begin{cases}n^{-1/2}+\sqrt{\log(1/\eta)/n},&m=1,\\n^{-1/2}\log^{3/2}(en)+\sqrt{\log(1/\eta)/n},&m=2,\\\{\log(en)/n\}^{1/m}+\sqrt{\log(1/\eta)/n},&m\ge3.\end{cases}
\]
Since $\log(en/\eta)\ge1$, each line is bounded by $C\rho_{n,m}(\eta)$ for $n\ge2$ and $\eta\in(0,1)$. This completes the proof.
\end{proof}

\subsubsection{Self-normalised confidence}
\begin{lemma}[Self-normalised concentration]\label[lemma]{l_s_normalise_lnu}
Suppose Assumptions \ref{a_noise} and \ref{a_beta} hold. For any $t\in[T]$ and $k\in\{1,2\}$, define
\begin{equation*}
S_{t,k}=\sum_{s=1}^{t-1}\widetilde X_s\xi_{s,k}\indc\{A_s=k\},\qquad V_{t,k}=\lambda I_{d+1}+\sum_{s=1}^{t-1}\widetilde X_s\widetilde X_s^\top\indc\{A_s=k\}.
\end{equation*}
Then, it holds that
\begin{equation*}
\mathbb P\Big[\|S_{t,k}\|_{V_{t,k}^{-1}}^2\le2\sigma^2\log\Big\{\frac{2\det(V_{t,k})^{1/2}\det(\lambda I_{d+1})^{-1/2}}{\eta}\Big\},\ \text{for all }t\in[T],\ k\in\{1,2\}\Big]\ge1-\eta.
\end{equation*}
\end{lemma}

\begin{proof}
The proof is constructed based on the proof of Theorem 1 in \citet{abbasi2011improved}. For any $t\in[T]$, denote the filtration $\mathcal F_t=\sigma(\{X_s,A_s,Y_s\}_{s\in[t]})$ and $\mathcal G_t=\sigma(\mathcal F_{t-1},X_t,A_t)$.

\noindent \textbf{Construction of Supermartingale.} For any $v\in\mathbb R^{d+1}$ and $t\in[T]\backslash\{1\}$, define
\begin{equation*}
M_{t,k}(v)=\exp\left\{\sigma^{-1}\langle v,S_{t,k}\rangle-\frac12 v^\top\left(\sum_{s=1}^{t-1}\indc\{A_s=k\}\widetilde X_s\widetilde X_s^\top\right)v\right\},
\end{equation*}
and set $M_{1,k}(v)=1$. Since
\begin{equation*}
S_{t+1,k}= \sum_{s =1}^{t}\widetilde X_s \xi_{s,k} \indc\{A_s =k\} =S_{t,k}+\widetilde X_t\xi_{t,k}\indc\{A_t=k\},
\end{equation*}
and 
    \[
    \sum_{s=1}^{t} \indc\{A_s =k\}\widetilde X_s\widetilde X_s^\top = \sum_{s=1}^{t-1} \indc\{A_s =k\}\widetilde X_s\widetilde X_s^\top + \indc\{A_t =k\}\widetilde X_t\widetilde X_t^\top,
    \]
we have
\begin{align*}
        M_{t+1,k}(v)
        =\;& \exp \Big\{\frac{1}{\sigma}\langle v, S_{t+1,k} \rangle -\frac{1}{2}v^\top \Big(\sum_{s=1}^{t} \indc\{A_s =k\}\widetilde X_s\widetilde X_s^\top\Big)v \Big\}\\
        =\;& M_{t,k}(v) \exp\Big\{\sigma^{-1}\big\langle v, \widetilde X_{t} \xi_{t,k} \indc\{A_t =k\}\big\rangle - \frac{\indc\{A_t =k\}}{2}v^\top\widetilde X_t\widetilde X_t^\top v\Big\}.
    \end{align*}
Since for $t\in[T]$, both $M_{t,k}(v)$ and $\widetilde X_t$ are $\mathcal G_t$-measurable, the sub-Gaussian condition gives
\begin{align*}
\mathbb E[M_{t+1,k}(v)|\mathcal G_t]&\le M_{t,k}(v)\exp\left\{-\frac{\indc\{A_t=k\}}{2}v^\top\widetilde X_t\widetilde X_t^\top v\right\}\exp\left\{\frac{\indc\{A_t=k\}}{2}v^\top\widetilde X_t\widetilde X_t^\top v\right\}\\
&=M_{t,k}(v).
\end{align*}
Hence $\{M_{t,k}\}_{t\in[T]}$ is a sequence of non-negative supermartingales.

\noindent \textbf{Method of mixtures.} Let $V\sim\mathcal N(0,\lambda^{-1}I_{d+1})$ be a Gaussian random variable which is independent of all other random variables and define the mixture process
\begin{equation*}
M^V_{t,k}=\int_{\mathbb R^{d+1}}M_{t,k}(v)\,\mathrm d\mathcal N(0,\lambda^{-1}I_{d+1})(v).
\end{equation*}
By Tonelli's theorem, we have that for any $t \in [T]$
    \begin{align*}
        \mathbb{E}[M^V_{t+1,k}|\mathcal{G}_t]
        =\;& \int_{\mathbb{R}^{d+1}} \mathbb{E}\{M_{t+1,k}(v)|\mathcal{G}_t\} \;\mathrm{d}\mathcal{N}(0,\lambda^{-1} I_{d+1})(v) \\
        \leq\;& \int_{\mathbb{R}^{d+1}} M_{t,k}(v)\;\mathrm{d}\mathcal{N}(0,\lambda^{-1} I_{d+1})(v) \\
        =\;& M^V_{t,k}.
    \end{align*}
 Hence, we can conclude that $\{M^V_{t,k}\}_{t\in[T]}$ is still a sequence of non-negative supermartingales with $\mathbb{E}[M^V_{1,k}] =1$. Next, we compute the integral explicitly. For notational simplicity, we let
\begin{equation*}
Q_{t,k}=\sum_{s=1}^{t-1}\indc\{A_s=k\}\widetilde X_s\widetilde X_s^\top.
\end{equation*}
 Then
    \begin{align*}
        M^{V}_{t,k}
        =\;& \frac{1}{(2\pi\lambda^{-1})^{(d+1)/2}}\int_{\mathbb{R}^{d+1}} \exp \Big\{\sigma^{-1}\langle v, S_{t,k} \rangle -\frac{v^\top Q_{t,k}v}{2} \Big\}\exp\Big(-\frac{\lambda\|v\|_2^2}{2}\Big)\;\mathrm{d}v\\
        =\;& \frac{1}{(2\pi\lambda^{-1})^{(d+1)/2}}\int_{\mathbb{R}^{d+1}} \exp \Big\{\sigma^{-1}\langle v, S_{t,k} \rangle -\frac{1}{2}v^\top (Q_{t,k} + \lambda I_{d+1})v \Big\}\;\mathrm{d}v.
    \end{align*}
    Note that 
    \begin{align*}
        &\sigma^{-1}v^\top S_{t,k} -\frac{v^\top (Q_{t,k} + \lambda I_{d+1})v}{2} \\
        =\;& -\frac{1}{2}\Big\{v- \sigma^{-1}(Q_{t,k}+\lambda I_{d+1})^{-1}S_{t,k}\Big\}^\top (Q_{t,k} + \lambda I_{d+1})\Big\{v- \sigma^{-1}(Q_{t,k}+\lambda I_{d+1})^{-1}S_{t,k}\Big\}\\
        & + \frac{1}{2\sigma^{2}} S_{t,k}^\top(Q_{t,k}+\lambda I_{d+1})^{-1}S_{t,k}.
    \end{align*}
    Then, the integral becomes
    \begin{align*}
        M^{V}_{t,k}
        =\;&  \frac{1}{(2\pi\lambda^{-1})^{(d+1)/2}} \exp\Big\{ \frac{1}{2\sigma^{2}} S_{t,k}^\top(Q_{t,k}+\lambda I_{d+1})^{-1}S_{t,k}\Big\} \\
        &\cdot\int_{\mathbb{R}^{d+1}}\exp\Big[-\frac{1}{2}\Big\{v- \sigma^{-1}(Q_{t,k}+\lambda I_{d+1})^{-1}S_{t,k}\Big\}^\top\\
        &\qquad\qquad (Q_{t,k} + \lambda I_{d+1})\Big\{v- \sigma^{-1}(Q_{t,k}+\lambda I_{d+1})^{-1}S_{t,k}\Big\}\Big]\;\mathrm d v\\
        =\;&   \frac{1}{(2\pi\lambda^{-1})^{(d+1)/2}} \exp\Big\{\frac{1}{2\sigma^{2}} S_{t,k}^\top(Q_{t,k}+\lambda I_{d+1})^{-1}S_{t,k}\Big\} (2\pi )^{(d+1)/2}\det(Q_{t,k}+\lambda I_{d+1})^{-1/2}\\
        =\;& \sqrt{\frac{\det(\lambda I_{d+1})}{\det(V_{t,k})}}\exp\Big\{\frac{1}{2\sigma^2} \|S_{t,k}\|^2_{V_{t,k}^{-1}}\Big\}.
    \end{align*}

\noindent \textbf{Tail bound.} By Ville's inequality, for any $\eta/2 >0$, we have with probability at least $1-\eta/2$ that
    \begin{align*}
       \sqrt{\frac{\det(\lambda I_{d+1})}{\det(V_{t,k})}}\exp\Big\{\frac{1}{2\sigma^2} \|S_{t,k}\|^2_{V_{t,k}^{-1}}\Big\} \leq \frac{2}{\eta}, \text{ for all } t \in [T].
    \end{align*}
    Thus, by another union bound over $k \in \{1,2\}$, we have that 
    \begin{align*}
        \mathbb{P}\Big[\|S_{t,k}\|^2_{V_{t,k}^{-1}} \leq 2\sigma^2 \log\Big\{\frac{2\det(V_{t,k})^{1/2}\det(\lambda I_{d+1})^{-1/2}}{\eta}\Big\}, \text{ for all } t \in [T], k \in \{1,2\}\Big] \geq 1- \eta.
    \end{align*}
\end{proof}

\begin{corollary}\label[corollary]{cor:self-normalised-deterministic-lnu}
Under the conditions of \Cref{l_s_normalise_lnu}, with probability at least $1-\eta$, it holds for all $t\in[T]$ and $k\in\{1,2\}$ that
\begin{equation*}
\|S_{t,k}\|_{V_{t,k}^{-1}}^2\le2\sigma^2\left\{\log\left(\frac{2}{\eta}\right)+\frac{d+1}{2}\log\left(1+\frac{T}{\lambda}\right)\right\}.
\end{equation*}
\end{corollary}

\begin{proof}
By \Cref{l_s_normalise_lnu}, with probability at least $1-\eta$, it holds for all $t\in[T]$ and $k\in\{1,2\}$ that
\begin{align*}
\|S_{t,k}\|_{V_{t,k}^{-1}}^2&\le2\sigma^2\log\left\{\frac{2\det(V_{t,k})^{1/2}\det(\lambda I_{d+1})^{-1/2}}{\eta}\right\}\\
&=2\sigma^2\left\{\log\left(\frac{2}{\eta}\right)+\frac12\log\left(\frac{\det(V_{t,k})}{\det(\lambda I_{d+1})}\right)\right\}.
\end{align*}
Since $\|\widetilde X_t\|_2\le\sqrt{d+1}$, Lemma 10 in \citet{abbasi2011improved} gives
\begin{equation*}
    \log\left(\frac{\det(V_{t,k})}{\det(\lambda I_{d+1})}\right)\le(d+1)\log\left(1+\frac{T}{\lambda}\right).
\end{equation*}
The result follows.
\end{proof}

\begin{lemma}[UCB confidence bound]\label[lemma]{l_ucb_lnu}
Suppose Assumptions \ref{a_noise} and \ref{a_beta} hold. For each $t\in[T]$ and $k\in\{1,2\}$, define
\begin{equation*}
V_{t,k}=\lambda I_{d+1}+\sum_{s=1}^{t-1}\widetilde X_s\widetilde X_s^\top\indc\{A_s=k\},\qquad \widehat\beta_{t,k}=V_{t,k}^{-1}\sum_{s=1}^{t-1}\widetilde X_sY_s\indc\{A_s=k\}.
\end{equation*}
Then, with probability at least $1-\eta$ that, for all $x\in[0,1]^d$, $k\in\{1,2\}$ and $t\in[T]$,
\begin{equation*}
    |\widetilde x^\top(\widehat\beta_{t,k}-\beta_k^*)|\le\left[\sqrt{2\sigma^2\left\{\log\left(\frac{2}{\eta}\right)+\frac{d+1}{2}\log\left(1+\frac{T}{\lambda}\right)\right\}}+C_\beta\sqrt\lambda\right]\sqrt{\widetilde x^\top V_{t,k}^{-1}\widetilde x}.
\end{equation*}
Consequently, taking
\begin{equation*}
\alpha_T=C\left[\sqrt{\sigma^2\left\{\log\left(\frac{1}{\eta}\right)+(d+1)\log\left(1+\frac{T}{\lambda}\right)\right\}}+C_\beta\sqrt\lambda\right],
\end{equation*}
it holds with probability at least $1-\eta$ that, for all $x\in[0,1]^d$, $k\in\{1,2\}$ and $t\in[T]$,
\begin{equation*}
    \widetilde x^\top\beta_k^*\le \widetilde x^\top\widehat\beta_{t,k}+\alpha_T\sqrt{\widetilde x^\top V_{t,k}^{-1}\widetilde x}.
\end{equation*}
\end{lemma}

\begin{proof}
Define the event for self-normalised norms
\begin{equation*}
\mathcal E_{\rm sn}=\left\{\|S_{t,k}\|_{V_{t,k}^{-1}}^2\le2\sigma^2\left\{\log\left(\frac{2}{\eta}\right)+\frac{d+1}{2}\log\left(1+\frac{T}{\lambda}\right)\right\},\ \text{for all }t\in[T]\text{ and }k\in\{1,2\}\right\}.
\end{equation*}
By \Cref{cor:self-normalised-deterministic-lnu}, it holds that $\mathbb P(\mathcal E_{\rm sn})\ge1-\eta$. In the rest of the proof, we condition on $\mathcal E_{\rm sn}$ happening.

For any $x\in[0,1]^d$, by the Cauchy--Schwarz inequality in the $V_{t,k}$-norm, we have that
\begin{align} \notag
|\widetilde x^\top(\widehat\beta_{t,k}-\beta_k^*)|&\le\|\widehat\beta_{t,k}-\beta_k^*\|_{V_{t,k}}\|\widetilde x\|_{V_{t,k}^{-1}}\\ \label{l_ucb_lnu_eq2}
&=\|\widehat\beta_{t,k}-\beta_k^*\|_{V_{t,k}}\sqrt{\widetilde x^\top V_{t,k}^{-1}\widetilde x}.
\end{align}
When $A_s=k$, we have $Y_s=\widetilde X_s^\top\beta_k^*+\xi_{s,k}$. Therefore,
\begin{align*}
\widehat\beta_{t,k}-\beta_k^*&=V_{t,k}^{-1}\sum_{s=1}^{t-1}\widetilde X_sY_s\indc\{A_s=k\}-\beta_k^*\\
&=V_{t,k}^{-1}\sum_{s=1}^{t-1}\widetilde X_s\widetilde X_s^\top\beta_k^*\indc\{A_s=k\}+V_{t,k}^{-1}S_{t,k}-\beta_k^*\\
&=V_{t,k}^{-1}(V_{t,k}-\lambda I_{d+1})\beta_k^*+V_{t,k}^{-1}S_{t,k}-\beta_k^*\\
&=V_{t,k}^{-1}S_{t,k}-\lambda V_{t,k}^{-1}\beta_k^*.
\end{align*}
Therefore,
\begin{equation} \label{l_ucb_lnu_eq1}
\|\widehat\beta_{t,k}-\beta_k^*\|_{V_{t,k}}\le\|S_{t,k}\|_{V_{t,k}^{-1}}+\sqrt\lambda\|\beta_k^*\|_2\le\|S_{t,k}\|_{V_{t,k}^{-1}}+C_\beta\sqrt\lambda.
\end{equation}
Combining \eqref{l_ucb_lnu_eq1} with \eqref{l_ucb_lnu_eq2} and by construction of $\mathcal E_{\rm sn}$, the first claim follows. The second claim follows directly by the choice of $\alpha_T$.
\end{proof}

\subsubsection{Lazy width summation and regret decomposition}
\begin{lemma}\label[lemma]{lem:lazy-determinant-comparison-lnu}
Let $t\ge s_0$ and $s=\tau(t)$ be the lazy update time used at round $t$. Then, for $k\in\{1,2\}$,
\begin{equation*}
    V_{s,k}^{-1}\preceq 2V_{t,k}^{-1}.
\end{equation*}
Consequently, for any $x\in[0,1]^d$,
\begin{equation*}
    \|\widetilde x\|_{V_{s,k}^{-1}}\le\sqrt2\|\widetilde x\|_{V_{t,k}^{-1}}.
\end{equation*}
\end{lemma}

\begin{proof}
At each lazy update time $s$, the algorithm sets $\ell_{s,k}:=\log\det(V_{s,k})$. Fix $t\ge s_0$ and let $s=\tau(t)$. If $t=s$, then $V_{t,k}=V_{s,k}$. If $t>s$, no update occurs at round $t$, so the determinant trigger, evaluated with $\ell_k=\log\det(V_{t,k})$ before observing $X_t$, has not fired. Thus, in both cases,
\begin{equation*}
    \log\det(V_{t,k})-\log\det(V_{s,k})<\log2.
\end{equation*}
Moreover, $V_{t,k}\succeq V_{s,k}$. Let the eigenvalues of $V_{s,k}^{-1/2}V_{t,k}V_{s,k}^{-1/2}$ be $\mu_1,\ldots,\mu_{d+1}$. Then $\mu_j\ge1$ for all $j$ and $\prod_{j=1}^{d+1}\mu_j=\det(V_{t,k})/\det(V_{s,k})<2$. Hence $\mu_j<2$ for all $j$, which gives $V_{t,k}\preceq2V_{s,k}$. Taking inverses yields $V_{s,k}^{-1}\preceq2V_{t,k}^{-1}$.
\end{proof}

\begin{lemma}[Realised lazy width sum]\label[lemma]{lem:lazy-width-sum-lnu}
For the sequence generated by \Cref{alg:lazy-global-model-ucb}, we have that
\begin{equation*}
    \sum_{t=s_0}^T\min\{1,\|\widetilde X_t\|_{V_{\tau(t),A_t}^{-1}}\}\le C\sqrt{T(d+1)\log\left(1+\frac{T}{\lambda}\right)},
\end{equation*}
where $C>0$ is an absolute constant.
\end{lemma}

\begin{proof}
By \Cref{lem:lazy-determinant-comparison-lnu}, we have that
\begin{equation*}
    \min\{1,\|\widetilde X_t\|_{V_{\tau(t),A_t}^{-1}}\}\le\sqrt2\min\{1,\|\widetilde X_t\|_{V_{t,A_t}^{-1}}\}.
\end{equation*}
Therefore, by the Cauchy--Schwarz inequality,
\begin{align*}
\sum_{t=s_0}^T\min\{1,\|\widetilde X_t\|_{V_{\tau(t),A_t}^{-1}}\}&\le\sqrt{2T\sum_{t=s_0}^T\min\{1,\|\widetilde X_t\|_{V_{t,A_t}^{-1}}^2\}}.
\end{align*}
For each arm $k$, Lemma 10 in \citet{abbasi2011improved} gives
\begin{equation*}
    \sum_{t=s_0}^T\min\{1,\|\widetilde X_t\|_{V_{t,k}^{-1}}^2\}\indc\{A_t=k\}\le 2\log\left(\frac{\det(V_{T+1,k})}{\det(\lambda I_{d+1})}\right).
\end{equation*}
Note that, since $V_{1,k}=\lambda I_{d+1}$ and $V_{t+1,k}=V_{t,k}+\indc\{A_t=k\}\widetilde X_t\widetilde X_t^\top$, the matrix determinant lemma \citep[e.g.~Theorem 18.1.1 in][]{harville1998matrix} gives
\begin{align*}
    \det(V_{t+1,k})
    &=\det(V_{t,k})\left\{1+\indc\{A_t=k\}\widetilde X_t^\top V_{t,k}^{-1}\widetilde X_t\right\}.
\end{align*}
Using the fact that $\min\{1,a\}\le2\log(1+a)$ for all $a\ge0$, we have
\begin{align*}
&\indc\{A_t=k\}\min\{1,\widetilde X_t^\top V_{t,k}^{-1}\widetilde X_t\} \le2\log\left\{1+\indc\{A_t=k\}\widetilde X_t^\top V_{t,k}^{-1}\widetilde X_t\right\} =2\log\left(\frac{\det(V_{t+1,k})}{\det(V_{t,k})}\right).
\end{align*}
Summing over $t=1,\ldots,T$ gives
\begin{align*}
    \sum_{t=1}^T\indc\{A_t=k\}\min\{1,\widetilde X_t^\top V_{t,k}^{-1}\widetilde X_t\}
    &\le2\sum_{t=1}^T\log\left(\frac{\det(V_{t+1,k})}{\det(V_{t,k})}\right) =2\log\left(\frac{\det(V_{T+1,k})}{\det(\lambda I_{d+1})}\right).
\end{align*}
Now let $N_k(T)=\sum_{t=1}^T\indc\{A_t=k\}$. Since $\|\widetilde X_t\|_2^2\le d+1$ and $N_k(T)\le T$, we have that
\begin{align*}
    \operatorname{tr}(V_{T+1,k})
    &=\lambda(d+1)+\sum_{t=1}^T\indc\{A_t=k\}\|\widetilde X_t\|_2^2 \le(d+1)(\lambda+T).
\end{align*}
Therefore, by Lemma 10 in \citet{abbasi2011improved},
\begin{equation*}
    \det(V_{T+1,k})\le\left\{\frac{\operatorname{tr}(V_{T+1,k})}{d+1}\right\}^{d+1}\le(\lambda+T)^{d+1}.
\end{equation*}
Since $\det(\lambda I_{d+1})=\lambda^{d+1}$, this implies
\begin{equation*}
    \log\left(\frac{\det(V_{T+1,k})}{\det(\lambda I_{d+1})}\right)\le(d+1)\log\left(1+\frac{T}{\lambda}\right).
\end{equation*}
Consequently,
\begin{equation*}
    \sum_{t=1}^T\indc\{A_t=k\}\min\{1,\widetilde X_t^\top V_{t,k}^{-1}\widetilde X_t\}
    \le2(d+1)\log\left(1+\frac{T}{\lambda}\right).
\end{equation*}
Summing the last display over $k\in\{1,2\}$, we obtain
\begin{align*}
\sum_{t=1}^T\min\{1,\|\widetilde X_t\|_{V_{t,A_t}^{-1}}^2\}
&=\sum_{k=1}^2\sum_{t=1}^T\indc\{A_t=k\}\min\{1,\widetilde X_t^\top V_{t,k}^{-1}\widetilde X_t\} \le4(d+1)\log\left(1+\frac{T}{\lambda}\right).
\end{align*}
Since the left-hand side is non-negative, this also implies
\begin{equation*}
    \sum_{t=s_0}^T\min\{1,\|\widetilde X_t\|_{V_{t,A_t}^{-1}}^2\}
    \le4(d+1)\log\left(1+\frac{T}{\lambda}\right).
\end{equation*}
Combining with the Cauchy--Schwarz inequality gives
\begin{align*}
\sum_{t=s_0}^T\min\{1,\|\widetilde X_t\|_{V_{\tau(t),A_t}^{-1}}\}
&\le\sqrt{2T\sum_{t=s_0}^T\min\{1,\|\widetilde X_t\|_{V_{t,A_t}^{-1}}^2\}} \le\sqrt{8T(d+1)\log\left(1+\frac{T}{\lambda}\right)}.
\end{align*}
The lemma thus follows.

\end{proof}


\begin{lemma}[Empirical optimism]\label[lemma]{lem:empirical-optimism}
On $\mathcal E_{\rm conf}$ in \eqref{eq:main-regret-confidence-event} and the event in \Cref{lem:initial-localisation} on which every feasible difference $b_1-b_2$, with $b_k\in\mathcal C_{s,k}$, belongs to $\widetilde\Theta_{j^*}$ for every update time $s$, the model pair $\bar\beta_s=(\bar\beta_{s,1},\bar\beta_{s,2})$ and policy $p_s$ returned by \Cref{alg:vertex-oracle} satisfy
\begin{equation}\label{eq:empirical-optimism}
\frac1{s-1}\sum_{i=1}^{s-1}r_{\bar\beta_s,p_s}(X_i)\ge\frac1{s-1}\sum_{i=1}^{s-1}r_{\beta^*,p^*}(X_i),
\end{equation}
where, for $b=(b_1,b_2)$ and a policy $p$, $r_{b,p}(x)=\widetilde x^\top b_1p(x)+\widetilde x^\top b_2\{1-p(x)\}$, $\beta^*=(\beta_1^*,\beta_2^*)$, and $p^*=p_{\theta^*,q^*}$ is chosen to be a canonical population-optimal policy as guaranteed by \Cref{thm:canonical_extension}.
\end{lemma}

\begin{proof}
On $\mathcal E_{\rm conf}$, we have $\beta_k^*\in\mathcal C_{s,k}$ for $k\in\{1,2\}$. Moreover, by \Cref{lem:initial-localisation}, every feasible difference $b_1-b_2$ belongs to $\widetilde\Theta_{j^*}\subseteq\widetilde\Theta$. Hence, \Cref{lem:vertex-optimisation} applies and gives
\begin{equation*}
\frac1{s-1}\sum_{i=1}^{s-1}r_{\bar\beta_s,p_s}(X_i)=\sup_{b_k\in\mathcal C_{s,k},\ q\in\mathcal Q_{b_1-b_2,L}}\frac1{s-1}\sum_{i=1}^{s-1}r_{b,p_{b_1-b_2,q}}(X_i).
\end{equation*}
Since $\theta^*=\beta_1^*-\beta_2^*$, $\beta_k^*\in\mathcal C_{s,k}$, and $q^*\in\mathcal Q_{\theta^*,L}$ with $p^*=p_{\theta^*,q^*}$, the true model-policy pair $(\beta^*,q^*)$ is feasible. Therefore,
\begin{equation*}
\sup_{b_k\in\mathcal C_{s,k},\ q\in\mathcal Q_{b_1-b_2,L}}\frac1{s-1}\sum_{i=1}^{s-1}r_{b,p_{b_1-b_2,q}}(X_i)\ge\frac1{s-1}\sum_{i=1}^{s-1}r_{\beta^*,p^*}(X_i).
\end{equation*}
Combining the preceding equality and inequality proves \eqref{eq:empirical-optimism}.
\end{proof}

\begin{lemma}[One-step regret]\label[lemma]{lem:one-step-regret}
On $\mathcal E_{\rm conf}\cap\mathcal E_{\rm emp}$, defined in \eqref{eq:main-regret-confidence-event} and \eqref{eq:main-regret-empirical-event}, and the event in \Cref{lem:initial-localisation} on which every feasible difference $b_1-b_2$, with $b_k\in\mathcal C_{s,k}$, belongs to $\widetilde\Theta_{j^*}$ and every candidate satisfies $\|b_k\|_2\le C_\beta+1$ for every update time $s$, for every $t\ge s_0$ and $s=\tau(t)$,
\begin{align}\label{eq:one-step-regret}
R(p^*)-R(p_s)\le&\ 2\Delta_{s-1}^{\mathrm{UCB}} +\{(1+\sqrt{d+1})\alpha_T+2B\}\mathbb E\left[\min\{1,\|\widetilde X_t\|_{V_{s,A_t}^{-1}}\}\mid\mathcal F_{t-1}\right].
\end{align}
where $B=\sqrt{d+1}(C_\beta+1)$.
\end{lemma}

\begin{proof}
Since $R(p)=\mathbb E\{r_{\beta^*,p}(X)\}$, where $X\sim P_X$, we decompose
\begin{align*}
R(p^*)-R(p_s)
={}&\left[\mathbb E\{r_{\beta^*,p^*}(X)\}-\frac1{s-1}\sum_{i=1}^{s-1}r_{\beta^*,p^*}(X_i)\right]\\
&+\left[\frac1{s-1}\sum_{i=1}^{s-1}r_{\beta^*,p^*}(X_i)-\frac1{s-1}\sum_{i=1}^{s-1}r_{\bar\beta_s,p_s}(X_i)\right]\\
&+\left[\frac1{s-1}\sum_{i=1}^{s-1}r_{\bar\beta_s,p_s}(X_i)-\mathbb E_X\{r_{\bar\beta_s,p_s}(X)\}\right]\\
&+\mathbb E_X\{r_{\bar\beta_s,p_s}(X)-r_{\beta^*,p_s}(X)\}.
\end{align*}
Since $p^*=p_{\theta^*,q^*}$, $\theta^*\in\widetilde\Theta$, and $\|\beta_k^*\|_2\le C_\beta$, we have $r_{\beta^*,p^*}\in\mathcal F_{\mathrm{UCB}}$. Moreover, on the event of \Cref{lem:initial-localisation}, $\theta_s=\bar\beta_{s,1}-\bar\beta_{s,2}\in\widetilde\Theta_{j^*}\subseteq\widetilde\Theta$, and the same lemma gives $\|\bar\beta_{s,k}\|_2\le C_\beta+1$, so it holds that $r_{\bar\beta_s,p_s}\in\mathcal F_{\mathrm{UCB}}$. Therefore, we can control the first and third terms of the order at most $\Delta_{s-1}^{\mathrm{UCB}}$ each on $\mathcal E_{\rm emp}$. Finally, \Cref{lem:empirical-optimism} shows that the second term is nonpositive.

It remains to control the final term. For $k\in\{1,2\}$, \Cref{lem:polyhedral-confidence} gives $\|\bar\beta_{s,k}-\beta_k^*\|_{V_{s,k}}\le(1+\sqrt{d+1})\alpha_T$ on $\mathcal E_{\rm conf}$. By \Cref{lem:initial-localisation} and \Cref{a_beta}, both $\bar\beta_{s,k}$ and $\beta_k^*$ have Euclidean norm at most $C_\beta+1$. Write $B=\sqrt{d+1}(C_\beta+1)$. Hence,
\begin{equation*}
|\widetilde x^\top(\bar\beta_{s,k}-\beta_k^*)|\le\min\{(1+\sqrt{d+1})\alpha_T\|\widetilde x\|_{V_{s,k}^{-1}},2B\}\le\{(1+\sqrt{d+1})\alpha_T+2B\}\min\{1,\|\widetilde x\|_{V_{s,k}^{-1}}\}.
\end{equation*}
Note that the update time $s=\tau(t)$ and the quantities $\bar\beta_s$, $p_s$, and $V_{s,k}$ used at round $t$ are $\mathcal F_{t-1}$-measurable. Since $X_t$ is independent of $\mathcal F_{t-1}$ and has law $P_X$, we have that
\begin{align}
&\mathbb E\{r_{\bar\beta_s,p_s}(X)-r_{\beta^*,p_s}(X)\}\notag\\
&=\mathbb E\left[\widetilde X_t^\top(\bar\beta_{s,1}-\beta_1^*)p_s(X_t)+\widetilde X_t^\top(\bar\beta_{s,2}-\beta_2^*)\{1-p_s(X_t)\}\mid\mathcal F_{t-1}\right]\notag\\
&\le\{(1+\sqrt{d+1})\alpha_T+2B\}\mathbb E\left[p_s(X_t)\min\{1,\|\widetilde X_t\|_{V_{s,1}^{-1}}\}+\{1-p_s(X_t)\}\min\{1,\|\widetilde X_t\|_{V_{s,2}^{-1}}\}\mid\mathcal F_{t-1}\right]\notag\\
&=\{(1+\sqrt{d+1})\alpha_T+2B\}\mathbb E\left[\sum_{k=1}^2\min\{1,\|\widetilde X_t\|_{V_{s,k}^{-1}}\}\mathbb P(A_t=k\mid X_t,\mathcal F_{t-1})\mid\mathcal F_{t-1}\right]\notag\\
&=\{(1+\sqrt{d+1})\alpha_T+2B\}\mathbb E\left[\mathbb E\left[\min\{1,\|\widetilde X_t\|_{V_{s,A_t}^{-1}}\}\mid X_t,\mathcal F_{t-1}\right]\mid\mathcal F_{t-1}\right]\notag\\
&=\{(1+\sqrt{d+1})\alpha_T+2B\}\mathbb E\left[\min\{1,\|\widetilde X_t\|_{V_{s,A_t}^{-1}}\}\mid\mathcal F_{t-1}\right]\notag\\
&=\{(1+\sqrt{d+1})\alpha_T+2B\}\mathbb E\left[\min\{1,\|\widetilde X_t\|_{V_{\tau(t),A_t}^{-1}}\}\mid\mathcal F_{t-1}\right], \notag
\end{align}
where the fourth equality follows from the tower property, and the last equality uses $s=\tau(t)$. Substituting this bound into the preceding decomposition proves \eqref{eq:one-step-regret}.

\end{proof}

\begin{lemma}[Summation of empirical errors]\label[lemma]{lem:error-sums}
For the lazy schedule in \Cref{alg:lazy-global-model-ucb}, $\tau(t)>t/2$ for every $t\ge s_0$. Consequently, with $u=\eta/(8T)$,
\begin{equation*}
\sum_{t=s_0}^T\rho_{\tau(t)-1,d-1}(u)\le C\begin{cases}\sqrt T\log(eT/u),&d=2,\\\sqrt T\log^{3/2}(eT/u),&d=3,\\T^{(d-2)/(d-1)}\log(eT/u),&d\ge4.\end{cases}
\end{equation*}
\end{lemma}

\begin{proof}
If an update occurs at $t$, then $\tau(t)=t$. Otherwise, the time trigger has not fired, so $t<2\tau(t)$. Thus $\tau(t)>t/2$. Monotonicity gives $\rho_{\tau(t)-1,m}(u)\le C\rho_{\lfloor t/2\rfloor,m}(u)$ after changing a universal constant. Summing the three displayed forms of $\rho$ and using $\sum_{t=1}^Tt^{-a}\le C_aT^{1-a}$ for $0<a<1$ proves the result.
\end{proof}

\section{Proofs of results in Section \ref{sec_lower}}\label[appendix]{sec_app_lower_results}

\subsection{Proof of Theorem \ref{thm:lower_bound}} \label[appendix]{sec_app_lower_bound}

Assume that $X_t\in\mathcal X=[0,1]^d$ is drawn independently from a density $g$ on $[0,1]^d$.
The learner observes $X_t$, chooses arm $A_t\in\{1,2\}$, and then observes
\[
    Y_t=r_{A_t}(X_t)+\xi_t,
    \qquad
    \xi_t\sim_{i.i.d.} N(0,1).
\]
The noises are independent of the contexts and the learner's internal randomization. We write
\[
    \Delta(x)=r_1(x)-r_2(x).
\]
The policy used at time $t$ is a measurable function
$p_t:[0,1]^d\to[0,1]$, where $p_t(x)$ is the conditional probability of pulling
arm $1$ at context $x$.  The policy may be random and adapted to the previous
history, but it is selected before seeing the current context and must satisfy the fairness constraint.

Define 
\[
    \mathcal{P}_{L,1} =  \left\{ p:[0,1]^d\to[0,1]: |p(x)-p(y)|\le L\|x-y\|_1 \text{ for all }x,y\in[0,1]^d \right\}.
\]
The population value of a policy $p$ under density $g$ is
\[
    J_g(p)  =  \int_{[0,1]^d} g(x)\Delta(x)p(x)\,\mathrm{d}x.
\]
We consider arbitrary non-anticipating learning algorithms.  Let $\mathcal F_{t-1}$ contain the observed history
$\sigma\{X_s,A_s,Y_s:1\le s\le t-1\}$ together with the learner's internal
randomization available when selecting $p_t$. All randomization kernels of
the learner are fixed across environments and independent of future contexts
and reward noise.  A learning algorithm $\pi$ is called non-anticipating if, at each round $t$, it selects an
$\mathcal F_{t-1}$-measurable random function
\[
p_t: [0,1]^d \to[0,1],\quad \text{such that }\quad 
p_t\in\mathcal{P}_{L,1} \quad \text{a.s.}.
\]
After observing the current covariate \(X_t\), the algorithm pulls arm
\(1\) with probability \(p_t(X_t)\) and arm \(2\) with probability
\(1-p_t(X_t)\). Thus
\(\pi\) induces a sequence of fair policies
\(\{p_t^{\pi}\}_{t=1}^T\subset\mathcal{P}_{L,1}\) almost surely.

Let $    p_g^\star\in\argmax_{q\in\mathcal{P}_{L,1}}J_g(q)
$ be the optimal policy. Then the regret of an algorithm $\{p_t^{\pi}\}_{t=1}^T$ is
\[
    \Reg_g  = \mathbb{E}_g\left[\sum_{t=1}^T  \{J_g(p_g^\star)-J_g(p_t^{\pi})\} \right].
\]
Since $J_g$ is linear in the policy, hence for $
    \bar p^{\pi}=\frac1T\sum_{t=1}^T p_t^{\pi}$, we can write
\[
    \Reg_g(\bar{p}^{\pi},T) =   T\,\mathbb{E}_g\{J_g(p_g^\star)-J_g(\bar p^{\pi})\}.
\]

\begin{theorem}[Lower bound]\label[theorem]{thm:LB}
For the reward functions $r_1(x)=x_1-1/2$ and $r_2(x)=0$, fix $L\in(0,1)$ and let $\mathcal{P}_{L,1}$ be the  collection of policies satisfying the $(L,1)$-Lipschitz constraint. Let $\mathcal{G}$ be the collection of bounded densities on $\mathcal{X}$. Then,  there is a constant
$c_{d,L}>0$ depending only on $d$ and $L$, such that for every integer $T\ge1$,
\[
   \inf_{\pi\in\Pi} \sup_{g\in\mathcal{G}}\Reg_{g}(\bar{p}^{\pi},T) \ge c_{d,L}
\begin{cases}
\sqrt T, & d=1,2,\\
T^{1-\frac{1}{d-1}}, & d\ge 3,
\end{cases}
\]
where $\Pi$ is the class of all the non-anticipating fair policies such that \[
\Pi = \Bigl\{\{p_t\}_{t=1}^\infty: p_t \text{ is } \mathcal F_{t-1}\text{-measurable and } p_t\in \mathcal{P}_{L,1} \text{ a.s. for every } t\Bigr\}.
\]
\end{theorem}

\begin{proof}
Fix an arbitrary learning algorithm $\pi\in\Pi$, common to all densities.
Write $\bar p=T^{-1}\sum_{t=1}^T p_t^\pi$. By convexity,
$\bar p\in\mathcal P_{L,1}$ almost surely.
All interaction laws below include the learner's internal randomization,
so $\bar p$ is the same measurable function of the full interaction under
every environment. Throughout the proof, the reward gap is 
\[
    \Delta(x)=x_1-\frac12,
    \qquad  r_1(x)=\Delta(x) \qquad \mbox{and} \qquad
    r_2=0, \qquad x \in [0, 1]^d.
\]
The following lower bound proof uses two different hard constructions.  When $d=1$,  we use a two-point density perturbation around the unique decision boundary $x_1=1/2$ and apply a Le Cam strategy.  The one-dimensional construction is embedded into $d=2$ below.  When $d\ge3$, we instead use an Assouad hypercube construction by encoding many independent bits in the $(d-1)$ tangential boundary coordinates $u=(x_2,x_3,\ldots,x_d)$.

\noindent\textbf{(1) Proof for  $d=1$ and $d=2$.}

\noindent\textbf{Step 1. Density construction.}
Note that we assume  $L\in (0,1)$. We set $
    c_0=(1-L)/4>0$,
and define
\[
    p^0(x)=\frac12+L\Delta(x),  \qquad  p^1(x)=p^0(x)+c_0.
\]
Clearly, both $p^0, p^1\in\mathcal{P}_{L,1}$.  Moreover, since
$p^0([0,1])\subset[(1-L)/2,(1+L)/2]$ and $c_0=(1-L)/4$, both policies take
values in $[0,1]$.

Let
\[
    s=\frac14,\qquad h=\frac{s}{16},\qquad
    x_+=\frac12+\frac{s}{2},\qquad
    x_-=\frac12-\frac{s}{2},
\]
and define
\[
    I_+=[x_+-h,x_++h],
    \qquad
    I_-=[x_--h,x_-+h].
\]
For  $\varepsilon\in(0,1/8]$, set
\[
    f_+(x) =  \left(\frac12+\varepsilon\right)\frac{\1_{I_+}(x)}{2h} +  \left(\frac12-\varepsilon\right)\frac{\1_{I_-}(x)}{2h},
\]
and
\[
    f_-(x) = \left(\frac12-\varepsilon\right)\frac{\1_{I_+}(x)}{2h} + \left(\frac12+\varepsilon\right)\frac{\1_{I_-}(x)}{2h}.
\]
To keep the densities bounded away from zero, we mix them with the uniform density:
\[
    g_+=(1+f_+)/2,
    \qquad
    g_-=(1+ f_-)/2.
\]
Then, by construction, we have  $g_\pm\ge1/2$ and $g_\pm\le1/2+1/(2h)$, where $h=1/64$ is fixed. Let  $f_0=({f_++f_-})/{2}$ and $g_0=({g_++g_-})/{2}$.

\noindent\textbf{Step 2. Optimality of $p^0$ for the baseline density $g_0$.}
We first prove that $p^0$ is optimal for $J_{f_0}$. Let $\tau(x)=1-x$, which maps $I_+$ onto $I_-$.  Since $f_0$ puts equal mass on the two intervals, it holds for any $q \in \mathcal{P}_{L,1}$ that
\begin{align*}
    \int_0^1 \Delta(x)q(x)f_0(x)\,\mathrm{d}x
    &= \frac{1}{4h}\left\{\int_{I_+}\Delta(x)q(x)\,\mathrm{d}x+\int_{I_-}\Delta(x)q(x)\,\mathrm{d}x\right\} \\
    &= \frac{1}{4h}\left\{\int_{I_+}\Delta(x)q(x)\,\mathrm{d}x+\int_{I_+}\Delta(\tau(x))q(\tau(x))\,\mathrm{d}x\right\} \\
    &= \frac{1}{4h}\left\{\int_{I_+}\Delta(x)q(x)\,\mathrm{d}x-\int_{I_+}\Delta(x)q(\tau(x))\,\mathrm{d}x\right\} \\
    &= \frac{1}{4h}\int_{I_+}\Delta(x)\{q(x)-q(\tau(x))\}\,\mathrm{d}x,
\end{align*}
where the third equality uses $\Delta(\tau(x))=-\Delta(x)$.
For $x\in I_+$, we have that $\Delta(x)\ge0$, and
\[
    q(x)-q(\tau(x))  \le  L|x-\tau(x)| =   2L\Delta(x),
\]
Therefore, we have that
\[
    \int_0^1 \Delta(x)q(x)f_0(x)\,\mathrm{d}x\le \frac{L}{2h}\int_{I_+}\Delta^2(x)\,\mathrm{d}x.
\]
Note that this upper bound is attained by $p^0$, because by construction
\[
    p^0(x)-p^0(\tau(x))=2L\Delta(x),\qquad x\in I_+.
\]
Hence $p^0$ maximizes $J_{f_0}$ over $\mathcal P_{L,1}$.

We then prove that $p^0$ is optimal for $J_{1}$, note that for any $q\in\mathcal{P}_{L,1}$, it holds that
\[
    \int_0^1 \Delta(x)q(x)\,\mathrm{d}x =  \int_{1/2}^1 \Delta(x)\{q(x)-q(1-x)\}\,\mathrm{d}x  \le \int_{1/2}^1 2L\Delta^2(x)\,\mathrm{d}x,
\]
and equality is attained by $p^0$. Indeed, by the definition of $p^0$, for every $x\in[1/2,1]$, we have that
\[
    p^0(x)-p^0(1-x)=2L\Delta(x).
\]
Therefore,
\[
    \int_0^1 \Delta(x)p^0(x)\,\mathrm{d}x=\int_{1/2}^1 \Delta(x)\{p^0(x)-p^0(1-x)\}\,\mathrm{d}x=\int_{1/2}^1 2L\Delta^2(x)\,\mathrm{d}x.
\]
Since by construction, we have that $g_0=(1+f_0)/2$, we have
\[
    J_{g_0}(q)=\frac12\int_0^1\Delta(x)q(x)\,\mathrm{d}x+\frac12\int_0^1\Delta(x)q(x)f_0(x)\,\mathrm{d}x.
\]
Therefore $p^0$ maximizes both components of $g_0$, and hence by linearity, maximizes $J_{g_0}$.

\noindent\textbf{Step 3. Separation of $p^1$ and $p^0$.} We next compute the separation between $p^1$ and $p^0$ under $g_+$.  Since
$p^1-p^0=c_0$ and $\int_0^1 \Delta(x)\,\mathrm{d}x=0$, we have 
\[
\begin{aligned}
    J_{g_+}(p^1)-J_{g_+}(p^0)
    &=  c_0\int_0^1 \Delta(x)g_+(x)\,\mathrm{d}x\\
    &=  \frac{c_0}{2}\int_0^1 \Delta(x)f_+(x)\,\mathrm{d}x.
\end{aligned}
\]
Note that
\[
    \frac{1}{2h}\int_{I_+}\Delta(x)\,\mathrm{d}x=\frac{1}{2h}\int_{x_+-h}^{x_++h}\left(x-\frac12\right)\,\mathrm{d}x=x_+-\frac12=\frac{s}{2},
\]
and similarly,
\[
    \frac{1}{2h}\int_{I_-}\Delta(x)\,\mathrm{d}x=\frac{1}{2h}\int_{x_--h}^{x_-+h}\left(x-\frac12\right)\,\mathrm{d}x=x_--\frac12=-\frac{s}{2}.
\]
Therefore,
\begin{align*}
    \int_0^1 \Delta(x)f_+(x)\,\mathrm{d}x
    &= \left(\frac12+\varepsilon\right)\frac{1}{2h}\int_{I_+}\Delta(x)\,\mathrm{d}x+\left(\frac12-\varepsilon\right)\frac{1}{2h}\int_{I_-}\Delta(x)\,\mathrm{d}x \\
    &= \left(\frac12+\varepsilon\right)\frac{s}{2}+\left(\frac12-\varepsilon\right)\left(-\frac{s}{2}\right) \\
    &= \frac{s}{4}+\frac{\varepsilon s}{2}-\frac{s}{4}+\frac{\varepsilon s}{2} \\
    &= \varepsilon s.
\end{align*}
Thus, we have that
\begin{equation}\label{gpdiff}
    J_{g_+}(p^1)-J_{g_+}(p^0)=c_0\varepsilon s/2.
\end{equation}
By a symmetry argument, we can also show 
\[
    J_{g_-}(p^0)-J_{g_-}(p^1)=c_0\varepsilon s/2.
\]

\noindent\textbf{Step 4. KL divergence}

Conditional on the history and the selected policy, the action kernel
and the conditional reward law given $(X_t,A_t)$ are the same under both
environments. Only the context density changes. The learner's internal
randomization also has the same conditional law. Therefore, the chain rule
for KL divergence gives
\[
    \KL(\mathbb P_{g_+}\|\mathbb P_{g_-})=T\,\KL(g_+\|g_-).
\]
Because $g_-\ge1/2$,
\[
\begin{aligned}
    \KL(g_+\|g_-) &\le \int_0^1 \left(\frac{g_+}{g_-} - 1\right) g_+ \mathrm{d}x \leq \int_0^1\frac{(g_+-g_-)^2}{g_-}\,\mathrm{d}x
    \le  \frac{1}{2} \int_0^1(f_+-f_-)^2\,\mathrm{d}x.
\end{aligned}
\]
On each of $I_+$ and $I_-$, note that $|f_+-f_-|=\varepsilon/h$, and the total length of the interval that they are different is $4h$.  Hence, we have that
\begin{equation}\label{ineq:KL}
      \KL(\mathbb P_{g_+}\|\mathbb P_{g_-})  =   T\,\KL(g_+\|g_-)\le 2T\varepsilon^2/h.
\end{equation}
\noindent\textbf{Step 5. Reduction to regret.}
Set
\[
    \mathcal A = \left\{ J_{g_+}(\bar p) \le \frac{J_{g_+}(p^0)+J_{g_+}(p^1)}2  \right\}.
\]
On $\mathcal A$, by the optimality of  $p_{g_+}^\star$ and \eqref{gpdiff}, we have 
\begin{equation}\label{opt_A}
    J_{g_+}(p_{g_+}^\star)-J_{g_+}(\bar p) \ge J_{g_+}(p^1)-J_{g_+}(\bar p)  \ge \frac{1}{2}[J_{g_+}(p^1)-J_{g_+}(p^0)] \geq   \frac{c_0\varepsilon s}{4}.
\end{equation}
On $\mathcal{A}^c$, we have  
\begin{equation}\label{eq:Ac-def}
J_{g_+}(\bar p)>\frac 1 2 \big(J_{g_+}(p^0)+J_{g_+}(p^1)\big).
\end{equation}
Moreover,  by \eqref{gpdiff}, we have 
\begin{equation}\label{eq:midpoint-lb}
\frac 1 2 \big[J_{g_+}(p^0)+J_{g_+}(p^1)\big] =J_{g_+}(p^0)+\tfrac12\big(J_{g_+}(p^1)-J_{g_+}(p^0)\big) = J_{g_+}(p^0)+\frac{c_0\varepsilon s}{4}.
\end{equation}
Combining \eqref{eq:Ac-def} and \eqref{eq:midpoint-lb} yields that,  on $\mathcal{A}^c,$
\begin{equation}\label{eq:gplus-diff-lb}
J_{g_+}(\bar p)-J_{g_+}(p^0)\;>\frac{c_0\varepsilon s}{4}.
\end{equation}
Next, using the decomposition that for any $p$,
\begin{equation}
J_{g_0}(p)\;=\;\frac12\Big(J_{g_+}(p)+J_{g_-}(p)\Big),
\end{equation}
we have 
\[
J_{g_+}(\bar p)-J_{g_+}(p^0) =\Big(J_{g_0}(\bar p)-J_{g_0}(p^0)\Big) +\frac12\Big( J_{g_+}(\bar p)-J_{g_-}(\bar p) - (J_{g_+}(p^0)-J_{g_-}(p^0))\Big).
\]
Since $p^0$ optimizes on $g_0$, we have $J_{g_0}(\bar p)-J_{g_0}(p^0)\le 0$, so
\[
\frac12\Big( J_{g_+}(\bar p)-J_{g_-}(\bar p) - (J_{g_+}(p^0)-J_{g_-}(p^0))\Big)
\;\ge\; J_{g_+}(\bar p)-J_{g_+}(p^0),
\]
i.e. 
\[
J_{g_-}(p^0)-J_{g_-}(\bar p)\ge J_{g_+}(\bar p)-J_{g_+}(p^0).
\]
Together with \eqref{eq:gplus-diff-lb}, this implies on $\mathcal{A}^c$,
\begin{equation}\label{eq:Ac-gap-gminus}
J_{g_-}(p^0)-J_{g_-}(\bar p)\;>\;\frac{c_0\varepsilon s}{4}.
\end{equation}
Thus, combining \eqref{opt_A} and \eqref{eq:Ac-gap-gminus}, 
\[
    \Reg_{g_+}(\bar{p},T)+\Reg_{g_-}(\bar{p},T) \ge \frac{Tc_0\varepsilon s}{4} \{\mathbb P_{g_+}(\mathcal A)+\mathbb P_{g_-}(\mathcal A^c)\}.
\]
The Bretagnolle--Huber inequality gives
\[
    \mathbb P_{g_+}(\mathcal A)+\mathbb P_{g_-}(\mathcal A^c) \ge \frac12  \exp\{-\KL(\mathbb P_{g_+}\|\mathbb P_{g_-})\}.
\]
Finally, choosing 
\[
    \varepsilon  =   \min\left\{\frac18,\sqrt{\frac{h}{2T}}\right\},
\]
since $h=1/64$ and $T\ge1$, we have
$\varepsilon=\sqrt{h/(2T)}$ and, by \eqref{ineq:KL},
$\KL(\mathbb P_{g_+}\|\mathbb P_{g_-})\le1$. Consequently,
\[
    \max_{g\in\{g_+,g_-\}}\Reg_g(\bar p,T)
    \ge \frac{c_0s}{16e}\sqrt{\frac h2}\,\sqrt T.
\]
For $d=2$, use the product densities $g_\pm(x_1,x_2)=g_\pm(x_1)$.
For every $q\in\mathcal P_{L,1}$, the function
\[
    \widetilde q(x_1)=\int_0^1q(x_1,x_2)\,\mathrm{d}x_2
\]
is a one-dimensional $L$-Lipschitz policy with the same population value.
Extend $p^0,p^1$ independently of $x_2$. Thus the baseline optimality and
separation arguments remain valid, and the independent uniform second
coordinate leaves the KL divergence unchanged. The same testing argument
proves the bound for $d=2$.

\noindent\textbf{(2) Proof for $d\geq 3$.}
We let $u=(x_2,x_3,\ldots,x_d)$, and 
\[
    p_0(x_1,u)=\frac12+L(x_1-1/2).
\]
Since $L<1$, $p_0$ takes values in
$[(1-L)/2,(1+L)/2]$ and belongs to $\mathcal{P}_{L,1}$.

\noindent\textbf{Step 1. Packing in the boundary coordinates}

Fix $h\in(0,1/64]$ and let $M = \lfloor 1/(32h)\rfloor$. For each multi-index
\[
    k=(k_1,\ldots,k_{d-1})\in\{0,1,\ldots,M-1\}^{d-1},
\]
define the grid point as 
\[
    a_k=\left(16hk_1,\ldots,16hk_{d-1}\right)+1/4,
\]
and the boundary cube as $ Q_k=a_k+[0,h]^{d-1}$. Denote $a_{k,r} = 1/4 + 16hk_{r}$, where $r \in \{1,\ldots,d-1\}$. Since
\[
    a_{k,r}+h \le \frac14+16h(M-1)+h \le \frac14+16hM \le \frac34,
\]
it holds that each $Q_k$ is contained in $[1/4,3/4]^{d-1}$. If $k\neq k'$, then for some coordinate $r$, we have that $|k_r-k'_r|\ge1$. Hence the distance between the two cubes in coordinate $r$ is at least $15h$. Since the $\ell_1$ distance between two points is at least the absolute difference in any one coordinate, we get $\dist_1(Q_k,Q_{k'})\ge 15h$.  In particular, the cubes are disjoint, and the number of cubes is $N=M^{d-1}$. In particular, for $h\le1/64$, we have that 
\[
    (64h)^{-(d-1)} \le N\le (32h)^{-(d-1)}.
\]
Inside each $Q_j$, $j\in[N]$, choose two subcubes $A_j$ and $B_j$ of side length $\ell=h/16$ as follows. Let $e_1=(1,0,\ldots,0)\in\mathbb R^{d-1}$ denote the first coordinate direction in the tangential variable $u\in[0,1]^{d-1}$, and set $A_j=a_j+[0,\ell]^{d-1}$ and $B_j = A_j + h e_1/2$. Thus $B_j$ is obtained from $A_j$ by shifting it by distance $h/2$ in the first tangential coordinate. Since $\ell=h/16$, both subcubes are contained in $Q_j=a_j+[0,h]^{d-1}$. Moreover, $ |A_j|=|B_j|=\ell^{d-1}$, and the separation between them in the first tangential coordinate is lower bounded by $h/4$, hence $\dist_1(A_j,B_j)\ge h/4$.

Define the signed bump
\[
    \varphi_j(u)=\indc_{A_j}(u)-\indc_{B_j}(u).
\]
Note that the supports of $\varphi_j$ are disjoint across $j$, and $ \int_{[0,1]^{d-1}}\varphi_j(u)\mathrm{d}u=0$.

\noindent\textbf{Step 2. Density Construction.}
Fix $\alpha=1/8$.  For each
$\theta=(\theta_1,\ldots,\theta_N)\in\{\pm1\}^N$, define
\[
    g_\theta(x_1,u) = 1+\alpha (x_1-1/2)\sum_{j=1}^N\theta_j\varphi_j(u).
\]
Then, it is easy to see that for each $\theta$, $g_{\theta}$ is a density on $\mathcal{X}$ that is upper and lower bounded.  Indeed, for every $u$, at most one $\varphi_j(u)$ is nonzero, and its value is in
$\{-1,0,1\}$. 
Since $|x_1-1/2|\le1/2$ and $\alpha=1/8$, we have
\begin{align*}
    |g_\theta(x_1,u)-1|
    &= \left|\alpha\left(x_1-\frac12\right)\sum_{j=1}^N\theta_j\varphi_j(u)\right| \le \alpha\left|x_1-\frac12\right|\left|\sum_{j=1}^N\theta_j\varphi_j(u)\right| \le \frac1{16}.
\end{align*}
Therefore, it holds that
\[
    \frac{15}{16}\le g_\theta(x_1,u)\le\frac{17}{16}.
\]
It remains to check that $g_\theta$ integrates to one. By the definition of $g_\theta$, we have that
\begin{align*}
    \int_{[0,1]^d}g_\theta(x_1,u)\,\mathrm{d}x_1\,\mathrm{d}u
    &= \int_{[0,1]^d}\left\{1+\alpha\left(x_1-\frac12\right)\sum_{j=1}^N\theta_j\varphi_j(u)\right\}\,\mathrm{d}x_1\,\mathrm{d}u \\
    &= 1+\alpha\sum_{j=1}^N\theta_j\int_0^1\left(x_1-\frac12\right)\,\mathrm{d}x_1\int_{[0,1]^{d-1}}\varphi_j(u)\,\mathrm{d}u \\
    &= 1.
\end{align*}
Thus each $g_\theta$ is a density on $[0,1]^d$.


\noindent\textbf{Step 3. Optimal Policies.}

Assume additionally that $h\le(1-L)/L$. 
For each $j\in[N]$, define $u_j:[0,1]^{d-1}\to\mathbb R$ as follows. For each $u=(x_2,x_3,\ldots,x_d)^{\top}$, if $\dist_1(u,Q_j)\ge h/2$, set $u_j(u)=0$; 
if $\dist_1(u,Q_j)<h/2$, set
\[
u_j(u)=
\begin{cases}
    L\left(h/2-\dist_1(u,Q_j)\right),
      & x_2<a_{j,1}+\dist_1(u,Q_j),\\
    L\left(\dist_1(u,Q_j)-h/2\right),
      & x_2>a_{j,1}+h-\dist_1(u,Q_j),\\
    L(a_{j,1}+h/2-x_2),
      & \text{otherwise.}
\end{cases}
\]
Here we recall that $a_{j,1} = 1/4 + 16hj_{1}$ is the first coordinate of the lower-left
corner \(a_j\) of the cube \(Q_j\).
\begin{figure}[H]
\centering
\begingroup
\definecolor{lbBump}{RGB}{18,112,103}
\definecolor{lbA}{RGB}{34,132,92}
\definecolor{lbB}{RGB}{46,104,175}
\definecolor{lbMuted}{RGB}{100,106,115}
\pgfmathdeclarefunction{lboutside}{1}{\pgfmathparse{max(0,-#1)+max(0,#1-1)}}
\pgfmathdeclarefunction{lbradius}{2}{\pgfmathparse{max(0,0.5-lboutside(#1)-lboutside(#2))}}
\pgfmathdeclarefunction{lbbump}{2}{\pgfmathparse{max(-lbradius(#1,#2),min(lbradius(#1,#2),0.5-#1))}}
\resizebox{\linewidth}{!}{%
\begin{tikzpicture}
\begin{axis}[
  width=8.5cm,height=6.8cm,
  title={(a) A slice through the paired subcubes},
  title style={font=\small,at={(0.5,1.035)}},
  xmin=-0.7,xmax=1.7,ymin=-0.7,ymax=0.82,
  axis lines=left,axis line style={draw=lbMuted},
  xlabel={$x_2$},ylabel={$u_j(u)$},
  xtick={-0.5,0,0.5,1,1.5},
  xticklabels={$a_{j,1}-h/2$,$a_{j,1}$,$a_{j,1}+h/2$,$a_{j,1}+h$,$a_{j,1}+3h/2$},
  xticklabel style={rotate=35,anchor=east},
  ytick={-0.5,0,0.5},yticklabels={$-Lh/2$,$0$,$Lh/2$},
  tick label style={font=\footnotesize},label style={font=\small},clip=false]
  \fill[black!3] (axis cs:0,-0.7) rectangle (axis cs:1,0.72);
  \fill[lbA!17] (axis cs:0,-0.7) rectangle (axis cs:0.0625,0.61);
  \fill[lbB!17] (axis cs:0.5,-0.7) rectangle (axis cs:0.5625,0.61);
  \draw[black!22] (axis cs:-0.7,0)--(axis cs:1.7,0);
  \draw[black!28,densely dashed] (axis cs:0,-0.7)--(axis cs:0,0.61);
  \draw[black!28,densely dashed] (axis cs:1,-0.7)--(axis cs:1,0.61);
  \addplot[thick,lbBump] coordinates {
    (-0.7,0) (-0.5,0) (0,0.5) (1,-0.5) (1.5,0) (1.7,0)
  };
  \node[lbA,anchor=south,font=\small] at (axis cs:0.03125,0.62) {$A_j$};
  \node[lbB,anchor=south,font=\small] at (axis cs:0.53125,0.62) {$B_j$};
  \draw[-{Stealth[length=1.1mm]},lbA] (axis cs:0.03125,0.60)--(axis cs:0.03125,0.51);
  \draw[-{Stealth[length=1.1mm]},lbB] (axis cs:0.53125,0.60)--(axis cs:0.53125,0.025);
  \addplot[only marks,mark=*,mark size=1.8pt,lbA] coordinates {(0.03125,0.46875)};
  \addplot[only marks,mark=*,mark size=1.8pt,lbB] coordinates {(0.53125,-0.03125)};
  \draw[lbMuted,densely dashed] (axis cs:0.03125,0.46875)--(axis cs:0.91,0.46875);
  \draw[lbMuted,densely dashed] (axis cs:0.53125,-0.03125)--(axis cs:0.91,-0.03125);
  \draw[{Stealth[length=1.2mm]}-{Stealth[length=1.2mm]},thick]
    (axis cs:0.87,-0.03125)--(axis cs:0.87,0.46875);
  \node[anchor=west,align=left,font=\footnotesize] at (axis cs:0.94,0.26)
    {difference\\$Lh/2$};
  \node[font=\footnotesize,lbBump] at (axis cs:0.29,-0.18) {slope $-L$};
  \node[font=\footnotesize,lbMuted] at (axis cs:-0.28,0.43) {taper};
  \node[font=\footnotesize,lbMuted] at (axis cs:1.30,-0.40) {taper};
  \draw[{Stealth[length=1mm]}-{Stealth[length=1mm]},black!60]
    (axis cs:0,-0.62)--(axis cs:1,-0.62);
  \node[fill=black!3,font=\footnotesize,inner sep=1pt] at (axis cs:0.5,-0.62) {$Q_j$};
\end{axis}

\begin{axis}[
  at={(9.35cm,0)},anchor=south west,width=6.8cm,height=6.8cm,
  title={(b) The tangential plane ($d=3$)},
  title style={font=\small,at={(0.5,1.035)}},
  view={0}{90},axis equal image,
  xmin=-0.7,xmax=1.7,ymin=-0.7,ymax=1.7,
  xlabel={$x_2$},ylabel={$x_3$},
  xtick={0,1},xticklabels={$a_{j,1}$,$a_{j,1}+h$},
  ytick={0,1},yticklabels={$a_{j,2}$,$a_{j,2}+h$},
  tick label style={font=\footnotesize},label style={font=\small},
  colormap={lbsigned}{rgb255(0cm)=(42,106,164);rgb255(1cm)=(249,249,249);rgb255(2cm)=(177,49,68)},
  point meta min=-0.5,point meta max=0.5,
  colorbar horizontal,
  colorbar style={at={(0.5,-0.24)},anchor=north,width=4.4cm,height=0.16cm,
    xtick={-0.5,0,0.5},xticklabels={$-Lh/2$,$0$,$Lh/2$},
    xlabel={$u_j(u)$},xlabel style={font=\footnotesize},
    tick label style={font=\footnotesize}},clip=false]
  \addplot3[surf,shader=interp,draw=none,domain=-0.7:1.7,y domain=-0.7:1.7,samples=49]
    {lbbump(x,y)};
  \draw[thick] (axis cs:0,0,0.6) rectangle (axis cs:1,1,0.6);
  \draw[densely dashed,black!65]
    (axis cs:-0.5,0,0.6)--(axis cs:-0.5,1,0.6)--(axis cs:0,1.5,0.6)
    --(axis cs:1,1.5,0.6)--(axis cs:1.5,1,0.6)--(axis cs:1.5,0,0.6)
    --(axis cs:1,-0.5,0.6)--(axis cs:0,-0.5,0.6)--cycle;
  \draw[black!55,dash pattern=on 2pt off 2pt]
    (axis cs:-0.6,0.03125,0.6)--(axis cs:1.6,0.03125,0.6);
  \draw[lbA,thick,fill=white] (axis cs:0,0,0.6) rectangle (axis cs:0.0625,0.0625,0.6);
  \draw[lbB,thick,fill=white] (axis cs:0.5,0,0.6) rectangle (axis cs:0.5625,0.0625,0.6);
  \node[font=\small,fill=white,fill opacity=0.85,text opacity=1,inner sep=1pt]
    at (axis cs:0.5,0.78,0.6) {$Q_j$};
  \node[lbA,font=\small] at (axis cs:0.02,-0.33,0.6) {$A_j$};
  \node[lbB,font=\small] at (axis cs:0.66,-0.33,0.6) {$B_j$};
  \draw[-{Stealth[length=1mm]},lbA] (axis cs:0.02,-0.23,0.6)--(axis cs:0.03,-0.02,0.6);
  \draw[-{Stealth[length=1mm]},lbB] (axis cs:0.66,-0.23,0.6)--(axis cs:0.54,-0.02,0.6);
  \node[font=\scriptsize,fill=white,inner sep=1pt] at (axis cs:0.5,1.61,0.6)
    {$\operatorname{dist}_1(u,Q_j)=h/2$};
\end{axis}
\end{tikzpicture}%
}
\endgroup
\caption{Localized perturbation used in the lower-bound construction.
Panel (a) shows a slice through the centres of $A_j$ and
$B_j=A_j+(h/2)e_1$. Within $Q_j$, the function $u_j$ is affine with slope
$-L$ in $x_2$, so every paired difference
$u_j(u)-u_j(u+(h/2)e_1)$, $u\in A_j$, equals $Lh/2$, attaining the
Lipschitz bound. Outside $Q_j$, the perturbation tapers to zero and
vanishes when $\operatorname{dist}_1(u,Q_j)\ge h/2$.
Panel (b) displays the perturbation for $d=3$, with red and blue denoting
positive and negative values. The solid square is $Q_j$, the outer
dashed outline bounds its $h/2$-neighbourhood in the $\ell_1$ metric,
and the horizontal dashed line indicates the slice in panel (a).}
\label{fig:lb-step3-construction}
\end{figure}
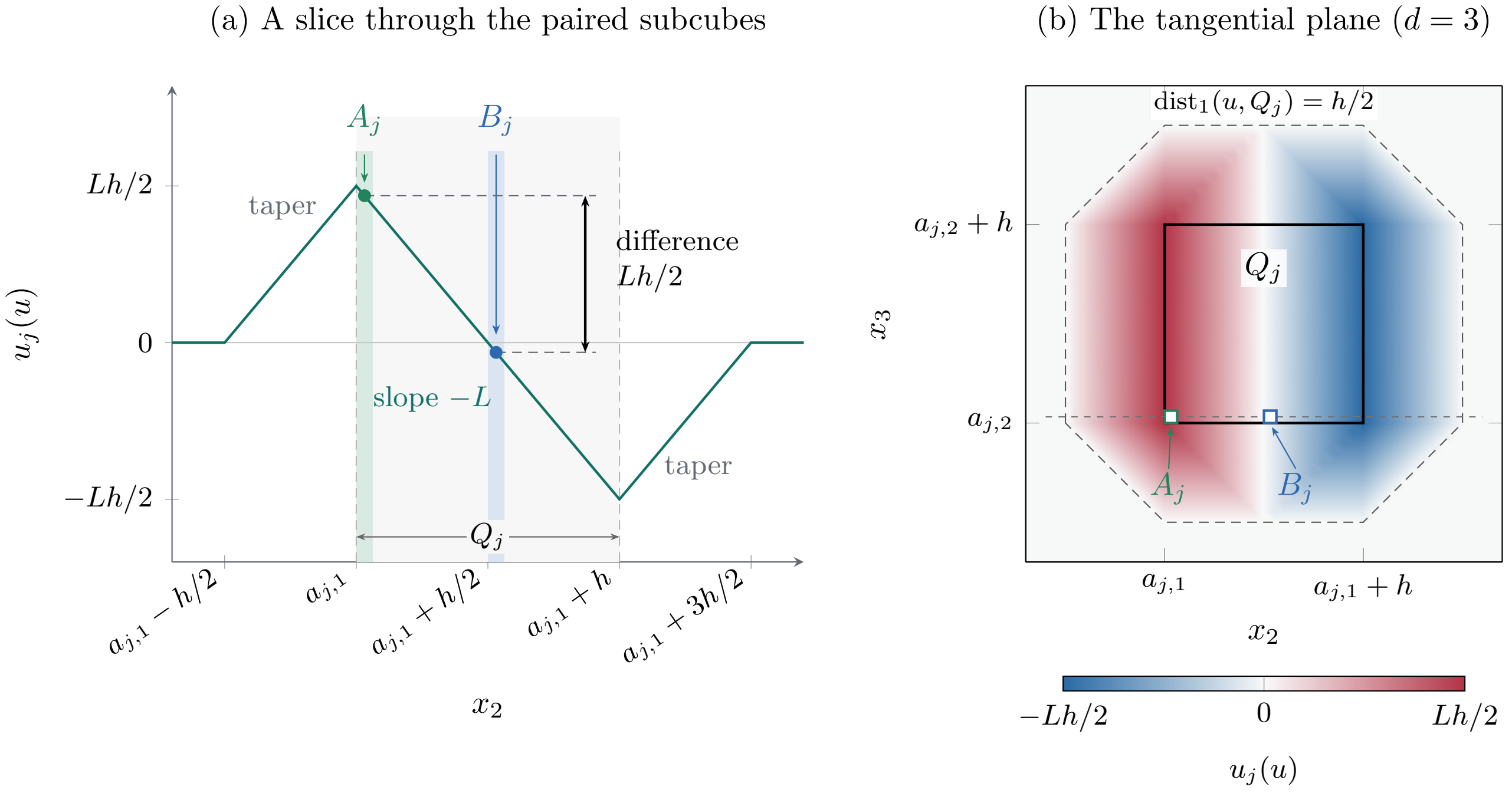

Intuitively, $u_j$ is a clipped affine function that retains slope $-L$ in the direction of $x_2$ on $Q_j$ and vanishes away from $Q_j$, thereby attaining the largest paired difference between $A_j$ and $B_j$ allowed by the Lipschitz constraint. Figure \ref{fig:lb-step3-construction} visualizes the construction.

Now set
\[
    U_\theta(u)=\sum_{j=1}^N\theta_j u_j(u),
    \qquad
    p_\theta^\circ(x_1,u)=p_0(x_1,u)+U_\theta(u).
\]
By Lemma \ref{lem:oracle-feasible}, $p_\theta^\circ$ is feasible and maximizes
$J_1$. By Lemma \ref{lem:bridge}, it also maximizes $J_{g_\theta}$.

\noindent\textbf{Step 4. Reduction of regret  to bit errors.}

For $q\in\mathcal{P}_{L,1}$ and $j \in [N]$, define
\[
    S_j(q)  = \int_{[0,1]^d} (x_1-1/2)^2\varphi_j(u)q(x_1,u)\,\mathrm{d}x_1\,du,
\]
and the plug-in bit estimator
\[
    \widehat\theta_j(q)
    =
    \begin{cases}
    1, & S_j(q)\ge0,\\
    -1, & S_j(q)<0.
    \end{cases}
\]

In Lemma \ref{lem:bridge}, we reduce the regret to classification error, and we have that for every $\theta \in \{\pm 1\}^N$ and every $q\in\mathcal{P}_{L,1}$,
\[
    J_{g_\theta}(p_\theta^\circ)-J_{g_\theta}(q) \ge  c_2\alpha L h^{d}
    \sum_{j=1}^N \1\{\widehat\theta_j(q)\ne\theta_j\}.
\]

\noindent\textbf{Step 5. Neighbouring KL divergence.}

Let $\theta^{(j)}$ be the vector obtained from $\theta$ by flipping the
$j$th coordinate.  Let $\mathbb P_\theta$ be the joint law of the full
interaction, including the learner's internal randomization, under $g_\theta$. Then by Lemma \ref{lem:neighbor-kl}, there exists $C_{\rm KL}>0$, depending only on $d$, such that for every
$\theta \in \{\pm 1\}^N$ and $j \in [N]$,
\[
    \KL(\mathbb P_\theta\|\mathbb P_{\theta^{(j)}}) \le C_{\rm KL}T\alpha^2h^{d-1}. 
\]

\noindent\textbf{Step 6. Bounding the classification error by Assouad's method}

For a fixed $\theta \in \{\pm 1\}^N$,
\[
\begin{aligned}
    \Reg_{g_\theta}(\bar{p},T)
    &=  T\,\mathbb{E}_\theta\{J_{g_\theta}(p_{g_\theta}^\star)-J_{g_\theta}(\bar p)\}\ge  T\,\mathbb{E}_\theta\{J_{g_\theta}(p_\theta^\circ)-J_{g_\theta}(\bar p)\}.
\end{aligned}
\]
Applying Lemma \ref{lem:bridge} with $q=\bar p$, we have 
\[
    \Reg_{g_\theta}(\bar{p},T) \ge  Tc_2\alpha L h^{d}\sum_{j=1}^N \mathbb{P}_\theta\{\widehat\theta_j(\bar p)\ne\theta_j\}.
\]
Average this inequality over $\theta$ gives
\begin{flalign}\label{ineq:avgregret}
    \frac1{2^N}\sum_{\theta \in \{\pm 1\}^N} \Reg_{g_\theta}(\bar{p},T) \ge Tc_2\alpha L h^{d} \sum_{j=1}^N \frac1{2^N}\sum_\theta  \mathbb{P}_\theta\{\widehat\theta_j(\bar p)\ne\theta_j\}.
\end{flalign}
Now, we choose
\[
    h_0=\min\left\{\frac1{64},\frac{1-L}{L}\right\},
    \qquad h=h_0T^{-1/(d-1)}.
\]
By Lemma \ref{lem:neighbor-kl}, with
$C_{\rm KL}=32/(45\cdot16^{d-1})$ and $\alpha=1/8$,
\[
    \KL(\mathbb P_\theta\|\mathbb P_{\theta^{(j)}})
    \le C_{\rm KL}\alpha^2h_0^{d-1}
    =\frac{h_0^{d-1}}{90\cdot16^{d-1}}\le1.
\]
For each $j$, pair every $\theta$ with $\theta_j=1$ with $\theta^{(j)}$,
and let $E_j=\{\widehat\theta_j(\bar p)=1\}$. For each such pair,
\[
\begin{aligned}
    &\mathbb P_\theta\{\widehat\theta_j(\bar p)\ne\theta_j\}
    +\mathbb P_{\theta^{(j)}}
        \{\widehat\theta_j(\bar p)\ne\theta_j^{(j)}\}\\
    &\qquad=\mathbb P_\theta(E_j^c)+\mathbb P_{\theta^{(j)}}(E_j)
    \ge\frac12\exp\{-\KL(\mathbb P_\theta\|\mathbb P_{\theta^{(j)}})\}
    \ge\frac1{2e},
\end{aligned}
\]
where the first inequality is the Bretagnolle--Huber inequality.
There are $2^{N-1}$ such pairs, so
\[
    \frac1{2^N}\sum_\theta
    \mathbb P_\theta\{\widehat\theta_j(\bar p)\ne\theta_j\}
    \ge\frac1{4e}>\frac1{12}.
\]
Substituting this into \eqref{ineq:avgregret} and using
$N\ge(64h)^{-(d-1)}$ gives
\[
\begin{aligned}
    \frac1{2^N}\sum_\theta\Reg_{g_\theta}(\bar p,T)
    &\ge\frac{Tc_2\alpha Lh^{d}N}{12}\ge\frac{c_2\alpha L}{12\cdot64^{d-1}}Th=\frac{c_2\alpha Lh_0}{12\cdot64^{d-1}}
      T^{(d-2)/(d-1)}.
\end{aligned}
\]
Thus, for every learner and every integer $T\ge1$, at least one density
in the construction has regret of the claimed order. Taking the infimum
over learners completes the proof.
 \end{proof}

\subsection{Auxiliary results for lower bound construction}


\begin{lemma}[Uniform-density optimizer]\label{lem:uniform-opt}
For every $q\in\mathcal{P}_{L,1}$, we have that
\[
    J_1(q):=\int_{[0,1]^d} (x_1-1/2)q(x_1,u)\,\mathrm{d}x_1\,du
    \le
    J_1(p_0).
\]
Moreover, if $U:[0,1]^{d-1}\mapsto \mathbb{R}$ is measurable and 
$p_0(x_1,u)+U(u)\in[0,1]$, then
\[
    J_1(p_0+U)=J_1(p_0).
\]
\end{lemma}

\begin{proof}
For each fixed $u\in[0,1]^{d-1}$, we pair every point $(x_1,u)$ with its reflection $(1-x_1,u)$ across the hyperplane $x_1=1/2$. Since $\Delta(x_1,u)=x_1-1/2$, this reflection changes the sign of the reward gap:
\[
    \Delta(1-x_1,u)=1-x_1-\frac12=-\left(x_1-\frac12\right)=-\Delta(x_1,u).
\]
Therefore, for any fixed $u$ and $q\in\mathcal{P}_{L,1}$, we have that
\begin{align*}
    &\int_0^1 \Delta(x)q(x_1,u)\,\mathrm{d}x_1\\
    =& \int_0^{1/2}\left(x_1-\frac12\right)q(x_1,u)\,\mathrm{d}x_1+\int_{1/2}^1\left(x_1-\frac12\right)q(x_1,u)\,\mathrm{d}x_1 \\
    =& -\int_{1/2}^1\left(x_1-\frac12\right)q(1-x_1,u)\,\mathrm{d}x_1+\int_{1/2}^1\left(x_1-\frac12\right)q(x_1,u)\,\mathrm{d}x_1 \\
    =& \int_{1/2}^1\left(x_1-\frac12\right)\{q(x_1,u)-q(1-x_1,u)\}\,\mathrm{d}x_1.
\end{align*}

    
For $x_1\ge1/2$, we have that $x_1-1/2\ge0$, and
\[
    q(x_1,u)-q(1-x_1,u) \le |q(x_1,u)-q(1-x_1,u)| \le L\|(x_1,u)-(1-x_1,u)\|_1  =  2L(x_1-1/2).
\]
Therefore we have that
\[
    \int_0^1 \Delta(x)q(x_1,u)\,\mathrm{d}x_1 \le \int_{1/2}^1 2L (x_1-1/2)^2\,\mathrm{d}x_1=\int_0^1L (x_1-1/2)^2\,\mathrm{d}x_1.
\]
Integrating over $u$ gives
\[
    J_1(q)\le L\int_0^1 (x_1-1/2)^2\,\mathrm{d}x_1=J_1(p_0),
\]
because $p_0(x_1,u)=1/2+L(x_1-1/2)$ and $\int_0^1(x_1-1/2)\,\mathrm{d}x_1=0$.

For the second claim, it is easy to see that 
\[
    J_1(p_0+U)-J_1(p_0) =  \left(\int_0^1 (x_1-1/2)\,\mathrm{d}x_1\right)
    \left(\int_{[0,1]^{d-1}}U(u)\,du\right) = 0.
\]
\end{proof}

\begin{lemma}[Feasibility and Optimality]\label{lem:oracle-feasible}
For every $\theta$, $p_\theta^\circ\in\mathcal{P}_{L,1}$ and
\[
    J_1(p_\theta^\circ)=J_1(p_0)=\sup_{q\in\mathcal{P}_{L,1}}J_1(q).
\]
\end{lemma}

\begin{proof}
We first prove that $U_\theta$ is $L$-Lipschitz. The piecewise definition
in Step 3 is equivalent to clipping $L(a_{j,1}+h/2-x_2)$ to the interval
\[
    \left[-\bigl(Lh/2-L\dist_1(u,Q_j)\bigr)_+,
    \bigl(Lh/2-L\dist_1(u,Q_j)\bigr)_+\right],
\]
where $t_+=\max\{t,0\}$. Since $u\mapsto\dist_1(u,Q_j)$ is
$1$-Lipschitz, the affine function and both clipping bounds are
$L$-Lipschitz. Clipping is a composition of a pointwise minimum and
maximum, each of which preserves a common Lipschitz bound.
Thus each $u_j$ is $L$-Lipschitz.

Moreover, $\|u_j\|_\infty\le Lh/2$, and its support is contained in the
$h/2$-neighbourhood of $Q_j$. Thus, for $j\ne j'$,
\[
    \dist_1\{\mathrm{supp}(u_j),\mathrm{supp}(u_{j'})\}
    \ge15h-h=14h.
\]
In particular, at most one $u_j$ is nonzero at any point.
If $u$ and $v$ belong to different supports, then
\[
    |U_\theta(u)-U_\theta(v)|\le Lh\le L\|u-v\|_1.
\]
If they belong to the same support, or one lies outside all supports,
the same inequality follows from the Lipschitz property of the
corresponding $u_j$. If both lie outside all supports, the difference
is zero. Therefore $U_\theta$ is $L$-Lipschitz.

Recall that $p_0$ is $L$-Lipschitz in the $x_1$ coordinate.
For any $(x_1,u),(x_1',v)$,
\[
    |p_\theta^\circ(x_1,u)-p_\theta^\circ(x_1',v)|
    \le L|x_1-x_1'|+L\|u-v\|_1
    =L\|(x_1,u)-(x_1',v)\|_1.
\]
Moreover, since the supports are disjoint,
\[
    \|U_\theta\|_\infty\le \frac{Lh}{2}\le\frac{1-L}{2}.
\]
Since $p_0$ is at distance at least $(1-L)/2$ from the endpoints
$0$ and $1$, $p_\theta^\circ$ takes values in $[0,1]$.
Hence $p_\theta^\circ\in\mathcal P_{L,1}$.

Finally, Lemma \ref{lem:uniform-opt} gives
$J_1(p_\theta^\circ)=J_1(p_0)=\sup_{q\in\mathcal P_{L,1}}J_1(q)$.
\end{proof}

\begin{lemma}[One-bit signal]\label{lem:one-bit}
There is a constant $c_1>0$, depending only on $d$, such that for every
$\theta \in \{\pm 1\}^N$ and every $j \in [N]$,
\[
    \theta_j S_j(p_\theta^\circ)
   =
    c_1\ L h^{d} .
\]
\end{lemma}

\begin{proof}
Since $\int\varphi_j(u)\,du=0$ and $p_0$ depends only on $x_1$,
we have $S_j(p_0)=0$. Therefore
\[
    S_j(p_\theta^\circ)
    =\int_{[0,1]^d}(x_1-1/2)^2\varphi_j(u)U_\theta(u)
      \,\mathrm{d}x_1\,\mathrm{d}u.
\]
For $u\in Q_j$, note that 
\[
    |L(a_{j,1}+h/2-x_2)|\le Lh/2.
\]
Consequently, $u_j(u)=L(a_{j,1}+h/2-x_2)$ throughout $Q_j$.
Also, $u_k(u)=0$ on $Q_j$ for $k\ne j$, by the separation of the cubes.
Since $B_j=A_j+(h/2)e_1$, for every $u\in A_j$ we have
\[
    U_\theta(u)-U_\theta(u+(h/2)e_1)=\theta_jLh/2.
\]
Changing variables in the integral over $B_j$ therefore gives
\[
\begin{aligned}
    \theta_jS_j(p_\theta^\circ)
    &=\int_0^1(x_1-1/2)^2\,\mathrm{d}x_1
      \int_{A_j}\theta_j
      \{U_\theta(u)-U_\theta(u+(h/2)e_1)\}\,\mathrm{d}u\\
    &=\frac1{12}\cdot\frac{Lh}{2}\cdot\ell^{d-1}
     =\frac{Lh^d}{24\cdot16^{d-1}}.
\end{aligned}
\]
This proves the claim with $c_1=1/(24\cdot16^{d-1})$.
\end{proof}

\begin{lemma}[Bridge inequality]\label{lem:bridge}
There exists $c_2>0$, depending only on $d$, such that for every $\theta \in \{\pm 1\}^N$ and
every $q\in\mathcal{P}_{L,1}$,
\[
    J_{g_\theta}(p_\theta^\circ)-J_{g_\theta}(q) \ge  c_2\alpha  L h^{d}
    \sum_{j=1}^N \1\{\widehat\theta_j(q)\ne\theta_j\}. 
\]
\end{lemma}

\begin{proof}
For every $q\in\mathcal P_{L,1}$, the definition of $g_\theta$ gives
\[
    J_{g_\theta}(q)=J_1(q)+\alpha\sum_{j=1}^N\theta_jS_j(q).
\]
Therefore
\[
    J_{g_\theta}(p_\theta^\circ)-J_{g_\theta}(q)
    =\{J_1(p_\theta^\circ)-J_1(q)\}
     +\alpha\sum_{j=1}^N\theta_j
       \{S_j(p_\theta^\circ)-S_j(q)\}.
\]
The first term is non-negative by Lemma \ref{lem:oracle-feasible}.
Note that $B_j=A_j+(h/2)e_1$, we can write
\[
    S_j(q)=\int_0^1(x_1-1/2)^2
    \int_{A_j}\{q(x_1,u)-q(x_1,u+(h/2)e_1)\}
    \,\mathrm{d}u\,\mathrm{d}x_1.
\]
The Lipschitz constraint and the proof of
Lemma \ref{lem:one-bit} imply
\[
    |S_j(q)|\le\frac1{12}\cdot\frac{Lh}{2}\cdot\ell^{d-1}
    =\theta_jS_j(p_\theta^\circ).
\]
Thus $\theta_j\{S_j(p_\theta^\circ)-S_j(q)\}\ge0$ for every $j$.
If $\widehat\theta_j(q)\ne\theta_j$, then $\theta_jS_j(q)\le0$.
Hence
\[
    \theta_j\{S_j(p_\theta^\circ)-S_j(q)\}
    \ge\theta_jS_j(p_\theta^\circ)
      \1\{\widehat\theta_j(q)\ne\theta_j\}.
\]
Applying Lemma \ref{lem:one-bit} and summing over $j$ proves the claim,
with $c_2=c_1=1/(24\cdot16^{d-1})$.
\end{proof}

\begin{lemma}[Neighbors KL]\label{lem:neighbor-kl}
There exists $C_{\rm KL}>0$, depending only on $d$, such that for every
$\theta \in \{\pm 1\}^N$ and $j \in [N]$,
\[
    \KL(\mathbb P_\theta\|\mathbb P_{\theta^{(j)}}) \le  C_{\rm KL}T\alpha^2h^{d-1}.
\]
One may take $C_{\rm KL}=32/(45\cdot16^{d-1})$.
\end{lemma}

\begin{proof}
The reward model is the same under all $\theta$; only the context density
changes.  Conditional on the history before round $t$, the learner uses the
same randomized decision rule under both environments, and the conditional
reward law given $(X_t,A_t)$ is also the same.  Therefore the chain rule for KL
gives
\[
    \KL(\mathbb P_\theta\|\mathbb P_{\theta^{(j)}}) = T\,\KL(g_\theta\|g_{\theta^{(j)}}).
\]
Moreover, by construction,
\[
    g_\theta(x_1,u)-g_{\theta^{(j)}}(x_1,u) = 2\alpha(x_1-1/2)\theta_j \varphi_j(u).
\]
Since $g_{\theta^{(j)}}\ge15/16$ by construction, 
\[
\begin{aligned}
    \KL(g_\theta\|g_{\theta^{(j)}})
    &\le
    \int_{[0,1]^d}
    \frac{(g_\theta-g_{\theta^{(j)}})^2}{g_{\theta^{(j)}}}\,\mathrm{d}x_1\,du\\
    &\le
    \frac{16}{15}
    \int_{[0,1]^d}4\alpha^2 (x_1-1/2)^2\varphi_j(u)^2\,\mathrm{d}x_1\,du\\
    &=
    \frac{64}{15}\alpha^2
    \left(\int_0^1(x_1-1/2)^2\,\mathrm{d}x_1\right)
    \{|A_j|+|B_j|\}\\
    &=
    C_{\rm KL}\alpha^2h^{d-1}.
\end{aligned}
\]
Multiplying by $T$ proves the result.
\end{proof}

\section{One-dimensional algorithms and theoretical guarantees}\label[appendix]{sec_app_one_dimensional}
This section collects the algorithms and theoretical guarantees for $d=1$. We first present a one-dimensional canonical extension for the fair oracle, and then present the greedy algorithm for continuous covariates and the UCB algorithm for general covariates separately, together with their regret guarantees and proofs.

\subsection{Common setup and canonical representation}\label[appendix]{sec_app_one_d_setup}
When $d=1$, the decision boundary is a single threshold point instead of a $(d-1)$-dimensional set. We take $\mathcal X=[0,1]$ and use the following one-dimensional analogue of \Cref{a_theta} for both algorithms.

\begin{assumption}[One-dimensional boundary parameters]\label[assumption]{a_theta_1d}
Suppose that there exist absolute constants $A_0>0$, $0<c_0<C_0<\infty$, and $\rho_0\in(0,1/2)$ such that
\[
\Theta_{\rm 1d}=\left\{(a,c)\in\mathbb R^2: |a|\le A_0,\ c_0\le c\le C_0,\ \rho_0\le -\frac ac\le1-\rho_0\right\}.
\]
We assume that $\theta^*=\beta_1^*-\beta_2^*=(a^*,c^*)\in\Theta_{\rm 1d}$.
\end{assumption}

Let $\widetilde A_0>A_0$, $0<\widetilde c_0<c_0<C_0<\widetilde C_0<\infty$, and $0<\widetilde\rho_0<\rho_0$. Define the enlarged parameter set
\[
\widetilde\Theta_{\rm 1d}=\left\{(a,c)\in\mathbb R^2: |a|\le \widetilde A_0,\ \widetilde c_0\le c\le \widetilde C_0,\ \widetilde\rho_0\le -\frac ac\le1-\widetilde\rho_0\right\}.
\]
By the strict enlargement, $\theta^*\in\widetilde\Theta_{\rm 1d}$ and there exists a constant $r_{\rm loc}^{(1)}>0$, depending only on $\theta^*$ and the constants above, such that $\{\theta:\|\theta-\theta^*\|_2\le2r_{\rm loc}^{(1)}\}\subseteq\widetilde\Theta_{\rm 1d}$.

For $\theta=(a,c)\in\widetilde\Theta_{\rm 1d}$, write $ x_\theta=-a/c$.
For $q\in[0,1]$, define the one-dimensional policy
\[
    p_{\theta,q}^{(1)}(x)=\left[q+L(x-x_\theta)\right]_{[0,1]},
    \qquad x\in[0,1].
\]
This is the one-dimensional version of the boundary-canonical policy, with the scalar $q$ playing the role of the value of the policy at the threshold $x_\theta$.

\begin{lemma}[One-dimensional canonical form and empirical optimisation]\label[lemma]{lem:one-d-canonical}
Suppose \Cref{a_theta_1d} holds. For every $\theta\in\widetilde\Theta_{\rm 1d}$ and $q\in[0,1]$, we have $p_{\theta,q}^{(1)}\in\mathcal P_{L,\nu}$. For any $p\in\mathcal P_{L,\nu}$, taking $q=p(x_\theta)$ gives
\begin{equation*}
(\widetilde x^\top\theta)p_{\theta,q}^{(1)}(x)\ge(\widetilde x^\top\theta)p(x),\qquad x\in[0,1].
\end{equation*}
Consequently, there exists $q^*\in[0,1]$ such that $p^*=p_{\theta^*,q^*}^{(1)}$ is a population-optimal fair policy. Moreover, with $H_s$ as in \eqref{eq:Hs} and $d=1$, for every $s\ge2$ and $\theta\in\widetilde\Theta_{\rm 1d}$, let $z^\theta$ be an optimiser in the definition of $H_s(\theta)$ and let $\overline p^\theta$ be its McShane extension. Then $q^\theta=\overline p^\theta(x_\theta)$ satisfies
\begin{equation}
q^\theta\in\argmax_{q\in[0,1]}\frac1{s-1}\sum_{i=1}^{s-1}(\widetilde X_i^\top\theta)p_{\theta,q}^{(1)}(X_i)
\end{equation}
and
\begin{equation*}
H_s(\theta)=\frac1{s-1}\sum_{i=1}^{s-1}(\widetilde X_i^\top\theta)p_{\theta,q^\theta}^{(1)}(X_i)=\max_{q\in[0,1]}\frac1{s-1}\sum_{i=1}^{s-1}(\widetilde X_i^\top\theta)p_{\theta,q}^{(1)}(X_i).
\end{equation*}
\end{lemma}

\begin{proof}
Since $\mathcal H_\theta=\{x_\theta\}$ and $c>0$, the canonical policy in \eqref{eq_canonical} is exactly $p_{\theta,q}^{(1)}$. Its feasibility, pointwise domination for every $\theta$, and population optimality at $\theta^*$ follow from \Cref{thm:canonical_extension}. Applying the general finite-LP argument in \Cref{sec_app_proof_empirical_optimisation} with this singleton boundary gives $q^\theta=\overline p^\theta(x_\theta)$ and the displayed equality with $H_s(\theta)$.
\end{proof}

\subsection{Greedy algorithm for continuous covariates}\label[appendix]{sec_app_one_d_greedy}
\subsubsection{Algorithm}\label[appendix]{sec_app_alg_d1}
\Cref{alg:lazy-one-dimensional-fair-greedy} is the one-dimensional version of \Cref{alg:lazy-adaptive-boundary-greedy}.
\begin{algorithm}[!htbp]
\caption{Lazy adaptive one-dimensional fair greedy algorithm ($d=1$)}
\label{alg:lazy-one-dimensional-fair-greedy}
\begin{algorithmic}
\Require{Time horizon $T$, fairness parameter $(L,\nu)$, ridge parameter $\lambda$, exploration period $T_o$}
\For{$k\in\{1,2\}$}
\State $V_k\gets\lambda I_2$, $B_k\gets0$.
\EndFor
\State $s\gets T_o$, $\widehat p(\cdot)=1/2$.
\For{$t=1,2,\ldots,2T_o$}
\State Observe $X_t$, set $\widetilde X_t=(1,X_t)^\top$.
\State Play arm $1$ with probability $1/2$ and arm $2$ otherwise; observe $Y_t$.
\State Update $V_{A_t}\gets V_{A_t}+\widetilde X_t\widetilde X_t^\top$ and $B_{A_t}\gets B_{A_t}+\widetilde X_tY_t$.
\EndFor
\For{$t=2T_o+1,\ldots,T$}
\If{$t\ge2s$}
\State $s\gets t$.
\For{$k\in\{1,2\}$}
\State Set $\widehat\beta_{k,s-1}\gets V_k^{-1}B_k$.
\EndFor
\State Set $\widehat\theta_s\gets\widehat\beta_{1,s-1}-\widehat\beta_{2,s-1}$ and write $\widehat\theta_s=(\widehat a_s,\widehat c_s)$.
\If{$\widehat c_s>0$ and $-\widehat a_s\in[0,\widehat c_s]$}
\State Set $\widehat x_s\gets-\widehat a_s/\widehat c_s$ and compute
\[
\widehat q_s\in\argmax_{q\in[0,1]}\frac1{s-1}\sum_{i=1}^{s-1}\widetilde X_i^\top\widehat\theta_s\left[q+L(X_i-\widehat x_s)\right]_{[0,1]}.
\]
\State Set $\widehat p(\cdot)\gets p_{\widehat\theta_s,\widehat q_s}^{(1)}(\cdot)$.
\Else
\State Set $\widehat p(\cdot)\gets1/2$.
\EndIf
\EndIf
\State Observe $X_t$ and set $\widetilde X_t=(1,X_t)^\top$.
\State Play arm $1$ with probability $\widehat p(X_t)$ and arm $2$ otherwise; observe $Y_t$.
\State Update $V_{A_t}\gets V_{A_t}+\widetilde X_t\widetilde X_t^\top$ and $B_{A_t}\gets B_{A_t}+\widetilde X_tY_t$.
\EndFor
\end{algorithmic}
\end{algorithm}

\subsubsection{Theoretical guarantee}
\begin{theorem}[One-dimensional regret upper bound]\label[theorem]{thm:one-d-regret}
Suppose Assumptions~\ref{a_contexts}-\ref{a_beta} and \ref{a_theta_1d} hold with $d=1$ and $0 < L \leq L_0'$. For $\eta\in(0,1/2)$, suppose that $2T_o\le T$ and $T_o$ satisfies
\[
    T_o\ge C_T\log(T/\eta)\vee\frac{C_{\rm loc}^{(1)}\log(T/\eta)}{\{r_{\rm loc}^{(1)}\}^2}.
\]
Then the policy played at every round in \Cref{alg:lazy-one-dimensional-fair-greedy} is $(L,\nu)$-fair. Moreover, with probability at least $1-\eta$,
\[
    \mathfrak R_T\le C_1T_o+C_2\sqrt{T\log(eT^2/\eta)},
\]
where $C_1,C_2<\infty$ depend only on $L,c_g,C_g,\sigma,C_\beta,\lambda$ and the constants defining $\widetilde\Theta_{\rm 1d}$. In particular, if $T\ge C_0\log(T/\eta)$ for a sufficiently large constant $C_0$ depending only on these parameters, $r_{\rm loc}^{(1)}$ and the fixed exploration multiplier, taking $T_o\asymp\sqrt{T\log(T/\eta)}$ gives
\[
    \mathfrak R_T\le C\sqrt{T\log(eT^2/\eta)}.
\]
\end{theorem}

\subsubsection{Proof of Theorem \ref{thm:one-d-regret}}
\begin{lemma}[Threshold perturbation]\label[lemma]{lem:one-d-threshold-perturbation}
For any $\theta=(a,c),\theta'=(a',c')\in\widetilde\Theta_{\rm 1d}$,
\[
    |x_\theta-x_{\theta'}|\le C_x\|\theta-\theta'\|_2,
\]
where $C_x=\widetilde c_0^{-1}+\widetilde A_0\widetilde c_0^{-2}$.
\end{lemma}

\begin{proof}
Since $x_\theta=-a/c$ and $x_{\theta'}=-a'/c'$, we have
\[
|x_\theta-x_{\theta'}|\le \frac{|a-a'|}{\widetilde c_0}+\frac{|a'|}{\widetilde c_0^2}|c-c'|\le C_x\|\theta-\theta'\|_2.
\]
\end{proof}

\begin{lemma}[One-dimensional policy diversity]\label[lemma]{lem:one-d-policy-diversity}
Suppose Assumptions~\ref{a_contexts} and \ref{a_theta_1d} hold with $d=1$ and $L > 0$. There exists a constant $\lambda_{\rm div}^{(1)}>0$, depending only on $L,c_g$ and $\widetilde\rho_0$, such that, for every $\theta\in\widetilde\Theta_{\rm 1d}$ and $q\in[0,1]$,
\[
    \mathbb E\left[\widetilde X\widetilde X^\top p_{\theta,q}^{(1)}(X)\right]\succeq\lambda_{\rm div}^{(1)}I_2
\]
and
\[
    \mathbb E\left[\widetilde X\widetilde X^\top\{1-p_{\theta,q}^{(1)}(X)\}\right]\succeq\lambda_{\rm div}^{(1)}I_2.
\]
\end{lemma}

\begin{proof}
Let $v_1=\widetilde\rho_0/2$ and $\varepsilon_1=\min\{1,L\widetilde\rho_0/2\}$. Since $x_\theta\in[\widetilde\rho_0,1-\widetilde\rho_0]$, for every $x\in[1-\widetilde\rho_0/2,1]$, we have $x-x_\theta\ge\widetilde\rho_0/2$. Therefore $p_{\theta,q}^{(1)}(x)\ge\varepsilon_1$ on $[1-\widetilde\rho_0/2,1]$.

Similarly, for every $x\in[0,\widetilde\rho_0/2]$, we have $x_\theta-x\ge\widetilde\rho_0/2$. If $q-L(x_\theta-x)<0$, then $p_{\theta,q}^{(1)}(x)=0$, while otherwise $1-p_{\theta,q}^{(1)}(x)=1-q+L(x_\theta-x)\ge L\widetilde\rho_0/2$. Hence $1-p_{\theta,q}^{(1)}(x)\ge\varepsilon_1$ on $[0,\widetilde\rho_0/2]$.

By \Cref{l_covariate_diverse} with $d=1$, any interval of length at least $v_1$ satisfies
\[
    \mathbb E\left[\widetilde X\widetilde X^\top\mathbf 1\{X\in A\}\right]\succeq c_g\frac{v_1^3}{128}I_2.
\]
Thus the two displayed lower bounds hold with $\lambda_{\rm div}^{(1)}=c_g\varepsilon_1v_1^3/128$.
\end{proof}

\begin{lemma}[One-dimensional estimation at lazy updates]\label[lemma]{lem:one-d-estimation}
Suppose Assumptions~\ref{a_contexts}-\ref{a_beta} and \ref{a_theta_1d} hold with $d=1$. Let $\mathcal S_T$ be the set of lazy update times up to time $T$. Suppose
\[
    T_o\ge C_T\log(T/\eta)\vee\frac{C_{\rm loc}^{(1)}\log(T/\eta)}{\{r_{\rm loc}^{(1)}\}^2}.
\]
Then, with probability at least $1-\eta$, for every $s\in\mathcal S_T$ and $k\in\{1,2\}$,
\[
    \|\widehat\beta_{k,s-1}-\beta_k^*\|_2\le C_{\rm est}^{(1)}\sqrt{\frac{\log(T/\eta)}{s-1}}.
\]
Moreover, for every $s\in\mathcal S_T$,
\[
    \widehat\theta_s\in\widetilde\Theta_{\rm 1d}
    \quad\text{and}\quad
    \|\widehat\theta_s-\theta^*\|_2\le2C_{\rm est}^{(1)}\sqrt{\frac{\log(T/\eta)}{s-1}}.
\]
\end{lemma}

\begin{proof}
The proof follows the same argument as the proof of \Cref{lem:localization-and-estimation}, so we only include the difference~here.

First, during the first $2T_o$ rounds, the algorithm uses a constant policy $p_t(\cdot)=1/2$. By the same allocation, conditional matrix Chernoff and noise concentration arguments as in the proof of \Cref{lem:initial-exploration-estimation}, after enlarging the constants in the lower bound on $T_o$, with probability at least $1-\eta/4$, we have $\|\widehat\theta_{s_0}-\theta^*\|_2\le r_{\rm loc}^{(1)}$ and
\[
    \lambda_{\min}(V_{k,2T_o})\ge \lambda+\frac{\lambda_X^{(1)}T_o}{4},\qquad k\in\{1,2\},
\]
where $s_0$ is the first lazy update time and $\lambda_X^{(1)}>0$ depends only on the covariate distribution.

Next, consider the post-initialisation rounds. Suppose that $\widehat\theta_r\in\widetilde\Theta_{\rm 1d}$ for every lazy update time $r\le t$. Then, for each round $2T_o<i\le t$, the policy used at round $i$ is generated by a parameter in $\widetilde\Theta_{\rm 1d}$. Therefore, by \Cref{lem:one-d-policy-diversity}, we have that
\[
    \mathbb E\left[\widetilde X_i\widetilde X_i^\top\mathbf 1\{A_i=k\}\mid\mathcal F_{i-1}\right]\succeq \lambda_{\rm div}^{(1)}I_2,\qquad k\in\{1,2\}.
\]
Moreover, for $i\le2T_o$, the independent arm assignments give
\[
    \mathbb E\left[\widetilde X_i\widetilde X_i^\top\mathbf 1\{A_i=k\}\mid\mathcal F_{i-1}\right]
    =\frac12\mathbb E\left[\widetilde X\widetilde X^\top\right]
    \succeq\frac{c_g}{256}I_2\succeq\lambda_{\rm div}^{(1)}I_2,\qquad k\in\{1,2\},
\]
where the first inequality follows from \Cref{l_covariate_diverse} with $A=[0,1]$, and the second follows from the value of $\lambda_{\rm div}^{(1)}$ in the proof of \Cref{lem:one-d-policy-diversity}.
Applying the same matrix concentration argument as in the proof of \Cref{lem:localization-and-estimation}, and taking a union bound over $t\le T$ and $k\in\{1,2\}$, we obtain that, with probability at least $1-\eta/2$, for every $t\in\{2T_o,\ldots,T\}$ and $k\in\{1,2\}$, if $\widehat\theta_r\in\widetilde\Theta_{\rm 1d}$ for every lazy update time $r\le t$, 
\[
    \lambda_{\min}(V_{k,t})\ge c_V^{(1)}t.
\]

On the same high-probability event, the martingale noise term is controlled uniformly over time by the same coordinatewise sub-Gaussian argument as in the proof of \Cref{lem:localization-and-estimation}. Thus, with probability at least $1-\eta/4$, for every $t\le T$ and $k\in\{1,2\}$,
\[
    \left\|\sum_{i=1}^t\widetilde X_i\xi_{i,k}\mathbf 1\{A_i=k\}\right\|_2
    \le C_\xi^{(1)}\sqrt{t\log(T/\eta)}.
\]
We work on the intersection of the preceding high-probability events, which by a union bound argument, has probability at least $1-\eta$.

It remains to prove that $\widehat\theta_s\in\widetilde\Theta_{\rm 1d}$ at every lazy update time. We prove this by \textbf{induction} over $s\in\mathcal S_T$. For the first lazy update time $s_0$, we already have $\|\widehat\theta_{s_0}-\theta^*\|_2\le r_{\rm loc}^{(1)}$. Since $ \{\theta:\|\theta-\theta^*\|_2\le2r_{\rm loc}^{(1)}\}\subseteq\widetilde\Theta_{\rm 1d}$, we have $\widehat\theta_{s_0}\in\widetilde\Theta_{\rm 1d}$.

Now fix a lazy update time $s>s_0$ and suppose that $\widehat\theta_r\in\widetilde\Theta_{\rm 1d}$ for every previous lazy update time $r<s$. Then the eigenvalue bound above gives $\lambda_{\min}(V_{k,s-1})\ge c_V^{(1)}(s-1), k\in\{1,2\}$. Writing $S_{k,s-1}=\sum_{i=1}^{s-1}\widetilde X_i\xi_{i,k}\mathbf 1\{A_i=k\}$, the ridge estimator satisfies
\[
    \widehat\beta_{k,s-1}-\beta_k^*=V_{k,s-1}^{-1}S_{k,s-1}-\lambda V_{k,s-1}^{-1}\beta_k^*.
\]
Therefore, we have that
\begin{align*}
\|\widehat\beta_{k,s-1}-\beta_k^*\|_2
&\le \frac{C_\xi^{(1)}\sqrt{(s-1)\log(T/\eta)}}{c_V^{(1)}(s-1)}+\frac{\lambda C_\beta}{c_V^{(1)}(s-1)} \le C_{\rm est}^{(1)}\sqrt{\frac{\log(T/\eta)}{s-1}}.
\end{align*}
Consequently, it follows that
\[
    \|\widehat\theta_s-\theta^*\|_2
    \le \|\widehat\beta_{1,s-1}-\beta_1^*\|_2+\|\widehat\beta_{2,s-1}-\beta_2^*\|_2
    \le 2C_{\rm est}^{(1)}\sqrt{\frac{\log(T/\eta)}{s-1}}.
\]
Since $s-1\ge2T_o$, the lower bound on $T_o$ implies
\[
    2C_{\rm est}^{(1)}\sqrt{\frac{\log(T/\eta)}{s-1}}\le r_{\rm loc}^{(1)}.
\]
Hence, we conclude that $\|\widehat\theta_s-\theta^*\|_2\le r_{\rm loc}^{(1)}$, and so $\widehat\theta_s\in\widetilde\Theta_{\rm 1d}$. This completes the induction and the lemma follows.
\end{proof}

For $n\ge1$, define
\[
\Delta_{n}^{L,1}=\sup_{\theta\in\widetilde\Theta_{\rm 1d},\ q\in[0,1]}\left|\frac1n\sum_{i=1}^n \widetilde X_i^\top\theta^*p_{\theta,q}^{(1)}(X_i)-\mathbb E\{\widetilde X^\top\theta^*p_{\theta,q}^{(1)}(X)\}\right|.
\]

\begin{lemma}[Control of $\Delta_n^{L,1}$]\label[lemma]{lem:one-d-control-delta}
Suppose Assumptions~\ref{a_contexts} and \ref{a_theta_1d} hold with $d=1$. There exists a constant $C<\infty$, depending only on $L,c_g,C_g,A_0,C_0$ and the constants defining $\widetilde\Theta_{\rm 1d}$, such that, for every $n\ge2$ and every $\eta\in(0,1)$, with probability at least $1-\eta$,
\[
    \Delta_n^{L,1}\le C\sqrt{\frac{\log(en/\eta)}{n}}.
\]
\end{lemma}

\begin{proof}
Let
\[
\mathcal F_{L,1}=\left\{x\mapsto \widetilde x^\top\theta^*p_{\theta,q}^{(1)}(x):\theta\in\widetilde\Theta_{\rm 1d},\ q\in[0,1]\right\}.
\]
Since $|\widetilde x^\top\theta^*|\le A_0+C_0=:B_1$ and $0\le p_{\theta,q}^{(1)}\le1$, the class is uniformly bounded by $B_1$. Moreover, $p_{\theta,q}^{(1)}$ depends on $\theta$ only through $x_\theta\in[\widetilde\rho_0,1-\widetilde\rho_0]$. For any $x_0,x_0'\in[\widetilde\rho_0,1-\widetilde\rho_0]$ and $q,q'\in[0,1]$,
\[
\sup_{x\in[0,1]}\left|[q+L(x-x_0)]_{[0,1]}-[q'+L(x-x_0')]_{[0,1]}\right|\le |q-q'|+L|x_0-x_0'|.
\]
Therefore, the class $\mathcal F_{L,1}$ admits brackets with
\[
    \log N_{[]}(\varepsilon,\mathcal F_{L,1},L_2(P))\le C_1\log(C_1/\varepsilon).
\]
Applying the same bracketing maximal inequality and bounded empirical-process concentration as in the proof of \Cref{lem:control-delta-L-n} gives the claim.
\end{proof}

\begin{proposition}[One-dimensional regret decomposition]\label[proposition]{prop:one-d-regret-decomposition}
Suppose Assumptions~\ref{a_contexts}-\ref{a_beta} and \ref{a_theta_1d} hold with $d=1$. Suppose that $T_o$ satisfies the lower bound in \Cref{lem:one-d-estimation}. Then, for every update time $s \in \mathcal{S}_T$, with probability at least $1-\eta$,
\[
    R_s(p^*)-R_s(\widehat p_s)\le C_{\theta,L}^{(1)}\|\widehat\theta_s-\theta^*\|_2+2\Delta_{s-1}^{L,1},
\]
where $p^*=p_{\theta^*,q^*}^{(1)}$, $q^*\in\argmax_{q\in[0,1]}R(p_{\theta^*,q}^{(1)})$, and $C_{\theta,L}^{(1)}<\infty$ depends only on $L,C_\beta$ and the constants defining $\widetilde\Theta_{\rm 1d}$.
\end{proposition}

\begin{proof}
Fix $t>2T_o$ and write $s=\tau(t)$. The following proof is constructed conditioning on the event described in \Cref{lem:one-d-estimation}, which has probability at least $1-\eta$. On this event, $\widehat\theta_s\in\widetilde\Theta_{\rm 1d}$. Writing $\widehat\theta_s=(\widehat a_s,\widehat c_s)$, this implies $\widehat c_s\ge\widetilde c_0>0$. Hence the fallback step in \Cref{alg:lazy-one-dimensional-fair-greedy} is not used, and the algorithm uses $\widehat p_s=p_{\widehat\theta_s,\widehat q_s}^{(1)}$. Define $p_s^\dagger=p_{\widehat\theta_s,q^*}^{(1)}$.

By \Cref{lem:one-d-threshold-perturbation} and the nonexpansiveness of the clipping map,
\[
    \sup_{x\in[0,1]}|p_{\theta^*,q^*}^{(1)}(x)-p_{\widehat\theta_s,q^*}^{(1)}(x)|\leq L|x_{\theta^*} - x_{\widehat\theta_s}| \le LC_x\|\widehat\theta_s-\theta^*\|_2.
\]
Since $|\widetilde X^\top\theta^*|\le\sqrt2(\|\beta_1^*\|_2+\|\beta_2^*\|_2)\le2\sqrt2 C_\beta$, it follows that
\[
    \mathbb E\left[\widetilde X^\top\theta^*\{p^*(X)-p_s^\dagger(X)\}\mid\mathcal F_{s-1}\right]\le C_1\|\widehat\theta_s-\theta^*\|_2.
\]
For the remaining term, adding and subtracting empirical averages gives the same decomposition as in the proof of \Cref{prop:lazy-regret-decomposition}. The two empirical-process terms are bounded by $\Delta_{s-1}^{L,1}$. Moreover, by empirical optimality of $\widehat q_s$ over $q\in[0,1]$,
\[
\frac1{s-1}\sum_{i=1}^{s-1}\widetilde X_i^\top\widehat\theta_s\{p_s^\dagger(X_i)-\widehat p_s(X_i)\}\le0.
\]
The remaining estimation term is bounded by
\[
\frac1{s-1}\sum_{i=1}^{s-1}|\widetilde X_i^\top(\theta^*-\widehat\theta_s)|\le\sqrt2\|\widehat\theta_s-\theta^*\|_2.
\]
Combining the preceding bounds gives
\[
    R_s(p^*)-R_s(\widehat p_s)\le C_1\|\widehat\theta_s-\theta^*\|_2+2\Delta_{s-1}^{L,1}+\sqrt2\|\widehat\theta_s-\theta^*\|_2.
\]
Therefore,
\[
    R_s(p^*)-R_s(\widehat p_s)\le C_{\theta,L}^{(1)}\|\widehat\theta_s-\theta^*\|_2+2\Delta_{s-1}^{L,1},
\]
where $C_{\theta,L}^{(1)}=C_1+\sqrt2$. This completes the proof.
\end{proof}

\begin{proof}[Proof of \Cref{thm:one-d-regret}]
The fairness guarantee follows because each policy is either a clipped function with slope $L$ or the constant policy $1/2$.

If $2T_o=T$, the claimed bound follows from the bounded one-round regret. Hence suppose $2T_o<T$. For the stated choice $T_o=\lceil K\sqrt{T\log(T/\eta)}\rceil$ with fixed $K>0$, write $D=C_T\vee C_{\rm loc}^{(1)}/\{r_{\rm loc}^{(1)}\}^2$. If
\[
    T\ge C_0\log(T/\eta),\qquad C_0\ge16\max\{1,K^2,D^2/K^2\},
\]
then $2T_o\le T$ and $T_o\ge2D\log(T/\eta)\ge D\log(2T/\eta)$. Thus the required lower bound on $T_o$ also holds with confidence parameter $\eta/2$.

Let $\mathcal S_T$ denote the set of lazy update times up to time $T$. By \Cref{lem:one-d-estimation}, after enlarging the constants in the assumed lower bound on $T_o$, with probability at least $1-\eta/2$, for every $s\in\mathcal S_T$,
\[
    \|\widehat\theta_s-\theta^*\|_2\le C\sqrt{\frac{\log(T/\eta)}{s-1}}.
\]
By \Cref{lem:one-d-control-delta}, together with a union bound over $n=2,\ldots,T$, with probability at least $1-\eta/2$, for every $n=2,\ldots,T$,
\[
    \Delta_n^{L,1}\le C\sqrt{\frac{\log(eT^2/\eta)}{n}}.
\]
We work on the intersection of these two events. By \Cref{prop:one-d-regret-decomposition}, for every post-initialization round $t>2T_o$ with $s=\tau(t)$,
\[
    R_s(p^*)-R_s(\widehat p_s)\le C\sqrt{\frac{\log(T/\eta)}{s-1}}+C\sqrt{\frac{\log(eT^2/\eta)}{s-1}}.
\]
The regret from the first $2T_o$ forced-exploration rounds is bounded by $C_1T_o$. Since the lazy update times grow geometrically and each policy computed at time $s$ is used for at most $s$ rounds,
\[
\mathfrak R_T\le C_1T_o+C_2\sqrt{\log(eT^2/\eta)}\sum_{s\in\mathcal S_T}s^{1/2}\le C_1T_o+C_3\sqrt{T\log(eT^2/\eta)}.
\]
This proves the theorem.
\end{proof}

\subsection{UCB algorithm for general covariates}\label[appendix]{sec:one-d-ucb}
\subsubsection{Algorithm}
\Cref{alg:lazy-one-dimensional-fair-ucb} is the one-dimensional version of \Cref{alg:lazy-global-model-ucb}. Its optimistic update uses the $16$ vertex pairs of the two confidence sets through \Cref{alg:vertex-oracle}.
\begin{algorithm}[!htbp]
\caption{Lazy one-dimensional polyhedral UCB algorithm ($d=1$)}
\label{alg:lazy-one-dimensional-fair-ucb}
\begin{algorithmic}
\Require{Time horizon $T$, fairness parameter $(L,\nu)$, ridge parameter $\lambda$, exploration period $T_o$, confidence radius $\alpha_T$ from \eqref{eq_ucb_alpha} with $d=1$}
\State Set $s_0\gets2T_o+1$.
\For{$k\in\{1,2\}$}
\State $V_k\gets\lambda I_2$, $B_k\gets0$.
\EndFor
\For{$t=1,2,\ldots,2T_o$}
\State Observe $X_t$ and set $\widetilde X_t=(1,X_t)^\top$.
\State Play arm $1$ with probability $1/2$ and arm $2$ otherwise; observe $Y_t$.
\State Update $V_{A_t}\gets V_{A_t}+\widetilde X_t\widetilde X_t^\top$ and $B_{A_t}\gets B_{A_t}+\widetilde X_tY_t$.
\EndFor
\State Set $s\gets s_0$ and, for $k\in\{1,2\}$, set $V_{s,k}\gets V_k$ and define $\ell_{s,k}:=\log\det(V_{s,k})$.
\For{$t=s_0,\ldots,T$}
\If{$t=s_0$ \textbf{or} $t\ge2s$ \textbf{or} $\exists k\in\{1,2\}:\log\det(V_k)-\ell_{s,k}\ge\log2$}
\State Set $s\gets t$.
\State For $k\in\{1,2\}$, set $V_{s,k}\gets V_k$, $\widehat\beta_{s,k}\gets V_k^{-1}B_k$, $\ell_{s,k}\gets\log\det(V_{s,k})$.
\State Apply the vertex enumeration and LP scoring steps of \Cref{alg:vertex-oracle} with $d=1$ to select $(\bar\beta_{s,1},\bar\beta_{s,2},z_s)$; set $\theta_s\gets \bar\beta_{s,1}-\bar\beta_{s,2}=(a_s,c_s)$.
\If{$c_s>0$ and $-a_s\in[0,c_s]$}
\State Set $x_{\theta_s}\gets-a_s/c_s$ and $q_s\gets\min\{1,\min_{1\le i\le s-1}[z_{s,i}+L|x_{\theta_s}-X_i|]\}$.
\State Set $p_s(\cdot)\gets p_{\theta_s,q_s}^{(1)}(\cdot)$.
\Else
\State Set $p_s(\cdot)\gets1/2$.
\EndIf
\State Set $\widehat p\gets p_s$.
\EndIf
\State Observe $X_t$ and set $\widetilde X_t=(1,X_t)^\top$.
\State Play arm $1$ with probability $\widehat p(X_t)$ and arm $2$ otherwise; observe $Y_t$.
\State Update $V_{A_t}\gets V_{A_t}+\widetilde X_t\widetilde X_t^\top$ and $B_{A_t}\gets B_{A_t}+\widetilde X_tY_t$.
\EndFor
\end{algorithmic}
\end{algorithm}

\subsubsection{Theoretical guarantee}
\begin{theorem}[One-dimensional regret upper bound]\label[theorem]{thm:one-d-ucb-regret}
Suppose Assumptions~\ref{a_contexts_non_cts}, \ref{a_noise}, \ref{a_beta} and \ref{a_theta_1d} hold with $d=1$. Implement \Cref{alg:lazy-one-dimensional-fair-ucb}, with $\lambda>0$ fixed, $T_o=\lceil\sqrt T\rceil$, and $\alpha_T$ in \eqref{eq_ucb_alpha} with $d=1$. Then the policy played at every round is $(L,\nu)$-fair. Moreover, for any $\eta\in(0,1/2)$, with probability at least $1-\eta$, we have that
\begin{equation}\label{eq:one-d-ucb-regret-bound}
\mathfrak R_T\le C_1(1+\sqrt2)\alpha_T\sqrt{T\left\{2\log\left(1+\frac T\lambda\right)+\log(8/\eta)\right\}}+C_2\sqrt{T\log(eT^2/\eta)},
\end{equation}
whenever $T\ge16\vee C_3(1+\sqrt2)^4\alpha_T^4/(\kappa_X^2\{r_{\rm loc}^{(1)}\}^4)$. The constants depend only on $L,\sigma,C_\beta,\lambda$ and the constants in \Cref{a_theta_1d}. In particular, $\mathfrak R_T=\widetilde O_{\mathrm{p}}(\sqrt T)$.
\end{theorem}

\subsubsection{Proof of Theorem \ref{thm:one-d-ucb-regret}}
We use the enlarged set $\widetilde\Theta_{\rm 1d}$ and policies $p_{\theta,q}^{(1)}$ defined in \Cref{sec_app_one_d_setup}, and take $r_{\rm loc}^{(1)}\le1$. For $d=1$, the ideal optimistic update is
\begin{equation}\label{eq:one-d-ideal-update}
    (\bar\beta_{s,1}^{\circ},\bar\beta_{s,2}^{\circ},q_s^{\circ})\in\argmax_{\substack{b_k\in\mathcal C_{s,k},\ k\in\{1,2\}\\q\in[0,1]}}\frac1{s-1}\sum_{i=1}^{s-1}r_{b,p_{b_1-b_2,q}^{(1)}}(X_i).
\end{equation}

\begin{lemma}[One-dimensional vertex optimisation]\label[lemma]{lem:one-d-vertex-optimisation}
Suppose every feasible difference belongs to $\widetilde\Theta_{\rm 1d}$. At each update time $s$ of \Cref{alg:lazy-one-dimensional-fair-ucb}, the selected model pair $\bar\beta_s=(\bar\beta_{s,1},\bar\beta_{s,2})$ and the policy $p_s$ satisfy
\begin{equation}\label{eq:one-d-vertex-opt}
\frac1{s-1}\sum_{i=1}^{s-1}r_{\bar\beta_s,p_s}(X_i)=\sup_{\substack{b_k\in\mathcal C_{s,k},\ k\in\{1,2\}\\q\in[0,1]}}\frac1{s-1}\sum_{i=1}^{s-1}r_{b,p_{b_1-b_2,q}^{(1)}}(X_i).
\end{equation}
\end{lemma}

\begin{proof}
By \Cref{lem:one-d-canonical}, the supremum in \eqref{eq:one-d-vertex-opt} equals $\max_{b_k\in\mathcal C_{s,k}}F_s(b_1,b_2)$. Similar as Lemma~\ref{lem:vertex-optimisation}, joint convexity of $F_s$ and the vertex representation in \Cref{lem:polyhedral-confidence} show that this maximum is attained at one of the $16$ vertex pairs. By the optimality of $z_v$, each score satisfies $S_v=F_s(v_1,v_2)$. Selecting the pair with the largest score and applying the canonical domination in \Cref{lem:one-d-canonical} gives the lower bound in \eqref{eq:one-d-vertex-opt}. The reverse bound follows because the returned model-policy pair is feasible in \eqref{eq:one-d-ideal-update}.
\end{proof}


Define the functional class
\begin{equation*}
\mathcal F_{\mathrm{UCB},1}=\left\{x\mapsto r_{b,p_{\theta,q}^{(1)}}(x):\|b_k\|_2\le C_\beta+1,\ k\in\{1,2\},\ \theta=b_1-b_2\in\widetilde\Theta_{\rm 1d},\ q\in[0,1]\right\}
\end{equation*}
and
\begin{equation}
\Delta_n^{\mathrm{UCB},1}=\sup_{f\in\mathcal F_{\mathrm{UCB},1}}\left|\frac1n\sum_{i=1}^nf(X_i)-\mathbb E\{f(X)\}\right|.
\end{equation}

\begin{lemma}[Control of $\Delta_n^{\mathrm{UCB},1}$]\label[lemma]{lem:one-d-ucb-control-delta}
Suppose Assumptions~\ref{a_beta} and \ref{a_theta_1d} hold. Let $X_1,\ldots,X_n$ be i.i.d. from an arbitrary Borel probability measure $P_X$ supported on $[0,1]$. There exists a constant $C<\infty$, depending only on $L,C_\beta$ and the constants defining $\widetilde\Theta_{\rm 1d}$, and independent of $P_X$, such that, for every $n\ge2$ and $\eta\in(0,1)$, with probability at least $1-\eta$,
\begin{equation*}
\Delta_n^{\mathrm{UCB},1}\le C\sqrt{\frac{\log(en/\eta)}{n}}.
\end{equation*}
\end{lemma}

\begin{proof}
The proof follows the same argument as in \Cref{lem:ucb-bracketing,lem:rademacher-policy,lem:ucb-control-delta}. We verify the entropy bound for the one-dimensional policy class and evaluate the resulting entropy integral.

For $x_0,x_0'\in[\widetilde\rho_0,1-\widetilde\rho_0]$ and $q,q'\in[0,1]$, the nonexpansiveness of the clipping map gives
\begin{equation*}
\sup_{x\in[0,1]}\left|[q+L(x-x_0)]_{[0,1]}-[q'+L(x-x_0')]_{[0,1]}\right|\le|q-q'|+L|x_0-x_0'|.
\end{equation*}
For $L>0$, discretising $q$ and the boundary location $x_0$ at resolutions of order $\varepsilon$ and $\varepsilon/L$, respectively, gives a sup-norm cover of the one-dimensional policy class with log covering number of order $\log(1/\varepsilon)$. If $L=0$, the policy is constant and discretising $q$ alone gives the same bound. Since we have that $r_{b,p_{\theta,q}^{(1)}}(x)=\widetilde x^\top b_2+(\widetilde x^\top\theta)p_{\theta,q}^{(1)}(x)$, and $b_2$ and $\theta$ range over bounded subsets of $\mathbb R^2$, discretising these finite-dimensional parameters as in \Cref{lem:ucb-bracketing} preserves the logarithmic entropy order. As in that proof, the resulting uniform approximations give finite brackets of pointwise width $\varepsilon$. Therefore, uniformly over all Borel probability measures $P$ on $[0,1]$, we have that
\begin{equation*}
\log N_{[]}(\varepsilon,\mathcal F_{\mathrm{UCB},1},L_2(P))\le C_2\log(C_2/\varepsilon),\qquad\varepsilon\in(0,1).
\end{equation*}
For any deterministic sample $x_1,\ldots,x_n\in[0,1]$, let $P_n=n^{-1}\sum_{i=1}^n\delta_{x_i}$. Selecting one class member from each nonempty bracket gives the same logarithmic bound for the $L_2(P_n)$ covering number. Moreover, every $f\in\mathcal F_{\mathrm{UCB},1}$ satisfies $\|f\|_\infty\le B_1:=\sqrt2(C_\beta+1)$. Applying the conditional sub-Gaussian increment calculation and truncated Dudley argument in the proof of \Cref{lem:rademacher-policy}, with this logarithmic entropy bound, gives
\[
\widehat{\mathfrak R}_n(\mathcal F_{\mathrm{UCB},1})
\le C\inf_{0<\alpha<1}\left\{\alpha+\frac1{\sqrt n}
\left[1+\int_\alpha^1\{\log(C_2/\varepsilon)\}^{1/2}\,\mathrm d\varepsilon\right]\right\}.
\]
The additional constant inside the brackets accounts for entropy scales above one. Since the integral remains finite as $\alpha\rightarrow 0$, we obtain $\widehat{\mathfrak R}_n(\mathcal F_{\mathrm{UCB},1})\le C_3n^{-1/2}$ for every sample. The finite brackets of pointwise width $\varepsilon$ imply separability in the sup-norm. Thus the empirical and Rademacher suprema are measurable, and the symmetrisation and bounded-differences argument in \Cref{lem:ucb-control-delta} applies and gives, with probability at least $1-\eta$,
\[
\Delta_n^{\mathrm{UCB},1}
\le2C_3n^{-1/2}+B_1\sqrt{\frac{2\log(1/\eta)}{n}}
\le C\sqrt{\frac{\log(en/\eta)}{n}}.
\]
\end{proof}

\begin{proof}[Proof of \Cref{thm:one-d-ucb-regret}]
The fairness guarantee follows because each policy is either a clipped function with slope $L$ or the constant policy $1/2$.

Let $B=\sqrt2(C_\beta+1)$ and choose $p^*=p_{\theta^*,q^*}^{(1)}$ as in \Cref{lem:one-d-canonical}. Similar as in the proof of Theorem~\ref{thm:main-regret}, the first $2T_o$ rounds contribute at most $4BT_o\le C\sqrt T$ to the regret.

We use the events $\mathcal E_{\rm conf}$, $\mathcal E_{\rm init}$, and $\mathcal E_{\rm mart}$ in \eqref{eq:main-regret-confidence-event}, \eqref{eq:initial-design-event}, and \eqref{eq:main-regret-martingale-event}, with $d=1$ and $\tau(t)$ the most recent update time of \Cref{alg:lazy-one-dimensional-fair-ucb}. The confidence and initial exploration arguments are similar as \Cref{cor:self-normalised-deterministic-lnu} and \Cref{lem:initial-localisation}. In particular, on $\mathcal E_{\rm conf}\cap\mathcal E_{\rm init}$, for every update time $s$ and every $b_k\in\mathcal C_{s,k}$, we have that
\begin{equation*}
\|(b_1-b_2)-\theta^*\|_2\le\frac{2(1+\sqrt2)\alpha_T}{\sqrt{\lambda+\kappa_XT_o/4}}\le\frac{4(1+\sqrt2)\alpha_T}{\sqrt{\kappa_X}T^{1/4}}\le2r_{\rm loc}^{(1)},
\end{equation*}
where the last inequality follows from the assumed lower bound on $T$, after choosing $C_3$ sufficiently large. The same argument gives $\|b_k-\beta_k^*\|_2\le r_{\rm loc}^{(1)}\le1$ and hence $\|b_k\|_2\le C_\beta+1$. Thus every feasible difference belongs to $\widetilde\Theta_{\rm 1d}$, so the boundary trace is well defined.

By \Cref{lem:one-d-ucb-control-delta} and a union bound over $n=2,\ldots,T$, the event
\begin{equation*}
\mathcal E_{\rm emp}^{(1)}=\left\{\Delta_n^{\mathrm{UCB},1}\le C\sqrt{\frac{\log(8enT/\eta)}{n}},\ \text{for every }n\in\{2,\ldots,T\}\right\}
\end{equation*}
has probability at least $1-\eta/8$. Similar as in the proof of \Cref{thm:main-regret}, the same lower bound on $T$ gives $P(\mathcal E_{\rm init}^c)\le\eta/4$, while $P(\mathcal E_{\rm conf}^c)\le\eta/8$ and $P(\mathcal E_{\rm mart}^c)\le\eta/8$. Hence, a union bound argument gives
\begin{equation*}
\mathbb{P}(\mathcal E_{\rm conf}\cap\mathcal E_{\rm init}\cap\mathcal E_{\rm mart}\cap\mathcal E_{\rm emp}^{(1)})\ge1-5\eta/8\ge1-\eta.
\end{equation*}
We work on this intersection for the rest of the proof.

Since the true model-policy pair is feasible in \eqref{eq:one-d-ideal-update}, \Cref{lem:one-d-vertex-optimisation} gives empirical optimism similar as \Cref{lem:empirical-optimism}. Note that both $r_{\beta^*,p^*}$ and $r_{\bar\beta_s,p_s}$ belong to $\mathcal F_{\mathrm{UCB},1}$, therefore, similar as \Cref{lem:one-step-regret}, for $t\ge s_0$ and $s=\tau(t)$, we have that
\begin{align*}
R(p^*)-R(p_s)\le{}&2\Delta_{s-1}^{\mathrm{UCB},1} +\{(1+\sqrt2)\alpha_T+2B\}\mathbb E\left[\min\{1,\|\widetilde X_t\|_{V_{s,A_t}^{-1}}\}\mid\mathcal F_{t-1}\right].
\end{align*}
Similar as \Cref{lem:error-sums}, we have that
\begin{equation*}
2\sum_{t=s_0}^T\Delta_{\tau(t)-1}^{\mathrm{UCB},1}\le C_1\sqrt{T\log(eT^2/\eta)}.
\end{equation*}
Finally, the determinant comparison and width summation are similar as Lemmas~\ref{lem:lazy-determinant-comparison-lnu} and \ref{lem:lazy-width-sum-lnu}. Combining them with $\mathcal E_{\rm mart}$ gives
\begin{align*}
\sum_{t=s_0}^T\mathbb E\left[\min\{1,\|\widetilde X_t\|_{V_{\tau(t),A_t}^{-1}}\}\mid\mathcal F_{t-1}\right]
\le C\sqrt{2T\log\left(1+\frac T\lambda\right)}+\sqrt{2T\log(8/\eta)}.
\end{align*}
Since $\alpha_T\ge C_\beta\sqrt\lambda$ and $\lambda$ is fixed, $(1+\sqrt2)\alpha_T+2B\le C(1+\sqrt2)\alpha_T$. Summing the one-step bounds and adding the exploration regret proves \eqref{eq:one-d-ucb-regret-bound}. The stated rate follows from $\alpha_T=O(\sqrt{\log(T/\eta)})$.
\end{proof}

\section{Structural results under arbitrary context laws}\label[appendix]{sec_app_structural_general}
Throughout this section, $P_X$ is an arbitrary Borel probability law on $\mathcal X=[0,1]^d$ and $\mathcal S=\supp(P_X)$.
\begin{lemma}[Equivalence between fairness over $\mathcal{S}$ and over $\cX$]\label[lemma]{lem:support-extension}
Let $\mathcal S=\supp(P_X)$. Suppose $p_\mathcal S:\mathcal S\to[0,1]$ satisfies $|p_\mathcal S(x)-p_\mathcal S(x')|\le L\|x-x'\|_\nu$ for all $x,x'\in\mathcal S$, then the extended policy
\begin{equation*}
\overline p(x)=\left[\inf_{z\in\mathcal S}\{p_\mathcal S(z)+L\|x-z\|_\nu\}\right]_{[0,1]},\qquad x\in[0,1]^d,
\end{equation*}
belongs to $\mathcal P_{L,\nu}$ and satisfies $\overline p=p_\mathcal S$ on $\mathcal S$. Consequently, maximising $R(p)$ over $\mathcal P_{L,\nu}$ is equivalent to maximising over $L$-Lipschitz maps defined only on $\mathcal S$.
\end{lemma}

\begin{proof}
The standard McShane argument \citep{mcshane1934extension} shows that $x\mapsto\inf_{z\in\mathcal S}\{p_\mathcal S(z)+L\|x-z\|_\nu\}$ is $L$-Lipschitz. Clipping to $[0,1]$ is nonexpansive, so it holds that $\overline p\in\mathcal P_{L,\nu}$.

Fix $x\in\mathcal S$. Taking $z=x$, by construction, we have that $\overline p(x)\le p_\mathcal S(x)$. Conversely, the Lipschitz property of $p_\mathcal S$ gives $p_\mathcal S(z)+L\|x-z\|_\nu\ge p_\mathcal S(x)$ for every $z\in\mathcal S$, so it holds that the infimum is at least $p_\mathcal S(x)$. Hence, for every $x\in\mathcal S$, we have that $\inf_{z\in\mathcal S}\{p_\mathcal S(z)+L\|x-z\|_\nu\}=p_\mathcal S(x)$. Since $p_\mathcal S(x)\in[0,1]$, clipping does not change this value, and therefore $\overline p(x)=p_\mathcal S(x)$ for every $x\in\mathcal S$.

Note that $X\in\mathcal S$ almost surely and the population reward depends only on the restriction of the policy to $\mathcal S$, therefore, the final claim holds.
\end{proof}

\begin{theorem}[Boundary randomisation]\label[theorem]{thm:support-mixing}
Suppose $P_X$ is an arbitrary Borel probability law on $[0,1]^d$, $L > 0$, and let $p^*\in\argmax_{p\in\mathcal P_{L,\nu}}R(p)$. Define the effective randomisation region by $\mathcal M_\mathcal S(p^*)=\{x\in\mathcal S:0<p^*(x)<1\}$. If $\mathcal H_{\theta^*}=\emptyset$, then $\Delta^*(\cdot)$ has constant sign on $[0,1]^d$ and a deterministic constant policy is optimal. If $\mathcal H_{\theta^*}\neq\emptyset$, then it holds that
\begin{equation}\label{eq:support-mixing}
\mathcal M_\mathcal S(p^*)\subseteq\left\{x\in\mathcal S:\dist_\nu(x,\mathcal H_{\theta^*})\le\frac1L\right\}.
\end{equation}
More precisely, for every $x\in\mathcal S$ with $\dist_\nu(x,\mathcal H_{\theta^*})>1/L$, one has $p^*(x)=1$ if $\Delta^*(x)>0$ and $p^*(x)=0$ if $\Delta^*(x)<0$.
\end{theorem}

\begin{proof}
For every policy $p$, $R(p)=\mathbb E[\widetilde X^\top\beta_2^*]+\mathbb E[\Delta^*(X)p(X)]$, so it suffices to consider the second term.

\noindent\textbf{Step 1. The case $\mathcal H_{\theta^*}=\emptyset$.} By continuity of the affine function $\Delta^*$ on the connected set $[0,1]^d$, its sign is constant. The objective is then pointwise increasing in $p$ when $\Delta^*>0$ and pointwise decreasing in $p$ when $\Delta^*<0$, so $p = 1$ or $p=0$, respectively, is optimal.

\noindent\textbf{Step 2. Saturation on the positive side.} Suppose $\mathcal H_{\theta^*}\neq\emptyset$. Fix $x_0\in\mathcal S$ such that $\dist_\nu(x_0,\mathcal H_{\theta^*})>1/L$ and $\Delta^*(x_0)>0$. Suppose for contradiction that $p^*(x_0)<1$. By continuity of $\Delta^*$ and $p^*$, and since $\dist_\nu(x_0,\mathcal H_{\theta^*})>1/L$, there exist $r>0$ and $\eta>0$ such that $A=\overline B_\nu(x_0,r)\cap[0,1]^d$ satisfies $\Delta^*(x)>0$ and $p^*(x)\le1-\eta$ for all $x\in A$, and $\dist_\nu(A,\mathcal H_{\theta^*})>1/L$. Define
\[
\phi_A(x)=\max\{0,1-L\dist_\nu(x,A)\}.
\]
Then $\phi_A$ is $L$-Lipschitz, takes values in $[0,1]$, and equals one on $A$. We claim that $\Delta^*(x)>0$ whenever $\phi_A(x)>0$. Indeed, if $\phi_A(x)>0$, then $\dist_\nu(x,A)<1/L$. If $\Delta^*(x)\le0$, choose $a\in A$ with $\|x-a\|_\nu<1/L$. By convexity of $\mathcal X$ and continuity of $\Delta^*$, the line segment between $a$ and $x$ intersects $\mathcal H_{\theta^*}$ at some point $h$. This gives $\dist_\nu(A,\mathcal H_{\theta^*})\le\|a-h\|_\nu\le\|a-x\|_\nu<1/L$, which is a contradiction.

Define $\widetilde p(x)=\max\{p^*(x),\phi_A(x)\}$. Since both $p^*$ and $\phi_A$ are $L$-Lipschitz and the pointwise maximum preserves this Lipschitz bound, we have $\widetilde p\in\mathcal P_{L,\nu}$. Since $\phi_A=1$ on $A$, we have $\widetilde p(x)-p^*(x)\ge\eta$ for all $x\in A$. Moreover, $\widetilde p-p^*$ can be positive only where $\phi_A>0$, so the claim above gives
\[
\Delta^*(x)\{\widetilde p(x)-p^*(x)\}\ge0,\qquad x\in\mathcal X.
\]
Since $x_0\in\mathcal S=\supp(P_X)$ and $A$ contains a relative open neighbourhood of $x_0$ in $\mathcal X$, we have $P_X(A)>0$. Therefore,
\begin{align*}
R(\widetilde p)-R(p^*)
&=\mathbb E[\Delta^*(X)\{\widetilde p(X)-p^*(X)\}]\\
&\ge\mathbb E[\Delta^*(X)\{\widetilde p(X)-p^*(X)\}\indc\{X\in A\}]\\
&\ge\eta\,\mathbb E[\Delta^*(X)\indc\{X\in A\}]>0,
\end{align*}
which contradicts the optimality of $p^*$. Hence $p^*(x_0)=1$.

\noindent\textbf{Step 3. Saturation on the negative side.} The case $\Delta^*(x_0)<0$ is analogous. Fix $x_0\in\mathcal S$ such that $\dist_\nu(x_0,\mathcal H_{\theta^*})>1/L$ and $\Delta^*(x_0)<0$. Suppose for contradiction that $p^*(x_0)>0$. Choose $r>0$ and $\eta>0$ such that $A=\overline B_\nu(x_0,r)\cap\mathcal X$ satisfies $\Delta^*(x)<0$ and $p^*(x)\ge\eta$ for all $x\in A$, and $\dist_\nu(A,\mathcal H_{\theta^*})>1/L$. Define $\phi_A(x)=\max\{0,1-L\dist_\nu(x,A)\}$. By the same line-segment argument as above, $\Delta^*(x)<0$ whenever $\phi_A(x)>0$. Set
\[
\widetilde p(x)=\min\{p^*(x),1-\phi_A(x)\}.
\]
The pointwise minimum of two $L$-Lipschitz functions is $L$-Lipschitz, so $\widetilde p\in\mathcal P_{L,\nu}$. Since $\phi_A=1$ on $A$, we have $\widetilde p(x)=0$ and $p^*(x)-\widetilde p(x)\ge\eta$ for all $x\in A$. Moreover, $\widetilde p-p^*$ can be negative only where $\phi_A>0$, so
\[
\Delta^*(x)\{\widetilde p(x)-p^*(x)\}\ge0,\qquad x\in\mathcal X.
\]
Again, since $x_0\in\mathcal S$ and $A$ contains a relative open neighbourhood of $x_0$, we have $P_X(A)>0$. Therefore,
\begin{align*}
R(\widetilde p)-R(p^*)
&=\mathbb E[\Delta^*(X)\{\widetilde p(X)-p^*(X)\}]\\
&\ge\mathbb E[\Delta^*(X)\{\widetilde p(X)-p^*(X)\}\indc\{X\in A\}]\\
&\ge\eta\,\mathbb E[-\Delta^*(X)\indc\{X\in A\}]>0,
\end{align*}
which contradicts the optimality of $p^*$. Hence $p^*(x_0)=0$.

Thus, whenever $x\in\mathcal S$ satisfies $\dist_\nu(x,\mathcal H_{\theta^*})>1/L$, the optimal policy is deterministic and selects the arm with the larger conditional mean reward. This proves \eqref{eq:support-mixing}.
\end{proof}

\section{Computation of the boundary-canonical policy}\label{sec:app-boundary-computation}

Throughout this section, we take $\nu=1$ in the fairness constraint.
At update time~$s$ in \Cref{alg:greedy-boundary-optimisation}, let $\widehat z_s$ be a solution of the linear program in \eqref{eq_alg_lp}.  We describe how to compute the resulting policy~$\widehat p_s(x)$ at any $x\in[0,1]^d$.  Write $n=s-1$, $\widehat\theta_s=(\widehat a_s,\widehat\gamma_s^\top)^\top$ and $\mathcal H=\mathcal H_{\widehat\theta_s}$. We assume $L>0$ and $\mathcal H\neq\varnothing$.

For each boundary point $h\in\mathcal H$, the McShane extension gives the probability of selecting arm~1 
\begin{equation}\label{eq:app-computed-trace}
\widehat q_s(h) =\widehat p_s^{\rm McExt}(h) =\min\Big\{1,\min_{i\in[n]} \bigl[\widehat z_{s,i}+L\|h-X_i\|_1\bigr]\Big\}.
\end{equation}
These boundary probabilities determine the canonical policy $\widehat p_s=p_{\widehat\theta_s,\widehat q_s}$ through \eqref{eq_canonical}. To compute $\widehat p_s$, we evaluate the upper envelope $U_{\widehat\theta_s,\widehat q_s}$ on $\mathcal{X}_{\widehat{\theta}_s}^+$, and the lower envelope $M_{\widehat\theta_s,\widehat q_s}$ on $\mathcal{X}_{\widehat{\theta}_s}^-$, both defined in \eqref{eq:envelope-definitions}.

\medskip
\noindent
\textbf{Evaluation on $\mathcal{X}_{\widehat{\theta}_s}^+$.}
For $x \in \mathcal{X}_{\widehat{\theta}_s}^+$, we compute
\begin{equation}\label{eq:app-positive-evaluation-1}
\widehat p_s(x)=\min\{1,U_{\widehat\theta_s,\widehat q_s}(x)\}, \quad \mbox{with }U_{\widehat\theta_s,\widehat q_s}(x) =\min_{h\in\mathcal H} \{\widehat q_s(h)+L\|h-x\|_1\}.
\end{equation}
Substituting \eqref {eq:app-computed-trace} into \eqref{eq:app-positive-evaluation-1} gives that 
\begin{equation}\label{eq:app-positive-evaluation}
\widehat p_s(x) =\min\Big\{1,\min_{i\in[n]} \big\{\widehat z_{s,i}+L \mathrm{d}_{\mathcal H}(X_i,x)\big\}\Big\}.
\end{equation}
Define the minimum total distance from $y$ to $x$ via a point on the boundary by
\[
\mathrm{d}_{\mathcal H}(y,x) =\min_{h\in\mathcal H} \{\|y-h\|_1+\|h-x\|_1\}.
\]
To compute $\mathrm{d}_{\mathcal H}(y,x)$, for $j \in [d]$, let
$\ell_j=\min\{x_j,y_j\}$, $u_j=\max\{x_j,y_j\}$ and
$B=\prod_{j=1}^d[\ell_j,u_j]$. Then
\begin{equation}\label{eq:app-two-hop-box}
\mathrm{d}_{\mathcal H}(y,x) =\|y-x\|_1+2\delta_B,
\quad \mbox{with } \delta_B=\min_{h\in\mathcal H,v\in B}\|h-v\|_1.
\end{equation}
It therefore remains to compute $\delta_B$, the distance between $B$ and $\mathcal H$.

If $B \cap \mathcal H \neq \emptyset$, then $\delta_B=0$.
Otherwise, choose a vertex of $B$ whose estimated reward gap is closest to zero. Let $r$ denote the magnitude of this gap, initially given by
\[
r=
\begin{cases}
\displaystyle
\widehat a_s+\sum_{j=1}^d
\min\{\widehat\gamma_{s,j}\ell_j,
      \widehat\gamma_{s,j}u_j\},
& B \subset \mathcal{X}_{\widehat{\theta}_s}^+,\\[4pt]
\displaystyle
-\widehat a_s-\sum_{j=1}^d
\max\{\widehat\gamma_{s,j}\ell_j,
      \widehat\gamma_{s,j}u_j\},
& B \subset \mathcal{X}_{\widehat{\theta}_s}^-.
\end{cases}
\]
For $j \in [d]$, let $c_j$ denote the maximum available movement of $B$ in coordinate $j$ toward the boundary while remaining in $[0,1]^d$. Then we have the following
\[
\begin{array}{c|cc}
&c_j\ (\widehat\gamma_{s,j}>0)
&c_j\ (\widehat\gamma_{s,j}<0)\\
\hline
B\text{ on the positive side}&\ell_j&1-u_j\\
B\text{ on the negative side}&1-u_j&\ell_j
\end{array}
\]
with $c_j=0$ when $\widehat\gamma_{s,j}=0$.
With $r$ fixed at its initial value, the required distance is
\begin{equation}\label{eq:app-fractional-knapsack}
\delta_B =\min\bigg\{ \sum_{j=1}^d t_j: \sum_{j=1}^d |\widehat\gamma_{s,j}|t_j=r, \, 0\leq t_j\leq c_j \bigg\}.
\end{equation}

To solve \eqref{eq:app-fractional-knapsack}, initialise $t_j=0$ for all $j$ and process coordinates with nonzero coefficients in nonincreasing order of $|\widehat\gamma_{s,j}|$. For each coordinate, let
\[
t_j=\min\bigg\{c_j,\frac{r}{|\widehat\gamma_{s,j}|}\bigg\}
\]
and update $r$ with $r-|\widehat\gamma_{s,j}|t_j$.
Each coordinate is processed at most once, and we stop when $r=0$, requiring at most $d$ updates. Since $\mathcal H$ is nonempty, the available movements suffice to reach the boundary. Substituting $\delta_B=\sum_j t_j$ into \eqref{eq:app-two-hop-box} and then \eqref{eq:app-positive-evaluation} gives $\widehat p_s(x)$.

\medskip
\noindent
\textbf{Evaluation on $\mathcal{X}_{\widehat{\theta}_s}^-$.}
For $x \in \mathcal{X}_{\widehat{\theta}}^-$, we compute
\begin{equation}\label{eq:app-negative-evaluation-1}
\widehat p_s(x)=\max\{0,M_{\widehat\theta_s,\widehat q_s}(x)\}, \quad \mbox{with }M_{\widehat\theta_s,\widehat q_s}(x) =\max_{h\in\mathcal H} \{\widehat q_s(h)-L\|h-x\|_1\}.
\end{equation}
To evaluate $M_{\widehat\theta_s,\widehat q_s}(x)$, we solve the following equivalent  problem:
\begin{equation}\label{eq:app-negative-regional}
\begin{aligned}
\max_{h,r,w,\epsilon}\quad
&r-L\sum_{j=1}^d w_j\\
\text{subject to}\quad
&h\in\mathcal H,\quad 0\leq r\leq1,\\
&\epsilon_{ij}\in\{-1,1\},\quad
 \epsilon_{ij}(h_j-X_{ij})\geq0,
&&i\in[n],\ j\in[d],\\
&r\leq\widehat z_{s,i}
 +L\sum_{j=1}^d\epsilon_{ij}(h_j-X_{ij}),
&&i\in[n],\\
&w_j\geq h_j-x_j,\quad w_j\geq x_j-h_j,
&&j\in[d].
\end{aligned}
\end{equation}
The sign constraints ensure that
$\epsilon_{ij}(h_j-X_{ij})=|h_j-X_{ij}|$.
 However, the formulation is nonlinear because of the products $\epsilon_{ij}h_j$.  Introducing binary variables $b_{ij}\in\{0,1\}$ and auxiliary variables $v_{ij}$ representing $|h_j-X_{ij}|$ gives the following equivalent mixed-integer linear programme:
\begin{equation}\label{eq:app-negative-regional-lp}
\begin{aligned}
\max_{h,r,w,v,b}\quad
&r-L\sum_{j=1}^d w_j\\
\text{subject to}\quad
&h\in\mathcal H,\quad 0\leq r\leq1,\\
&b_{ij}\in\{0,1\},
&&i\in[n],\ j\in[d],\\
&v_{ij}\geq h_j-X_{ij},\quad
 v_{ij}\geq X_{ij}-h_j,
&&i\in[n],\ j\in[d],\\
&v_{ij}\leq h_j-X_{ij}+2X_{ij}(1-b_{ij}),
&&i\in[n],\ j\in[d],\\
&v_{ij}\leq X_{ij}-h_j+2(1-X_{ij})b_{ij},
&&i\in[n],\ j\in[d],\\
&r\leq\widehat z_{s,i}+L\sum_{j=1}^d v_{ij},
&&i\in[n],\\
&w_j\geq h_j-x_j,\quad w_j\geq x_j-h_j,
&&j\in[d].
\end{aligned}
\end{equation}
The constraints enforce $v_{ij}=|h_j-X_{ij}|$, and maximisation sets $r=\widehat q_s(h)$ and $w_j=|h_j-x_j|$. The optimal value therefore equals $M_{\widehat\theta_s,\widehat q_s}(x)$, which is clipped at zero to obtain $\widehat p_s(x)$.

\section{Additional numerical results}\label{sec:app-numerical}

In this section, we present additional simulation details and results in \Cref{sec:app-simulation} and describe the preprocessing procedures for the real-data analysis in \Cref{sec:app-real-data}.

\subsection{Additional results for the simulation studies}\label{sec:app-simulation}

In this subsection, we first give the oracle policies for the Arm 1 small and Arm 2 small settings and then present additional simulation results to complement \Cref{sec:numerical-simulation}.

\medskip
\noindent
\textbf{Explicit oracle policies for the Arm 1 small and Arm 2 small settings.} Both settings have uniform covariates on $[0,1]^2$. The following policies are optimal under the $\ell_1$-Lipschitz fairness constraint for every $L\geq0$.

In the Arm 1 small setting, $\Delta^*(x)=x_2-0.82+0.27x_1$. An optimal policy is
\[
p_{L,1}^*(x)= \bigg[\min\bigg\{ L(x_2-0.64+0.54x_1), \frac12+L(x_2-0.82+0.27x_1) \bigg\}\bigg]_{[0,1]}.
\]

In the Arm 2 small setting, $\Delta^*(x)=x_1+x_2-0.2$. An optimal policy is
\[
p_{L,1}^*(x)= \begin{cases}
[1+L(x_1+x_2-0.3)]_{[0,1]}, & 0\leq L\leq10/3,\\[4pt]
\displaystyle
\bigg[\frac12+L\bigg\{ x_1+x_2-\frac{0.2+\sqrt{0.04-1/(3L^2)}}{2} \bigg\}\bigg]_{[0,1]}, & L>10/3.
\end{cases}
\]

\begin{figure}[ht]
 \centering
  \begin{subfigure}{.49\linewidth}
  \includegraphics[width=\linewidth]{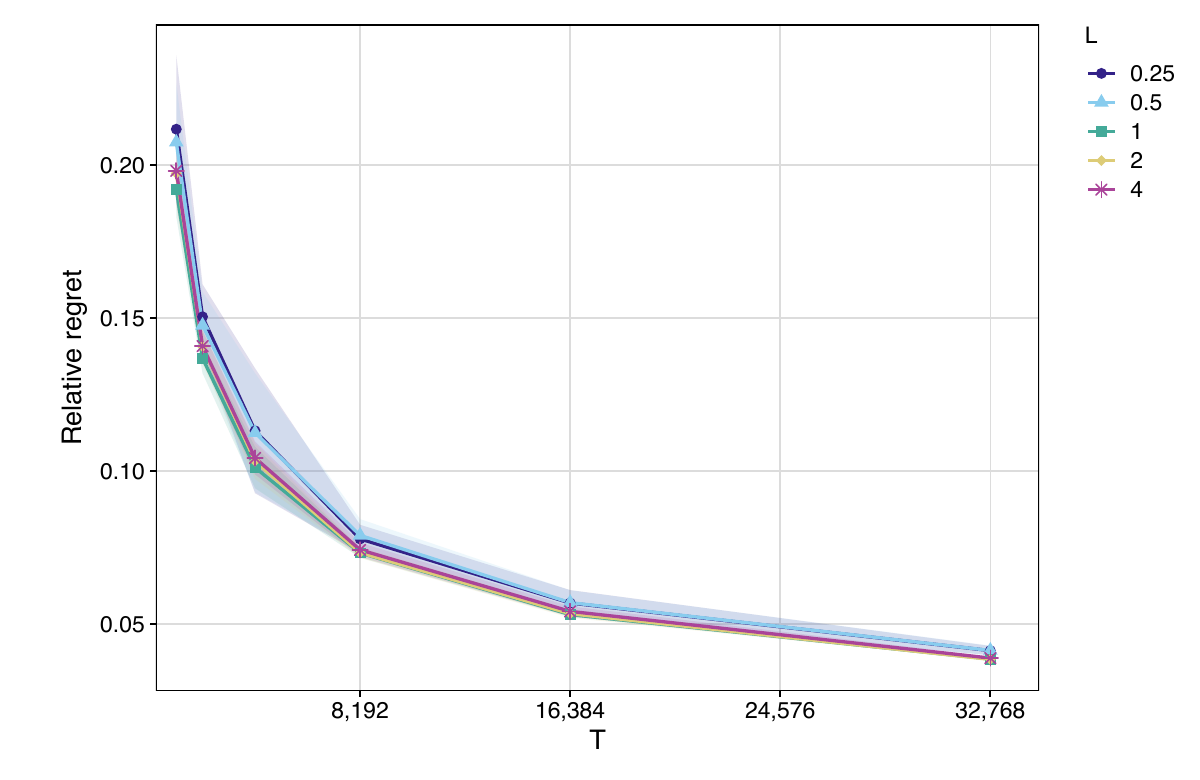}
  \caption{Relative regret on linear axes.}
 \end{subfigure}
 \begin{subfigure}{.49\linewidth}
\includegraphics[width=\linewidth]{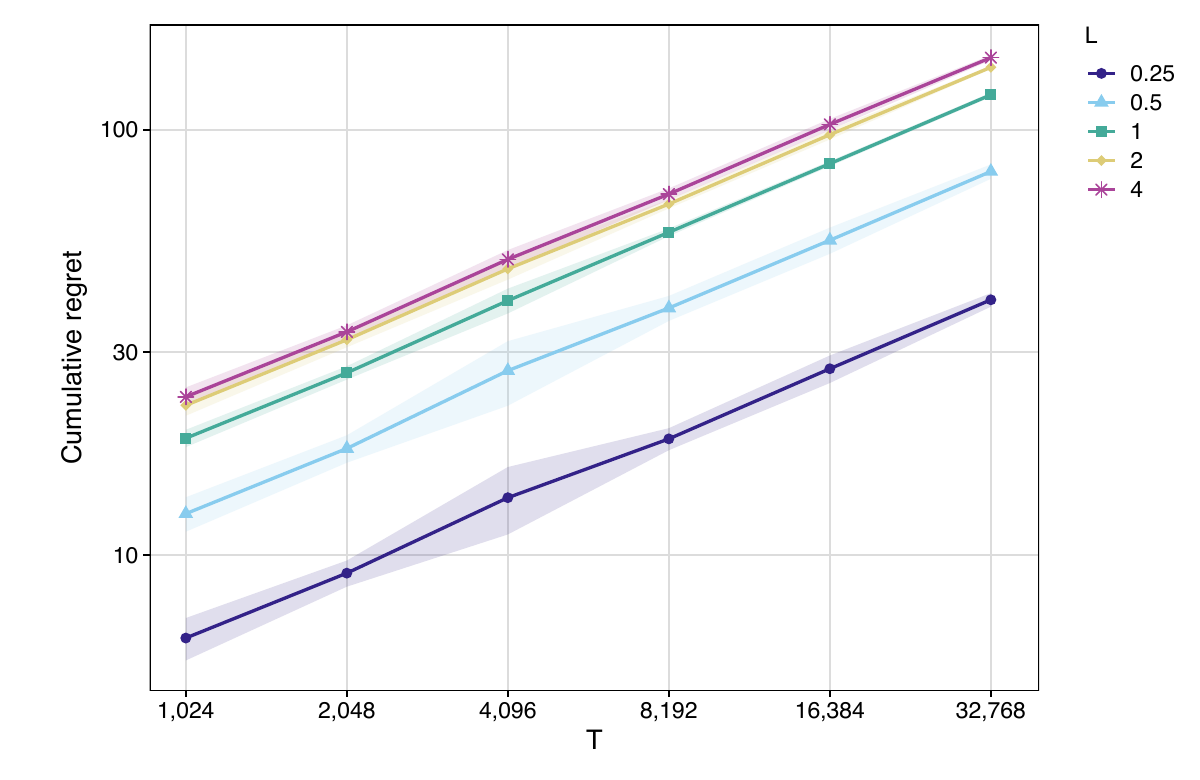}
  \caption{Cumulative regret on log--log axes.}
 \end{subfigure}\hfill
\caption{Relative regret (left) and cumulative regret (right) against $T$ for $L\in\{0.25,0.5,1,2,4\}$,  under uniform covariates and homogeneous rewards with $d=2$. Lines show means over $50$ independent Monte Carlo trials; shaded bands show the mean $\pm$ one standard deviation.}
 \label{fig:app-vary-L}
\end{figure}

\medskip
\noindent
\textbf{Varying the fairness parameter.}
We vary $L\in\{0.25,0.5,1,2,4\}$ under uniform covariates and homogeneous rewards with $d=2$, using the six horizons $T\in\{1024,2048,4096,8192,16384,32768\}$. \Cref{fig:app-vary-L} reports relative and cumulative regret, averaged over $50$ independent Monte Carlo trials. Relative regret decreases with $T$ for every value of $L$, consistent with the sublinear cumulative regret bound for fixed $L$ in \Cref{thm:fair_regret_control}. The relative regret curves remain close across the values of $L$ considered. Mean cumulative regret increases with $L$, while the curves are approximately linear on log--log axes, with similar slopes of about $0.53$.

\begin{figure}[ht]
\centering
\includegraphics[width=.8\linewidth]{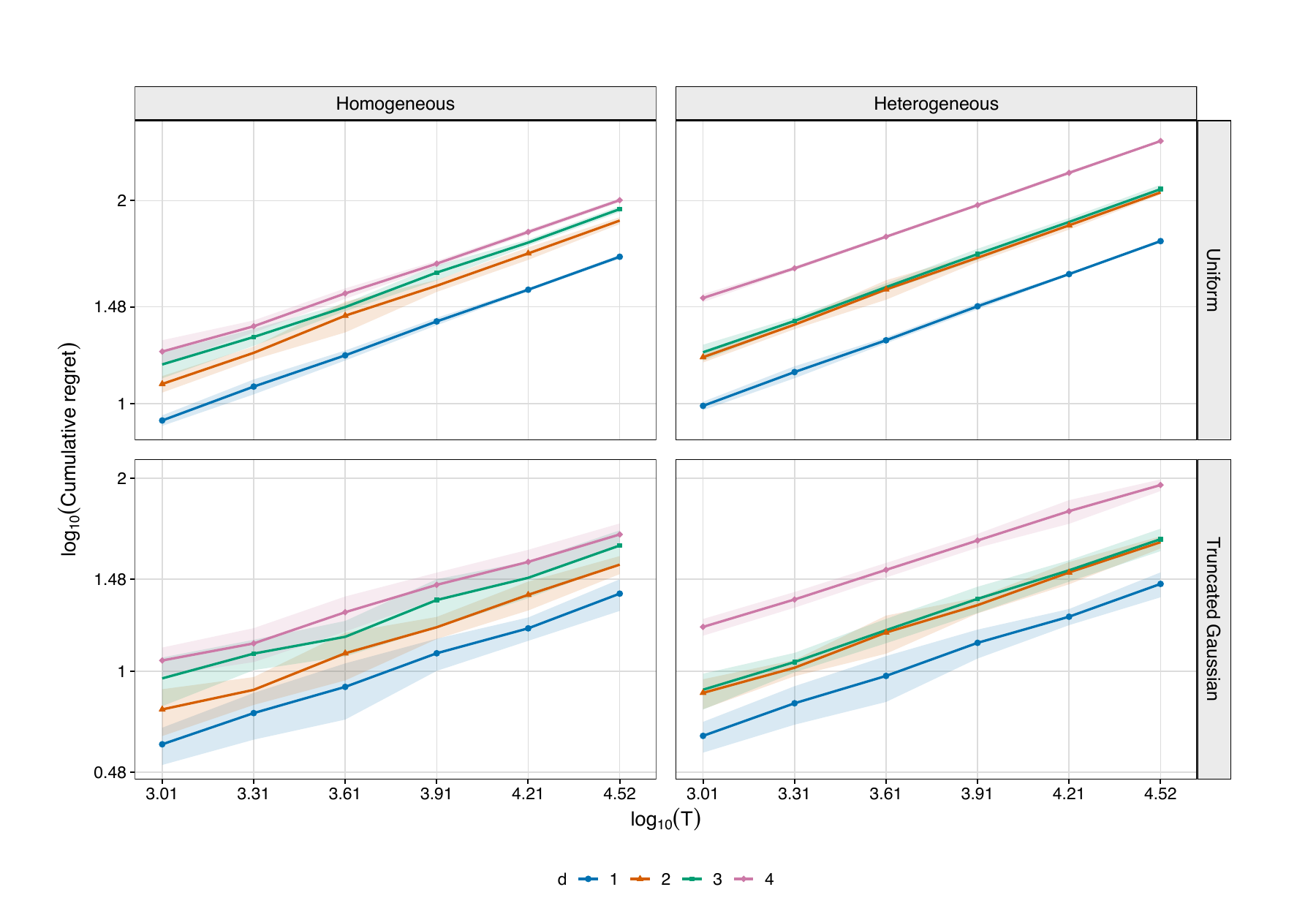}
\caption{Cumulative regret against $T$ on log--log axes for $d=1,2,3,4$ at $L=0.5$. Columns show homogeneous (left) and heterogeneous (right) rewards; rows show uniform (top) and truncated-Gaussian (bottom) covariates. Lines show means over $50$ independent Monte Carlo trials; shaded bands show the mean $\pm$ one standard deviation.}
 \label{fig:app-ordinary-loglog}
\end{figure}

\medskip
\noindent
\textbf{Cumulative regret in the baseline settings.}
\Cref{fig:app-ordinary-loglog} plots cumulative regret against $T$ on log--log axes for the four baseline settings and $d=1,2,3,4$, with $L=0.5$. The curves are approximately linear, with slopes generally close to $1/2$ over the observed horizons.

\subsubsection{An additional simulation setting}\label{sec:app-mixed-design}

We next consider an additional simulation setting with $d=4$, modifying both the covariate distribution and the reward functions. We first define the covariate distribution. Let
\[
X=(Z,U),\quad Z=(G-1)/5,\quad G\in\{1,\ldots,6\},
\quad \mbox{and} \quad U\sim\mathrm{Unif}([0,1]^3).
\]
Here, $Z$ encodes a six-category covariate whose conditional distribution depends on $U$. 

To specify how the category probabilities vary with $U$, partition $[0,1]^3$ into $m^3$ equal cubic cells with $m=4$. For a point $u$, define its cell indices and its rescaled coordinates within that cell by
\[
k_j(u)=\min\{\lfloor mu_j\rfloor,m-1\},
\quad r_j(u)=mu_j-k_j(u), 
\quad \forall j \in \{1,2,3\}.
\]
Write $k(u)=(k_1(u),k_2(u),k_3(u))$ for the cell containing $u$. Within each cell, we construct two equal-volume subregions that favour opposite groups of categories. Independently draw a sign $\eta_k\in\{-1,1\}$ with equal probabilities for each cell $k$, and keep these signs fixed across horizons and replications. Set
\begin{equation}\label{eq:mixed-H}
H(u)=\eta_{k(u)}
\Big[\mathbbm{1}\{r_1(u)\in[1/4,1/2)\}
-\mathbbm{1}\{r_1(u)\in[1/2,3/4)\}\Big]
\prod_{j=2}^{3}\mathbbm{1}\{r_j(u)\in[1/4,3/4)\}.
\end{equation}
Thus, $H(u)$ takes opposite signs on the two subregions and is zero elsewhere.

We use $H(u)$ to define the conditional category probabilities:
\begin{equation}\label{eq:mixed-weights}
w_g(u)=\mathbb P(G=g\mid U=u)
=\frac{1+ s_gH(u)/2}{6},
\quad \mbox{where }
s_g=\begin{cases}
-1,&g\leq3,\\
1,&g\geq4.
\end{cases}
\end{equation}
When $H(u)=1$, categories $4,5,6$ are favoured; when $H(u)=-1$, categories $1,2,3$ are favoured. When $H(u)=0$, all six categories are equally likely. These conditional probabilities change discontinuously across subregions. Since the positive and negative contributions of $H$ cancel within each cell, $\int_{[0,1]^3}H(u)\mathrm du=0$. Each category therefore has marginal probability $1/6$, although $G$ and $U$ are dependent.

We next specify the mean rewards:
\[
\mu_1(g,u)=\frac{g-1}{5}-\frac12
+0.03\sum_{j=1}^{3}(u_j-1/2)
\quad \mbox{and} \quad
\mu_2(g,u)=0.
\]
We generate rewards by adding independent $\mathcal N(0,0.3^2)$ noise to these means. We use the lazy fair greedy algorithm in \Cref{alg:lazy-adaptive-boundary-greedy} with $L=0.5$ and $\nu=1$, following the same exploration schedule as in the baseline experiments.

\medskip
\noindent
\textbf{Numerical oracle.}
We approximate the fair oracle on a regular grid $\mathcal V=\{0,h,\ldots,1\}^3$, with $h=1/64$, for each of the six categories.  Let $p_{g,v}\in[0,1]$ denote the probability of selecting arm~1 for category $g$ at grid point $v\in\mathcal V$.  We solve a linear programme in these $6\times65^3=1647750$ variables to maximise population reward subject to the $\ell_1$-Lipschitz fairness constraint at the grid points. We evaluate population reward for this linear programme using two-point Gaussian quadrature in each coordinate, splitting the integration domain at grid boundaries and discontinuities of $H$.

We then extend the resulting probabilities to all $u\in[0,1]^3$ by trilinear interpolation within each category. Fo $u$ in a grid cell $v+[0,h]^3$, define
\[
\alpha_j=\frac{u_j-v_j}{h}, \quad \forall  j=1,2,3.
\]
The interpolated policy is
\[
p^{\mathrm{num}}\bigg(\frac{g-1}{5},u\bigg) = \sum_{b\in\{0,1\}^3} \bigg\{\prod_{j=1}^3\omega_{b_j}(\alpha_j)\bigg\}
p_{g,v+hb}, \quad \mbox{where } \omega_0(a)=1-a \mbox{ and } \omega_1(a)=a.
\]
Thus, the probability at $u$ is a weighted average of the probabilities at the cell's eight vertices. The weights are nonnegative and sum to one, so the interpolated policy takes values in $[0,1]$. It agrees with the LP solution at grid points and preserves the $\ell_1$-Lipschitz fairness constraint.

\medskip
\noindent
\textbf{Reward evaluation and results.}
We evaluate population reward by integrating over $U$ using $32768$ Halton points in $[0,1]^3$ after skipping the first $128$, and summing over the six categories with weights $w_g(u)$. We denote the resulting reward estimate for a policy $p$ by $\widehat R(p)$. Within each replication, the numerical oracle and learned policies use the same integration points and weights. We compute estimated cumulative and relative regret as
\[
\widehat{\mathfrak R}_T =\sum_{t=1}^{T} \big\{\widehat R(p^{\mathrm{num}}) -\widehat R(\widehat p_t)\big\} \quad \mbox{and} \quad
\widehat{\mathrm{RelReg}}_T
=\frac{\widehat{\mathfrak R}_T}{T\widehat R(p^{\mathrm{num}})},
\]
respectively, including initial exploration.

\begin{figure}[ht]
\centering
\begin{subfigure}{.45\linewidth}
    \centering
    \includegraphics[width=\linewidth]{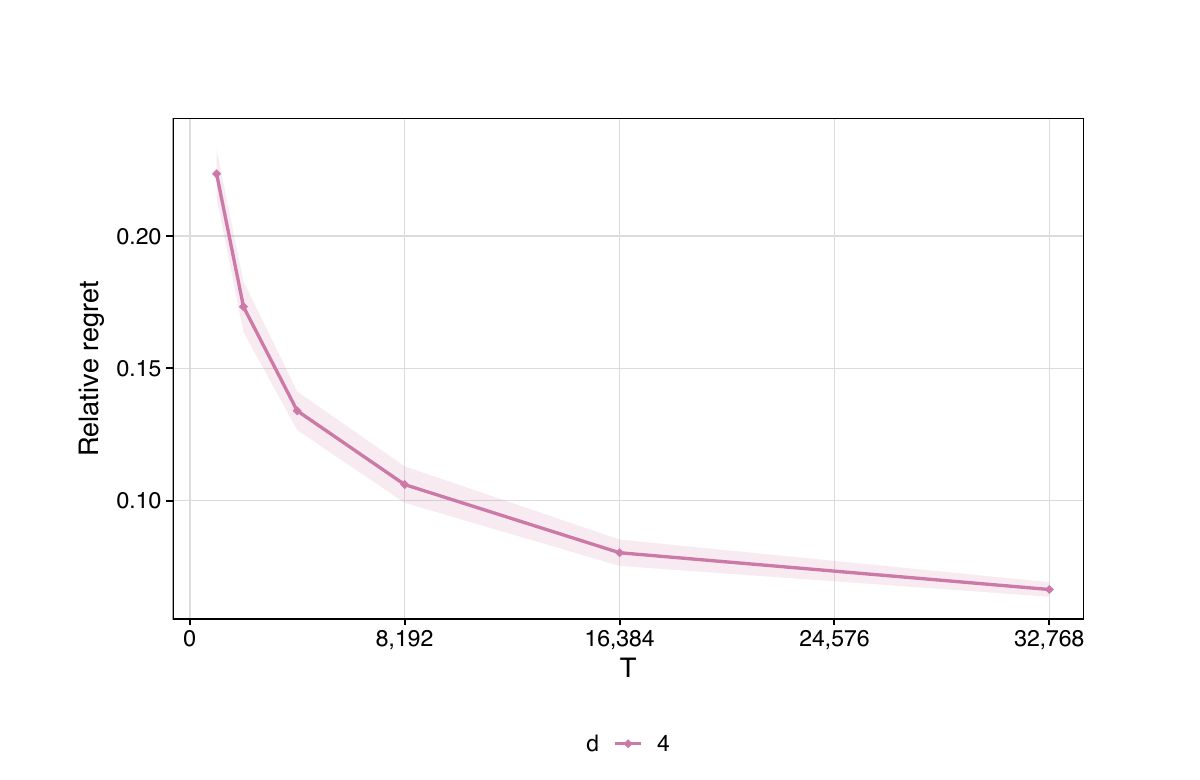}
    \caption{Relative regret on linear axes.}
\end{subfigure}\hfill
\begin{subfigure}{.45\linewidth}
    \centering
    \includegraphics[width=\linewidth]{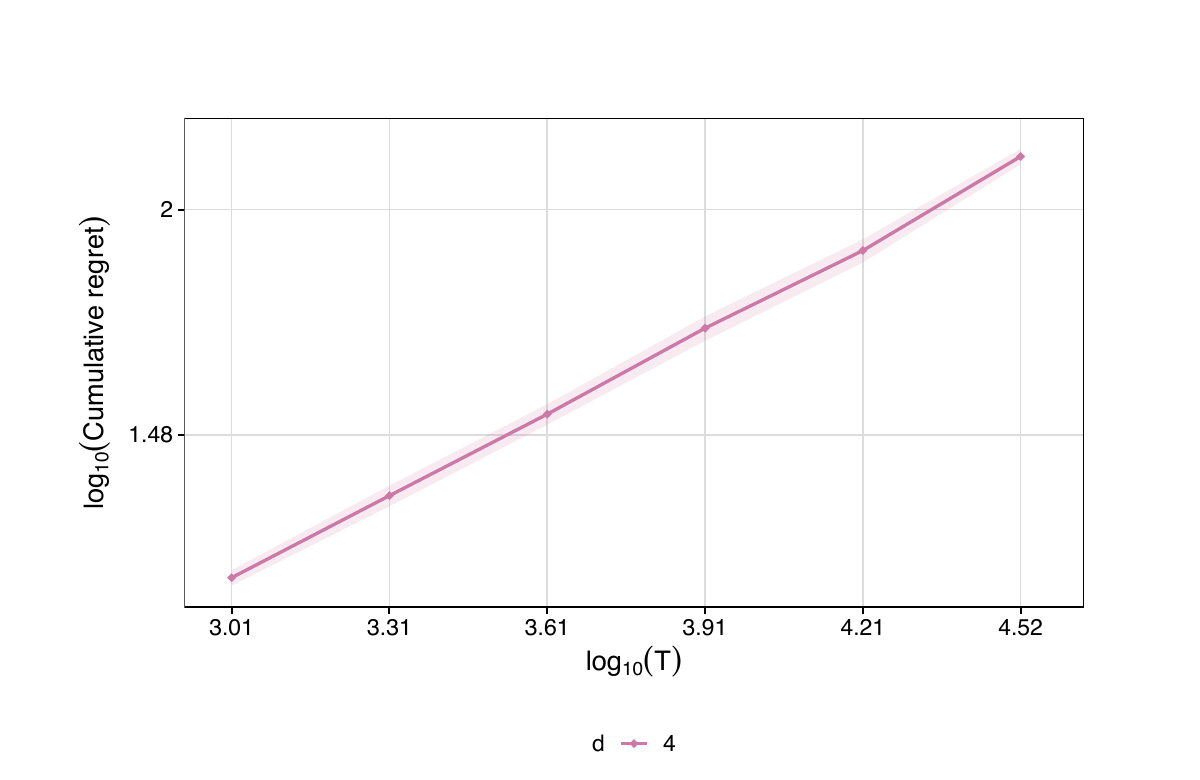}
    \caption{Cumulative regret on log--log axes.}
\end{subfigure}
\caption{Estimated relative regret (left) and cumulative regret (right) against $T$ for the additional simulation setting with $d=4$ and $L=0.5$. Lines show means over $50$ independent Monte Carlo trials; shaded bands show the mean $\pm$ one standard deviation.}
\label{fig:app-mixed}
\end{figure}

\Cref{fig:app-mixed} reports mean estimated relative and cumulative regret at $L =0.5$, based on $50$ independent Monte Carlo trials, at $T=1024,2048,4096,8192,16384,32768$, respectively. Relative regret decreases with $T$. Mean cumulative regret plotted against $T$ is approximately linear on log--log axes, with an empirical slope of $0.645$, numerically close to the exponent $2/3$ in \Cref{thm:fair_regret_control}.

\subsection{Additional details for the real data analysis}\label{sec:app-real-data}

In this subsection, we describe the sample selection and covariate processing used in \Cref{sec:numerical-real-data}.

\medskip
\noindent
\textbf{Datasets and sample selection.}
For Adult \citep{becker1996adult}, we remove records with missing values and retain records from Canada, El Salvador, Germany, India, Mexico, the Philippines, Puerto Rico and the United States. We set the response to $B=1$ when annual income exceeds $\$50000$ and $B=0$ otherwise. The $12$ covariates are age, work class, final weight, numeric education level, marital status, occupation, relationship, race, sex, capital gain, capital loss and weekly working hours. The resulting sample contains $43303$ records.

For Jobs \citep{tankha2023jobs}, we remove the index column and the detailed list of technologies used, then discard records with missing values. We use the binary employed field as the response. The $11$ covariates are age group, accessibility indicator, education level, employment status at the time of application, gender, a mental health indicator, professional developer status, years of coding experience, years of professional coding experience, previous salary and the number of technologies used. Country is excluded from the covariates, and all observed gender categories are retained. The resulting sample contains $73462$ records.

\medskip
\noindent
\textbf{Covariate processing.}
For each dataset with $n$ eligible records, we randomly select $\lfloor0.1n\rfloor$ records without replacement to form a preprocessing subset. We fit the covariate transformation on this subset without using response labels and apply it unchanged to all covariates used for online learning. The preprocessing subset is excluded from the online sampling pool, and both remain fixed across horizons, fairness parameters and replications.

For Adult, the categorical covariates are work class, marital status, occupation, relationship, race and sex; the remaining covariates are numeric. For Jobs, the categorical covariates are age group, accessibility, education level, gender, mental health and professional developer status; the remaining covariates are numeric. We one-hot encode categorical covariates using one indicator for each category observed in the preprocessing subset. An unseen category is represented by zeros in all indicator columns for the corresponding covariate. After encoding, Adult and Jobs have $46$ and $21$ covariates, respectively.

We standardise each encoded covariate using its mean and sample standard deviation in the preprocessing subset. We then fit principal component analysis (PCA)  to the standardised covariates in this subset and retain the first three components for policy learning \citep{jolliffe2016principal}.

Let $S_{ij}$ denote the score of record $i$ on component $j$, and let $a_j$ and $b_j$ be the minimum and maximum scores in the preprocessing subset. We rescale and clip the scores as
\[
X_{ij} =\bigg[\frac{S_{ij}-a_j}{b_j-a_j}\bigg]_{[0,1]},
\quad j=1,2,3.
\]
Applying this fixed transformation gives $X_i\in[0,1]^3$ for every record. 
\end{document}